\documentclass[a4paper, reqno, 11pt, notitlepage]{amsart}
\usepackage{fullpage}
\usepackage[utf8]{inputenc}

\usepackage{latexsym}
\usepackage{amsmath}
\usepackage{amssymb}
\usepackage{amsthm}
\usepackage{amscd}
\usepackage{mathrsfs}
\usepackage[all]{xy}
\usepackage{graphicx}
\usepackage{comment}

\usepackage{mathtools}
\usepackage{dsfont}
\usepackage{adjustbox}
\usepackage{multicol}
\usepackage{quiver}
\usepackage{changepage} 
\usepackage[activate={true,nocompatibility},final,tracking=true,kerning=true,spacing=true]{microtype}
\microtypecontext{spacing=nonfrench}
\usepackage{eucal}
\usepackage{yfonts}
\DeclareFontFamily{U}{mathx}{}
\DeclareFontShape{U}{mathx}{m}{n}{<-> mathx10}{}
\DeclareSymbolFont{mathx}{U}{mathx}{m}{n}
\DeclareMathAccent{\widecheck}{0}{mathx}{"71}
\usepackage{accents}

\numberwithin{equation}{subsection}

\usepackage{array}
\usepackage{booktabs}
\usepackage{tabularx}
\usepackage{longtable}
\usepackage{arydshln}

\makeatletter

\renewcommand{\tocsection}[3]{%
  \indentlabel{\@ifnotempty{#2}{\bfseries\ignorespaces#1 #2\quad}}\bfseries#3}

\renewcommand{\tocsubsection}[3]{%
  \indentlabel{\@ifnotempty{#2}{\ignorespaces#1 #2\quad}}#3}

\newcommand\@dotsep{4.5}
\def\@tocline#1#2#3#4#5#6#7{\relax
  \ifnum #1>\c@tocdepth %
  \else
    \par \addpenalty\@secpenalty\addvspace{#2}%
    \begingroup \hyphenpenalty\@M
    \@ifempty{#4}{%
      \@tempdima\csname r@tocindent\number#1\endcsname\relax
    }{%
      \@tempdima#4\relax
    }%
    \parindent\z@ \leftskip#3\relax \advance\leftskip\@tempdima\relax
    \rightskip\@pnumwidth plus1em \parfillskip-\@pnumwidth
    #5\leavevmode\hskip-\@tempdima{#6}\nobreak
    \leaders\hbox{$\m@th\mkern \@dotsep mu\hbox{.}\mkern \@dotsep mu$}\hfill
    \nobreak
    \hbox to\@pnumwidth{\@tocpagenum{\ifnum#1=1\bfseries\fi#7}}\par%
    \nobreak
    \endgroup
  \fi}
\AtBeginDocument{%
\expandafter\renewcommand\csname r@tocindent0\endcsname{0pt}
}
\def\l@subsection{\@tocline{2}{0pt}{2.5pc}{5pc}{}}
\makeatother

\makeatletter
\renewcommand{\paragraph}{%
  \@startsection{paragraph}{4}%
    {\z@}%
    {6pt}%
    {-\fontdimen2\font}%
    {\normalfont\bfseries}%
}

\makeatother

\makeatletter
\def\subsubsection{\@startsection{subsubsection}{3}
  \z@{.5\linespacing\@plus.7\linespacing}{-.5em}
  {\normalfont\bfseries}}
\makeatother

\usepackage[T1]{fontenc}

\usepackage{tikz-cd}
\usepackage{nicefrac}
\usepackage{authblk}
\usepackage{textcomp}
\usepackage{graphicx}
\usepackage{framed} \usepackage{color}  

\usepackage[pagebackref]{hyperref} 
    \definecolor{darkblue}{rgb}{0,0,.85} 
    \definecolor{darkred}{rgb}{0.84,0,0}
    \hypersetup{colorlinks = true,
            linkcolor = darkblue,
            urlcolor  = darkblue,
            citecolor = darkred,
            anchorcolor = darkblue}

\usepackage{accents}

\usepackage[shortlabels]{enumitem}

\usepackage[foot]{amsaddr}

\usepackage{bm}

\usepackage{stmaryrd}

\usepackage{graphics}

\usepackage{mathtools}

\setlist[enumerate]{font=\upshape}
\setlist[itemize]{font=\upshape}

\newtheorem{thm}{Theorem}[section]
\newtheorem{prop}[thm]{Proposition}
\newtheorem{propnota}[thm]{Proposition/Notation}

\newtheorem{thmi}{Theorem}
\newtheorem{propi}[thmi]{Proposition}
\newtheorem{cori}[thmi]{Corollary}

\newtheorem{lem}[thm]{Lemma}

\newtheorem{cor}[thm]{Corollary}
\newtheorem{conj}[thm]{Conjecture}
\newtheorem{conji}[thmi]{Conjecture}

\theoremstyle{definition}

\newtheorem{eg}[thm]{Example}

\newtheorem{obs}[thm]{Observation}

\newtheorem{obsnota}[thm]{Observation/Notation}

\newtheorem{rem}[thm]{Remark}
\newtheorem{remdef}[thm]{Remark/Definition}
\newtheorem{remi}{Remark}

\newtheorem{defn}[thm]{Definition}
\newtheorem{defni}[thmi]{Definition}
\newtheorem{nota}[thm]{Notation}
\newtheorem{setup}[thm]{Setup}
\newtheorem{assumption}[thm]{Assumption}

\newtheorem*{defni*}{Definition}

\newtheorem{constr}[thm]{Construction}
\newtheorem{construction}[thm]{Construction}

\DeclareMathOperator{\Spec}{Spec}

\DeclareMathOperator{\Gal}{Gal}
\DeclareMathOperator{\id}{id}
\DeclareMathOperator{\tr}{tr}

\DeclareMathOperator{\End}{End}
\DeclareMathOperator{\Hom}{Hom}

\DeclareMathOperator{\Spf}{Spf}

\DeclareMathOperator{\Aut}{Aut} 
\DeclareMathOperator{\Fil}{Fil}

\DeclareMathOperator{\Nr}{Nr}

\DeclareMathOperator{\Spa}{Spa}

\DeclareMathOperator{\Art}{Art}

\DeclareMathOperator{\pr}{pr}

\DeclareMathOperator{\Spd}{Spd}

\DeclareMathOperator{\GL}{GL}

\DeclareMathOperator{\Int}{Int}

\newcommand{\dotimes}{\mathbin{\otimes^{\scriptscriptstyle\mathbf L}}}

\newcommand{\dquot}{{/^{\scriptscriptstyle\mathbf L}}\,}

\renewcommand{\sp}{\mathrm{sp}}

    \DeclareFontFamily{U}{wncy}{}
    \DeclareFontShape{U}{wncy}{m}{n}{<->wncyr10}{}
    \DeclareSymbolFont{mcy}{U}{wncy}{m}{n}
    \DeclareMathSymbol{\Sha}{\mathord}{mcy}{"58}

\newcommand{\A}{\mathbb{A}}

\newcommand{\F}{\mathbb{F}}

\newcommand{\bT}{\mathbb{T}}

\newcommand{\bb}[1]{\mathbb{#1}}

\newcommand{\mc}[1]{\mathcal{#1}}

\newcommand{\fD}{\mathfrak{D}}
\newcommand{\fE}{\mathfrak{E}}

\newcommand{\mf}[1]{\mathfrak{#1}}

\usepackage{relsize}

\newcommand{\syn}{\mathrm{syn}}

\newcommand{\wh}{\widehat}
\newcommand{\wt}{\widetilde}

\renewcommand{\sc}{\mathrm{sc}}
\newcommand{\mb}[1]{\mathbf{#1}}
\newcommand{\Z}{\mathbb{Z}}

\newcommand{\Q}{\mathbb{Q}}

\newcommand{\C}{\mathbb{C}}
\newcommand{\univ}{\mathrm{univ}}

\renewcommand{\ll}{\llbracket}
\newcommand{\rr}{\rrbracket}

\renewcommand{\j}{\jmath}

\makeatletter
\renewcommand{\email}[2][]{%
  \ifx\emails\@empty\relax\else{\g@addto@macro\emails{,\space}}\fi%
  \@ifnotempty{#1}{\g@addto@macro\emails{\textrm{(#1)}\space}}%
  \g@addto@macro\emails{#2}%
}
\makeatother

\newcommand{\triv}{\mathrm{triv}}

\newcommand{\der}{\mathrm{der}}

\newenvironment{sectionappendices}
{%
  \setcounter{subsection}{0}%
  \renewcommand{\thesubsection}{\thesection.\Alph{subsection}}%
}
{}

\newcommand{\sectionappendix}[1]{%
  \subsection{#1}%
  \setcounter{thm}{0}%
  \setcounter{question}{0}%
  \setcounter{thmi}{0}%
  \setcounter{remi}{0}%
}

\newcommand{\stacks}[1]{\cite[\href{https://stacks.math.columbia.edu/tag/#1}{Tag~#1}]{StacksProject}}

\newcommand{\an}{\mathrm{an}}

\newcommand{\cat}[1]{\mathbf{#1}}

\newcommand{\Sh}{\mathrm{Sh}}

\newcommand{\ms}[1]{\mathscr{#1}}

\newcommand{\R}{\mathbb{R}}

\def\Item(#1){\item[\llap{(}\refstepcounter{enumi}$\bullet$] #1)}

\newcommand{\ovc}[1]{\accentset{\circ}{#1}}

\DeclareMathOperator*{\twolim}{2-lim}

\DeclareMathAccent{\wtilde}{\mathord}{largesymbols}{"65}

\newcommand*\isomto{%
        \xrightarrow{\raisebox{-0.2 em}{\smash{\ensuremath{\sim}}}}%
    }

    \newcommand*\isomfrom{%
        \xleftarrow{\raisebox{-0.2 em}{\smash{\ensuremath{\sim}}}}%
    }

    \newcommand{\ov}[1]{\overline{#1}}
    \newcommand{\et}{\mathrm{\acute{e}t}}

    \DeclareMathOperator{\Res}{Res}
    
   \newcommand{\defeq}{\vcentcolon=}
    
    \newcommand{\be}{\begin{equation*}}
    \newcommand{\ee}{\end{equation*}}
    \newcommand{\bx}{\begin{equation*}\xymatrix}
    \newcommand{\ex}{\end{equation*}}

\usepackage{relsize}
\usepackage[bbgreekl]{mathbbol}
\usepackage{amsfonts}

\DeclareSymbolFontAlphabet{\mathbbl}{bbold}
\newcommand{\Prism}{{\mathlarger{\mathbbl{\Delta}}}}

\newcommand{\smallprism}{{{\mathsmaller{\Prism}}}}
\newcommand{\smallN}{{{\mathsmaller{\mc{N}}}}}

\usepackage{scalerel}

\NewDocumentCommand\limder{e{_}}{\mathchoice
{\lim  \IfValueT{#1}{_{\mathclap{#1}}}{}^{\!1\!}\mathop{}}
{\lim^1  \IfValueT{#1}{_{#1}}}
{\lim^1  \IfValueT{#1}{_{#1}}}
{\lim^1  \IfValueT{#1}{_{#1}}}}

\newcommand{\crys}{\mathrm{crys}}

\newcommand{\mr}{\mathrm}
\newcommand{\ad}{\mathrm{ad}}

\title{On the Langlands--Kottwitz--Scholze method}
\author{Alex Youcis}
\address{\fontsize{7.5}{8.4}\selectfont Department of Mathematics, Bahen Centre, University of Toronto, Toronto, ON, M5S 2E4, Canada}
\email{alex.youcis@gmail.com}

\date{\today}

\begin{document}

\begin{abstract} Using recent advances in the integral canonical models of Shimura varieties using syntomic methods, we give an extension of the Langlands--Kottwitz--Scholze method from \cite{ScholzeLK}. Namely, we define analogues $\phi^{\mc{G},\mu}_{\tau,h}$ of the local test functions from op.\@ cit.\@ in full generality, show they satisfy reasonable harmonic-analytic properties, and give rise to trace formulae for the Galois-Hecke action on the cohomology of Shimura varieties of bad reduction. Additionally, we use these new local test functions and the work of Fargues--Scholze to state a more robust and unconditional version of the Scholze--Shin conjecture from \cite{ScholzeShin}. Along the way we study the $(\mc{G},\mu)$-apertures of \cite{GMM}, in particular showing that such $(\mc{G},\mu)$-apertures carry no more information than a Tannakian version of Fontaine--Laffaille theory, at least over nice bases, and give a very explicit description of the universal deformations of $(\mc{G},\mu)$-apertures.
\end{abstract}

\maketitle

\tableofcontents

\section*{Introduction}

\setcounter{subsection}{0}
\renewcommand{\thesubsection}{\arabic{subsection}}
\addtocontents{toc}{\protect\setcounter{tocdepth}{0}}

In the series of papers \cite{ScholzeLKModularCurve,ScholzeLKSimple,ScholzeLK},  there are defined functions
\begin{equation*}
\phi_{\tau,h}^{\mc{G},\mu}\colon G(\Q_{p^r})\to \Q,
\end{equation*}
associated to certain local data $(\mc{G},\mu,\tau,h)$ which we soon recall, called \emph{local test functions}. These local test functions fulfill two, ostensibly disparate but important, roles:

\medskip

\noindent\textbf{The Langlands--Kottwitz--Scholze method} Let $(\mb{G},\mb{X})$ be a Shimura datum,\footnote{For the sake of simplicity, in this introduction we assume that $\mb{G}=\mb{G}^c$; see \S\ref{sss:etale-realization}. We do not make this assumption in the body of the paper.} and $\Sh_\mathsf{K}(\mb{G},\mb{X})$ its canonical model  at level $\mathsf{K}=\mathsf{K}_p\mathsf{K}^p$ over its reflex field $\mathbf{E}$. Also fix $\mc{F}_{\xi,\mathsf{K}}$ an $\ell$-adic automorphic \'etale sheaf attached to an $\ell$-adic representation $\xi$ of $\mb{G}$ for some prime $\ell\ne p$; see \S\ref{sss:etale-realization}.

In \cite{KottwitzAnnArbor} one finds a formula of the form
\begin{equation*}
\tr\left(\tau\times f^p\mathds{1}_{\mathsf{K}_p}\mid H^\ast_\et(\Sh_\mathsf{K}(\mb{G},\mb{X})_{\ov{E}},\mc{F}_{\xi,\mathsf{K}})\right)=\sum_{(\gamma_0;\gamma,\delta)}c(\gamma_0;\gamma,\delta)\tr(\xi(\gamma_0))\mr{O}_\gamma(f^p)\mr{TO}_\delta(\mathds{1}_{S_j})
\end{equation*}
where $\mathsf{K}_p$ is \emph{hyperspecial}, and where
\begin{itemize}
\item $H^\ast_\et =\sum_{i\geqslant 0}(-1)^i H^i_\et$;
\item $\mc{G}$ is a reductive model of $G=\mb{G}_{\Q_p}$ with $\mathsf{K}_p=\mc{G}(\Z_p)$;
\item $v$ is a place of $\mb{E}$ lying over $p$ and where $E=\mb{E}_v$;
\item $\tau$ is an element of $W_E$ projecting to a $j^\text{th}$ power of geometric Frobenius;
\item $S_j=\mc{G}(\mc{O}_{E_j})\sigma(-\mu_h)(p)\mc{G}(\mc{O}_{E_j})$ where $\mu_h$ is the Hodge cocharacter;
\item $\mathds{1}_S$ means the indicator function of $S$;
\item and $(\gamma_0;\gamma,\delta)$  is a Kottwitz triple; see Construction \ref{constr:classical-KP}.
\end{itemize}
We call this formula and its method of proof, the \emph{Langlands--Kottwitz method}, which is quite important for verifying the decomposition of the cohomology of Shimura varieties as in op.\@ cit.\@

The Langlands--Kottwitz method is sufficient from a global perspective, but is locally quite restrictive at a given $p$ as it only considers cases where the cohomology is unramified at $p$. 

To remedy this, Scholze developed the \emph{Langlands--Kottwitz--Scholze} method, which allows $\mathds{1}_{\mathsf{K}_p}$ to be replaced by any locally constant compactly supported $\Q$-valued function $h$ on $\mc{G}(\Z_p)$ in exchange for replacing $\mathds{1}_{S_j}$ with the local test function $\phi^{\mc{G},-\mu_h}_{\tau,h}$ mentioned above. This gives rise to the following formula:
\begin{equation*}
T(\tau,h,f^p)\defeq \tr\left(\tau\times f^p h\mid H^\ast_\et(\Sh_\mathsf{K}(\mb{G},\mb{X})_{\ov{E}},\mc{F}_{\xi,\mathsf{K}})\right)=\sum_{(\gamma_0;\gamma,\delta)}c(\gamma_0;\gamma,\delta)\tr(\xi(\gamma_0))\mr{O}_\gamma(f^p)\mr{TO}_\delta(\phi_{\tau,h}^{\mc{G},-\mu_h}).
\end{equation*}
That said, note that Scholze constructed these local test functions only for so-called EL- and PEL-data of type $A$ or $C$ \cite[\S\S3--4]{ScholzeLK}. In the PEL setting, he established the above trace formula under a properness hypothesis
\cite[Theorem 5.7]{ScholzeLK}; this properness assumption was then later removed by Lan--Stroh in \cite[Theorem 6.36]{LanStrohI}.

\medskip

\paragraph*{The Scholze--Shin conjecture} On the other hand, the functions $\phi^{\mc{G},\mu}_{\tau,h}$ appear when trying to understand the relationship between a tempered representation $\pi$ of $G(\Q_p)$ and the semisimple $L$-parameter $\varphi_\pi\colon W_{\Q_p}\to {}^L G$ associated to $\pi$ by a hypothetical local Langlands correspondence. 

Namely, it is a conjecture of Scholze--Shin (see \cite[Conjecture 7.2]{ScholzeShin}) that one has the equality
\begin{equation*}
S\Theta_\pi(f_{\tau,h}^{\mc{G},\mu})=\tr\left(\tau\mid r_{\mu}\circ \varphi_\pi\right)S\Theta_\pi(h)q^{j\langle \rho_G,\mu\rangle},
\end{equation*}
where 
\begin{itemize}
\item $S\Theta_\pi$ is the stable distribution associated to (the $A$-packet of) $\pi$;
\item $r_{\mu}$ is the representation of $\wh{G}\rtimes W_E$ from Construction \ref{constr:Kottwitz};
\item $f_{\tau,h}^{\mc{G},\mu}$ is the twisted transfer of $\phi_{\tau,h}^{\mc{G},\mu}$ from a function on $G(E_j)$ to one on $G(\Q_p)$;
\item $\rho_G$ is the half-sum of positive roots;
\item and $q=p^r$ is the size of the residue field $k$ of $E$.
\end{itemize}
In other words, $f_{\tau,h}^{\mc{G},\mu}$ is a function which allows one (up to the cut-off function $h$) to relate the (stable) trace characters of $\pi$ and $\varphi_\pi$.

There are better versions of this conjecture allowing both for endoscopic groups of $G$ as well as more general stable characters; see \S\ref{s:Scholze-Shin}. In fact, using these more refined versions, this matching of traces between $\pi$ and $\varphi_\pi$ frequently does what one might hope: It characterizes the local Langlands correspondence for $G$; see \cite{BMYCharacterization}.

\medskip

\paragraph*{The local test functions}  Now that their importance has been made clear, let us explain the definition of the local test functions $\phi_{\tau,h}^{\mc{G},\mu}$ constructed by Scholze. The local data $(\mc{G},\mu,\tau,h)$ is:
\begin{itemize}
\item $\mc{G}$ is a reductive group $\Z_p$-scheme;
\item $\mu$ is a minuscule cocharacter of $\mc{G}$ defined over $W=W(k)$ for $k/\bb{F}_p$ a finite extension;
\item $j\geqslant 1$ an integer;
\item $k_j/k$ is the extension of degree $j$, $W_j=W(k_j)$, and $K_j=W_j[\nicefrac{1}{p}]$ (we omit $j$ when it's $1$); 
\item $\tau$ is an element $W_K$ projecting to a $j^\text{th}$-power of geometric Frobenius;
\item $h\colon \mc{G}(\Z_p)\to \Q$ is a locally constant compactly supported function.
\end{itemize}
The output should then be a locally constant and compactly supported function $\phi_{\tau,h}^{\mc{G},\mu}\colon G(K_j)\to \Q$.

In \cite{ScholzeLK}, Scholze defined the functions $\phi^{\mc{G},\mu}_{\tau,h}$ when $(\mc{G},\mu)$ is of so-called EL- or PEL-type $A$ or $C$. Namely, for $b$ in $\mc{G}(W_j)\sigma(\mu)(p)\mc{G}(W_j)$ one attaches a $p$-divisible group with EL- or PEL-structure $\bb{X}_b$ (see \cite[Definition 3.3]{ScholzeLK}) over $k_j$. If $\bb{X}_b^\univ\to \mf{D}(\mc{G},b,\mu)$ is the universal deformation of this object, then $T_p(\bb{X}_b^\univ)$ produces a $\mc{G}(\Z_p)$-local system over $\mc{D}(\mc{G},b,\mu)\defeq \mf{D}(\mc{G},b,\mu)_\eta$. The frame space $\mc{D}_\infty(\mc{G},b,\mu)$ for $T_p(\bb{X}_b^\univ)$ has a smooth $\mc{G}(\Z_p)$-action by modifying the level, thus giving an action of $\tau\times h$ on the $\ell$-adic cohomology of $\mc{D}_\infty(\mc{G},b,\mu)_C$ where $C=\wh{\ov{K}}$. Scholze then defines:
\begin{equation*}
\phi_{\tau,h}^{\mc{G},\mu}(b)\defeq \tr\left(\tau\times h\mid H^\ast_\et(\mc{D}_\infty(\mc{G},b,\mu)_C,\Q_\ell)\right),
\end{equation*}
and shows that it's locally constant, compactly supported, $\Q$-valued, and independent of $\ell\ne p$.

\medskip

\paragraph*{The goal of this article} Scholze's construction of local test functions naturally extends to the Hodge-type setting; see \cite{Youcis}. But beyond the Hodge-type case, the lack of a good notion of `$p$-divisible groups with $\mc{G}$-structure' has posed a serious obstruction.

 In this article we use recent advances in integral $p$-adic Hodge theory, namely the theory of $(\mc{G},\mu)$-apertures as in \cite{DrinfeldShimurian,GMM}, to overcome this difficulty, establish an (essentially fully general) version of the Langlands--Kottwitz--Scholze method, and formulate a more robust version of the Scholze--Shin conjectures using the recent work of Fargues--Scholze.

In the rest of the introduction, we explain some of these points more precisely. For now though, we would like to state our version of the Langlands--Kottwitz--Scholze method for Shimura varieties of abelian-type, where it is most robust.

\begin{thmi}[{see Corollary \ref{cor:Langlands-Kottwitz-Scholze-abelian-type}}]\label{thmi:abelian-type-trace-formula}
Suppose that $(\mb{G},\mb{X})$ is of abelian type. Then there exists an integer $j_0\geqslant 1$, depending only on $f^p$, such that for every $j\geqslant j_0$ and $\tau$ projecting to $\mr{Frob}_k^j$:
\begin{equation*}
T(\tau,h,f^p)=
\sum_{\gamma_0\in\Sigma_{\mathsf{K}^p}}\sum_{\substack{c\in\mf{KP}(\gamma_0)\cap\mf{KP}_a(p^n)\\ \alpha(c)=0}}
\ov{\iota}_{\mb{G}}(\gamma_0)^{-1}
 c_2(\gamma_0)\tr(\xi(\gamma_0))
 c_1(c,\mathsf{K}^p)
 \mr{O}_\gamma(f^p)
 \mr{TO}_\delta\left(\phi_{\tau,h}^{\mc{G},-\mu_h}\right),
\end{equation*}
where $n=rj$. If $\mb{G}^{\mr{ad}}$ is $\Q$-anisotropic, then one may take $j_0=1$.
\end{thmi}

\begin{remi}[Some previous work] There is considerable work done concerning local test functions and the cohomology of Shimura varieties at bad level both inside and outside the Langlands--Kottwitz--Scholze method. We give a brief list of them below:
\begin{itemize}[leftmargin=.5cm]
\item In \cite{HainesDrinfeld,HaistBC}, Haines gives explicit test functions in the Drinfeld case and develops the stable Bernstein-center formulation of the Haines--Kottwitz conjecture for general level.
\item  In Haines--Ng\^o \cite{HainesNgo}, Haines--Richarz \cite{HRParahoric,HRWeil}, and Ansch\"utz--Gleason--Louren\c{c}o--Richarz \cite{AGLR} nearby-cycle descriptions of test functions are
given at parahoric level.
\item In \cite{HainesRapoport}, Haines--Rapoport explicitly determine the local test functions for $\Gamma_1(p)$-level in the Drinfeld case.
\item In \cite{ChiHaines}, Chi--Haines establish local test-function identities and cohomological applications in certain non-quasi-split EL settings.
\item In \cite{ShenUniformized,ShenBad}, Shen obtains related identities and cohomological descriptions for certain quaternionic
and unitary Shimura varieties, including $p$-adically uniformized cases.
\end{itemize}
In particular, it would be quite interesting to understand whether the methods of \cite{ChiHaines} can be extended to work in the more general situations considered here.
\end{remi}

\subsection{Deformation spaces of $(\mc{G},\mu)$-apertures} 

We now more precisely discuss the deformation spaces of $(\mc{G},\mu)$-apertures used in the formulation of our general local test functions. In particular, we give a very concrete description of the universal deformation over these spaces. In the following, we maintain the meaning of $(\mc{G},\mu)$ and $j$ from above, and write $\mr{BT}^{\mc{G},\mu}_\infty$ for the formal Artin stack over $W$ parameterizing $(\mc{G},\mu)$-apertures; see \cite{GMM}.

The mechanism for our simpler description of the universal deformation of a $(\mc{G},\mu)$-aperture is a general comparison of $(\mc{G},\mu)$-apertures to certain crystalline data over reasonable bases. Namely, let $\mf{X}$ be a base formal $W$-scheme as in Remark \ref{rem:base-algebras}; for example $\mf{X}$ is smooth over $W$ or $\Spf(W\ll t_1,\ldots,t_d\rr)$ for some $d$. We then define $\cat{FFCrys}^{\mr{sd},\mu}_\mc{G}(\mf{X})$ to be the category of strongly divisible filtered $F$-crystals with $\mc{G}$-structure which are of type $\mu$; see \S\ref{ss:filtered-F-crystals}. In essence, one can think of this as a Tannakian version of the category of Fontaine--Laffaille modules.

\begin{thmi}[{see Theorem \ref{thm:BT-over-base-ring}}]\label{thmi:equiv} Suppose $p>2$ and that $\mf{X}$ is a base formal $W$-scheme. Then, there is an equivalence of categories
\begin{equation*}
\bb{D}_\mr{crys}\colon \mr{BT}^{\mc{G},\mu}_\infty(\mf{X})\to \cat{FFCrys}^{\mr{sd},\mu}_\mc{G}(\mf{X}).
\end{equation*}
\end{thmi}

\begin{remi} One may see this as a Tannakian version of \cite[Theorem A]{IKY3} or of \cite{TVX} in the case of a point. That said, interestingly enough, as we vary the representation, we are considering here $F$-gauges of arbitrarily high Hodge--Tate weights. In some sense, the minisculeness condition on $\mu$ shows that in aggregate these $F$-gauges act as if they have Hodge--Tate weights in the Fontaine--Laffaille range.
\end{remi}

\begin{remi} Using this, one can give an alternative construction of the syntomic realization functors from \cite{IKY2} and \cite{MY}, which is more in line with the methods used in this article, e.g., doesn't use any hard facts from \cite{GuoReinecke}; see \S\ref{sss:alternative-proofs}. 
\end{remi}

Using Theorem \ref{thm:BT-over-base-ring}, one can give a very concrete description of universal deformations of $(\mc{G},\mu)$-apertures by describing the output of $\bb{D}_\crys$ applied to their universal deformation. 

By Proposition \ref{prop:isom-classes-in-BT} an isomorphism class of a $(\mc{G},\mu)$-aperture over $k_j$ is equivalent to the data of an element  $b$ of $\mc{G}(W_j)\sigma(\mu)(p)\mc{G}(W_j)$. We then further write:
\begin{itemize}
\item $\mc{U}_{\mu,j}\subseteq \mc{G}_{W_j}$ for the unipotent defined by $\mu$;
\item $\mc{U}^{-}_{\mu,j}$ for its opposite unipotent;
\item $\wh{\mc{U}}^-_{\mu,j}$ for the completion of $\mc{U}^{-}_{\mu,j}$ at the origin; 
\item and $R_{\mc{G},\mu,j}$ for the ring of global sections of $\wh{\mc{U}}^-_{\mu,j}$. 
\end{itemize}
Finally, for reasons that will be soon clear we alternatively write $\mf{D}(\mc{G},b,\mu)=\wh{\mc{U}}^{-}_{\mu,j}$.

Then, \cite[Theorem 4.4.2]{ItoDeformation}, \cite[Proposition 3.32]{IKY2}, and \cite[Proposition 10.2.9]{GMM} explain that there is a $(\mc{G},\mu)$-aperture $\mf{Q}_b^\univ$ on $\mf{D}(\mc{G},b,\mu)$ which is the universal deformation of $b$. We would then like to describe the equivalent piece of data $\bb{D}_\mr{crys}(\mf{Q}_b^\univ)$. In fact, we would like to describe an even more concrete equivalent piece of data.

Fix an isomorphism $R_{\mc{G},\mu,j}\simeq W_j\ll t_1,\ldots,t_d\rr$, which induces a Frobenius lift $\sigma$ on $R_{\mc{G},\mu,j}$ by transporting that on $W_j\ll t_1,\ldots,t_d\rr$ which is the usual Frobenius on $W_j$ and raising the coordinates to the $p^\text{th}$-power. By Proposition \ref{prop:torsors-with-connection-equiv-filtered}, the data of $\bb{D}_\crys(\mf{Q}_b^\univ)$ is equivalent to  a quadruple of the form $(Q_b^\univ,\varphi_b^\univ,\nabla_b^\univ,\Fil^\bullet Q_b^\univ)$ where:
\begin{itemize} 
\item $Q_b^\univ$ is a $\mc{G}$-torsor on $R_{\mc{G},\mu,j}$; 
\item $\nabla_b^\univ$ is an {integrable and topologically quasi-nilpotent} $\mc{G}$-connection; see \S\ref{sss:G-conn};
\item  $\varphi_b^\univ\colon \sigma^\ast Q_b^\univ[\nicefrac{1}{p}]\isomto Q_b^\univ[\nicefrac{1}{p}]$ is $\nabla_b^\univ$-horizontal and type $\sigma(\mu)$; see Proposition \ref{prop:type-mu-equiv-char-p};
\item and $\Fil^\bullet Q_b^\univ$ is a filtration on $Q_b^\univ$ (as in \S\ref{sec:filtered-Tannakian-appendix}) satisfying Griffiths transversality in the sense of Definition \ref{defn:Griffiths}, and for which $\varphi_b^\univ$ is strongly divisible.
\end{itemize}
Let us then write:
\begin{itemize}
\item $P_\mu\subseteq \mc{G}_{W_j}$ for the parabolic associated to $\mu$;
\item $U_{\mu,j}^{\univ}$ for the tautological element of $\wh{\mc{U}}^-_{\mu,j}(R_{\mc{G},\mu,j})$;
\item $\eta_b^\univ$ for the Maurer--Cartan element $(U_{\mu,j}^\univ)^{-1}dU_{\mu,j}^\univ$; see Notation \ref{nota:Maurer-Cartan};
\item $\Delta_b$ for the operator on $\mf{g}\otimes_{\Z_p}\wh{\Omega}^1_{R_{\mc{G},\mu,j}/W_j}$ given by $\Delta_b\left(x\otimes \eta\right)=\mr{Ad}\left(({U}_{\mu,j}^\univ)^{-1} b\right)(x)\otimes \,\sigma^\ast(\eta)$.
\end{itemize}

\begin{thmi}[{see Theorem \ref{thm:explicit-deformation}}]\label{thmi:explicit-deformation}  Write $b=k_1\sigma(\mu)(p)$ for $k_1$ in $\mc{G}(\Z_p)$,\footnote{This restricted setting is for simplicity here in the intro, and is not assumed in the main body of the paper.}. Then, 
\begin{equation*}
\mc{Q}_b^\mr{univ}= \mc{G}_{R_{\mc{G},\mu,j}},\quad \varphi_b^\univ=({U}^\univ_{\mu,j})^{-1} b,\quad \Fil^\bullet \mc{Q}_b^\mr{univ}=P_\mu\in (\mc{G}/P_\mu)(R_{\mc{G},\mu,j}),\quad \nabla_b^\mr{univ}
=\sum_{n\geqslant 0}\Delta_b^n(\eta_b^\univ).
\end{equation*}
\end{thmi}

\begin{remi} Versions of Theorem \ref{thmi:explicit-deformation} are relatively well-known in Hodge-type situations, e.g., see \cite[Theorem 3.6]{KimRZ}. What is actually new here is the description beyond the Hodge-type setting (although compare with \cite{ItoDeformation}), and an explicit description of the universal connection (even in the Hodge-type setting).
\end{remi}

\subsection{The local test functions $\phi_{\tau,h}^{\mc{G},\mu}$} We now give a precise description of the local test functions $\phi_{\tau,h}^{\mc{G},\mu}$ studied in this article. Let the notation concerning $\mc{G},\mu,\tau,h$, and $j$ be as above.

To formulate the definition of $\phi_{\tau,h}^{\mc{G},\mu}$, let us write:
\begin{itemize}
\item $\mc{D}(\mc{G},b,\mu)$ for the Berthelot generic fiber of $\mf{D}(\mc{G},b,\mu)$.\footnote{This is non-canonically an open polydisk over $K_j$.}
\end{itemize}
The universal deformation $\mf{Q}_b^\univ$ over $\mf{D}(\mc{G},b,\mu)$ gives rise to a $\mc{G}(\Z_p)$-local system on $\mc{D}(\mc{G},b,\mu)$ by considering \'etale realization; see Notation \ref{nota:etale-realization}. We then write
 \begin{itemize}
 \item $\mc{D}_\infty(\mc{G},b,\mu)\to \mc{D}(\mc{G},b,\mu)$ for the frame space for this $\mc{G}(\Z_p)$-local system. 
\end{itemize}
Modifying the $\mc{G}(\Z_p)$-level structure gives a $\mc{G}(\Z_p)$-action on $\mc{D}_\infty(\mc{G},b,\mu)$. Using this, we see that if $C$ is a completed algebraic closure of $K_j$ we have an admissible action of $W_{K_j}\times \mc{G}(\Z_p)$ on $H^i_\et(\mc{D}_\infty(\mc{G},b,\mu)_C,\Q_\ell)$ for any prime $\ell\ne p$.  

\begin{defni} We define the \emph{local test function} associated to $(\mc{G},\mu,\tau,h)$ to be:
\begin{equation*}
\phi_{\tau,h}^{\mc{G},\mu}\colon G(K_j)\to \C,\qquad b\mapsto \begin{cases} \tr\left(\tau\times h| H_\et^\ast(\mc{D}_\infty(\mc{G},b,\mu)_C,\Q_\ell)\right) & \mbox{if}\quad b\in \mc{G}(W_j)\sigma(\mu)(p)\mc{G}(W_j)\\ 0 & \mbox{if}\quad \text{otherwise.}\end{cases}
\end{equation*}
\end{defni}

To actually be useful in applications, we need to know that $\phi^{\mc{G},\mu}_{\tau,h}$ has reasonable harmonic-analytic properties. This is one of the main results of this paper.

\begin{thmi}[{see Theorem \ref{thm:phi-nice-properties} }]\label{thmi:test-functions} The function $\phi_{\tau,h}^{\mc{G},\mu}$ is locally constant, compactly supported, takes values in $\Q$, and is independent of $\ell\ne p$.
\end{thmi}

The key to proving the most difficult part of Theorem \ref{thmi:test-functions}, the local constancy, is a fine study of the action of $J_b^\mr{int}=\mr{Aut}(b)$ on the finite-level spaces $\mc{D}_\mathsf{K}(\mc{G},b,\mu)\defeq \mc{D}_\infty(\mc{G},b,\mu)/\mathsf{K}$ where $\mathsf{K}\subseteq\mc{G}(\Z_p)$ is a compact open subgroup. Namely, in Proposition \ref{prop:uniform-continuity-K-level}, we establish that this action is uniformly continuous in an appropriate sense. We also make use of some, potentially independently interesting, independence of $\ell$-results as described in Appendix \ref{s:l-independence-appendix}.

\begin{remi} While the function $\phi^{\mc{G},\mu}_{\tau,h}$ is difficult to compute in general, it is feasible in the simplest possible case: When the cocharacter $\mu$ is central. This is particularly simple when $\mc{G}=\mc{T}$ is a torus. Namely, if $b$ belongs to $\mc{T}(W_j)\sigma(\mu)(p)\mc{T}(W_j)$ then 
\begin{equation*}
\phi_{\tau,h}^{\mc{T},\mu}(b)=h\left(r_\mu(p^{-j}\mr{Art}_\mathsf{K}(\tau))\cdot \mr{Nr}_{K_j/\Q_p}(b\sigma(\mu)(p)^{-1})\right),
\end{equation*}
where $r_\mu$ is as in Construction \ref{constr:Kottwitz}, $\mr{Art}_\mathsf{K}$ is the (appropriately normalized) Artin reciprocity map, and $\mr{Nr}_{K_j/\Q_p}\colon \mc{T}(K_j)\to\mc{T}(\Q_p)$ is the norm map. See Proposition \ref{prop:central-character-calc}  for the general central-character calculation, and compare with \cite[Proposition 4.10]{ScholzeLK} which contains the special case when $\mc{G}=\bb{G}_{m,\Z_p}$ and $\mu$ is the identity map.
\end{remi}

\begin{remi} As mentioned before, there have previously been cases of the local test functions $\phi_{\tau,h}^{\mc{G},\mu}$ studied: in the EL- or PEL-type $A$ or $C$ case in \cite{ScholzeLK}, and in the Hodge-type (and some abelian-type) cases in \cite{Youcis}. In \S\ref{ss:comparison}, we give a precise comparison between our local test functions and the ones from these articles.
\end{remi}
\begin{remi} In the type $A$ or $C$ cases studied in \cite{ScholzeLK}, the local test functions can be described using infinite-level tubes inside of Rapoport--Zink space. In \S\ref{ss:shtuka-description-test-functions} we discuss the generalization of this observation, showing that $\phi_{\tau,h}^{\mc{G},\mu}$ can be understood as computing the traces of the cohomology of infinite-levels tubes inside of moduli spaces of integral shtukas for $\mc{G}$. As discussed in Remark \ref{rem:parahoric-extension} this gives a natural guess for how to define $\phi^{\mc{G},\mu}_{\tau,h}$ when $\mc{G}$ is a parahoric group scheme. That said, this is not so useful as the finer properties of the local test functions in the reductive case (e.g., Theorem \ref{thmi:test-functions}) require serious deformation theory unavailable in the setting of shtukas.
\end{remi}

\subsection{The primitive trace formula} Now that we have described local test functions in full generality, it is natural to ask whether there exists an equally general version of the Langlands--Kottwitz--Scholze method, giving a trace formula for the cohomology of Shimura varieties at bad (specifically subhyperspecial) level. We now discuss the first step in the realization of this idea. We keep our notation concerning Shimura varieties, and $\mc{G},\mu=-\mu_h,\tau,h$, and $j$ as above.

As the local test functions are defined in terms of $(\mc{G},\mu)$-apertures, we make use of the notion of \emph{integral canonical models} from \cite{IKY2,MY}. Namely, suppose that $\mathsf{K}=\mathsf{K}_0\mathsf{K}^p$ is a neat level with $\mathsf{K}_0=\mc{G}(\Z_p)$. Then, as in \S\ref{sss:etale-realization}, we can build a $\mc{G}(\Z_p)$-local system 
\begin{equation*}
\mr{Et}_{\mathsf{K},p}^\circ=\varprojlim_{K_p\subseteq \mathsf{K}_0}\mr{Sh}_{\mathsf{K}_p\mathsf{K}^p}(\mb{G},\mb{X})_E\to \mr{Sh}_\mathsf{K}(\mb{G},\mb{X})_E.
\end{equation*}
In \cite{MY} there is constructed a (pro-)algebraization $\mr{BT}^{\mc{G},-\mu_h,\mr{alg}}_\infty$ of $\mr{BT}^{\mc{G},-\mu_h}_\infty$ and for any $W$-scheme $X$, an \'etale realization functor $T_\et\colon \mr{BT}^{\mc{G},-\mu_h,\mr{alg}}_\infty(X)\to B\mc{G}(\Z_p)(X[\nicefrac{1}{p}])$. Moreover, it is shown in \cite[Theorem 4.1.5]{MY} that this \'etale realization functor is fully faithful if $X$ is normal.

Then, as in \cite[Definition 3.39]{IKY2} and \cite[Definition 6.5.2]{MY}, a separated $W$-model $\ms{S}_{\mathsf{K}}(\mb{G},\mb{X})$ of $\Sh_\mathsf{K}(\mb{G},\mb{X})_E$ is an \emph{integral canonical model} if:
\begin{enumerate}
\item The $\mc{G}(\bb{Z}_p)$-local system $\mr{Et}^\circ_{\mathsf{K},p}$ is the \'etale realization of a $(\mc{G},-\mu_h)$-aperture on $\ms{S}_{\mathsf{K}}(\mb{G},\mb{X})$ such that the corresponding map $\ms{S}_\mathsf{K}(\mb{G},\mb{X})\to\mr{BT}^{\mc{G},-\mu_h,\mr{alg}}_\infty$ is $p$-adically formally \'etale.
\item One has an equality of classical points
\begin{equation*}
|(\wh{\ms{S}}_{\mathsf{K}})_\eta|^\mr{cl}=\left\{x\in |\Sh_\mathsf{K}^\mr{an}|^\mr{cl}:\mb{Et}_{\mathsf{K},p}^\circ\text{ is potentially crystalline at }x\right\}.
\end{equation*}
\end{enumerate}
An integral canonical model $\ms{S}_\mathsf{K}$ is furthermore called \emph{limpid} if $\mr{Et}_{\mathsf{K},\ell}$ extends to $\ms{S}_{\mathsf{K}}$ for all $\ell\ne p$.

\begin{remi} Integral canonical models, if they exist, are unique; see Proposition \ref{prop:mapping-prop}. The notion of integral canonical models characterizes good integral models by conditions generalizing the Serre--Tate deformation theorem and the Coleman--Iovita theorem from the Siegel-type case. Limpid integral canonical models are known to exist for essentially all possible cases, and agree with Kisin's integral models in the abelian-type case; see Theorem \ref{thm:ICMs-exist} and Remark \ref{rem:ICMs-exist}.
\end{remi}

Given the fact that integral canonical models are designed to be compatible with $(\mc{G},-\mu_h)$-apertures, it is then perhaps not surprising that they are the proper setting to consider for generalizations of the Langlands--Kottwitz--Scholze formula. And, in fact, with no extra assumptions, they admit a weak form of the formula from the Langlands--Kottwitz--Scholze method, which we call the \emph{primitive trace formula}.

\begin{thmi}[{see Theorem \ref{thm:prim-trace-formula-II}}]\label{thmi:primitive-trace-formula} Let $g=g_pg^p\in\mathsf{K}_0\mathsf{K}^p$ and $\mathsf{K}'=\mathsf{K}_p\mathsf{K}^p$ with $\mathsf{K}_p\subseteq \mathsf{K}_0$. Then, assuming that $\mathscr{S}_\mathsf{K}(\mb{G},\mb{X})$ is a proper\footnote{This is not necessary; see the precise statement of Theorem \ref{thm:prim-trace-formula-II}.} limpid integral canonical model of $\Sh_\mathsf{K}(\mb{G},\mb{X})_E$, one has
 \begin{equation*}
\tr\left(\tau\times e(g^{-1},\mathsf{K}')\mid {H}^\ast_\et(\Sh_{\mathsf{K}'}(\mb{G},\mb{X})_{\ov{E}},\mc{F}_{\xi,\mathsf{K}}\right)=\sum_{y\in\mr{Fix}(c^{(j)})(\ov{k})}\tr\left((u^{(j)}_{c_2({y})})^\dashv\mid ({\mc{F}}_{\xi,\mathsf{K}^p})_{c_2({y})}\right) \phi_{\tau,e((g_p)^{-1},\mathsf{K}_p)}^{\mc{G},-\mu_h}(\delta(y)).
\end{equation*}
\end{thmi}

Let us explain the precise meaning of the notation used above:
\begin{itemize}
\item $e(h,\mathsf{L})=\tfrac{1}{\mr{vol}(\mathsf{L})}\mathds{1}_{\mathsf{L}h\mathsf{L}}$ in all usages,
\item $c=c(g^p,\mathsf{K})$ is the integral Hecke correspondence on $\ms{S}_{\mathsf{K}}(\mb{G},\mb{X})$; see Definition \ref{defn:integral-Hecke-corr},
\item $u=u(g^p,\mathsf{K})$ is the natural $c$-cohomological Hecke correspondence on $\mc{F}_{\xi,\mathsf{K}}$; see loc.\@ cit.\@,
\item $c^{(j)}$ is the $j^\text{th}$-Frobenius twist of $c$; see Notation \ref{nota:Frobenius-twist-corr},
\item $\mr{Fix}(c^{(j)})(\ov{k})$ is the fixed-point set of $c^{(j)}$; see Definition \ref{defn:fixed-scheme},
\item $(u^{(j)}_{c_2({y})})^\dashv$ is the naive local term; see Definition \ref{defn:naive-local-term},
\item $\delta\colon \mr{Fix}(c^{(j)})(\ov{k})\to \mc{G}(W_j)\sigma(-\mu_h)(p)\mc{G}(W_j)$ is the map from Construction \ref{constr:corr-delta}.
\end{itemize}
In words, the primitive trace formula says that the trace on the left-hand side of the formula can be computed as a sum over fixed-points of a Frobenius-twisted integral Hecke correspondence on a limpid integral canonical model. The summand at each fixed point has two contributions: 
\medskip

\noindent\textbf{An away-from-$p$ contribution:} This is computed entirely in terms of $\mc{F}_{\xi,\mathsf{K}}$ and the away-from-$p$ parts $g^p$ and $\mathsf{K}^p$ of $g$ and $\mathsf{K}$, respectively.

\medskip 

\noindent\textbf{The at-$p$ contribution:} This is computed in terms of the local test function $\phi^{\mc{G},-\mu_h}_{\tau,h}(\delta(y))$ where $h=e((g_p)^{-1},\mathsf{K}_p)$ depends only on the at-$p$ parts $g_p$ and $\mathsf{K}_p$ of $g$ and $\mathsf{K}$, respectively, and $\delta(y)$ is a $(\mc{G},-\mu_h)$-aperture obtained from using the fact that $\ms{S}_{\mathsf{K}}(\mb{G},\mb{X})$ is an integral canonical model.

\medskip

The proof of the primitive trace formula combines a detailed study of integral canonical models, together with a general trace formula relating the Fujiwara--Varshavsky trace formula for correspondences on the special fiber of a formal scheme to traces of correspondences of certain covers of its generic fiber. This general rigid-analytic trace formula is discussed in \S\ref{s:prelims}, and generalizes previous work of Fargues in \cite{Fargues}.

\subsection{The Langlands--Kottwitz--Scholze method}

For global applications to the cohomology of Shimura varieties, discussed below, the primitive trace formula is not sufficient; we need a version of the Langlands--Kottwitz--Scholze method as discussed at the beginning of this introduction. We again maintain our notation concerning Shimura varieties, and $\mc{G},\mu=-\mu_h,\tau,h$, and $j$.

Passing from the primitive trace formula to the desired Langlands--Kottwitz--Scholze-type trace formula requires a detailed understanding of the fixed point sets $\mr{Fix}(c^{(j)})(\ov{k})$. This is exactly the purview of the \emph{Langlands--Rapoport conjecture}, which posits a decomposition
\begin{equation*}
\ms{S}_{\mathsf{K}_0}(\ov{k})=\bigsqcup_{[\varphi]} \varprojlim_{\mathsf{L}^p}
I_\varphi(\Q)\backslash\left(X^p(\varphi)\times X_p(\varphi)\right)/\mathsf{L}^p,\qquad \ms{S}_{\mathsf{K}_0}=\varprojlim_{\mathsf{L}^p\subseteq \mathsf{K}^p}\mathscr{S}_{\mathsf{K}_0\mathsf{L}^p}(\mb{G},\mb{X}),
\end{equation*}
which are equivariant for the prime-to-$p$ Hecke action and Frobenius. Morally, one should think:
\begin{itemize}
\item $[\varphi]$ is like a motive $M$ with $\mb{G}$-structure over $\ov{k}$;
 \item $X^p(\varphi)$ is parameterizing $\A_f^p$-lattices inside $M\otimes \A_f^p$;
\item $X_p(\varphi)$ is parameterizing $\Z_p$-lattices inside of $M\otimes \Q_p$;
\item $I_\varphi(\Q)$ is the actual automorphism group of $M$ that provides coherency between these away-from-$p$ and at-$p$ realizations.
\end{itemize}
The usual double-coset description of $\mr{Sh}_\mathsf{K}(\mb{G},\mb{X})(\C)$ should be viewed as a motivic decomposition, which can be made precise assuming the Hodge conjecture. The Langlands--Rapoport conjecture is a similar motivic decomposition over $\ov{\bb{F}}_p$ (sans the archimedean place which is not relevant).

While it was Langlands--Rapoport that had the original vision for this conjecture and its applications to the Langlands--Kottwitz method, it was Kisin--Shin--Zhu in \cite{KSZ} that gave this very strong rigorous foundations. Namely, they very precisely axiomatized not only the Langlands--Rapoport conjecture, but certain weaker \emph{twisted variants} called the \emph{Langlands--Rapoport-$\tau$-conjecture}. They moreover showed that a Langlands--Kottwitz-type trace formula follows from such a Langlands--Rapoport-$\tau$ conjecture.

We may then directly imitate their method and almost show that, under the assumption of a Langlands--Rapoport-$\tau$ conjecture one can upgrade the primitive trace formula to a Langlands--Kottwitz--Scholze-type trace formula. The `almost' here is because we require that this Langlands--Rapoport-$\tau$ conjecture is compatible with the map $\delta$ from the primitive trace formula, something we refer to as \emph{compatibility with syntomic realization}; see Definition \ref{defn:LR-realization-compatible}.

\begin{thmi}[{The Langlands--Kottwitz--Scholze formula, see Theorem \ref{thm:Langlands-Kottwitz-Scholze-formula}}]
Suppose that $\mathscr{S}_{\mathsf{K}}(\mb{G},\mb{X})$ is proper,\footnote{Again, this is not actually necessary; see Theorem \ref{thm:Langlands-Kottwitz-Scholze-formula} for the precise general assumptions.} and that there exists a Langlands--Rapoport-$\tau$ conjecture compatible with syntomic realization. Then, with $n=rj$, one has:
\begin{equation*}
T(\tau,h,f^p)=
\sum_{\gamma_0\in\Sigma_{\mathsf{K}^p}}\sum_{\substack{c\in\mf{KP}(\gamma_0)\cap\mf{KP}_a(p^n)\\ \alpha(c)=0}}
\ov{\iota}_{\mb{G}}(\gamma_0)^{-1}
 c_2(\gamma_0)\tr(\xi(\gamma_0))
 c_1(c,\mathsf{K}^p)
 \mr{O}_\gamma(f^p)
 \mr{TO}_\delta\left(\phi_{\tau,h}^{\mc{G},-\mu_h}\right),
\end{equation*}
where for each $c$ in the inner sum, $(\gamma_0,\gamma,\delta)$ denotes any representative of its associated classical Kottwitz parameter.
\end{thmi}

\begin{remi} When $(\mb{G},\mb{X})$ is of abelian type, Kisin--Shin--Zhu verify that the Langlands--Rapoport-$\tau$ conjecture holds. Using the construction of the syntomic realization in \cite{IKY2} one can show that Kisin--Shin--Zhu's Langlands--Rapoport-$\tau$ bijections are compatible with syntomic realization. Thus, one arrives at Theorem \ref{thmi:abelian-type-trace-formula}, at least after relaxing the unnecessary properness assumptions.
\end{remi}

Finally, we mention that imitating \cite[\S8]{KSZ}, which provides an (Arthur--Selberg) stabilization of a Langlands--Kottwitz-type trace formula, we can deduce a stabilization of such Langlands--Kottwitz--Scholze-type trace formulae (see the main-body theorem for notation):

\begin{thmi}[{see Theorem \ref{thm:geometric-stabilization}}]\label{thmi:stabilization}
Suppose that $\mathscr{S}_{\mathsf{K}}(\mb{G},\mb{X})$ is proper,\footnote{Again, this is not actually necessary; see Theorem \ref{thm:geometric-stabilization} for the precise general assumptions.} and that there exists a Langlands--Rapoport-$\tau$ conjecture compatible with syntomic realization. Then,
\begin{equation}\label{eq:geometric-stabilization-intro}
T(\tau,h,f^p)=
\sum_{\mf{e}\in\mc{E}_{\mr{ell}}(\mb{G})}
\iota(\mb{G},\mf{e})
\mr{ST}_{\mr{ell},\Omega_{\mf{e}}}^{\mb{H}_1}
(f_{\tau,h,\xi}^{\mf{e}}).
\end{equation}
\end{thmi}

\subsection{The Scholze--Shin conjecture and applications to Shimura varieties}\label{ss:Scholze-Shin-intro} Finally, we explain that our general local test functions, together with recent deep work of Fargues--Scholze and Hansen, allows one to give a more robust and precise version of the Scholze--Shin conjecture. We fix our notation concerning $\mc{G},\mu,\tau,h$ and $j$ as above. We furthermore write $G=\mc{G}_{\Q_p}$.

To begin, recall that in \cite{FarguesScholze} there is constructed a $\C$-algebra map
\begin{equation*}
\Psi_G^\mr{FS}\colon \mc{Z}^\mr{spec}(G)\to \mc{Z}(G),
\end{equation*} 
where
\begin{itemize}
\item $\mc{Z}^\mr{spec}(G)$ is the spectral Bernstein center (aka the stable Bernstein center); see \S\ref{sss:spectral-Bernstein};
\item $\mc{Z}(G)$ is the Bernstein center; see \S\ref{ss:Bernstein-center}.
\end{itemize}
One may think of $\mc{Z}^\mr{spec}(G)$ and $\mc{Z}(G)$ as the \emph{global sections} of the \emph{coarse} moduli spaces of $L$-parameters and smooth irreducible representations of $G(\Q_p)$, respectively. Moreover, a homomorphism like $\Psi_G^\mr{FS}$ can be seen as an incarnation of a hypothetical local Langlands correspondence for $G$; see \cite[Proposition 4.23]{HaistBC}. In fact, from $\Psi_G^\mr{FS}$ one may build an association of a semisimple $L$-parameter $\varphi_\pi^\mr{FS}\colon W_{\Q_p}\to {}^L G$ for any smooth irreducible representation $\pi$ of $G(\Q_p)$.

Now, from the data $\mu$ and $\tau$ one can cook up a natural element $z_{\mu,\tau}^\mr{spec}$ of $\mc{Z}^\mr{spec}(G)$. Namely, one may think of an element of the coarse moduli space of $L$-parameters as a semisimple $L$-parameter $\varphi\colon W_{\Q_p}\to  {}^L G$ and one then defines 
\begin{equation*}
z_{\mu,\tau}^\mr{spec}(\varphi)=\tr \left(\varphi(\tau)\mid r_{\mu}\right) q^{j\langle\rho_G,\mu\rangle},
\end{equation*}
where
\begin{itemize}
\item $r_{\mu}\colon \wh{G}\rtimes W_K\to\GL(V_\mu)$ is the representation as in Construction \ref{constr:Kottwitz};
\item $q$ is the size of $k$;
\item $\langle \rho_G,\mu\rangle$ is the pairing of $\mu$ against the half-sum of a set of positive roots.
\end{itemize}
The Scholze--Shin conjecture can then be seen as describing an explicit description of the images of $z_{\mu,\tau}^\mr{spec}$ under $\Psi_G^\mr{FS}$, which we denote $z_{\mu,\tau}=\Psi_G^\mr{FS}(z_{\mu,\tau}^\mr{spec})$.

Now, $z_{\mu,\tau}$ being an element of the Bernstein center $\mc{Z}(G)$ is not an element of the Hecke algebra of $G(\Q_p)$, but has the property that it can be convolved against any such Hecke operator $h$ to produce a Hecke operator $z_{\mu,\tau}\ast h$. Moreover, as the Langlands conjecture only concerns the stable parts of distributions, we will only be interested in such functions up to stable equivalence $=_\mr{st}$; i.e., where we disregard functions with vanishing stable orbital integrals.

\begin{conji}[{see Conjecture \ref{conj:Scholze-Shin}}] One has $z_{\mu,\tau}\ast h=_\mr{st} f_{\tau,h}^{\mc{G},\mu}$, where $f_{\tau,h}^{\mc{G},\mu}$ is a twisted transfer of $\phi_{\tau,h}^{\mc{G},\mu}$ to a function on $G(\Q_p)$.
\end{conji}

\begin{remi} Implicit in this statement being well-posed, we are using the result from \cite{HansenStableBernstein} which implies that $\Psi_G^\mr{FS}$ takes values in the very-stable elements of the Bernstein center; see Definition \ref{defn:stable-very-stable-center}. Hansen's result also allows us to state a more familiar implication of the Scholze--Shin conjecture. Let $\Theta=\sum_i a_i \Theta_{\pi_i}$ be an atomically stable virtual character of $G(\Q_p)$ (see Definition \ref{defn:atomically-stable-character}). Then, if the Scholze--Shin conjecture \ref{conj:Scholze-Shin} holds, we have the following identity
\begin{equation*}
\Theta\left(f_{\tau,h}\right)= z_{\mu,\tau}^{\mr{spec}}\left(\varphi_\Theta^{\mr{FS}}\right)
\Theta\left(h\right).
\end{equation*}
Here $\varphi_\Theta^\mr{FS}=\varphi_{\pi_i}^\mr{FS}$ which is independent of the choice of $\pi_i$ by \cite{HansenStableBernstein}; see Proposition \ref{prop:FS-atomically-stable-character}.
\end{remi}

Finally following Scholze--Shin in \cite{ScholzeShin}, and Haines in \cite{HaistBC}, we are able to apply the Scholze--Shin conjecture, in tandem with the stabilized version of the Langlands--Kottwitz--Scholze formula from Theorem \ref{thmi:stabilization}, to obtain concrete consequences about the cohomology of Shimura varieties.  

In the following, we use the following notation:

\begin{itemize}
\item $d=\dim_\C\mb{X}$;
\item $(-)^\mr{ss}$ is the semisimplification of Weil representations from Notation \ref{nota:ss}; 
\item Let $\pi_p$ be an irreducible smooth $G(\Q_p)$-representation. Write
\begin{equation*}
R_{\pi_p}\defeq \left(r_{-\mu_h}\circ
\varphi_{\pi_p}^{\mr{FS}}|_{W_E}\right)^{\mr{ss}}
\otimes|\cdot|_E^{-\tfrac{d}{2}},
\end{equation*}
with $|\tau|_E\defeq |\mr{Art}_E(\tau)|_E$ normalized so that geometric Frobenius has absolute value $q^{-1}$;
\item $a_\xi(\pi_f)$ is as in Notation \ref{nota:Haines-coefficients}.
\end{itemize}

\begin{propi}[{see Proposition \ref{prop:cohomology-decomp-no-end}}]\label{propi:coh-decomp}
Assume that $(\mb{G},\mb{X})$ is a Shimura datum where $\mb{G}$ is $\Q$-anisotropic and with no endoscopy \emph{(}in the sense of \cite[Theorem 6.3.2]{HaistBC}\emph{)}, \eqref{eq:geometric-stabilization-intro} holds, and the Scholze--Shin conjecture holds for all $\tau$ and $h$. Then, there is a decomposition of virtual
$\mb{G}(\A_f^p)\times\mathsf{K}_0\times W_E$-representations
\begin{equation*}
H_c^*(\mb{G},\mb{X},\xi)^\mr{ss}
=\sum_{\pi_f=\pi^p\otimes\pi_p}a_\xi(\pi_f)
\bigg[\pi^p\boxtimes\pi_p|_{\mathsf{K}_0}\boxtimes R_{\pi_p}\bigg],
\end{equation*}
where $\pi_f$ ranges over irreducible admissible $\mb{G}(\A_f)$-representations.
\end{propi}

For a Weil representation $V$, we denote by $\zeta_S^\mr{ss}(V)$ the partial semisimple Hasse--Weil $\zeta$-functions in the sense of Definition \ref{nota:ss-local-factor}. Then, as in \cite{HaistBC} the above decomposition has concrete implications for the partial semisimple Hasse--Weil $\zeta$-functions of Shimura varieties.

\begin{defni} Fix a neat level $\mathsf{K}$ and write $S(\mathsf{K})$ for the set of primes $v$ of $\mb{E}$ not dividing $\ell$ and lying over a prime $p$ where $\mathsf{K}_p$ is contained in a hyperspecial subgroup of $\mb{G}(\Q_p)$. We then write
\begin{equation*}
\zeta^\mr{ss}_{S(\mathsf{K})}(\mb{G},\mb{X},\mathsf{K},\xi,s)\defeq \zeta_{S(\mathsf{K})}^\mr{ss}\left(H^\ast_c(\mb{G},\mb{X},\mathsf{K},\xi),s\right).
\end{equation*}
\end{defni}

\begin{cori}Assume the setup of \emph{Proposition \ref{propi:coh-decomp}} holds for every place $v$ in $S(\mathsf{K})$, and write $R_{\pi_p,v}$ for $R_{\pi_p}$ formed with $E=\mb{E}_v$. Then, we have an equality
\begin{equation*}
\zeta^\mr{ss}_{S(\mathsf{K})}(\mb{G},\mb{X},\mathsf{K},\xi,s)=\prod_{v\in S(\mathsf{K})}\prod_{\pi_f=\pi^p\otimes\pi_p}\det\left(1-q_v^{-s}\Phi_v\mid R_{\pi_p,v}^{I_v}\right)^{
-a_{\xi}(\pi_f)\dim\pi_f^{\mathsf{K}}}.
\end{equation*}
\end{cori}

\subsection*{Notation and Conventions}
\label{conv:notation-and-conventions}
In the entirety of this paper we make use of the following notations and conventions (with more specific notations/conventions defined locally below):
\begin{itemize}[leftmargin=.5cm,itemsep=2pt,topsep=3pt]
\item For a functor of groupoids $F\colon\ms{G}_1\to\ms{G}_2$ and an
object $x$ of $\ms{G}_2$, the fiber
$\mr{fib}(\ms{G}_1\to\ms{G}_2;x)$ has objects
$(y,\iota\colon F(y)\isomto x)$, with $y$ an object of $\ms{G}_1$.
A morphism $(y,\iota)\to(y',\iota')$ is an isomorphism
$a\colon y\isomto y'$ satisfying $\iota'\circ F(a)=\iota$.
\item For a groupoid $\ms{G}$, we denote by $\pi_0(\ms{G})$ the set of isomorphism classes in $\ms{G}$.
\item Closed subsets of (formal) schemes carry their reduced scheme
structures.
\item A non-archimedean field $K$ is complete for a non-trivial rank $1$
absolute value $|\cdot|$ inducing its topology. We write
$\mc{O}_K=\{x\in K:|x|\leqslant1\}$ for its valuation ring.
\item We use the terminology formally of finite type as in \cite[Definition 2.2]{ALY2}.
\item For a Huber ring $A$, write $A^\circ$ and $A^{\circ\circ}$ for its
power-bounded and topologically nilpotent elements, respectively.
\item An adic ring is a topological ring $R$ whose topology is induced by the powers of some ideal $J\subseteq R$, which we call an ideal of definition. We will write $\Spf(R,J)$ for the formal spectrum of this topological ring, or $\Spf(R)$ when the topology/ideal is implicit.
\item A Tate ring is a Huber ring containing a pseudouniformizer, i.e., a topologically nilpotent unit.
\item For a Noetherian adic ring $R$ with ideal of definition $J$, write
$R_\mr{red}=R/R^{\circ\circ}=R/\sqrt{J}$ for its topological reduction.
\item We write $\mr{Int}(g)(x)=gxg^{-1}$.
\item For an algebraic group $L$ over a ring $R$ and a subset $S\subseteq L(R)$ we write $\mr{Cent}(S,L)$ for the centralizer of $S$ in $L$. We further write $L_\gamma=\mr{Cent}(\gamma,L)$ and $I_\gamma=L_\gamma^\circ$.
\item For the Weil group $W_K$ of a local field $K$, we write $\mathsf{v}\colon W_K\to \Z$ for the valuation taking the value $1$ on (a lift of) geometric Frobenius.
\item The symbol $\sigma$ will often (but not always) denote a Frobenius lift on a ring, most often the Witt vectors of a perfect field of characteristic $p$.
\item Our conventions concerning cocharacters, filtrations, and their associated parabolic subgroups follow \S\ref{sec:filtered-Tannakian-appendix} and more specifically Remark \ref{rem:sign-conventions}.
\item For a finite-dimensional representation $\rho$ on $V$, we write
$\tr(x\mid\rho)$ or $\tr(x\mid V)$ for $\tr(\rho(x))$.
\item For an endomorphism $u$ of a perfect complex $M$ over a commutative ring, we write $\tr(u)$ for $\sum_i(-1)^i\tr(u^i)$ for any bounded finite-projective representative $(P^\bullet,u^\bullet)$ of $M$ and $u$.
\item For any cohomology theory $H^i$, we will denote $H^\ast=\sum_{i\geqslant 0}(-1)^i H^i$.
\item Our conventions concerning derived algebraic geometry are as in \cite{GMM}.
\end{itemize}

\paragraph*{Acknowledgements}

I would like to thank Piotr Achinger, Alexander Bertoloni Meli, Patrick Daniels, David Hansen, Hiroki Kato, Naoki Imai, Keerthi Madapusi, Emile Okada, Arghya Sadhukhan, Peter Scholze, Sug Woo Shin, Peihang Wu, and Bogdan Zavyalov for helpful conversations and suggestions.

\medskip

\paragraph*{AI usage} Most of the ideas for this paper were written between 2024 and early 2026 before the AI models reached the point they are currently at in Fall 2026. That said, since Summer 2026 they have become an indispensable tool for literature searching, TeX help, and proof checking/correction/streamlining; I have used them very extensively for this purpose (mainly ChatGPT 5.6 Sol, and ChatGPT 6 Astra). 

The proof checking/correcting/streamlining was the most impactful part of AI's usage here. The vast majority of the issues that it found were important to fix, but fairly routine: normalization errors, missing hypotheses used in a proof, incorrect bounds made in estimation proofs, etc. While I cannot list all instances of these suggestions here, there are two which incur particular intellectual debt. In my original proofs of Theorem \ref{thm:l-independence-formally-of-finite-type} and Proposition \ref{prop:uniform-continuity-K-level}, I had lengthy, multi-step arguments that used the skeleton of the current proofs throughout. ChatGPT pointed out that I was greatly overcomplicating things, leading to the current streamlined proofs.

Despite the above, I would like to emphasize the obvious: All remaining errors in the paper are entirely my responsibility.

\addtocontents{toc}{\protect\setcounter{tocdepth}{2}}
\renewcommand{\thesubsection}{\thesection.\arabic{subsection}}

\section{A rigid-analytic trace formula}
\label{s:prelims}

In this section we establish a formula which computes traces of an analytic correspondence on a rigid $K$-space $Y$ in terms of cohomology of Berthelot tubes of a formal model. The reader should compare the material of this section with \cite{Fargues}, which establishes many of these results for the algebraizable situation.

\begin{nota}\label{nota:rigid-analytic-general} We shall freely use standard terminology, notation, and results about rigid spaces, formal schemes, their generic fibers, Berthelot tubes, etc. For a reminder on these concepts see \cite[\S2]{ALY2} and \cite[\S2.1]{AchingerYoucis}, and the references therein. 

That said, for the convenience of the reader we do recall the following notation and conventions:
\begin{itemize}
\item For a Huber ring $A$, we shorten the notation $\Spa(A,A^\circ)$ to $\Spa(A)$.
\item A \emph{rigid $K$-space} is an adic space locally of finite type over $\Spa(K)$.
\item For a locally formally of finite type formal $\mc{O}_K$-scheme $\mf{X}$, we denote by $\mf{X}_\eta$ its generic fiber in the sense of Berthelot. Equationally, if $\mf{X}^\mr{ad}$ is the adic space associated to $\mf{X}$, then $\mf{X}_\eta=\mf{X}^\mr{ad}\times_{\Spa(\mc{O}_K)}\Spa(K)$.
\item For a closed point $x$ of $\mf{X}$, we define the \emph{Berthelot tube} to be $\mf{X}(x)\defeq \mr{sp}_\mf{X}^{-1}(x)^\circ$, where $\mr{sp}_\mf{X}\colon |\mf{X}_\eta|\to |\mf{X}|$ is the specialization map.
\end{itemize} 
\end{nota}

\subsection{The Fujiwara--Varshavsky trace formula}\label{ss:FV-trace-formula} In this subsection we recall Fujiwara's solution to Deligne's conjecture (see \cite{Fujiwara}), later improved by Varshavsky in \cite{Varshavsky}. The output of this is a formula relating the trace of a cohomological correspondence on a variety over a finite field to a weighted point count of the fixed points of that correspondence.

\subsubsection{Correspondences} We first recall the basic theory of correspondences.

\begin{setup} Throughout we set:
\begin{itemize}
\item $S$ to be a Noetherian (formal) scheme;
\item $\Lambda$ to be a finite ring with $|\Lambda|$ in $\mc{O}_S^\times$, or a finite extension of $\Z_\ell$ or $\Q_\ell$ (with $\ell$ in $\mc{O}_S^\times$);
\item $D_\mr{ctf}(X,\Lambda)$ to be the category defined in \stacks{0F4M}.
\end{itemize}
\end{setup}

\begin{defn} We recall the following definitions.
\begin{enumerate} 
\item A \emph{correspondence} on (topologically of) finite type separated (formal) $S$-schemes is a pair of arrows of (formal) schemes (not necessarily over $S$):
\begin{equation*}
    c\colon \qquad X_1\xleftarrow{c_1}X_0\xrightarrow{c_2}X_2.
\end{equation*}
For a property $P$ of morphisms of (formal) schemes, we say that $c$ is $P$ if each $c_i$ is. If $X_1=X_2$ we may abbreviate the notation for a correspondence to $c=(c_1,c_2)\colon X_0\rightrightarrows X$.
\item A \emph{morphism of correspondences} $d\to c$ is a triple $f=(f_1,f_0,f_2)$ with $f_i\colon Y_i\to X_i$ a morphism of (formal) $S$-schemes such that the following diagram commutes:
\begin{equation}\label{eq:map-of-corr}
    \begin{tikzcd}
	{Y_1} & {Y_0} & {Y_2} \\
	{X_1} & {X_0} & {X_2.}
	\arrow["{d_1}"', from=1-2, to=1-1]
	\arrow["{d_2}", from=1-2, to=1-3]
	\arrow["{c_1}", from=2-2, to=2-1]
	\arrow["{c_2}"', from=2-2, to=2-3]
	\arrow["{f_1}"', from=1-1, to=2-1]
	\arrow["{f_0}"', from=1-2, to=2-2]
	\arrow["{f_2}"', from=1-3, to=2-3]
\end{tikzcd}
\end{equation}
For a property $P$ of morphisms of (formal) schemes, we say that $f$ is $P$ if each $f_i$ is.

\item For a correspondence $c$ and sheaves $\mc{F}_i$ in $D_\mr{ctf}(X_i,\Lambda)$, a \emph{$c$-cohomological correspondence} is a morphism $u\colon c_{2!}c_1^\ast\mc{F}_1\to\mc{F}_2$ in $D_\mr{ctf}(X_2,\Lambda)$. We denote the $\Lambda$-module of $c$-cohomological correspondences from $\mc{F}_1$ to $\mc{F}_2$ by $\mr{Coh}_c(\mc{F}_1,\mc{F}_2)$.
\end{enumerate}
\end{defn}

We often make use of the following alternative presentation of a cohomological correspondence. 

\begin{nota}\label{nota:adjoint-correspondence} As $c_2^!$ is right adjoint to $c_{2!}$ we see that a $c$-cohomological correspondence $u\colon c_{2!}c_1^\ast\mc{F}_1\to \mc{F}_2$ is equivalent to a morphism $c_1^\ast\mc{F}_1\to c_2^!\mc{F}_2$ which we denote $u^{\dashv}$.
\end{nota}

\begin{rem}Suppose that $c$ is finite \'etale so that, in particular, $c_2^!=c_2^\ast$. In this case, we can describe for a geometric point $x$ of $X_2$ the map $u_x\colon (c_{2!}c_1^\ast\mc{F}_1)_x\to \mc{F}_{2x}$ easily in terms of $u^\dashv$. In particular, by \stacks{0F5F} we have that
\begin{equation*}
    (c_{2!}c_1^\ast\mc{F}_1)_x=\bigoplus_{y\in c_2^{-1}(x)}\mc{F}_{1c_1(y)},
\end{equation*}
and we then have
\begin{equation}\label{eq:stalk-map}
    u_x=\sum_{y\in c_2^{-1}(x)}u^\dashv_y,
\end{equation}
where for each point $y$ of $c_2^{-1}(x)$ we have the map $u^\dashv_y\colon \mc{F}_{1c_1(y)}\to \mc{F}_{2c_2(y)}=\mc{F}_{2x}$.
\end{rem}

\begin{defn} Suppose that $c$ is a finite \'etale correspondence. We define the \emph{tensor product}
\begin{equation*}
    \otimes\colon \mr{Coh}_c(\mc{F}_1,\mc{F}_2)\times 
    \mr{Coh}_c(\mc{G}_1,\mc{G}_2)\to \mr{Coh}_c(\mc{F}_1\otimes_\Lambda\mc{G}_1,\mc{F}_2\otimes_\Lambda\mc{G}_2)
\end{equation*}
uniquely characterized so that $(u\otimes v)^\dashv=u^\dashv\otimes v^\dashv$. 
\end{defn}

\begin{rem} It's easy to see that one has the formula
\begin{equation}\label{eq:stalk-tensor}
    (u\otimes v)_x=\sum_{y\in c_2^{-1}(x)}(u_y^\dashv\otimes v_y^\dashv),
\end{equation}
for any geometric point $x$ of $X_2$, with notation as in \eqref{eq:stalk-map}.
\end{rem}

\subsubsection{Functoriality of cohomological correspondences} While pushforwards and pullbacks of cohomological correspondences exist in much greater generality (e.g., see \cite[Part 1]{FYZ}) we will only need them in the basic case of finite \'etale morphisms of finite \'etale correspondences.

\begin{construction}
Let $f\colon d\to c$ be a finite \'etale morphism of finite \'etale correspondences. In this case we have that $f_{i!}=f_{i\ast}$, $c_{i!}=c_{i\ast}$, and $d_{i!}=d_{i\ast}$, a fact which we use below. That said, for readability we use both notations.

\begin{enumerate}
    \item Let $u\colon d_{2!}d_1^\ast\mc{F}_1\to\mc{F}_2$ be a $d$-cohomological correspondence. Define 
    \begin{equation*}
        f_\ast(u)\colon c_{2!}c_1^{\ast}f_{1\ast}\mc{F}_1\to f_{2\ast}\mc{F}_2
    \end{equation*} 
    to be the $c$-cohomological correspondence given by the following composition
    \begin{equation*}
        c_{2!}c_1^{\ast}(f_{1\ast}\mc{F}_1)\xrightarrow{\mr{BC}_\ast}c_{2!}f_{0\ast}d_1^\ast\mc{F}_1=f_{2\ast}d_{2!}d_1^\ast\mc{F}_1\xrightarrow{f_{2\ast}(u)}f_{2\ast}\mc{F}_2,
    \end{equation*}
    where $\mr{BC}_\ast$ is as in \stacks{07A7}.
    \item Suppose that $v\colon c_{2!}c_1^{\ast}\mc{F}_1\to \mc{F}_2$ is a $c$-cohomological correspondence. We define the $d$-cohomological correspondence 
    \begin{equation*}
        f^\ast(v)\colon d_{2!}d_1^\ast f_1^\ast\mc{F}_1\to f_2^\ast\mc{F}_2
    \end{equation*}
    to be such that $f^\ast(v)^\dashv$ agrees with
    \begin{equation*}
        d_1^\ast f_1^\ast\mc{F}_1=f_0^\ast c_1^\ast\mc{F}_1\xrightarrow{f^\ast_0(v^\dashv)} f_0^\ast c_2^\ast\mc{F}_2=d_2^\ast f_2^\ast\mc{F}_2.
    \end{equation*}

\end{enumerate}
\end{construction}

\begin{rem}As $c$ and $d$ are finite \'etale, one can describe the maps $f_\ast(u)$ and $f^\ast(v)$ explicitly on stalks using \stacks{0F5F} as follows:
\begin{enumerate}
    \item Fix $x$ a geometric point of $X_2$. Then, there is a functorial identification 
    \begin{equation}\label{eq:pushforward-stalk}
        f_\ast(u)_x=\bigg(a\colon \bigoplus_{y\in c_2^{-1}(x)}\bigoplus_{\ell\in f_1^{-1}(c_1(y))}\mc{F}_{1\ell}\to \bigoplus_{w\in f_2^{-1}(x)}\mc{F}_{2w}\bigg)
    \end{equation}
    where the $(y,\ell,w)$-matrix entry of $a$ is given by
    \begin{equation}\label{eq:pushforward-matrix-entry}
        a(y,\ell,w)=\sum_{q\in X(y,\ell,w)}u_q^\dashv,
    \end{equation}
    where $X(y,\ell,w)=f_0^{-1}(y)\cap d_1^{-1}(\ell)\cap d_2^{-1}(w)$.
     \item Fix $y$ a geometric point of $Y_2$. Then, there is a functorial identification
    \begin{equation}\label{eq:pullback-stalk}
        f^\ast(v)_y=\sum_{z\in d_2^{-1}(y)}f^\ast(v)^\dashv_z=\sum_{z\in d_2^{-1}(y)}v_{f_0(z)}^\dashv.
    \end{equation}
\end{enumerate}
\end{rem}

Finally, we establish an analogue of the projection formula for cohomological correspondences. 

\begin{defn} Let $d$ be a finite \'etale correspondence. We define
\begin{equation*} \mathrm{tr}_d\colon d_{2!}d_1^\ast\Lambda\to\Lambda,
\end{equation*}
to be the unique element of $\mr{Coh}_d(\Lambda,\Lambda)$ so that $\tr_d^\dashv$ is the identity.
\end{defn}

\begin{prop}\label{prop:proj-formula}Suppose that $f\colon d\to c$ is a finite \'etale morphism of finite \'etale correspondences. Then, for a $c$-cohomological correspondence $u\colon c_{2!}c_1^\ast\mc{F}_1\to\mc{F}_2$ there is an identification
\begin{equation*}
    \bigg[f_\ast(f^\ast(u))\colon c_{2!}c^{\ast}_1 f_{1\ast}f_1^\ast\mc{F}_1\to f_{2\ast}f_2^\ast\mc{F}_2\bigg]\simeq \bigg[u\otimes f_\ast(\tr_d)\colon c_{2!}c^{\ast}_1(\mc{F}_1\otimes_\Lambda f_{1\ast}\Lambda)\to \mc{F}_2\otimes_\Lambda f_{2\ast}\Lambda\bigg]
\end{equation*}
\end{prop}
\begin{proof} More precisely, we claim that the diagram
\begin{equation}\label{eq:proj-formula-proof}
    \begin{tikzcd}[column sep=4em, row sep=2em]
	{c_{2!}c^{\ast}_1(\mc{F}_1\otimes_\Lambda f_{1\ast}\Lambda)} & {\mc{F}_2\otimes_\Lambda f_{2\ast}\Lambda} \\
	{c_{2!}c^{\ast}_1 f_{1\ast}f_1^\ast\mc{F}_1} & {f_{2\ast}f_2^\ast\mc{F}_2,}
	\arrow["{f_\ast(f^\ast(u))}", from=2-1, to=2-2]
	\arrow["{c_{2!}c_1^\ast(\dagger)}"', from=1-1, to=2-1]
	\arrow["{u\otimes f_\ast(\tr_d)}", from=1-1, to=1-2]
	\arrow["\dagger", from=1-2, to=2-2]
\end{tikzcd}
\end{equation}
commutes, where $\dagger$ is the isomorphism from \stacks{0B54}. To check this, it suffices to check it at the stalk of an (arbitrary) geometric point $x$ of $X_2$. Using \stacks{0F5F} we may identify the diagram obtained by passing to the stalk at $x$ of \eqref{eq:proj-formula-proof} with

\begin{equation*}\begin{tikzcd}[column sep=5em,row sep=2.25em]
	{\displaystyle \bigoplus_{y\in c_2^{-1}(x)}\left(\mathcal{F}_{1c_1(y)}\otimes_\Lambda \bigoplus_{\ell\in f_1^{-1}(c_1(y))}\Lambda\right)} & {\displaystyle \mathcal{F}_{2x}\otimes_\Lambda\bigoplus_{w\in f_2^{-1}(x)}\Lambda} \\
	{\displaystyle \bigoplus_{y\in c_2^{-1}(x)}\bigoplus_{\ell\in f_1^{-1}(c_1(y))}\mc{F}_{1f_1(\ell)}} & {\displaystyle \bigoplus_{w\in f_2^{-1}(x)}\mc{F}_{2f_2(w)},}
	\arrow["{(m(y,\ell,w))}", from=1-1, to=1-2]
	\arrow["{(n(y,\ell,w))}", from=2-1, to=2-2]
	\arrow[from=1-2, to=2-2]
	\arrow[from=1-1, to=2-1]
\end{tikzcd}
\end{equation*}
where the vertical maps are the obvious ones (using that $f_1(\ell)=c_1(y)$ and $f_2(w)=x$), and the $m(y,\ell,w)$ and $n(y,\ell,w)$ are the matrix entries of their respective maps. Using \eqref{eq:pushforward-stalk} and \eqref{eq:pullback-stalk} one sees that
\begin{equation*}
    m(y,\ell,w)=\sum_{q\in X(y,\ell,w)}(u_y^\dashv\otimes \mathrm{id}).
\end{equation*}
On the other hand, using \eqref{eq:stalk-tensor} and \eqref{eq:pushforward-stalk} we see that
\begin{equation*}
    n(y,\ell,w)=\sum_{q\in X(y,\ell,w)}u_{f_0(q)}^\dashv,
\end{equation*}
which are equal as for $q$ in $X(y,\ell,w)$ we have that $f_0(q)=y$. 
\end{proof}

\subsubsection{Induced morphisms on cohomology}
We are most interested in cohomological correspondences over perfect fields where cohomological correspondences induce maps on cohomology. 

\begin{setup}\label{setup:corr-coh-field} Throughout we fix $k$ to be a perfect field with absolute Galois group $\Gamma_k$.
\end{setup}

\begin{construction} Suppose that $c\colon X_1\xleftarrow{c_1}X_0\xrightarrow{c_2} X_2$ is a correspondence with $c_1$ and $c_2$ proper, and that $u\colon c_{2!}c_1^\ast\mc{F}_1\to\mc{F}_2$ is a $c$-cohomological correspondence. We can define morphisms
\begin{equation}\label{eq:maps-on-coh}
    R\Gamma_c(u)\colon R\Gamma_c(X_1,\mc{F}_1)\to R\Gamma_c(X_2,\mc{F}_2),\qquad R\Gamma(u)\colon R\Gamma(X_1,\mc{F}_1)\to R\Gamma(X_2,\mc{F}_2)
\end{equation}
in $D^b_\mr{c}(\Spec(k),\Lambda)$. The first map in \eqref{eq:maps-on-coh} is given as the composition
\begin{equation*}
    R\Gamma_c(X_1,\mc{F}_1)\xrightarrow{c_1^\ast}R\Gamma_c(X_0,c_1^\ast\mc{F}_1)\xrightarrow{\mr{BC}_!}R\Gamma_c(X_2,c_{2!}c_1^\ast\mc{F}_1)\xrightarrow{u}R\Gamma_c(X_2,\mc{F}_2),
\end{equation*}
where $\mr{BC}_!$ denotes the base change morphism (see \stacks{0F7L}), and the second map in \eqref{eq:maps-on-coh} is obtained similarly.
\end{construction} 

\begin{rem}The maps in \eqref{eq:maps-on-coh} are equivalent to the $\Gal(\ov{k}/k)$-equivariant maps
\begin{equation*}
    R\Gamma_c(\ov{u})\colon R\Gamma_c(X_{1,\ov{k}},\mc{F}_{1,\ov{k}})\to R\Gamma_c(X_{2,\ov{k}},\mc{F}_{2,\ov{k}}),\qquad R\Gamma(\ov{u})\colon R\Gamma(X_{1,\ov{k}},\mc{F}_{1,\ov{k}})\to R\Gamma(X_{2,\ov{k}},\mc{F}_{2,\ov{k}})
\end{equation*}
of $\Lambda$-modules, where $\ov{u}\colon c_{2!,\ov{k}}c_{1,\ov{k}}^\ast\mc{F}_{1,\ov{k}}\to\mc{F}_{2,\ov{k}}$ is the correspondence obtained by base change. 
\end{rem}

\begin{prop}\label{prop:fin-morphism-pushforward} Let $f\colon d\to c$ be a finite \'etale morphism of finite \'etale correspondences. For a $d$-cohomological correspondence $u\colon d_{2!}d_1^\ast\mc{F}_1\to\mc{F}_2$ there are identifications of complexes
\begin{equation*}
         \bigg[R\Gamma(u)\colon R\Gamma(Y_1,\mc{F}_1)\to R\Gamma(Y_2,\mc{F}_2)\bigg]\simeq \bigg[R\Gamma(f_\ast(u))\colon R\Gamma(X_1,f_{1\ast}\mc{F}_1)\to R\Gamma(X_2,f_{2\ast}\mc{F}_2)\bigg],
    \end{equation*}
and
    \begin{equation*}
         \bigg[R\Gamma_c(u)\colon R\Gamma_c(Y_1,\mc{F}_1)\to R\Gamma_c(Y_2,\mc{F}_2)\bigg]\simeq \bigg[R\Gamma_c(f_\ast(u))\colon R\Gamma_c(X_1,f_{1\ast}\mc{F}_1)\to R\Gamma_c(X_2,f_{2\ast}\mc{F}_2)\bigg].
    \end{equation*}
\end{prop}
\begin{proof}For the first claim, we more precisely claim that we have a commutative diagram
\begin{equation*}
    \begin{tikzcd}[sep=large]
	{R\Gamma(X_1,f_\ast\mathcal{F}_1)} & {R\Gamma(X_2,f_\ast\mathcal{F}_2)} \\
	{R\Gamma(Y_1,\mathcal{F})} & {R\Gamma(Y_2,\mathcal{F}_2),}
	\arrow["\wr"', from=1-1, to=2-1]
	\arrow["{R\Gamma(f_\ast(u))}", from=1-1, to=1-2]
	\arrow["\wr", from=1-2, to=2-2]
	\arrow["{R\Gamma(u)}", from=2-1, to=2-2]
\end{tikzcd}
\end{equation*}
where the left vertical map is that obtained by the composition
\begin{equation*}
    R\Gamma(X_1,f_\ast\mathcal{F}_1)\to R\Gamma(Y_1,f^\ast f_\ast\mathcal{F}_1)\to R\Gamma(Y_1,\mathcal{F}_1),
\end{equation*}
with the last map obtained by the counit map $f^\ast f_\ast\mathcal{F}_1\to\mc{F}_1$, and the right vertical arrow is obtained similarly; these are both isomorphisms by \cite[Proposition 5.7.4 and Corollary 5.6.9 (i)]{LeiFu}. The commutativity follows because $R\Gamma(u)$ can be interpreted as $g_\ast$, where $g$ is the map from $d$ to the trivial correspondence $c_\mr{tr}$ (e.g., see \cite[\S1.1]{Varshavsky}), and the base change map commutes with compositions; cf.\@ \cite[\S1.1.7]{Varshavsky}. The proof of the second claim follows similarly.
\end{proof}

\subsubsection{Twisting correspondences}\label{ss:Galois-twisted-correspondence} We now discuss the procedure of twisting correspondences by elements of the absolute Galois group. We are still in the setting of Setup \ref{setup:corr-coh-field}.

\begin{nota} Fix $\tau$ in $\Gamma_k$. Then, for any $k$-scheme $X$, we have a natural morphism
\begin{equation*}
    \tau_X\colon X_{\ov{k}}\simeq \tau^\ast X_{\ov{k}}\to X_{\ov{k}}.
\end{equation*}
\end{nota}

\begin{rem} If $f\colon Y\to X$ is a morphism of $k$-schemes, then $f_{\ov{k}}\circ \tau_Y=\tau_X\circ f_{\ov{k}}$.
\end{rem}

\begin{defn} For a correspondence
\begin{equation*}
    c\colon \quad X_1\xleftarrow{c_1}X_0\xrightarrow{c_2}X_2,
\end{equation*}
 over $k$ and an element $\tau$ of $\Gal(\ov{k}/k)$ we define the \emph{$\tau$-twist} of $c$ to be the correspondence
\begin{equation*}
    c^\tau\colon\quad  X_{1,\ov{k}}\xleftarrow{c_{1,\ov{k}}\circ \tau_{X_0}}X_{0,\ov{k}}\xrightarrow{c_{2,\ov{k}}}X_{2,\ov{k}}.
\end{equation*}
\end{defn}

\begin{constr}For a $c$-cohomological correspondence $u\colon c_{2!}c_1^\ast\mc{F}_1\to\mc{F}_2$ we obtain a $c^\tau$-cohomological correspondence $u^\tau$ induced by $\ov{u}$ and the natural identifications
\begin{equation}\label{eq:tau-pullback-identifications}
    (c_{1,\ov{k}}\circ \tau_{X_0})^\ast\mc{F}_{1,\ov{k}}=(\tau_{X_1}\circ c_{1,\ov{k}})^\ast\mc{F}_{1,\ov{k}}\simeq c_{1,\ov{k}}^\ast\tau_{X_1}^\ast\mc{F}_{1,\ov{k}}\simeq c_{1,\ov{k}}^\ast\mc{F}_{1,\ov{k}},
\end{equation}
where the last isomorphism is the natural one coming from the descent datum of $\mc{F}_{1,\ov{k}}$ to $\mc{F}_1$. 
\end{constr}

\begin{rem}\label{rem:twisted-coh-corr-on-coh} Unraveling the definitions shows that
\begin{equation*}
    R\Gamma_c(u^\tau)=\tau\times R\Gamma_c(\ov{u}),\qquad R\Gamma(u^\tau)=\tau\times R\Gamma(\ov{u}),
\end{equation*}
where $\tau\times R\Gamma_c(\ov{u})$ denotes the composition in either order (which is unambiguous as they commute), and where we define $\tau\times R\Gamma(\ov{u})$ similarly.
\end{rem} 

The following is an immediate corollary of Proposition \ref{prop:proj-formula} and Proposition \ref{prop:fin-morphism-pushforward} as twisting morphisms commutes with pushforwards, pullbacks, and tensor products in the evident way.

\begin{prop} Let $f\colon d\to c$ be a finite \'etale morphism of finite \'etale correspondences over $k$ and let $\tau$ be an element of $\Gal(\ov{k}/k)$. For a $d$-cohomological correspondence $v\colon d_{2!}d_1^\ast\mc{F}_1\to\mc{F}_2$ there are identifications of complexes
\begin{equation*}
         \resizebox{\hsize}{!}{$\bigg[ R\Gamma(v^\tau)\colon R\Gamma(Y_{1,\ov{k}},\mc{F}_{1,\ov{k}})\to R\Gamma(Y_{2,\ov{k}},\mc{F}_{2,\ov{k}})\bigg]\simeq \bigg[R\Gamma(f_\ast(v)^\tau)\colon R\Gamma(X_{1,\ov{k}},f_{1\ast}\mc{F}_{1,\ov{k}})\to R\Gamma(X_{2,\ov{k}},f_{2\ast}\mc{F}_{2,\ov{k}})\bigg],$}
    \end{equation*}
and
    \begin{equation*}
       \resizebox{\hsize}{!}{$\bigg[R\Gamma_c(v^\tau)\colon R\Gamma_c(Y_{1,\ov{k}},\mc{F}_{1,\ov{k}})\to R\Gamma_c(Y_{2,\ov{k}},\mc{F}_{2,\ov{k}})\bigg]\simeq \bigg[R\Gamma_c(f_\ast(v)^\tau)\colon R\Gamma_c(X_{1,\ov{k}},f_{1\ast}\mc{F}_{1,\ov{k}})\to R\Gamma_c(X_{2,\ov{k}},f_{2\ast}\mc{F}_{2,\ov{k}})\bigg],$}
    \end{equation*}
If, moreover, $v\simeq f^\ast(u)$ for a $c$-cohomological correspondence $u\colon c_{2!}c_1^\ast\mc{F}_1\to\mc{F}_2$ then there are further quasi-isomorphisms of complexes
\begin{equation*}
         \resizebox{\hsize}{!}{$\bigg[ R\Gamma(v^\tau)\colon R\Gamma(Y_{1,\ov{k}},\mc{F}_{1,\ov{k}})\to R\Gamma(Y_{2,\ov{k}},\mc{F}_{2,\ov{k}})\bigg]\simeq \bigg[R\Gamma(u^\tau\times f_\ast(\tr_d)^\tau)\colon R\Gamma(X_{1,\ov{k}},f_{1\ast}\mc{F}_{1,\ov{k}})\to R\Gamma(X_{2,\ov{k}},f_{2\ast}\mc{F}_{2,\ov{k}})\bigg],$}
    \end{equation*}
and
    \begin{equation*}
      \resizebox{\hsize}{!}{$\bigg[R\Gamma_c(v^\tau)\colon R\Gamma_c(Y_{1,\ov{k}},\mc{F}_{1,\ov{k}})\to R\Gamma_c(Y_{2,\ov{k}},\mc{F}_{2,\ov{k}})\bigg]\simeq \bigg[R\Gamma_c(u^\tau\times f_\ast(\tr_d)^\tau)\colon R\Gamma_c(X_{1,\ov{k}},f_{1\ast}\mc{F}_{1,\ov{k}})\to R\Gamma_c(X_{2,\ov{k}},f_{2\ast}\mc{F}_{2,\ov{k}})\bigg].$}
    \end{equation*}
\end{prop}

\subsubsection{The Fujiwara--Varshavsky trace formula} We are now ready to discuss the Fujiwara--Varshavsky trace formula which computes the trace of a cohomological correspondence in terms of so-called naive local terms at fixed points.

\begin{setup} Fix a perfect field $k$ and a correspondence $c=(c_1,c_2)\colon X_0\rightrightarrows X$.
\end{setup}

\begin{defn}\label{defn:fixed-scheme} Define the \emph{fixed scheme} of $c$ to be $\mr{Fix}(c)=\mr{Eq}(c_1,c_2)$. We consider $\mr{Fix}(c)$ as a closed subscheme of $X_0$. 
\end{defn}

\begin{defn}\label{defn:naive-local-term} Suppose that $c_2$ is finite \'etale, and let $w\colon c_{2!}c_1^\ast\mc{F}\to\mc{F}$ be a $c$-cohomological correspondence. For a $\ov{k}$-point $y$ of $\mr{Fix}(c)$, we write $x$ for the common point $c_1(y)=c_2(y)$. Thus, $w_y^\dashv$ is a morphism $\mc{F}_x\to\mc{F}_x$ of complexes of $\Lambda$-modules, and we call the element $\mr{tr}(w_y^\dashv)$ of $\Lambda$ the \emph{naive local term \emph{(}at $y$\emph{)}}.
\end{defn}

\begin{nota}\label{nota:Frobenius-twist-corr}
Let $k$ be a field of order $p^r$, and let $n\geqslant 0$.
For a $k$-scheme $Y$, write $\Phi_Y=(F_Y^r)_{\ov{k}}$, with
$F_Y$ the absolute Frobenius.
For a correspondence of $k$-schemes $c\colon X_1\xleftarrow{c_1}X_0\xrightarrow{c_2}X_2$, define
\begin{equation*}
c^{(n)}\colon\quad X_{1,\ov{k}}\xleftarrow{c_1^{(n)}=\Phi_{X_1}^{n}\circ c_{1,\ov{k}}}
X_{0,\ov{k}}
\xrightarrow{c_{2,\ov{k}}}X_{2,\ov{k}}.
\end{equation*}
\end{nota}

\begin{rem}\label{rem:frobenius-galois-pullback}
Let $\sigma$ in $\Gamma_k$ denote geometric Frobenius, and set $\sigma_Y\defeq \sigma\otimes 1$ viewed as a morphism on $Y_{\ov{k}}$. Then, the identity $\Phi_Y^n=F_{Y_{\ov{k}}}^{rn}\circ\sigma_Y^n$ gives rise to a canonical identification $(\Phi_Y^n)^\ast\simeq(\sigma_Y^n)^\ast$; see \stacks{03SN} and \stacks{03SV}.

Thus, for a $c$-cohomological correspondence
$u\colon c_{2!}c_1^\ast\mc{F}_1\to\mc{F}_2$ defined over $k$,
define the $c^{(n)}$-cohomological correspondence $u^{(n)}$ by
\begin{equation*}
 c_{2,\ov{k}!}(c_1^{(n)})^\ast\mc{F}_{1,\ov{k}}\xrightarrow{\sim}
 c_{2,\ov{k}!}(c_{1,\ov{k}}\circ\sigma_{X_0}^n)^\ast
       \mc{F}_{1,\ov{k}}\xrightarrow{u^{\sigma^n}}\mc{F}_{2,\ov{k}}.
\end{equation*}
Thus, from the above discussion, $R\Gamma_c(u^{(n)})=R\Gamma_c(u^{\sigma^n})$, and similarly for ordinary cohomology.
\end{rem}

\begin{thm}[{cf.\@ \cite[Theorem 2.3.2]{Varshavsky} and \cite{Fujiwara}}]\label{thm:Fujiwara--Varshavsky-trace-formula} Suppose that $k$ is a finite field of order $q$ and that $c$ is finite \'etale. Then, for every $n\geqslant 1$ the $\ov{k}$-scheme $\mr{Fix}(c^{(n)})$ is finite. Moreover, there exists an $m\geqslant 1$ such that if $q^n>m$ and $u$ is any $c^{(n)}$-cohomological correspondence, then
\begin{equation*}
    \tr\left(R\Gamma_c({u})\right)=\sum_{y\in\mr{Fix}(c^{(n)})(\ov{k})}\tr\left(u_y^\dashv\right).
\end{equation*}
If $c$ admits a finite compactification $\widetilde{c}$,\footnote{A finite compactification of $c$ is a finite correspondence $\widetilde c\colon\widetilde X_0\rightrightarrows\widetilde X$ of proper $k$-schemes, together with dense open immersions $j_0\colon X_0\hookrightarrow\widetilde X_0$ and $j\colon X\hookrightarrow\widetilde X$ satisfying $\widetilde c_i\circ j_0=j\circ c_i$ for $i=1,2$. The ramification index is that of \cite[Notation 2.2.2(b)]{Varshavsky}.} then $m$ may be taken to be the ramification index of $\widetilde{c}_2$.
\end{thm}

\subsection{A rigid-analytic trace formula}

In this subsection we explain how Theorem \ref{thm:Fujiwara--Varshavsky-trace-formula} allows one to give a formula computing the trace of a Galois-twisted correspondence on the generic fiber of a formal scheme in terms of a weighted point count on the special fiber.

\begin{nota}\label{nota:rigid-analytic-trace-formula} We fix the following notation:
\begin{multicols}{2}
\begin{itemize}
    \item $K$ is a finite extension of $\Q_p$;
    \item $\mc{O}_K$ is the valuation ring of $K$;
    \item $k$ is the residue field of $K$ and $q=|k|$;
    \item $C\defeq \wh{\ov{K}}$;
    \item $\eta=\Spa(K)$ and $\Gamma_\eta=\Gal(\ov{K}/K)$; 
    \item $s=\Spec(k)$ and $\Gamma_s=\Gal(\ov{k}/k)$; 
    \item $r\colon \Gamma_\eta\to \Gamma_s$ is the reduction map;
    \item $W_K$ is the Weil group of $K$;
    \item $\mathsf{v}\colon W_K\to \Z$ is the usual (geometrically-normalized) valuation;
    \item $\Lambda$ is a finite ring with order invertible in $k$, or a finite extension of $\Z_\ell$ or $\Q_\ell$ (with $\ell$ invertible in $k$).
\end{itemize}
\end{multicols}
\end{nota}

\subsubsection{Analytic correspondences} We now discuss the theory of analytic correspondences. This theory is obtained, mutatis mutandis, from the algebraic theory (once the appropriate assumptions/definitions are made), and so we omit details analogous to those in \S\ref{ss:FV-trace-formula}.

\begin{defn}[{cf.\@ \cite{BhattHansen}}] A map $f\colon Y\to X$ of rigid $K$-spaces is (in decreasing generality):
\begin{itemize}
    \item \emph{reasonable} if it is separated, taut (see \cite[Definition 5.1.2]{HuberEC}), and finite-dimensional (i.e., that $\mr{dim.tr}(f)<\infty$ with notation as in \cite[Definition 1.8.3]{HuberEC});
    \item \emph{Zariski-compactifiable} if there exists a Zariski open immersion $Y\hookrightarrow \ov{Y}$ over $X$, where $\ov{Y}$ is a proper adic $X$-space.
\end{itemize}
A rigid $K$-space $X$ is \emph{reasonable} (resp.\@ \emph{Zariski-compactifiable}) if $X\to\Spa(K)$ is.
\end{defn}

\begin{nota} Denote by $D_\mr{zc}(X,\Lambda)$ the category of Zariski-constructible $\Lambda$-modules in the sense of \cite[Definition 3.1]{BhattHansen}. 
\end{nota}

\begin{rem}
The following hold for
locally bounded complexes and reasonable $f$:
\begin{itemize}[leftmargin=.5cm]
    \item $f^\ast$ preserves Zariski-constructibility;
    see \cite[Proposition 3.4(3) and Theorem 3.36(1)]{BhattHansen};

    \item $Rf^!$ preserves Zariski-constructibility;
    see \cite[Remark 3.13 and Theorem 3.36(5)]{BhattHansen};

    \item if $\mc{F}$ and $\mc{G}$ are Zariski-constructible and
    $\mc{F}$ has finite Tor-dimension, then
    $\mc{F}\otimes_\Lambda^L\mc{G}$ and
    $R\mathscr{H}om_\Lambda(\mc{F},\mc{G})$
    are Zariski-constructible;
    see \cite[Corollary 3.14 and Theorem 3.36(1)]{BhattHansen};

    \item if $f$ is proper, $Rf_\ast$ preserves
    Zariski-constructibility;
    see \cite[Theorems 3.10 and 3.36(3)]{BhattHansen};

    \item if $f$ is Zariski-compactifiable, both $Rf_\ast$ and $Rf_!$
    send lisse complexes to Zariski-constructible complexes;
    see \cite[Corollary 3.11 and Theorem 3.36(4)]{BhattHansen}.
\end{itemize}
\end{rem}

\begin{defn}\label{defn:analytic-corr} We make the following definitions:
\begin{enumerate} 
\item An \emph{analytic correspondence} of locally of finite type separated rigid $K$-spaces is a pair of arrows (not necessarily over $K$):
\begin{equation*}
    c\colon \qquad X_1\xleftarrow{c_1}X_0\xrightarrow{c_2}X_2.
\end{equation*}
We always assume that $c_2$ is Zariski-compactifiable. For a property $P$ of morphisms of adic spaces, we say that $c$ is $P$ if each $c_i$ is. If $X_1=X_2$ we may abbreviate the notation to $c=(c_1,c_2)\colon X_0\rightrightarrows X$.
\item A \emph{morphism of analytic correspondences} $d\to c$ is a triple $f=(f_1,f_0,f_2)$ where $f_i\colon Y_i\to X_i$ is a morphism of rigid $K$-spaces where the diagram as in \eqref{eq:map-of-corr} commutes.
\item For an analytic correspondence $c$ and sheaves $\mc{F}_i$ in $D_\mr{zc}(X_i,\Lambda)$, a \emph{$c$-cohomological correspondence} is a morphism $u\colon c_{2!}c_1^\ast\mc{F}_1\to\mc{F}_2$ in $D(X_2,\Lambda)$. We denote the $\Lambda$-module of such objects by $\mr{Coh}_c(\mc{F}_1,\mc{F}_2)$.
\end{enumerate}
\end{defn}

\begin{nota} If $c$ is finite \'etale, then a $c$-cohomological correspondence $u\colon c_{2!}c_1^\ast\mc{F}_1\to \mc{F}_2$ is equivalent to a morphism $c_1^\ast\mc{F}_1\to c_2^\ast\mc{F}_2$ which we denote $u^{\dashv}$.
\end{nota}

\begin{defn} A geometric point of a rigid $K$-space $X$ is a morphism $x\colon \Spa(K,K^+)\to X$ with $K$ separably closed. 
\end{defn}

\begin{rem} Let notation be as in Definition \ref{defn:analytic-corr}. For a geometric point $x$ of $X$ we have that
\begin{equation*}
    (c_{2!}c_1^\ast\mc{F}_1)_x=\bigoplus_{y\in c_2^{-1}(x)}\mc{F}_{1c_1(y)},
\end{equation*}
as follows from \cite[Corollary 5.2.3]{HuberEC}, and moreover we have that 
\begin{equation}\label{eq:stalk-map-analytic}
    u_x=\sum_{y\in c_2^{-1}(x)}u^\dashv_y,
\end{equation}
where for each point $y$ of $c_2^{-1}(x)$ we have the map $u^\dashv_y\colon \mc{F}_{1c_1(y)}\to \mc{F}_{2c_2(y)}=\mc{F}_{2x}$.
\end{rem}

\begin{defn} Suppose that $c$ is a finite \'etale correspondence. We define the \emph{tensor product}
\begin{equation*}
    \otimes\colon \mr{Coh}_c(\mc{F}_1,\mc{F}_2)\times 
    \mr{Coh}_c(\mc{G}_1,\mc{G}_2)\to \mr{Coh}_c(\mc{F}_1\otimes_\Lambda\mc{G}_1,\mc{F}_2\otimes_\Lambda\mc{G}_2)
\end{equation*}
uniquely characterized so that $(u\otimes v)^\dashv=u^\dashv\otimes v^\dashv$. 
\end{defn}

\subsubsection{Functoriality of analytic correspondences}
We again only need basic functorialities for finite \'etale morphisms of analytic correspondences.

\begin{construction}
Let $f\colon d\to c$ be a finite \'etale morphism of finite \'etale analytic correspondences; we again make use of the fact that $f_{i!}=f_{i\ast}$, $c_{i!}=c_{i\ast}$, and $d_{i!}=d_{i\ast}$.
\begin{enumerate}
    \item Let $u\colon d_{2!}d_1^\ast\mc{F}_1\to\mc{F}_2$ be a $d$-cohomological correspondence. Define 
    \begin{equation*}
        f_\ast(u)\colon c_{2!}c_1^{\ast}f_{1\ast}\mc{F}_1\to f_{2\ast}\mc{F}_2
    \end{equation*} 
    to be the $c$-cohomological correspondence given by the following composition
    \begin{equation*}
        c_{2!}c_1^{\ast}(f_{1\ast}\mc{F}_1)\xrightarrow{\mr{BC}_\ast}c_{2!}f_{0\ast}d_1^\ast\mc{F}_1=f_{2\ast}d_{2!}d_1^\ast\mc{F}_1\xrightarrow{f_{2\ast}(u)}f_{2\ast}\mc{F}_2.
    \end{equation*}
    
    \item Suppose that $v\colon c_{2!}c_1^{\ast}\mc{F}_1\to \mc{F}_2$ is a $c$-cohomological correspondence. We define the $d$-cohomological correspondence 
    \begin{equation*}
        f^\ast(v)\colon d_{2!}d_1^\ast f_1^\ast\mc{F}_1\to f_2^\ast\mc{F}_2
    \end{equation*}
     so that $f^\ast(v)^\dashv$ agrees with
    \begin{equation*}
        d_1^\ast f_1^\ast\mc{F}_1=f_0^\ast c_1^\ast\mc{F}_1\xrightarrow{f_0^\ast(v^\dashv)} f_0^\ast c_2^\ast\mc{F}_2=d_2^\ast f_2^\ast\mc{F}_2.
    \end{equation*}

\end{enumerate}
\end{construction}

The following is obtained mutatis mutandis from Proposition \ref{prop:proj-formula}.

\begin{defn} Let $d$ be a finite \'etale analytic correspondence. We define
\begin{equation*} \mathrm{tr}_d\colon d_{2!}d_1^\ast\Lambda\to\Lambda,
\end{equation*}
to be the unique element of $\mr{Coh}_d(\Lambda,\Lambda)$ so that $\tr_d^\dashv$ is the identity.
\end{defn}

\begin{prop}\label{prop:proj-formula-analytic} Suppose that $f\colon d\to c$ is a finite \'etale morphism of finite \'etale analytic correspondences. For a $c$-cohomological correspondence $u\colon c_{2!}c_1^\ast\mc{F}_1\to\mc{F}_2$ there is an identification
\begin{equation*}
    \bigg[f_\ast(f^\ast(u))\colon c_{2!}c^{\ast}_1 f_{1\ast}f_1^\ast\mc{F}_1\to f_{2\ast}f_2^\ast\mc{F}_2\bigg]\simeq \bigg[u\otimes f_\ast(\tr_d)\colon c_{2!}c^{\ast}_1(\mc{F}_1\otimes_\Lambda f_{1\ast}\Lambda)\to \mc{F}_2\otimes_\Lambda f_{2\ast}\Lambda\bigg]
\end{equation*}
\end{prop}

\subsubsection{Induced morphisms on cohomology}
Again cohomological correspondences induce morphisms on cohomology groups. 

\begin{construction}\label{constr:analytic-corr-coh} Let $c\colon X_1\xleftarrow{c_1}X_0\xrightarrow{c_2}X_2$ be an analytic correspondence of reasonable rigid $K$-spaces with $c_1$  and $c_2$ proper, and $u\colon c_{2!}c_1^\ast\mc{F}_1\to\mc{F}_2$ a $c$-cohomological correspondence. Define
\begin{equation}\label{eq:maps-on-coh-an}
    R\Gamma_c(u)\colon R\Gamma_c(X_1,\mc{F}_1)\to R\Gamma_c(X_2,\mc{F}_2),\qquad R\Gamma(u)\colon R\Gamma(X_1,\mc{F}_1)\to R\Gamma(X_2,\mc{F}_2)
\end{equation}
to be the morphisms in $D^b(\Spec(K),\Lambda)$, with the former map given as the composition
\begin{equation*}
    R\Gamma_c(X_1,\mc{F}_1)\xrightarrow{c_1^\ast}R\Gamma_c(X_0,c_1^\ast\mc{F}_1)\xrightarrow{\mr{BC}_!}R\Gamma_c(X_2,c_{2!}c_1^\ast\mc{F}_1)\xrightarrow{u}R\Gamma_c(X_2,\mc{F}_2),
\end{equation*}
with $\mr{BC}_!$ the base change morphism (see \cite[Theorem 5.5.9 i)]{HuberEC}), and similarly for the latter.
\end{construction} 

\begin{rem}The maps in \eqref{eq:maps-on-coh-an} are equivalent to the $\Gamma_\eta$-equivariant maps
\begin{equation*}
    R\Gamma_c(\ov{u})\colon R\Gamma_c(X_{1,C},\mc{F}_{1,C})\to R\Gamma_c(X_{2,C},\mc{F}_{2,C}),\qquad R\Gamma(\ov{u})\colon R\Gamma(X_{1,C},\mc{F}_{1,C})\to R\Gamma(X_{2,C},\mc{F}_{2,C}),
\end{equation*}
with $\ov{u}\colon c_{2!,C}c_{1,C}^\ast\mc{F}_{1,C}\to\mc{F}_{2,C}$ the map obtained by base change. If the $X_i$ are Zariski-compactifiable and the $\mc{F}_i$
are bounded lisse complexes, these complexes have finitely
generated cohomology over $\Lambda$, and vanish outside a
a specific range; see \cite[Corollary 3.11 and Theorem 3.36]{BhattHansen}. 
\end{rem}

The following is obtained, mutatis mutandis, from the proof of Proposition \ref{prop:fin-morphism-pushforward}.

\begin{prop}\label{prop:fin-morphism-pushforward-analytic} Let $f\colon d\to c$ be a finite \'etale morphism of finite \'etale correspondences. For a $d$-cohomological correspondence $u\colon d_{2!}d_1^\ast\mc{F}_1\to\mc{F}_2$ there is an identification of complexes
\begin{equation*}
         \bigg[R\Gamma(u)\colon R\Gamma(Y_1,\mc{F}_1)\to R\Gamma(Y_2,\mc{F}_2)\bigg]\simeq \bigg[R\Gamma(f_\ast(u))\colon R\Gamma(X_1,f_{1\ast}\mc{F}_1)\to R\Gamma(X_2,f_{2\ast}\mc{F}_2)\bigg],
    \end{equation*}
and
    \begin{equation*}
         \bigg[R\Gamma_c(u)\colon R\Gamma_c(Y_1,\mc{F}_1)\to R\Gamma_c(Y_2,\mc{F}_2)\bigg]\simeq \bigg[R\Gamma_c(f_\ast(u))\colon R\Gamma_c(X_1,f_{1\ast}\mc{F}_1)\to R\Gamma_c(X_2,f_{2\ast}\mc{F}_2)\bigg]
    \end{equation*}
\end{prop}

\subsubsection{Twisting correspondences} We now discuss Galois twisting for analytic correspondences.

\begin{nota} Fix $\tau$ in $\Gamma_\eta$. Then, for any rigid $K$-space $X$, we have a natural morphism
\begin{equation*}
    \tau_X\colon X_{C}\simeq \tau^\ast X_{C}\to X_{C},
\end{equation*}
\end{nota}
\begin{rem}  If $f\colon Y\to X$ is a morphism of rigid $K$-spaces, then $f_{C}\circ \tau_Y=\tau_X\circ f_{C}$.
\end{rem}
\begin{defn} For an analytic correspondence
\begin{equation*}
    c\colon \quad X_1\xleftarrow{c_1}X_0\xrightarrow{c_2}X_2,
\end{equation*}
and an element $\tau$ of $\Gamma_\eta$ we define the \emph{$\tau$-twist} of $c$ to be the correspondence
\begin{equation*}
    c^\tau\colon \quad X_{1,C}\xleftarrow{c_{1,C}\circ \tau_{X_0}}X_{0,C}\xrightarrow{c_{2,C}}X_{2,C}.
\end{equation*}
\end{defn}

\begin{nota}For a $c$-cohomological correspondence $u\colon c_{2!}c_1^\ast\mc{F}_1\to\mc{F}_2$ we obtain a $c^\tau$-cohomological correspondence $u^\tau$ using the obvious analogues of the identifications in \eqref{eq:tau-pullback-identifications}. 
\end{nota}

\begin{rem} Unraveling the definitions shows that
\begin{equation*}
    R\Gamma_c(u^\tau)=\tau\times R\Gamma_c(\ov{u}),\qquad R\Gamma(u^\tau)=\tau\times R\Gamma(\ov{u}),
\end{equation*}
where the right-hand sides have the same meaning as in Remark \ref{rem:twisted-coh-corr-on-coh}. 
\end{rem}

The following is an immediate corollary of Proposition \ref{prop:proj-formula-analytic} and Proposition \ref{prop:fin-morphism-pushforward-analytic}, as twisting morphisms commutes with pushforwards, pullbacks, and tensor products in the evident way.

\begin{prop}\label{prop:analytic-pushforward} Let $f\colon d\to c$ be a finite \'etale morphism of finite \'etale analytic correspondences, and let $\tau$ be an element of $\Gamma_\eta$. For a $d$-cohomological correspondence $v\colon d_{2!}d_1^\ast\mc{F}_1\to\mc{F}_2$ there is an identification of complexes
\begin{equation*}
         \resizebox{\hsize}{!}{$\bigg[ R\Gamma(v^\tau)\colon R\Gamma(Y_{1,C},\mc{F}_{1,C})\to R\Gamma(Y_{2,C},\mc{F}_{2,C})\bigg]\simeq \bigg[R\Gamma(f_\ast(v)^\tau)\colon R\Gamma(X_{1,C},f_{1\ast}\mc{F}_{1,C})\to R\Gamma(X_{2,C},f_{2\ast}\mc{F}_{2,C})\bigg],$}
    \end{equation*}
and
    \begin{equation*}
       \resizebox{\hsize}{!}{$\bigg[R\Gamma_c(v^\tau)\colon R\Gamma_c(Y_{1,C},\mc{F}_{1,C})\to R\Gamma_c(Y_{2,C},\mc{F}_{2,C})\bigg]\simeq \bigg[R\Gamma_c(f_\ast(v)^\tau)\colon R\Gamma_c(X_{1,C},f_{1\ast}\mc{F}_{1,C})\to R\Gamma_c(X_{2,C},f_{2\ast}\mc{F}_{2,C})\bigg],$}
    \end{equation*}
If, moreover, $v\simeq f^\ast(u)$ for a $c$-cohomological correspondence $u\colon c_{2!}c_1^\ast\mc{F}_1\to\mc{F}_2$ then there are further identifications
\begin{equation*}
         \resizebox{\hsize}{!}{$\bigg[ R\Gamma(v^\tau)\colon R\Gamma(Y_{1,C},\mc{F}_{1,C})\to R\Gamma(Y_{2,C},\mc{F}_{2,C})\bigg]\simeq \bigg[R\Gamma(u^\tau\times f_\ast(\tr_d)^\tau)\colon R\Gamma(X_{1,C},f_{1\ast}\mc{F}_{1,C})\to R\Gamma(X_{2,C},f_{2\ast}\mc{F}_{2,C})\bigg],$}
    \end{equation*}
and
    \begin{equation*}
      \resizebox{\hsize}{!}{$\bigg[R\Gamma_c(v^\tau)\colon R\Gamma_c(Y_{1,C},\mc{F}_{1,C})\to R\Gamma_c(Y_{2,C},\mc{F}_{2,C})\bigg]\simeq \bigg[R\Gamma_c(u^\tau\times f_\ast(\tr_d)^\tau)\colon R\Gamma_c(X_{1,C},f_{1\ast}\mc{F}_{1,C})\to R\Gamma_c(X_{2,C},f_{2\ast}\mc{F}_{2,C})\bigg].$}
    \end{equation*}
\end{prop}

\subsubsection{Nearby cycle functors}

We now recall the setup for nearby cycle functors in the formal scheme setting, and their interaction with correspondences. 

\begin{nota} We abuse notation by using $\eta$ and $s$ to also denote the \'etale topoi of $\eta$ and $s$, respectively. These may be identified with the topoi of discrete sets with a continuous $\Gamma_\eta$-action (resp.\@ $\Gamma_s$-action), whose covers are surjective maps of $\Gamma_\eta$-sets (resp.\@ $\Gamma_s$-sets). We have a natural map of topoi $\eta\to s$ corresponding to the reduction map $r$.
\end{nota}

\begin{defn} Let $\mf{X}$ be a finite type flat formal $\mc{O}_K$-scheme. The \emph{Deligne topos} is the fiber product in the $2$-category of topoi (see \cite[Proposition 3.4]{GiraudTopos})
\begin{equation*}
    \mf{X}_s\times_s\eta \defeq \mf{X}_{s,\et}\times_s \eta,
\end{equation*}
where $\mf{X}_{s,\et}\to s$ is the obvious map of topoi. We let $p_\mf{X}\colon \mf{X}_s\times_s \eta \to \mf{X}_s$ be the natural projection.
\end{defn}

\begin{rem}\label{rem:topos-ident} We identify $\mf{X}_s\times_s \eta$ with the topos $T_{\Gamma_\eta}(\mf{X}_{\ov{s}})$ of sheaves on $\mf{X}_{\ov{s}}$ with a continuous action of the topological group $\Gamma_\eta$; see \cite[Definitions A.1.1 and A.1.2, and Lemma A.1.4]{HansenZavyalov}. For an object $\mc{G}$ of $D(\mf{X}_s\times_s \eta,\Lambda)$ we denote by $\mc{G}_{\ov{s}}$ the induced object of $D(\mf{X}_{\ov{s}},\Lambda)$.
\end{rem}

\begin{rem} Under the identification of Remark~\ref{rem:topos-ident}, $p_\mf{X}^\ast(\mc{F})$ corresponds to $\mc{F}_{\ov{s}}$ endowed with the action of $\Gamma_\eta$ obtained from the natural $\Gamma_s$-action via $r\colon\Gamma_\eta\to\Gamma_s$. We often use $p_\mf{X}^\ast$ to regard an object of $D(\mf{X}_s,\Lambda)$ as an object of $D(\mf{X}_s\times_s\eta,\Lambda)$ with this action.
\end{rem}

\begin{constr}[{\cite[\S3.5]{HuberEC}}]\label{constr:specialization-morphism} Let $\mf{X}$ be a finite type flat formal $\mc{O}_K$-scheme. Then, there is a morphism of topoi $\lambda_\mf{X}\colon\mf{X}_{\eta,\et}\to \mf{X}_{s,\et}$ with corresponding morphism of sites sending an \'etale morphism $\mf{Y}_s\to \mf{X}_s$ to $\mf{Y}_\eta\to\mf{X}_\eta$ with $\mf{Y}\to\mf{X}$ the unique \'etale deformation of $\mf{Y}_s\to \mf{X}_s$.
\end{constr}

\begin{nota} The $2$-functorial equivalence of \'etale topoi $\mf{X}_\et\simeq \mf{X}_{s,\et}$ induces an equivalence $D(\mf{X},\Lambda)\simeq D(\mf{X}_s,\Lambda)$. For an object $\mc{F}$ of $D(\mf{X},\Lambda)$ we denote by $\mc{F}_s$ its corresponding object in $D(\mf{X}_s,\Lambda)$, and we furthermore write $\mc{F}_\eta\defeq \lambda_\mf{X}^\ast(\mc{F}_s)$.\end{nota}

\begin{rem} If $\mf{X}=\wh{\ms{X}}$ for a separated finite type flat $\mc{O}_K$-scheme $\ms{X}$, and $\mc{F}$ is a locally constant constructible $\Lambda$-module on $\ms{X}$, then for the induced object $\wh{\mc{F}}$ of $D(\mf{X},\Lambda)$, we have an identification $\wh{\mc{F}}_\eta=j_\ms{X}^\ast \mc{F}_\eta^\mr{an}$, where $\mc{F}_\eta^\mr{an}$ is the natural object of $D(\ms{X}_\eta^\mr{an},\Lambda)$ (see \cite[\S3.8]{HuberEC}), and $j_\ms{X}$ is the natural open embedding $\wh{\ms{X}}_\eta\hookrightarrow \ms{X}_\eta^\mr{an}$.
\end{rem}

\begin{defn}\label{defn:nearby-cycles} We define the \emph{nearby cycles functor} to be
\begin{equation*}
    R\Psi_\mf{X}\colon D(\mf{X}_\eta,\Lambda)\to D(\mf{X}_s\times_s\eta,\Lambda)
\end{equation*}
where $\Psi_\mf{X}\colon \mf{X}_{\eta,\et}\to \mf{X}_s\times_s\eta$ is obtained from the $2$-commutative diagram
\begin{equation*}
    \begin{tikzcd}
	{\mf{X}_{\eta,\et}} & {\mf{X}_{s,\et}} \\
	\eta & s.
	\arrow["{\lambda_\mf{X}}", from=1-1, to=1-2]
	\arrow[from=1-2, to=2-2]
	\arrow["r",from=2-1, to=2-2]
	\arrow[from=1-1, to=2-1]
\end{tikzcd}
\end{equation*}
\end{defn}

\begin{nota} We write
\begin{itemize}
\item $D_\mr{ctf}(\mf{X}_s\times_s\eta,\Lambda)$ for the subcategory of $D(\mf{X}_s\times_s\eta,\Lambda)$ spanned by $\mc{G}$ with $\mc{G}_{\ov{s}}$ in $D_\mr{ctf}(\mf{X}_{\ov{s}},\Lambda)$;
\item  $D_{\mr{zc},\mr{ftd}}(\mf{X}_\eta,\Lambda)$ for the category of Zariski-constructible sheaves with finite Tor-dimension.
\end{itemize}
\end{nota}

 We record the following omnibus result concerning the nearby cycles functor. 

\begin{prop}\label{prop:nearby-omnibus} The following statements are true.
\begin{enumerate}[itemsep=0.4em]
    \item The nearby cycles functor commutes in a $2$--functorial way with 
    \begin{itemize}
        \item pushforward; see \cite[Lemma A.3.4]{HansenZavyalov};
        \item \'etale pullback;
        \item exceptional pushforward on the category of finite Tor-dimension Zariski-constructible sheaves; see \cite[Theorem A.3.9]{HansenZavyalov};
        \item extension of ground field; see \cite[Theorem A.3.4]{HansenZavyalov}.
    \end{itemize}
    \item There is an identification $R\Psi_\mf{X}(\mc{F})_{\ov{s}}\simeq R\lambda_{\mf{X}_{\mc{O}_C}\ast}(\mc{F}_{\wh{\ov{\eta}}})$ in $D(\mf{X}_{\ov{s}},\Lambda)$; see \cite[Lemma A.3.4]{HansenZavyalov}.
    \item The functor $R\Psi_\mf{X}$ maps $D_{\mr{zc},\mr{ftd}}(\mf{X}_\eta,\Lambda)$ to $D_\mr{ctf}(\mf{X}_s\times_s\eta,\Lambda)$; see \cite[Lemma A.3.4 (5)]{HansenZavyalov}.
    \item Let $\ms{X}$ be a flat separated finite type $\mc{O}_K$-scheme. For $\mc{G}$ an object of  $D(\ms{X}_\eta,\Lambda)$, there is a functorial identification
\begin{equation*}
 R\Psi_{\ms{X}}(\mc{G})
 \xrightarrow{\sim}
 R\Psi_{\wh{\ms{X}}}\bigl(j_{\ms{X}}^\ast\mc{G}^{\mr{an}}\bigr),
\end{equation*}
where the left-hand side denotes the algebraic nearby cycles functor from \cite[Expos\'e XIII]{SGA7-2}; see \cite[Theorem A.4.4 and Remark B.1]{HansenZavyalov}.
    \item For any locally constant constructible $\Lambda$-module
$\mc{F}$ on $\mf{X}$ and any object $\mc{G}$ of $D(\mf{X}_\eta,\Lambda)$, there is a natural isomorphism
    \begin{equation*}
        \mc{F}_s\otimes_\Lambda R\Psi_\mf{X}(\mc{G})\isomto R\Psi_\mf{X}(\mc{F}_\eta\otimes_\Lambda\mc{G});
    \end{equation*}
    see \cite[Theorem 2.5]{HansenNearby}.
    \item For $\mf{X}/\mc{O}_K$ smooth and $\mc{F}$ locally constant, there is a natural identification $R\Psi_\mf{X}(\mc{F}_\eta)\simeq \mc{F}_s$.
\end{enumerate}
\end{prop}
\begin{proof} It remains only to show (6). There is a natural identification $\mc{F}_\eta=\lambda_\mf{X}^\ast(\mc{F})=\Psi_\mf{X}^\ast(\mc{F})$, and thus, there is a natural unit map $\mc{F}_s\to R\Psi_\mf{X}\Psi_\mf{X}^\ast\mc{F}_s=R\Psi_\mf{X}\mc{F}_\eta$; we show that this is an isomorphism. As this can be checked \'etale locally, we may assume that $\mc{F}= \underline{M}$ for a finite $\Lambda$-module $M$. Now, it suffices to check that this unit map is an isomorphism on the stalks over the algebraic closure $x$ of a closed point of $\mf{X}$ (see the argument in the proof of \cite[Proposition 2.2.3.2]{GaitsgoryLurie}). But, then we have an identification $(R\Psi_\mf{X}\mc{F}_\eta)_x\simeq R\Gamma(\mf{X}(x)_C,M)$. As $\mf{X}/\mc{O}_K$ is smooth we have $\mf{X}(x)_C\simeq \bb{D}^d_C$, the $d$-dimensional open polydisk over $C$. But, by our assumption on $\Lambda$, one has that $R\Gamma(\bb{D}^d_C,M)\simeq M[0]$ via the unit map, as desired.
\end{proof}

\subsubsection{Nearby cycles, correspondences, and tubes} We now wish to compare cohomologies of correspondences between the generic and special fibers of a formal scheme using nearby cycles.

\begin{setup}\label{setup:fml-sch-corr}Let $c\colon \mf{X}_1\xleftarrow{c_1}\mf{X}_0\xrightarrow{c_2}\mf{X}_2$ be a finite \'etale correspondence of finite type flat formal $\mc{O}_K$-schemes, and notate its special and generic fibers as follows:
\begin{equation*}
    c_s\colon \qquad \mf{X}_{1s}\xleftarrow{c_{1s}}\mf{X}_{0s}\xrightarrow{c_{2s}}\mf{X}_{2s},\qquad c_\eta\colon\quad  \mf{X}_{1\eta}\xleftarrow{c_{1\eta}}\mf{X}_{0\eta}\xrightarrow{c_{2\eta}}\mf{X}_{2\eta}.
\end{equation*}
Furthermore, let $u\colon c_{2\eta!}c_{1\eta}^\ast\mc{F}_1\to\mc{F}_2$ be a $c_\eta$-cohomological correspondence with $\mc{F}_i$ an object of $D_{\mr{zc},\mr{ftd}}(\mf{X}_{i\eta},\Lambda)$. Finally, let $\tau$ be an element of $W_K$ with $\mathsf{v}(\tau)=n\geqslant 0$.
\end{setup}

\begin{construction}
Define the $c_s^{(n)}$-cohomological correspondence $R\Psi(u^\tau)$
by the composite
\begin{equation*}
 c_{2,\ov{s}!}(c_{1s}^{(n)})^\ast R\Psi(\mc{F}_1)_{\ov{s}} \xrightarrow{\sim}
 c_{2,\ov{s}!}c_{1,\ov{s}}^\ast R\Psi(\mc{F}_1)_{\ov{s}} \xrightarrow{\sim}
 R\Psi(c_{2\eta!}c_{1\eta}^\ast\mc{F}_1)_{\ov{s}} \xrightarrow{R\Psi(u)_{\ov{s}}}
 R\Psi(\mc{F}_2)_{\ov{s}},
\end{equation*}
The first arrow uses Remark~\ref{rem:frobenius-galois-pullback}
and the $\tau$-equivariant structure on nearby cycles;
the second uses compatibility with \'etale pullback and exceptional pushforward
as in Proposition~\ref{prop:nearby-omnibus}.
\end{construction}

From Proposition \ref{prop:nearby-omnibus} the following proposition is clear.

\begin{prop}\label{prop:nearby-cycles-of-correspondence} There is an identification of complexes of $\Lambda$-modules between
\begin{equation*}
    \left[R\Gamma_c(u^\tau)\colon R\Gamma_c(\mf{X}_{1C},\mc{F}_{1C})\to R\Gamma_c(\mf{X}_{2C},\mc{F}_{2C})\right]
\end{equation*}
and
\begin{equation*}\left[R\Gamma_c(R\Psi(u^\tau))\colon R\Gamma_c(\mf{X}_{1\ov{s}},R\Psi(\mc{F}_{1})_{\ov{s}})\to R\Gamma_c(\mf{X}_{2\ov{s}},R\Psi(\mc{F}_{2})_{\ov{s}})\right].
\end{equation*}
\end{prop}

\begin{setup}\label{setup:fml-sch-corr-2} We continue as in Setup \ref{setup:fml-sch-corr}, except we now require that $\mf{X}_1=\mf{X}_2$ and $\mc{F}_1=\mc{F}_2$. We denote these common objects by $\mf{X}$ and $\mc{F}$, respectively.
\end{setup}

\begin{defn}Fix a point $y$ of $\mr{Fix}(c_s^{(n)})(\ov{k})$. Set 
\begin{equation*}
    x=\Phi^n_{\mf{X}_{s}}(c_{1s}(y))=c_{2s}(y).
\end{equation*}
We may then restrict $c_\eta^\tau$ to get an analytic correspondence on the Berthelot tubes:
\begin{equation*}
    c^\tau(y)\colon \quad \mf{X}_{C}(x)\xleftarrow{c^\tau(y)_1}\mf{X}_{0C}(y)\xrightarrow{c^\tau(y)_2}\mf{X}_C(x),
\end{equation*}
and we may restrict $u^\tau$ to a $c^\tau(y)$-cohomological correspondence $u^\tau(y)$. 
\end{defn}

\begin{prop}[{cf.\@ \cite[Theorem 3.5.8]{HuberEC} and \cite[Th\'eor\`eme 5.10.1]{Fargues}}]\label{prop:stalks-of-nearby-cycles} There is an identification of complexes of $\Lambda$-modules
\begin{equation*}
    \bigg[R\Psi(u^\tau)^\dashv_y\colon R\Psi(\mc{F})_x\to R\Psi(\mc{F})_x\bigg]\simeq\bigg[R\Gamma(u^\tau(y))\colon R\Gamma(\mf{X}_C(x),\mc{F}_C)\to R\Gamma(\mf{X}_C(x),\mc{F}_C)\bigg].
\end{equation*}
\end{prop}

\subsubsection{A rigid-analytic trace formula} We are now ready to state the aforementioned rigid-analytic trace formula. We maintain the setup of Setup \ref{setup:fml-sch-corr-2}.

We first have the following general version of a rigid-analytic analogue of Theorem \ref{thm:Fujiwara--Varshavsky-trace-formula}, which follows immediately by combining Theorem \ref{thm:Fujiwara--Varshavsky-trace-formula}, Proposition \ref{prop:nearby-cycles-of-correspondence}, and Proposition \ref{prop:stalks-of-nearby-cycles}.

\begin{prop}\label{prop:general-rigid-analytic-trace-formula}There exists an $m\geqslant 1$ such that if $q^n>m$ then we have the equality
\begin{equation*}
    \tr\left(\tau\times R\Gamma_c(\ov{u})\right)=\sum_{y\in \mr{Fix}(c_{s}^{(n)})(\ov{k})}\tr\left(R\Gamma(u^\tau(y))\right).
\end{equation*}
If $c_s$ admits a finite compactification $\widetilde c_s$,
then $m$ may be taken to be the ramification index of
$\widetilde c_{2s}$.
\end{prop}

For the applications to Shimura varieties we specialize to an even finer situation.

\begin{setup}\label{setup:bad-level-corr} We retain the setup from Setup \ref{setup:fml-sch-corr-2} but with several changes:

\begin{enumerate}
\item Assume that $\mf{X}$ is a smooth formal $\mc{O}_K$-scheme. 
\item Assume that the $c_\eta$-cohomological correspondence is of the form $u_\eta$ for a $c$-cohomological correspondence $u$, for some locally constant $\Lambda$-module $\mc{F}$ on $\mf{X}$.

\end{enumerate}
Additionally, we let $d\colon Y_0\rightrightarrows Y$ be a finite \'etale analytic correspondence, and $f\colon d\to c_\eta$ a finite \'etale morphism of analytic correspondences with $f_1=f_2$. Finally: 
\begin{enumerate} 
\item[(3)] Assume that we have a $d$-cohomological correspondence $v\colon d_{2!}d_1^\ast\mc{G}\to \mc{G}$, with a fixed isomorphism $v\simeq f^\ast(u_\eta)$.
\end{enumerate}
\end{setup}

\begin{rem}\label{tau-twisted-generic-corr}From Setup \ref{setup:bad-level-corr} we get a correspondence
\begin{equation*}
    d^\tau(y)\colon f_0^{-1}(\mf{X}_{0C}(y))\rightrightarrows f_1^{-1}(\mf{X}_{C}(x)),
\end{equation*}
with the $d^\tau(y)$-cohomological correspondence $\tr_{d^\tau(y)}$. Evidently $f_\ast(\tr_{d^\tau(y)})$ is equal to $f_\ast(\tr_d)^\tau(y)$.
\end{rem}

\begin{prop}\label{prop:specific-rigid-analytic-trace-formula} There exists an $m\geqslant 1$ such that if $q^n>m$ then we have the equality
\begin{equation*}
    \tr\left(\tau\times R\Gamma_c(\ov{v})\right)=\sum_{y\in \mr{Fix}(c_{s}^{(n)})}\tr\left((u_s^{(n)})_y^\dashv\right)\tr\left( R\Gamma(\tr_{d^\tau(y)})\right).
\end{equation*}
If $\widetilde{c}_s$ is a finite compactification of $c_s$, then $m$ may be taken to be the ramification index of $\widetilde{c}_{2s}$.
\end{prop}

\begin{lem}\label{lem:proj-formula-corr-nearby} There is a canonical identification of complexes of $\Lambda$-modules between
\begin{equation*}
    \left[R\Psi_\mf{X}(u_\eta^\tau\otimes f_\ast(\tr_d)^\tau)\colon c_{2,\ov{s}!}(c_{1,s}^{(n)})^\ast R\Psi_\mf{X}(\mc{F}_\eta\otimes f_\ast(\Lambda))\to R\Psi_\mf{X}(\mc{F}_\eta\otimes f_\ast(\Lambda))\right],
\end{equation*}
and
\begin{equation*}
\left[u_s^{(n)}\otimes R\Psi_\mf{X}(f_\ast(\tr_d)^\tau)\colon c_{2,\ov{s}!}(c_{1,s}^{(n)})^\ast\left(\mc{F}_s\otimes_\Lambda R\Psi_\mf{X}(f_\ast(\Lambda))\right)\to \mc{F}_s\otimes R\Psi_\mf{X}(f_\ast(\Lambda))\right].
\end{equation*}
\end{lem}
\begin{proof} More precisely, we claim that the following diagram is commutative 
\begin{equation*}
    \begin{tikzcd}[sep=large]
	{c_{2,\ov{s}!}(c_{1,s}^{(n)})^\ast\left(\mc{F}_s\otimes_\Lambda R\Psi_\mf{X}(f_\ast(\Lambda))\right)} & {c_{2,\ov{s}!}(c_{1,s}^{(n)})^\ast R\Psi_\mf{X}(\mc{F}_\eta\otimes_\Lambda f_\ast(\Lambda))} \\
	{\mc{F}_s\otimes_\Lambda R\Psi_\mf{X}(f_\ast(\Lambda))} & {R\Psi(\mc{F}_\eta\otimes_\Lambda f_\ast(\Lambda)),}
	\arrow["\sim", from=1-1, to=1-2]
	\arrow["{u_s^{(n)}\otimes R\Psi_\mf{X}(f_\ast(\tr_d)^\tau)}"', from=1-1, to=2-1]
	\arrow["\sim"', from=2-1, to=2-2]
	\arrow["{R\Psi(u_\eta^\tau\otimes f_\ast(\tr_d)^\tau)}", from=1-2, to=2-2]
\end{tikzcd}
\end{equation*}
where the horizontal isomorphisms are as in (5) of Proposition \ref{prop:nearby-omnibus}. This follows from the naturality of the isomorphism in loc.\@ cit.\@\end{proof}

\begin{proof}[Proof of Proposition \ref{prop:specific-rigid-analytic-trace-formula}] By Proposition \ref {prop:analytic-pushforward} the left-hand side is equal to $\tr\left(R\Gamma_c(u_\eta^\tau\otimes f_\ast(\tr_d)^\tau)\right)$. But, by Proposition \ref{prop:nearby-cycles-of-correspondence} this is equal to $\tr\left(R\Gamma_c(R\Psi_\mf{X}(u_\eta^\tau \otimes f_\ast(\tr_d)^\tau)\right)$. By Lemma \ref{lem:proj-formula-corr-nearby} this is, in turn, equal to $\tr\left(R\Gamma_c\left( u_s^{(n)}\otimes R\Psi_\mf{X}(f_\ast(\tr_d)^\tau)\right)\right)$. If $q^n>m$ then by Theorem \ref{thm:Fujiwara--Varshavsky-trace-formula} this then is equal to the following sum
\begin{equation*}
    \sum_{y\in\mr{Fix}(c_s^{(n)})(\ov{k})}\tr\left(\left(u_s^{(n)}\otimes R\Psi_\mf{X}(f_\ast(\tr_d)^\tau)\right)_y^\dashv\right).
\end{equation*}
But, 
\begin{equation*}
    \left(u_s^{(n)}\otimes R\Psi_\mf{X}(f_\ast(\tr_d)^\tau)\right)_y^\dashv\simeq (u^{(n)}_s)_y^\dashv \otimes \left(R\Psi_\mf{X}(f_\ast(\tr_d)^\tau)\right)_y^\dashv,
\end{equation*}
and so this term is equal to
\begin{equation*}
    \sum_{y\in\mr{Fix}(c_s^{(n)})(\ov{k})}\tr\left((u^{(n)}_s)_y^\dashv\right) \tr\left(R\Psi_\mf{X}(f_\ast(\tr_d)^\tau)\right)_y^\dashv.
\end{equation*}
But, by Proposition \ref{prop:stalks-of-nearby-cycles} we have that 
\begin{equation*}
    \tr\left(R\Psi_\mf{X}(f_\ast(\tr_d)^\tau)\right)_y^\dashv=\tr\left(R\Gamma(f_\ast(\tr_d)^\tau(y))\right)=\tr\left(R\Gamma(f_\ast(\tr_{d^\tau(y)})\right).
\end{equation*}
We are then done by Proposition \ref{prop:fin-morphism-pushforward-analytic} which implies this last term is equal to $\tr\left(R\Gamma(\tr_{d^\tau(y)})\right)$.
\end{proof}

\section{\texorpdfstring{$(\mathcal{G},\mu)$}{(G,mu)}-apertures and crystalline theory}\label{s:apertures-and-crystalline-theory}

Key to having a robust generalization of the deformation spaces from \cite{ScholzeLK}, especially with no restrictions on the unramified group involved, is the recently developed theory of \emph{$(\mc{G},\mu)$-apertures} as in \cite{DrinfeldShimurian,GMM}. In this section we very briefly recall the theory of $(\mc{G},\mu)$-apertures, compare this theory to $F$-crystals with $\mc{G}$ structure, and describe their deformation theory. 

\begin{rem} We will take for granted that the reader is familiar with the basic setup of the stacks $\mf{X}^\smallprism$, $\mf{X}^\smallN$, $\mf{X}^\mr{syn}$, and the associated notions of prismatic $F$-crystals and prismatic $F$-gauges; we refer the interested reader to \cite[\S2-3]{MY}, and the references therein, for further details. We will furthermore be using the theory of filtered torsors freely throughout; see \S\ref{sec:filtered-Tannakian-appendix} below for a quick overview of the theory.
\end{rem}

\begin{nota}\label{nota:deformation-spaces-setup} We fix the following notation:

\begin{multicols}{2}
\begin{itemize}[leftmargin=.2cm]
\item $k$ is a perfect extension of $\bb{F}_p$;
\item $W=W(k)$ is the ring of Witt vectors; 
\item $\sigma$ is the natural Frobenius lift on $W$;
\item $K=W[\nicefrac{1}{p}]$,
\item $\mc{G}$ is a reductive group $\mathbb{Z}_p$-scheme;
\item $\mf{g}$ is the Lie algebra of $\mc{G}$;
\item $G$ is the generic fiber of $\mc{G}$;
\item $\mu\colon \bb{G}_{m,W}\to \mc{G}_W$ is a minuscule cocharacter; 
\item $\sigma^j(\mu)$ is the cocharacter of $\mc{G}$ obtained as $(\sigma^j)^\ast \mu$ using the identification $(\sigma^j)^\ast\mc{G}_W\simeq \mc{G}_W$.\footnote{For example, one has $\sigma(\mu)(p)=\sigma(\mu(p))$.}
\end{itemize}
\end{multicols}
Finally, we write 
\begin{equation*}
P_\mu=\left\{g\in\mc{G}_W: \lim_{t\to0}\mu(t)g\mu(t)^{-1}\text{ exists}\right\}
\end{equation*}
for the parabolic subgroup associated to $\mu$; see \cite[Theorem 4.1.7]{ConradReductive}. We refer the reader to Remark \ref{rem:sign-conventions} for our sign conventions (e.g., $\mu$ vs. $\mu^{-1}$) and how they compare to other sources.
\end{nota}

\subsection{Definition and general theory} We first recall the feature of $\mf{X}^\mr{syn}$ that will allow us to formulate the correct analogue of being of `relative position $\mu$' from the theory of shtukas. 
\begin{defn}\label{defn:syntomic-de-Rham-point} Let $\mf{X}$ be a derived $p$-adic formal scheme. The \emph{mod $p^n$ filtered de Rham point} 
\begin{equation*}
x_{\mr{dR},n}^\mr{syn}\colon (\mf{X}\dotimes \bb{Z}/p^n)\times \bb{A}^1/\bb{G}_m\to \mf{X}^\syn\dotimes\bb{Z}/p^n,
\end{equation*}
is obtained as the composition of the map $x^\smallN_\mr{dR}\colon \mf{X}\times \bb{A}^1/\bb{G}_m\to \mf{X}^\smallN$ from \cite[\S6.7]{GMM} and the projection map $\mf{X}^\smallN\to \mf{X}^\syn$. When $n=\infty$ we omit $n$ from the notation/terminology.
\end{defn}

\begin{defn}\label{defn:BTGmu} Let $\mf{X}$ be a derived $p$-adic formal scheme and $n$ in $\bb{N}\cup\{\infty\}$.
\begin{enumerate}
\item For a $\mc{G}$-bundle $\mf{Q}$ on $\mf{X}^\mr{syn}\dotimes \bb{Z}/p^n$ we define the \emph{Hodge-filtered de Rham realization}, denoted $\Fil^\bullet_\mr{Hdg} T_\mr{dR}(\mf{Q})$, to be the filtered $\mc{G}$-bundle on $\mf{X}\dotimes\bb{Z}/p^n$ given by $(x^\mr{syn}_{\mr{dR},n})^\ast\mf{Q}$.
\item A $\mc{G}$-bundle $\mf{Q}$ on $\mf{X}^\mr{syn}\dotimes\bb{Z}/p^n$ is of \emph{type $\mu$} if $\Fil^\bullet_\mr{Hdg} T_\mr{dR}(\mf{Q})$ is of type $\mu$.
\item An \emph{$n$-truncated $(\mc{G},\mu)$-aperture} on $\mf{X}$ (or just a \emph{$(\mc{G},\mu)$-aperture} if $n=\infty$) is a $\mc{G}$-bundle $\mf{Q}$ on $\mf{X}^\mr{syn}\dotimes \bb{Z}/p^n$ of type $\mu$. 
\end{enumerate}
The $\infty$-groupoid of $n$-truncated $(\mc{G},\mu)$-apertures $\mf{Q}$ over $\mf{X}$ is denoted by $\mr{BT}^{\mc{G},\mu}_n(\mf{X})$.
\end{defn}

\begin{rem} When $n<\infty$, the stack $\mf{X}\dotimes \bb{Z}/p^n$ need not be classical and so we are implicitly using Remark \ref{rem:derived-prestack} in the above definition. But, one can avoid all derived geometry using \cite[Proposition 5.5.2]{GMM} which implies that $\mf{Q}$ is of type $\mu$ if and only if $\Fil^\bullet_\mr{Hdg} T_\mr{dR}(\mf{Q})_\kappa$ is of type $\mu$ for every geometric point $\Spec(\kappa)\to \mf{X}$ (or even just one if $\mf{X}$ is connected).
\end{rem}

\begin{rem}[Independence of topology]\label{rem:bounded-deformation-space-BT}\label{rem:integrability} Suppose that $R$ is a $W$-algebra complete with respect to the $J$-adic topology for some ideal $J\subseteq R$ such that $p^k$ is in $J$ for some $k$, i.e., such that we have a (necessarily unique) map $\Spf(R,J)\to \Spf(\bb{Z}_p)$. Then, there is potential ambiguity about the notation $\mr{BT}^{\mc{G},\mu}_\infty(R)$: Do we mean maps $\Spf(R,J)\to \mr{BT}^{\mc{G},\mu}_\infty$ or maps $\Spf(R,(p))\to \mr{BT}^{\mc{G},\infty}$?

But, this ambiguity is illusory by the integrability of each $\mathrm{BT}_n^{\mathcal{G},\mu}$ (i.e., condition (3) in \cite[Theorem 7.1.6]{DAG}) for finite $n$. Indeed, for finite $n$ this integrability gives natural equivalences
\begin{equation*}
\begin{aligned} \mr{BT}^{\mc{G},\mu}_n\big(\Spf(R,(p))\big) &=\varprojlim_k \mr{BT}^{\mc{G},\mu}_n\big(\Spec(R/p^k)\big)\\ &=\varprojlim_r\varprojlim_k \mr{BT}^{\mc{G},\mu}_n\big(\Spec(R/(p^k,J^r)\big)\\ &= \mr{BT}^{\mc{G},\mu}_n\big(\Spf(R,J)\big),
\end{aligned}
\end{equation*}
from which the same equality follows for $n=\infty$ by passing to the limit. Thus, for an adic $W$-algebra $R$ we may unambiguously write $\mr{BT}^{\mc{G},\mu}_\infty(R)$.
\end{rem}

Concerning the geometry of $\mr{BT}^{\mc{G},\mu}_n$, we have the following omnibus theorem of Gardner--Madapusi. In particular, it shows that the derived prestack $\mr{BT}^{\mc{G},\mu}_n$ is $1$-truncated on discrete inputs, i.e., is just a usual stack.

\begin{thm}[{\cite[Theorems D and G]{GMM}}]\label{thm:GMM}Fix $n$ in $\bb{N}\cup\{\infty\}$.
\begin{enumerate} 
\item If $n$ is finite, the derived prestack $\mr{BT}^{\mc{G},\mu}_n$ is a quasi-compact smooth $p$-adic formal Artin stack over $\Spf(W)$ of dimension $0$ with affine diagonal. 
\item If $n$ is finite, the natural truncation map $\mr{BT}^{\mc{G},\mu}_{n+1}\to \mr{BT}^{\mc{G},\mu}_n$ is smooth and surjective.
\item{\emph{[Grothendieck--Messing theory]}} For $(R'\to R,\gamma)$ a \emph{(}pro-nilpotent\emph{)} divided power thickening of $p$-complete animated $W$-algebras, there is a Cartesian diagram of groupoids
\begin{equation}\label{eq:GR-1}
    \begin{tikzcd}
	{\mr{BT}^{\mc{G},\mu}_n(R')} & {BP_{\mu}(R'\dquot p^n)} \\
	{\mr{BT}^{\mc{G},\mu}_n(R)} & B\mc{G}(R'\dquot p^n)\times_{B\mc{G}(R\dquot p^n)} BP_{\mu}(R\dquot p^n).
	\arrow[from=1-1, to=1-2]
	\arrow[from=1-1, to=2-1]
	\arrow[from=1-2, to=2-2]
	\arrow[from=2-1, to=2-2]
\end{tikzcd}
\end{equation}
\end{enumerate}
\end{thm}

\begin{rem}\label{rem:notation-for-GM} For the convenience of the reader, we explicate some of the maps appearing in the statement of Grothendieck--Messing theory:
\begin{itemize}
\item The natural morphism $\mr{BT}^{\mc{G},\mu}_n(R)\to BP_{\mu}(R\dquot p^n)$ is the one classifying the Hodge-filtered de Rham realization using the equivalence from Proposition \ref{prop:type-nu-classification}.
\item The map $\mr{BT}^{\mc{G},\mu}_n(R)\to B\mc{G}(R'\dquot p^n)$ is given by pullback along the map $\Spec(R')\to R^\syn$ described in \cite[\S6.8]{GMM}.
\end{itemize}

\end{rem}

We now discuss the important fact that $(\mc{G},\mu)$-apertures on $\mf{X}$ give rise to local systems on the generic fiber $\mf{X}_\eta$. We then use this to define an algebraization of $\mr{BT}^{\mc{G},\mu}_n$ (for $n$ finite).

\begin{nota} For $1\leqslant n\leqslant \infty$ we denote by $B\mc{G}(\bb{Z}/p^n)$ (or $B\mc{G}(\bb{Z}_p)$ when $n=\infty$) the classifying stack for the (pro)-finite group $\mc{G}(\bb{Z}/p^n)$ on the pro-\'etale site of a (formal) scheme or adic space.
\end{nota}

\begin{rem} When $\mc{G}=\mr{GL}_{n,\bb{Z}_p}$ one has that $B\mc{G}(\bb{Z}/p^n)$ is the moduli stack of $\bb{Z}/p^n$-local systems of rank $n$. More generally, using the material from \S\ref{sec:filtered-Tannakian-appendix} one may interpret $B\mc{G}(\bb{Z}/p^n)(X)$ as the category of $\mc{G}$-objects $\cat{Rep}_{\bb{Z}_p}(\mc{G})\to \cat{Loc}_{\bb{Z}/p^n}(X)$ of the category $\cat{Loc}_{\bb{Z}/p^n}(X)$ of $\bb{Z}/p^n$-local systems on $X$; see \cite[Proposition 2.3 and Remark 2.8]{IKY1}. 
\end{rem} 

\begin{nota}\label{nota:etale-realization} For a $p$-adic formal $W$-scheme $\mf{X}$ and an integer $1\leqslant n\leqslant \infty$, we denote by 
\begin{equation}\label{eq:etale-realization}
T_\et\colon \mr{BT}^{\mc{G},\mu}_n(\mf{X})\to B\mc{G}(\bb{Z}/p^n)(\mf{X}_\eta),
\end{equation}
the \'etale realization functor as in \cite[\S3.8]{MY}.
\end{nota}

We then obtain an algebraization of $\mr{BT}^{\mc{G},\mu}_n$ by gluing it to $B\mc{G}(\bb{Z}/p^n)$.

\begin{defn}[{see \cite[Construction 4.1.1]{MY}}]\label{defn:algebraic-apertures} For an animated $W$-scheme $Y$, we define 
\begin{equation*}
\mr{BT}^{\mc{G},\mu,\mr{alg}}_n(Y)\defeq B\mc{G}(\bb{Z}/p^n)(Y_\eta)\times_{B\mc{G}(\bb{Z}/p^n)(\wh{Y}_\eta)}\mr{BT}^{\mc{G},\mu}_n(\wh{Y}),
\end{equation*}
where the left-structure map is via pullback and the right-structure map is $T_\et$.
\end{defn}

\begin{rem} The notion of gluing here is meant roughly in the Beauville--Laszlo sense as in \cite{AchingerYoucis}. That said, op.\@ cit.\@ applies only to separated algebraic spaces, and it would be interesting to see a generalization applied to more directly provide algebraic structure on $\mr{BT}^{\mc{G},\mu,\mr{alg}}_n$.
\end{rem}

We then have the following omnibus result.

\begin{thm}[{\cite[Theorem 4.1.3. and Lemma 4.1.10]{MY}}]Fix $n$ in $\bb{N}\cup\{\infty\}$.
\begin{enumerate} 
\item If $n$ is finite, the prestack $\mr{BT}^{\mc{G},\mu,\mr{alg}}_n$ is a quasi-compact smooth Artin stack over $\Spec(W)$ of dimension $0$ with affine diagonal. 
\item For finite $n$, the natural map $\mr{BT}^{\mc{G},\mu,\mr{alg}}_{n+1}\to \mr{BT}^{\mc{G},\mu,\mr{alg}}_n$ is smooth and surjective.
\item{\emph{[Grothendieck--Messing theory]}} For $(R'\to R,\gamma)$ a \emph{(}pro-nilpotent\emph{)} divided power thickening of animated $W$-algebras, there is a Cartesian diagram of groupoids
\begin{equation}\label{eq:alg-GR-1}
    \begin{tikzcd}
	{\mr{BT}^{\mc{G},\mu,\mr{alg}}_n(R')} & {BP_{\mu}(R'\dquot  p^n)} \\
	{\mr{BT}^{\mc{G},\mu,\mr{alg}}_n(R)} & B\mc{G}(R'\dquot p^n)\times_{B\mc{G}(R\dquot p^n)} BP_{\mu}(R\dquot  p^n).
	\arrow[from=1-1, to=1-2]
	\arrow[from=1-1, to=2-1]
	\arrow[from=1-2, to=2-2]
	\arrow[from=2-1, to=2-2]
\end{tikzcd}
\end{equation}
\end{enumerate}
\end{thm}

\begin{rem} Suppose that $R$ is a $p$-completed $W$-algebra and that $1\leqslant n\leqslant \infty$. Then, we can consider $\mr{BT}^{\mc{G},\mu,\mr{alg}}_n(R):=\mr{BT}^{\mc{G},\mu,\mr{alg}}_n(\Spec(R))$ or $\mr{BT}^{\mc{G},\mu}_n(R)=\mr{BT}^{\mc{G},\mu}_n(\Spf(R))$. These are, in fact, the same. Indeed, by definition we have an identification of groupoids
\begin{equation*}
\mr{BT}^{\mc{G},\mu,\mr{alg}}_n(R)=B\mc{G}(\bb{Z}/p^n)(\Spec(R[\nicefrac{1}{p}]))\times_{B\mc{G}(\bb{Z}/p^n)(\Spa(R[\nicefrac{1}{p}]))}\mr{BT}^{\mc{G},\mu}_n(R).
\end{equation*}
Thus, the claim will follow if we know that the natural map 
\begin{equation*}
B\mc{G}(\bb{Z}/p^n)(\Spec(R[\nicefrac{1}{p}]))\to B\mc{G}(\bb{Z}/p^n)(\Spa(R[\nicefrac{1}{p}])),
\end{equation*}
is an equivalence. But, using the Tannakian perspective on $\mc{G}(\bb{Z}/p^n)$-local systems, and working component-by-component, we are reduced to knowing that $\pi_1^\et(\Spa(R[\nicefrac{1}{p}]))\to \pi_1^\et(\Spec(R[\nicefrac{1}{p}]))$ is an isomorphism. But, this result is standard, e.g., see \cite[Example 1.6.6]{HuberEC}.
\end{rem}

\subsection{$(\mc{G},\mu)$-apertures and  $F$-crystals} We now explain that over nice characteristic $p$ schemes, the theory of $(\mc{G},\mu)$-apertures coincides with the more classical (Tannakian) theory of $F$-crystals.

\begin{nota}We continue with the notation from Notation \ref{nota:deformation-spaces-setup}, but further assume $p>2$ and fix $S$ to be a regular quasi-compact $k$-scheme.\footnote{In what follows, the quasi-compactness is easily removed, but the regularity is essential (see Remark \ref{rem:why-regular}).} As $S$ is regular and quasi-compact, there exists a flat (in fact quasi-syntomic) cover $\Spec(R)\to S$ with $R$ perfect, which we fix.\footnote{In fact, they are equivalent conditions; see \cite[Theorem 2.1]{Kunz}.}
\end{nota}

\begin{nota} Write $S_\mr{crys}$ for the absolute (big fppf) crystalline site of $S$, endowed with the usual structure sheaf, and $F_S$ for the Frobenius endomorphism on the associated ringed topos.
\end{nota}

\begin{rem} For other notation and terminology concerning this site and its attendant objects (e.g., (iso)crystals) see \cite[\S2.3.1]{IKY1}. In particular, we shall often heavily abbreviate notation for objects of $S_\mr{crys}$ when the meaning is unambiguous. For example, for the object of $\Spec(k)_\mr{crys}$ that would precisely be written as the diagram $W(k)\twoheadrightarrow k\xleftarrow{\mr{id}}k$ we may just write $W(k)$.

Moreover, for any object $T$ of $S_\mr{crys}$, and a sheaf $\mc{F}$ on $S_\mr{crys}$, we shall write $\mc{F}_T$ for the restriction of $\mc{F}$ to the Zariski site of $T$ considered as a subsite of $S_\mr{crys}$ in the obvious way.
\end{rem}

\subsubsection{$F$-crystals with $\mc{G}$-structure of type $\mu$} We begin by describing the crystalline objects that we shall ultimately show recover the theory of $(\mc{G},\mu)$-apertures.

\begin{defn} An \emph{$F$-crystal \emph{(}in vector bundles\emph{)}} on $S$ is a pair $(\mc{E},\varphi)$ where $\mc{E}$ is a locally finite free crystal on $S$ and $\varphi$ is an isomorphism of isocrystals
\begin{equation*}
\varphi \colon F_S^\ast\mc{E}[\nicefrac{1}{p}]\isomto\mc{E}[\nicefrac{1}{p}].
\end{equation*}
With the obvious notion of morphism, we denote by $\mb{FCrys}(S)$ the category of $F$-crystals on $S$.
\end{defn}

\begin{obs} Evidently the category $\mb{FCrys}(S)$ carries the natural structure of an exact $\bb{Z}_p$-linear $\otimes$-category in the sense of \cite[Appendix A]{IKY1}, and thus it makes sense to consider the category $\mb{FCrys}_\mc{G}(S)$ of $F$-crystals on $S$ with $\mc{G}$-structure as in loc.\@ cit. Such objects are of the form $(\mc{Q},\varphi)$ where $\mc{Q}$ is a $\mc{G}$-torsor on $S_\crys$ and $\varphi\colon F_S^\ast\mc{Q}[\nicefrac{1}{p}]\isomto\mc{Q}[\nicefrac{1}{p}]$ is an isomorphism.
\end{obs}

\begin{rem}\label{rem:crys-prism} There is a bi-exact $\bb{Z}_p$-linear $\otimes$-equivalence
\begin{equation*}
(-)^\mr{crys}\colon \cat{Vect}^\varphi(S_\smallprism)\isomto \mb{FCrys}(S)
\end{equation*}
where the source is the category of prismatic $F$-crystals on $S$, even if $S$ is just quasi-syntomic; see \cite[\S1.1.1]{IKY3}. We denote the natural quasi-inverse of this equivalence by $(-)^\smallprism$. 

In fact, one can understand this equivalence as coming from an isomorphism of stacks over $W$:
\begin{equation}\label{eq:isom-crys-prism-stacks}
(S/W)^\mr{crys}\simeq \sigma_W^\ast(S^\smallprism);
\end{equation}
see \cite[Lemma 1.6]{IKY3} for this fact and the meaning of this notation. That said, ultimately this equivalence can be understood from something quite simple: qrsp objects $\Spec(T)$ form a basis of $S_{\mr{qsyn}}$, and $\mr{A}_\mr{crys}(T)$ is both an initial object of $T_\smallprism$ and $T_\mr{crys}$. 

As this equivalence is $\bb{Z}_p$-linear, monoidal, and exact, we deduce an equivalence of associated Tannakian categories, i.e., $(-)^\crys\colon \cat{Tors}^\varphi_\mc{G}(S_\smallprism)\simeq \cat{FCrys}_\mc{G}(S)$ with quasi-inverse $(-)^\smallprism$; see \cite[Definition 1.28]{IKY1} for the source category. We shall often abusively use this equivalence below without comment when emphasizing it is not important.
\end{rem}

The following proposition establishes that, under the equivalence $\cat{Tors}^\varphi_\mc{G}(S_\smallprism)\simeq \cat{FCrys}_\mc{G}(S)$ from Remark \ref{rem:crys-prism}, various notions of `bounded by $\nu$' coincide, at least up to Frobenius twist.

\begin{prop}\label{prop:type-mu-equiv-char-p} Let $(\mc{Q},\varphi)$ be an object of $\mb{FCrys}_\mc{G}(S)$, and let $\Spec(R)\to S$ be a flat cover where $R$ is perfect. Then, the following are equivalent for a minuscule cocharacter $\nu$:
\begin{enumerate}
\item For every geometric point $\Spec(\kappa)\to S$ there exists a trivialization of $\mc{Q}|_{W(\kappa)}$ which transports the Frobenius $\varphi|_{W(\kappa)}$ to an element of $\mc{G}(W(\kappa))\nu(p)\mc{G}(W(\kappa))$.
\item There exists an \'etale cover $R\to R'$ and a trivialization of $\mc{Q}|_{W(R')}$ which transports the Frobenius $\varphi|_{W(R')}$ to an element of $\mc{G}(W(R'))\nu(p)\mc{G}(W(R'))$.
\item There exists a $p$-adically flat cover $(W(R),(p))\to (A,(p))$ in $S_\smallprism$ such that $\mc{Q}^\smallprism|_{\Spec(A)}$ admits a trivialization which transports the Frobenius $\varphi|_{A}$ to an element of $\mc{G}(A)\sigma^{-1}(\nu)(p)\mc{G}(A)$.
\end{enumerate}
Assume further that $S$ is a smooth $k$-scheme. Then, the above are equivalent to:
\begin{itemize}
\item[(4)] For every \'etale cover of $S$ of the form $\mf{X}'_k\to S$ where $\mf{X}'/W$ is a smooth formal scheme equipped with a Frobenius lift $F_{\mf{X}'}$, there exists an \'etale cover $\Spf(A)=\mf{X}''\to \mf{X}'$ such that $\mc{Q}_{\mf{X}''}$ has a trivialization which transports the Frobenius $\varphi$ to an element of $\mc{G}(A)\nu(p)\mc{G}(A)$.
\end{itemize}
\end{prop}
\begin{proof} For (1) implies (2), observe that there exists some \'etale cover $R\to R'$ such that $\mc{Q}_{\Spec(W(R'))}$ is trivializable.\footnote{Indeed, this follows as $\mc{G}$ is smooth and affine as $H^1_\et(\Spec(R),\mc{G})\simeq H^1_\mr{fl}(\Spec(R),\mc{G})\simeq H^1_\mr{fl}(\Spec(W(R)),\mc{G})$ for any perfect $R$; e.g., see \cite[Theorem 2.1.6]{CesnaviciusBouthier}.\label{footnote:Hensel}} Choosing a trivialization $\mc{Q}_{\Spec(W(R'))}\simeq \mc{G}$ allows one to view $\varphi$ as giving a map $\Spec(R')\to \mr{Gr}_\mc{G}$, where the target is the Witt vector affine Grassmannian as in \cite[\S1]{ZhuAffGmix}. But, as $\mr{Gr}_{\mc{G},\nu}=\mr{Gr}_{\mc{G},\leqslant \nu}$ (see \cite[Corollary 1.24]{ZhuAffGmix}), we see that (1) implies that the composition $\Spec(\kappa)\to \Spec(R')\to \mr{Gr}_\mc{G}$ factorizes through $\mr{Gr}_{\mc{G},\nu}$ for all geometric points $\kappa$ of $R'$. As $\mr{Gr}_{\mc{G},\nu}$ is a closed subfunctor and $\Spec(R')$ is reduced, we deduce that $\Spec(R')\to\mr{Gr}_\mc{G}$ factorizes through $\mr{Gr}_{\mc{G},\nu}$. By definition, this implies that after replacing $R'$ with a further \'etale cover $R''$, under the chosen trivializations the Frobenius defines an element of $\mc{G}(W(R''))\nu(p)\mc{G}(W(R''))$.

That (2) implies (3) is simple: Indeed, with notation as in (2), $(W(R),(p))\to (W(R'),(p))$ is a $p$-adically faithfully flat cover in $S_\smallprism$. The introduction of the inverse Frobenius twist comes from the Frobenius twist in the equivalence in Equation \eqref{eq:isom-crys-prism-stacks}.

For (3) implies (1), first assume that $\Spec(\kappa)\to S$ factorizes through $\Spec(A/p)\to S$. By the universal property of the Witt vectors we have a map of prisms $(A,(p))\to (W(\kappa),(p))$ and thus the coset condition on Frobenius holds for $W(\kappa)$ as desired. In general, there is an extension of algebraically closed fields $\kappa'/\kappa$, such that $\Spec(\kappa')\to \Spec(\kappa)\to S$ does lift to $\Spec(A/p)$. By the same argument in the first paragraph, we may define a map $\Spec(\kappa)\to \mr{Gr}_\mc{G}$ and, as the claim in (1) holds for $\kappa'$, we deduce that the composition $\Spec(\kappa')\to\Spec(\kappa)\to \mr{Gr}_\mc{G}$ factors through $\mr{Gr}_{\mc{G},\nu}$ and thus so must also the map $\Spec(\kappa)\to \mr{Gr}_\mc{G}$ implying the claim.

Finally, to see that the statements (2) and (4) are equivalent one may apply mutatis mutandis the argument given in \cite[Proposition 3.17]{IKY1}.
\end{proof}

With this we can, unambiguously relative to our various identifications of categories, define when an object of $\cat{FCrys}_\mc{G}(S)$ is `bounded by $\nu$'.

\begin{defn}\label{defn:crystalline-type} An object $\mc{Q}$ of $\mb{FCrys}_\mc{G}(S)$ satisfying any of the conditions from Proposition \ref{prop:type-mu-equiv-char-p} is called \emph{of type $\nu$}. We denote by $\cat{FCrys}_\mc{G}^{\nu}(S)$ the full subcategory of such objects.
\end{defn}

\subsubsection{The crystalline-aperture equivalence} We can now precisely describe our claim that $(\mc{G},\mu)$-gauges on $S$ are the same as objects of $\cat{FCrys}_\mc{G}^{\sigma(\mu)}(S)$; the twist on the cocharacter here is (ultimately) a consequence of the twist in \eqref{eq:isom-crys-prism-stacks}; cf.\@  \cite[Proposition 1.1]{IKY3}.

\begin{constr} In \cite[Construction 1.21]{IKY2} there is constructed a functor
\begin{equation*}
\mr{R}_{\mf{X}}\colon \mr{BT}^{\mc{G},\mu}_\infty(\mf{X})\to \cat{Tors}^{\varphi,\mu}_\mc{G}(\mf{X}_\smallprism),
\end{equation*}
called the \emph{restriction functor}, for any quasi-syntomic $p$-adic formal scheme $\mf{X}$, with the target the category of prismatic $\mc{G}$-torsors with $F$-structure bounded by $\mu$ as in op.\@ cit. 
\end{constr}

\begin{rem}If $\mf{X}$ is now our regular $k$-scheme $S$, then from Proposition \ref{prop:type-mu-equiv-char-p} we know that the equivalence from Remark \ref{rem:crys-prism} restricts to an equivalence $\cat{Tors}^{\varphi,\mu}_\mc{G}(S_\smallprism)\simeq \cat{FCrys}^{\sigma(\mu)}_\mc{G}(S)$. Thus, we obtain
\begin{equation*}
\mr{R}_{S}\colon \mr{BT}^{\mc{G},\mu}_\infty(S)\to \cat{FCrys}_\mc{G}^{\sigma(\mu)}(S),
\end{equation*}
which we have notated abusively, and which we also call the restriction functor.
\end{rem}

\begin{prop}\label{prop:BT-in-char-p} Suppose that $S$ is a regular $k$-scheme. Then, the restriction functor 
\begin{equation*}
\mr{R}_S\colon \mr{BT}^{\mc{G},\mu}_\infty(S)\to \mb{FCrys}_\mc{G}^{\sigma(\mu)}(S),
\end{equation*}
is a natural equivalence.
\end{prop}
\begin{proof}This follows from \cite[Proposition 1.39]{IKY2}: the proof of (2) in loc.\@ cit.\@ only requires the existence of a quasi-syntomic cover by a perfectoid.
\end{proof}

\begin{rem}[Why regular?]\label{rem:why-regular} Proposition \ref{prop:BT-in-char-p} really requires that $S$ is regular. For general quasi-syntomic $k$-schemes $S$ one still has an equivalence $\cat{Tors}^{\varphi,\mu}_\mc{G}(S_\smallprism)\simeq \cat{FCrys}_\mc{G}^{\sigma(\mu)}(S)$ (see \cite[\S1.1.1]{IKY3}), and the restriction functor $\mr{R}_S$ is fully faithful (see \cite[Proposition 1.39]{IKY2}), but it can fail to be essentially surjective.
\end{rem}

\subsubsection{A crystalline description of the Hodge-filtered de Rham realization} It will be useful below to understand how to describe the Hodge-filtered $\mc{G}$-bundle $\Fil^\bullet_\mr{Hdg} T_\mr{dR}(\mf{Q})$ associated to a $(\mc{G},\mu)$-aperture $\mf{Q}$ on $S$ in terms of the equivalence from \ref{prop:BT-in-char-p}, at least when $S$ is `nice'.

\begin{defn}\label{defn:base-k-algebras} Let $A$ be a $k$-algebra with $\Spec(A)$ connected. We call $A$ a \emph{base $k$-algebra} if $A$ is $J_A$-adically complete for an ideal $J_A\subseteq A$ such that the pair $(A,J_A)=(A_n,I_n)$ is obtained from the following iterative procedure. For some $d\geqslant 0$, set
\begin{equation*}
A_0=\bT_{d,k}\defeq k[t_1^{\pm 1},\ldots,t_d^{\pm 1}],\qquad I_0=(0).
\end{equation*}
For each $i=0,\ldots,n-1$, form the pair $(A_{i+1},I_{i+1})$ by one of the following operations:
\begin{itemize}
\item $A_{i+1}$ is an \'etale $A_i$-algebra and $I_{i+1}=I_iA_{i+1}$;
\item $A_{i+1}$ is a localization $(A_i)_{\mf{p}}$ at a prime $\mf{p}\subseteq A_i$ and $I_{i+1}=I_iA_{i+1}$;
\item $A_{i+1}$ is the $I$-adic completion of $A_i$ for an ideal $I\subseteq A_i$ and $I_{i+1}=(I_i,I)A_{i+1}$.
\end{itemize}
A \emph{base $k$-scheme} is a $k$-scheme which admits an affine open cover by spectra of base $k$-algebras.
\end{defn}

\begin{eg} The two most natural examples of base $k$-schemes are:
\begin{itemize}
\item smooth $k$-schemes;
\item the $k$-scheme $\Spec(k\llbracket t_1,\ldots,t_d\rrbracket)$ for some $d\geqslant 0$.
\end{itemize}
\end{eg}

\begin{defn} Suppose that $A$ is a base $k$-algebra. If $J_A$ is an ideal arising from a presentation, then a $J_A$-formally \'etale map $w\colon \bT_{d,k}\to A$ is called a \emph{formal framing} of $A$.
\end{defn}

\begin{rem}\label{rem:base-algebras} In the following we shall freely use the notions of \emph{base $W$-algebras}, \emph{base formal $W$-schemes}, and \emph{formal framings} as in \cite[\S1.1.5]{IKY1}; these are straightforward analogues of Definition \ref{defn:base-k-algebras} in the formal scheme setting where we use the model object $\bb{T}_d=W\langle t_1^{\pm1},\ldots,t_d^{\pm 1}\rangle$.

We always regard a base $W$-algebra as being equipped with its $p$-adic topology. For a formal framing $w$ we write $\phi_w$ for the associated Frobenius lift as in loc.\@ cit.\@
\end{rem}

\begin{rem}\label{rem:base-algebra-lifts} We make some simple remarks here concerning base $k$-schemes that we use in the sequel without comment:

\begin{itemize}
\item Reduction modulo $p$ carries base formal $W$-schemes to base $k$-schemes. 
\item If $w\colon \bb{T}_{d,k}\to A$ is a framing, then there exists a framed base formal $W$-algebra $\wt{w}\colon \bb{T}_d\to \wt{A}$ with reduction $w$. In particular, every base $k$-scheme is Zariski locally the special fiber of a base formal $W$-scheme.
\item A base $k$-scheme is regular, excellent, and has a finite $p$-basis; this follows from the last bullet and \cite[Proposition 1.12]{IKY1}.
\end{itemize}
\end{rem}

\begin{constr}[{cf.\@ \cite{OgusRemarks}}] Suppose that $S$ is a base $k$-scheme and $(\mc{E},\varphi)$ an $F$-crystal on $S$. We will describe the construction of a filtration $\Fil^\bullet_\mr{MNO} \mc{E}_S$ in two steps.

\medskip 

\noindent\textbf{Step 1:} Assume that $(\mc{E},\varphi)$ is an $F$-crystal on $S=\mf{X}_k$ where $\mf{X}$ is a base formal $W$-scheme equipped with a framing $w$. The choice of $w$ gives the Frobenius lift $\phi_w$ and thus allows one to endow $\mc{E}_\mf{X}$ with an isomorphism $\varphi\colon \phi_w^\ast\mc{E}_\mf{X}[\nicefrac{1}{p}]\to \mc{E}_\mf{X}[\nicefrac{1}{p}]$. Let us define the \emph{Nygaard filtration}
\begin{equation*}
\Fil^r_\mr{Nyg}\phi_w^\ast \mc{E}_\mf{X}\defeq \left\{x\in \phi_w^\ast\mc{E}_\mf{X}: \varphi(x)\in p^r\mc{E}_\mf{X}\right\}.
\end{equation*}
We then define the \emph{Mazur--Nygaard filtration} on $F_S^\ast\mc{E}_S$ as follows:
\begin{equation*}
\Fil^r_\mr{MN}F_S^\ast\mc{E}_S=\mathrm{im}\left(\Fil^r_\mr{Nyg}\phi_w^\ast\mc{E}_\mf{X}\to F_S^\ast\mc{E}_S\right).
\end{equation*}
It is then not hard to show (the argument of \cite[Lemma 2.6]{KatoOgus} works mutatis mutandis) that $\Fil^\bullet_{\mr{MN}}F_S^\ast\mc{E}_S$ is independent of the choice of $\mf{X}$ and $w$ and commutes with \'etale localization. Thus, it canonically glues to a filtration $\Fil^\bullet_{\mr{MN}}F_S^\ast\mc{E}_S$ on $F_S^\ast \mc{E}_S$ for general base $k$-schemes $S$.

\medskip

\noindent\textbf{Step 2:} Continue with notation as in Step 1. As $F_S^\ast\mc{E}_S=F_{S/k}^\ast \mc{E}_{S^{(1)}}$, where $F_{S/k}$ is the relative Frobenius and $S^{(1)}$ the Frobenius twist, it naturally carries the usual Cartier connection
\begin{equation*}
\nabla_{\mr{Car},\mc{E}}\colon F_S^\ast\mc{E}_S\to F_S^\ast\mc{E}_S\otimes_{\mc{O}_S}\Omega^1_{S/k}
\end{equation*}
given in local coordinates by sending $a\otimes v$ in $F_S^\ast\mc{E}_S$ to $da\otimes v$, which has $p$-curvature $0$.

A simple calculation (e.g., using a framed base formal $W$-scheme lift $\mf{X}$) shows that the Mazur--Nygaard filtration $\Fil^\bullet_\mr{MN}F_S^\ast\mc{E}_S$ is stable under $\nabla_\mr{Car}$. Thus, by Cartier descent (see \cite[Theorem 5.1]{KatzNilpotent} and \cite[\S2.5.2]{deJongCrystalline}), it comes from a unique filtration $\Fil^\bullet_{\mr{MNO}} \mc{E}_{S^{(1)}}$, which we call the \emph{Mazur--Nygaard--Ogus filtration}. As the map $S^{(1)}\to S$ is an isomorphism, we deduce the existence of a unique filtration $\Fil^\bullet_{\mr{MNO}}\mc{E}_S$ which we also call the Mazur--Nygaard--Ogus filtration. It can be easily shown that the Mazur--Nygaard--Ogus filtration commutes with \'etale base change.

As can be deduced mutatis mutandis from \cite[Proposition 3.11]{KatoOgus}, one can give a more concrete local description of the Mazur--Nygaard--Ogus filtration. Namely, suppose that $S=\mf{X}_k$ for a base formal $W$-scheme $\mf{X}$ with framing $w$. Then, one has that
\begin{equation}\label{eq:MNO-filtration-on-a-lift}
\Fil^r_\mr{MNO}\mc{E}_S=\mr{im}\bigg(\big\{x\in \mc{E}_\mf{X}:\varphi(1\otimes x)\in p^r\mc{E}_\mf{X}\big\}\to \mc{E}_S\bigg),
\end{equation}
where $1\otimes x$ denotes the image of $x$ under the natural map $\mc{E}_\mf{X}\to \phi_w^\ast\mc{E}_\mf{X}$.
\end{constr}

We now observe that the Mazur--Nygaard and Mazur--Nygaard--Ogus filtrations admit natural Tannakian analogues.

\begin{prop} Suppose that $S$ is a base $k$-scheme, and let $(\mc{Q},\varphi)$ be an object of $\cat{FCrys}^{\sigma(\mu)}_\mc{G}(S)$. Then, the assignments
\begin{equation*}
\mr{Fil}^\bullet_\mr{MN} \,F_S^\ast\,\omega_Q|_S\colon \cat{Rep}_{\bb{Z}_p}(\mc{G})\to \cat{FilVect}(S),\qquad \Lambda\mapsto \Fil^\bullet_\mr{MN} F_S^\ast\mc{Q}^{(\Lambda)}_S,
\end{equation*}
and
\begin{equation*}
\mr{Fil}^\bullet_\mr{MNO} \,\omega_Q|_S\colon \cat{Rep}_{\bb{Z}_p}(\mc{G})\to \cat{FilVect}(S),\qquad \Lambda\mapsto \Fil^\bullet_\mr{MNO} \mc{Q}^{(\Lambda)}_S
\end{equation*}
are objects of $\mc{G}\text{-}\cat{FilVect}^{\sigma(\mu)}(S)$ and $\mc{G}\text{-}\cat{FilVect}^{\mu}(S)$, respectively.
\end{prop}
\begin{proof} Being a $\mc{G}$-object of the appropriate type is an \'etale local property. Thus, we are free by Proposition \ref{prop:type-mu-equiv-char-p} to assume that $S=\Spf(A)_k$ for some base formal $W$-scheme $\mf{X}=\Spf(A)$ equipped with a framing $w$, and that $\varphi$ is in the double coset $\mc{G}(A)\sigma(\mu)(p)\mc{G}(A)$. Given this, it is trivial to check that $\Fil^\bullet_\mr{MN} \,F_S^\ast\,\omega_Q|_S\simeq \Fil^\bullet_{\sigma(\mu)}\mc{G}_S$ and $\Fil^\bullet_\mr{MNO}\,\omega_Q|_S\simeq \Fil^\bullet_{\mu}\mc{G}_S$.\end{proof}

\begin{defn}\label{defn:MNO-filtration} Suppose that $S$ is a base $k$-scheme and let $(\mc{Q},\varphi)$ be an object of $\cat{FCrys}_\mc{G}^{\sigma(\mu)}(S)$. Then, we define the \emph{Mazur--Nygaard--Ogus filtration} $\Fil^\bullet_\mr{MNO}\mc{Q}_S$ on $\mc{Q}_S$ to be the object of $\cat{FilTors}_\mc{G}^\mu(S)$ corresponding to $\Fil^\bullet_\mr{MNO}\, \omega_Q|_S$ under Proposition \ref{prop:Filtered-Rees-equiv}. 
\end{defn}

We then have the following purely crystal-theoretic description of the Hodge-filtered de Rham realization $\Fil^\bullet_\mr{Hdg} T_\mr{dR}(\mf{Q})$ of a $(\mc{G},\mu)$-aperture.

\begin{prop}\label{prop:Hodge-filtered-bundle-descp-char-p} Suppose that $S$ is a base $k$-scheme, and let $\mf{Q}$ be a $(\mc{G},\mu)$-aperture on $S$ and set $(\mc{Q},\varphi)=\mr{R}_S(\mf{Q})$. Then, there is a natural identification of objects of $\cat{FilTors}_\mc{G}^\mu(S)$:
\begin{equation*}
\Fil^\bullet_\mr{Hdg} T_\mr{dR}(\mf{Q})=\Fil^\bullet_\mr{MNO}\mc{Q}_S.
\end{equation*}
 \end{prop}
\begin{proof} Let us first observe that there is a natural identification $T_\mr{dR}(\mf{Q})=\mc{Q}_S$. This is a crystalline-de Rham comparison and follows mutatis mutandis from the proof of \cite[Theorem 1.19]{IKY3}. Thus, it suffices to show that under this comparison the two filtrations are equal. This can be shown representation-by-representation and so fix such a representation $\Lambda$ and write $(\mc{E},\varphi)$ for the evaluation of $(\mc{Q},\varphi)$ at $\Lambda$. 

As $F_S$ is faithfully flat (as $S$ is regular), it suffices to show that
\begin{equation*}
F_S^\ast \Fil^\bullet_\mr{MNO} \mc{E}_S=\Fil^\bullet_\mr{MN} F_S^\ast \mc{E}_S=F_S^\ast \Fil^\bullet_\mr{Hdg} T_\mr{dR}(\mf{Q}).
\end{equation*}
Finally, as the above comparison is natural in $S$, it suffices to \'etale localize and thus to assume that $S=\mf{X}_k$ for $\mf{X}=\Spf(A)$ a base formal $W$-scheme with framing $w$. 

Now, \'etale localizing further, we may choose a Frobenius-equivariant embedding $A\to W(R)$; see \cite[p.\@ 12]{IKY1}. We then see by the crystal property that
\begin{equation*}
\Fil^\bullet_\mr{MN}\, F_S^\ast\,\mc{E}_S\otimes_{A_k}R=\mathrm{im}\bigg(\big\{x\in \phi_{W(R)}^\ast \mc{E}_{W(R)}: \varphi(x)\in p^r\mc{E}_{W(R)}\big)\to \mc{E}_R\bigg).
\end{equation*}
On the other hand, it follows from \cite[Claim 3.40]{GuoLi} (which is applicable here given the equivalence in \cite[Proposition 1.39]{IKY2}) that the pullback $\mf{Q}_{R^\smallN}$ of $\mf{Q}$ along $R^\smallN\to R^\mr{syn}$ may be understood as the filtered $W(R)$-module corresponding to the saturated Nygaard filtration. From this, one deduces that 
\begin{equation*}
\begin{aligned}
F_R^\ast \Fil^\bullet_\mr{Hdg}T_\mr{dR}(\mf{Q})_R &=F_R^\ast \Fil^\bullet_\mr{Hdg}T_\mr{dR}(\mf{Q}_R)\\ &=\mathrm{im}\bigg(\big\{x\in \phi_{W(R)}^\ast \mc{E}_{W(R)}: \varphi(x)\in p^r\mc{E}_{W(R)}\big)\to \mc{E}_R\bigg).
\end{aligned}
\end{equation*}
Thus, 
\begin{equation*}
F_S^\ast \Fil^\bullet_\mr{MNO}\mc{E}_S=F_S^\ast \Fil^\bullet_\mr{Hdg}T_\mr{dR}(\mf{Q})
\end{equation*} after base change along the faithfully flat cover $\Spec(R)\to S$. The conclusion follows.
\end{proof}

\subsubsection{$\mc{G}$-connections and quadruples}\label{sss:G-conn} In this section we would like to recall a concrete description of objects $\cat{FCrys}_\mc{G}^{\sigma(\mu)}(S)$ in terms of certain quadruples of data depending on the choice of a base formal $W$-scheme $\mf{X}$ with special fiber $S$.

\begin{nota} Throughout this subsubsection, $\mf{X}$ is a base formal $W$-scheme with $\mf{X}_k=S$. We moreover fix a framing $w$ of $\mf{X}$.
\end{nota}

To describe these quadruples, it is helpful to first recall the basic theory of $\mc{G}$-connections; see \cite[Appendix A]{Wakabayashi} for more details.

\begin{defn}\label{defn:Atiyah-sequence} Let ${Q}$ be a $\mc{G}$-bundle on $\mf{X}$. We then define the \emph{Atiyah bundle} to be
\begin{equation*}
\mr{At}_{\mf{X}/W}(Q):=(f_\ast T_{Q/W})^\mc{G},
\end{equation*}
where $f\colon Q\to\mf{X}$ is the structure map, and $T_{{Q}/W}$ is the tangent bundle of the map of formal schemes $Q\to\Spf(W)$. The \emph{Atiyah sequence} is then defined to be
\begin{equation}\label{eq:Atiyah-sequence}
0\to \mr{ad}(\mc{Q})\to \mr{At}_{\mf{X}/W}({Q})\to T_{\mf{X}/W}\to 0,
\end{equation}
where $\mr{ad}({Q})$ is the adjoint bundle ${Q}^{(\mf{g})}$, the first map is the natural one, and the second map is given by the differential of $f$.
\end{defn}

\begin{defn} Let notation be as in Definition \ref{defn:Atiyah-sequence}. A \emph{$\mc{G}$-connection} on ${Q}$ is then a section $\nabla\colon T_{\mf{X}/W}\to \mr{At}_{\mf{X}/W}({Q})$ of \eqref{eq:Atiyah-sequence}. We denote the set of $\mc{G}$-connections on $Q$ by $\mc{G}\text{-}\mr{Conn}(Q)$.

Observe that $\mr{At}_{\mf{X}/W}(Q)$ and $T_{\mf{X}/W}$ both carry natural Lie brackets. We call a $\mc{G}$-connection $\nabla$ \emph{integrable} if it is a map of Lie algebras. We call $\nabla$ \emph{topologically quasi-nilpotent} if the induced connection on $Q^{(\Lambda)}$ (see Remark \ref{rem:G-connections-and-normal-connections} below) is topologically (for the $p$-adic topology) quasi-nilpotent for every $\Lambda$ in $\cat{Rep}_{\Z_p}(\mc{G})$.
\end{defn}

\begin{eg}\label{eg:trivial-torsor-connections} If $\mc{G}$ denotes the trivial $\mc{G}$-torsor on $\mf{X}$, then one has a natural identification
\begin{equation*}
\mr{At}_{\mf{X}/W}(\mc{G})\simeq T_{\mf{X}/W}\oplus (\mf{g}\otimes_{\Z_p}\mc{O}_\mf{X}).
\end{equation*}
From this one obtains a natural identification
\begin{equation}\label{eg:connections-on-trivial-torsor}
\mc{G}\text{-}\mr{Conn}(\mc{G})\simeq \mf{g}\otimes_{\Z_p}\Omega^1_{\mf{X}/W},
\end{equation}
where if $\nu$ is an element of $\mf{g}\otimes\Omega^1_{\mf{X}/W}$ we obtain the $\mc{G}$-connection
\begin{equation*}
\nabla_\nu\colon T_{\mf{X}/W}\to \mr{At}_{\mf{X}/S}(\mc{G})=T_{\mf{X}/W}\oplus(\mf{g}\otimes_{\Z_p}\mc{O}_\mf{X}),\qquad D\mapsto (D,\nu(D)).
\end{equation*}
As per usual, $\nabla_\nu$ is integrable if and only if $\nu$ satisfies the Maurer--Cartan equation $\tfrac{1}{2}[\nu,\nu]=-d\nu$. For a representation $\rho\colon\mc{G}\to\GL(\Lambda)$, the induced
connection is $d+d\rho(\nu)$.
Writing $\nu=\sum_i A_i\,dx_i$, the connection $\nabla_\nu$ is
topologically quasi-nilpotent if and only if, for every $\rho$,
every local section $v$ of $\Lambda\otimes_{\Z_p}\mc{O}_\mf{X}$,
and every $i$, one has the following limit in the $p$-adic topology:
\begin{equation*}
 \lim_{s\to \infty}\left[ \left(\frac{\partial}{\partial x_i}+d\rho(A_i)\right)^s(v)\right]=0.
\end{equation*}
\end{eg}

\begin{rem}\label{rem:G-connections-and-normal-connections} Let us spell out the precise relationship between Example \ref{eg:trivial-torsor-connections} when $\mc{G}=\GL(\Lambda)$ and usual connections on $\Lambda\otimes\mc{O}_\mf{X}$. Namely, in this case we see that if $A$ belongs to $\mf{gl}(\Lambda)\otimes_{\Z_p}\Omega^1_{\mf{X}/W}=\End_{\Z_p}(\Lambda)\otimes \Omega^1_{\mf{X}/W}$ then we have the actual connection 
\begin{equation*}
\nabla_A\colon \Lambda_{\Z_p}\mc{O}_\mf{X}\to \Lambda\otimes_{\Z_p}\Omega^1_{\mf{X}/W},\qquad s\mapsto ds+As
\end{equation*}
where $As$ has the obvious meaning, e.g., if $A=T\otimes \varepsilon$ for $T$ in $\End_{\Z_p}(\Lambda)$ and $\varepsilon$ in $\Omega^1_{\mf{X}/W}$, and $s=\lambda\otimes f$ for $\lambda$ in $\Lambda$ and $f$ in $\mc{O}_\mf{X}$, then $As=T(\lambda)\otimes f\varepsilon$ in $\Lambda\otimes_{\Z_p}\Omega^1_{\mf{X}/W}$. In particular, choosing a basis of $\Lambda$ identifies $A$ with a matrix of differential $1$-forms on $\mf{X}/W$.
\end{rem}

The reason that $\mc{G}$-connections are desirable is a Tannakian version of the well-known description of crystals on $S$ as modules with integrable connection. To state this, we set the following notation.

\begin{nota}We denote by $\cat{Tors}^{\varphi,\nabla}_\mc{G}(\mf{X})$ the category of triples $(Q,\varphi,\nabla)$ where $Q$ is a $\mc{G}$-torsor on $\mf{X}$, $\nabla$ is an {integrable and topologically quasi-nilpotent} $\mc{G}$-connection, and $\varphi\colon \phi_w^\ast Q[\nicefrac{1}{p}]\to Q[\nicefrac{1}{p}]$ is a $\nabla$-horizontal quasi-isogeny of $\mc{G}$-torsors on $\mf{X}$.\footnote{Technically this category depends on the choice of formal framing $w$, but we suppress this from the notation.}
\end{nota}

The following proposition is proved, mutatis mutandis, from the classical case of vector bundles; see \cite[Proposition A.6.2]{Wakabayashi}.

\begin{prop}\label{prop:torsors-with-connection-equiv} There is a natural equivalence
\begin{equation*}
 \cat{FCrys}_\mc{G}(S)\isomto \cat{Tors}^{\varphi,\nabla}_\mc{G}(\mf{X}),\qquad \mc{Q}\mapsto (\mc{Q}_\mf{X},\varphi_\mc{Q},\nabla_\mc{Q}).
\end{equation*}
\end{prop}

\begin{rem}\label{rem:tannakian-connection}Let us say that a \emph{connection} on an object 
\begin{equation*}
\omega\colon\cat{Rep}_{\Z_p}(\mc{G})\to\cat{Vect}(\mf{X})
\end{equation*}
of $\mc{G}\text{-}\cat{Vect}(\mf{X})$ 
is a connection for each representation $\Lambda$ of $\mc{G}$
\begin{equation*}
\nabla(\Lambda)\colon \omega(\Lambda)\to \omega(\Lambda)\otimes_{\mc{O}_\mf{X}}\Omega^1_{\mf{X}/W}
\end{equation*}
compatible in $\Lambda$ in the obvious way. We say that $\nabla=(\nabla(\Lambda))$ is \emph{integrable} or \emph{topologically quasi-nilpotent} if $\nabla(\Lambda)$ is for each $\Lambda$. 

It is evident that if $Q$ is a $\mc{G}$ torsor on $\mf{X}$ and $\omega_Q$ is its associated $\mc{G}$-object as in Proposition \ref{prop:Tannakian-equiv}, giving a connection on $Q$ is equivalent to giving a connection on $\omega_Q$, and moreover that the notions of integrability and topological quasi-nilpotence match.
\end{rem}

We now define the aforementioned quadruples which will give a concrete description of an object of $\cat{FCrys}^{\sigma(\mu)}_\mc{G}(S)$. To do this, we fix some auxiliary data.

\begin{nota}\label{nota:tensor-package} We fix a {tensor package} for $\mathcal{G}$: a $\bb{Z}_p$-module $\Lambda$ and a finite set of tensors $\{s_{\alpha,0}\}\subseteq \Lambda^\otimes$ such that $\mc{G}=\mr{Fix}(\{s_{\alpha,0}\})$. See \cite[Appendix A]{IKY1} for this notion, the notation we've used here, and the existence of tensor packages.
\end{nota}

\begin{defn}\label{defn:quad} Consider quadruples $\big(M,\nabla,\varphi,\{s_\alpha\}\big)$ where:
\begin{itemize}
\item $M$ is a vector bundle on $\mf{X}$;
\item $\varphi$ is an isomorphism $\varphi\colon \phi_w^\ast M[\nicefrac{1}{p}]\to M[\nicefrac{1}{p}]$ satisfying condition (1) from Proposition \ref{prop:type-mu-equiv-char-p} for $\nu=\sigma(\mu)$;
\item $\nabla\colon M\to M\otimes \Omega^1_{\mf{X}/W}$ is a topologically quasi-nilpotent integrable connection;
\item $\{s_\alpha\}\subseteq M^\otimes$ is a set of tensors;
\end{itemize}
subject to the conditions:
\begin{enumerate}
\item $\varphi$ and $\{s_\alpha\}$ are $\nabla$-horizontal;
\item the tensors $\{s_\alpha\}$ are $\varphi$-fixed;\footnote{More precisely, the map on tensor constructions induced by $\varphi$ carries $s_\alpha\otimes1$ to $s_\alpha$, where $s_\alpha\otimes1$ is viewed in $(\phi_w^\ast M[\nicefrac1p])^\otimes$.}
\item and there exists an \'etale cover $\mf{X}'\to\mf{X}$ and an isomorphism $M_{\mf{X}'}\isomto \Lambda\otimes\mc{O}_{\mf{X}'}$ carrying the tensors $\{s_\alpha\}$ to the tensors $\{s_{\alpha,0}\otimes 1\}$. 
\end{enumerate}
With the obvious notion of morphism, such quadruples form a category $\mb{Quad}^\mu_\mc{G}(\mf{X})$.\footnote{Again we are suppressing the dependence on $w$ from the notation.}
\end{defn}

The following follows easily from Proposition \ref{prop:torsors-with-connection-equiv} and Remark \ref{rem:tannakian-connection}.

\begin{prop}\label{prop:F-crystal-quad-equiv} There is a natural equivalence
\begin{equation*}
\cat{FCrys}_\mc{G}^{\sigma(\mu)}(S)\isomto \cat{Quad}_\mc{G}^\mu(\mf{X}).
\end{equation*}
\end{prop}

\subsubsection{The case of finite fields}\label{ss:case-of-finite-fields} Because it will play such a prominent role below, we now specially analyze the case of $(\mc{G},\mu)$-apertures over finite fields, giving a simple group-theoretic description of this groupoid.

\begin{nota}\label{nota:finite-field}
Assume that $k$ is finite, and write:
\begin{itemize}
\item $r=[k:\F_p]=[K:\Q_p]$ and $q=|k|=p^r$;
\item $k_j$ for the degree-$j$ extension of $k$ in $\ov{k}$,
for $1\leqslant j<\infty$, and $k_\infty=\ov{k}$;
\item $W_j=W(k_j)$ and $K_j=W_j[\nicefrac1p]$,
for $1\leqslant j\leqslant\infty$;
\item and $\sigma$ for the $p$-Frobenius lift on $W_j$.
\end{itemize}
\end{nota}

\begin{setup}\label{setup:simplifying-over-finite-fields} Using Propositions \ref{prop:BT-in-char-p} and \ref{prop:torsors-with-connection-equiv}, we see that there are natural equivalences
\begin{equation*}
\mr{BT}^{\mc{G},\mu}_\infty(k_j)\isomto \cat{FCrys}^{\sigma(\mu)}_\mc{G}(k_j)\isomto \cat{Tors}_\mc{G}^{\varphi,\sigma(\mu)}(W_j),
\end{equation*}
where this last category denotes the category of pairs $(Q,\varphi)$ where $Q$ is a $\mc{G}$-torsor on $W_j$ and $\varphi\colon \sigma^\ast Q[\nicefrac{1}{p}]\to Q[\nicefrac{1}{p}]$ is an isomorphism satisfying (4) of Proposition \ref{prop:type-mu-equiv-char-p}.

We denote by $\cat{Tors}^{\varphi,\sigma(\mu),\mr{neut}}_\mc{G}(W_j)$ the full subcategory of $\cat{Tors}_\mc{G}^{\varphi,\sigma(\mu)}(W_j)$ consisting of those objects $(Q,\varphi)$ with $Q=\mc{G}$. We denote by $\mr{BT}^{\mc{G},\mu,\mr{neut}}_\infty(k_j)$ and $\cat{FCrys}_\mc{G}^{\sigma(\mu),\mr{neut}}(k_j)$ the essential preimages of $\cat{Tors}^{\varphi,\sigma(\mu),\mr{neut}}_\mc{G}(W_j)$ in $\mr{BT}^{\mc{G},\mu}_\infty(k_j)$ and $\cat{FCrys}^{\sigma(\mu)}_\mc{G}(k_j)$, respectively. We call these the \emph{neutral subcategories} of their respective ambient categories.
\end{setup}

\begin{obs} All $\mc{G}$-torsors on $W_j$ are trivial by combining Footnote \ref{footnote:Hensel} with \cite{LangGroups}. Thus, the inclusion of each neutral subcategory into its ambient category is an equivalence.
\end{obs}

\begin{nota}\label{nota:sigma-conjugacy-classes} Define the action groupoid
\begin{equation*} 
\mc{C}_j(\mc{G})\defeq \big[G(K_j)/_{\mr{Ad}_\sigma} \mc{G}(W_j)\big]
\end{equation*}
where we have the action
\begin{equation*}
\mr{Ad}_\sigma\colon G(K_j)\times \mc{G}(W_j)\to G(K_j),\qquad (b,c)\mapsto cb\sigma(c)^{-1}.
\end{equation*}
We denote by
\begin{equation*}
C_j(\mc{G})\defeq \pi_0\left(\mc{C}_j(\mc{G})\right)\quad \text{and}\quad \langle-\rangle\colon G(K_j)\to C_j(\mc{G}),
\end{equation*}
the set of isomorphism classes and the natural map taking $b$ to its isomorphism class $\langle b\rangle$. Write
\begin{equation*}
C_j(\mc{G},\sigma(\mu))=\left\langle \mc{G}(W_j)\sigma(\mu)(p)\mc{G}(W_j)\right\rangle \subseteq C_j(\mc{G}),
\end{equation*}
 and denote by $\mc{C}_j(\mc{G},\sigma(\mu))$ the associated full subgroupoid of $\mc{C}_j(\mc{G})$.
\end{nota}

We can then give the following group-theoretic classification of $\mr{BT}^{\mc{G},\mu}_\infty(k_j)$.

\begin{prop}\label{prop:isom-classes-in-BT} The process above defines equivalences of groupoids
\begin{equation*}
\mr{BT}^{\mc{G},\mu}_\infty(k_j)\isomfrom \mr{BT}^{\mc{G},\mu,\mr{neut}}_\infty(k_j)\isomto \mc{C}_j(\mc{G},\sigma(\mu))\hookrightarrow \mc{C}_j(\mc{G}),
\end{equation*}
and thus a bijection
\begin{equation*}
\pi_0\left(\mr{BT}^{\mc{G},\mu}_\infty(k_j)\right)\isomto C_j(\mc{G},\sigma(\mu)).
\end{equation*}
\end{prop}
\begin{proof} It suffices to prove the second claim. By definition, the image of $\pi_0\left(\mr{BT}^{\mc{G},\mu}_\infty(k_j)\right)$ in $C_j(\mc{G})$ is the image in $C_j(\mc{G})$ of $G(K_j)\cap \mc{G}(W(\overline{k}))\sigma(\mu)(p)\mc{G}(W(\overline{k}))$, where this intersection is taken inside of $G(W(\overline{k})[\nicefrac{1}{p}])$. But, this is $\mc{G}(W_j)\sigma(\mu)(p)\mc{G}(W_j)$: apply flat descent for $\mr{Gr}_{\mc{G},\sigma(\mu)}$ from \cite{ZhuAffGmix}.
\end{proof}

\subsection{$(\mc{G},\mu)$-apertures and filtered $F$-crystals}\label{ss:filtered-F-crystals} We now give a mixed characteristic version of Proposition \ref{prop:BT-in-char-p} for {base formal $W$-schemes} as in Remark \ref{rem:base-algebras}. The key difference is that instead of plain $F$-crystals with $\mc{G}$-structure, we will use \emph{filtered} $F$-crystals with $\mc{G}$-structure. 

\begin{nota}Return to Notation \ref{nota:deformation-spaces-setup}, and fix a base formal $W$-scheme $\mf{X}$.
\end{nota}

\subsubsection{Various categories of filtered $F$-crystals with $\mc{G}$-structure} We begin by describing the type of filtered $F$-crystals with $\mc{G}$-structure appearing in our comparison to $(\mc{G},\mu)$-apertures.

\begin{defn}[{see \cite[\S2.1]{IKY3}}]\label{defn:FFcrys} The category $\cat{FFCrys}(\mf{X})$ of \emph{filtered $F$-crystals} on $\mf{X}$ is the category of triples $(\mc{E},\varphi,\Fil^\bullet \mc{E})$ where $(\mc{E},\varphi)$ is an object of $\cat{FCrys}(\mf{X}_k)$ and $\Fil^\bullet \mc{E}\subseteq \mc{E}_{\mf{X}}$ is a decreasing, separated, and exhaustive filtration satisfying the following conditions:
\begin{itemize}
\item $\Fil^\bullet \mc{E}$ satisfies Griffiths transversality with respect to the natural connection $\nabla_\mc{E}$ on $\mc{E}_\mf{X}$;\footnote{Recall that this means that $\nabla_\mc{E}(\Fil^i \mc{E})\subseteq \Fil^{i-1} \mc{E}\otimes_{\mc{O}_\mf{X}}\Omega^1_{\mf{X}/W}$.}
\item $\Fil^\bullet \mc{E}$ is locally split.\footnote{Recall that this means that the graded pieces $\mr{Gr}^i\mr{Fil}^\bullet \mc{E}$ are vector bundles on $\mf{X}$ for all $i$.}
\end{itemize}
\end{defn}

In practice that category of filtered $F$-crystals will not be sufficient for our comparison to $(\mc{G},\mu)$-apertures, and we will need the following condition similar to that which appears in Fontaine--Laffaille theory; cf.\@ \cite{Faltings89}.

\begin{defn}[{see \cite[\S2.1.4]{IKY3}}]\label{defn:strongly-divisible} A filtered $F$-crystal $(\mc{E},\varphi,\Fil^\bullet \mc{E})$ is \emph{strongly divisible} if for any affine open $\Spf(R)\subseteq \mf{X}$ and any formal framing $w$ of $R$ one has that
\begin{equation*}
\varphi\left(\sum_{r\in\bb{Z}}p^{-r}\phi_w^\ast\Fil^r\mc{E}(R)\right)=\mc{E}_\mf{X}(R).
\end{equation*}
\end{defn}

Of course, we are interested in Tannakian versions of the above categories.

\begin{defn} The category $\cat{FFCrys}(\mf{X})$ naturally carries the structure of an exact $\bb{Z}_p$-linear $\otimes$-category as in \cite[\S2.1.1--2.1.2]{IKY3}. We then write $\cat{FFCrys}_\mc{G}(\mf{X})$ for its category of $\mc{G}$-objects, and $\cat{FFCrys}^{\mr{sd}}_\mc{G}(\mf{X})$ for the full subcategory of $\cat{FFCrys}_\mc{G}(\mf{X})$ taking values in $\cat{FFCrys}^\mr{sd}(\mf{X})$. 
\end{defn}

We will often switch between the literal interpretation of objects of $\cat{FFCrys}_\mc{G}(\mf{X})$ as exact $\Z_p$-linear $\otimes$-functors $\omega\colon \cat{Rep}_{\Z_p}(\mc{G})\to \cat{FFCrys}(\mf{X})$ and the following more intrinsic definition.

\begin{rem} Using Propositions \ref{prop:Tannakian-equiv} and \ref {prop:Filtered-Rees-equiv}, the category $\cat{FFCrys}_\mc{G}(\mf{X})$ is naturally equivalent to the category of triples $(\mc{Q},\varphi,\Fil^\bullet \mc{Q}_\mf{X})$ where $(\mc{Q},\varphi)$ is an object of $\cat{FCrys}_\mc{G}(\mf{X}_k)$ and $\Fil^\bullet \mc{Q}_\mf{X}$ is a locally split filtration on $\mc{Q}_\mf{X}$ satisfying Griffiths transversality on each representation.
\end{rem}

We next impose the necessary type-$\mu$ condition on our filtered $F$-crystals with $\mc{G}$-structure.

\begin{defn}\label{defn:filtered-crystalline-type} An object $(\mc{Q},\varphi,\Fil^\bullet \mc{Q}_\mf{X})$ of $\cat{FFCrys}_\mc{G}(\mf{X})$ is \emph{of type $\mu$} if $\Fil^\bullet \mc{Q}_\mf{X}$ is. We denote by $\cat{FFCrys}^{\mu}_\mc{G}(\mf{X})$ and $\cat{FFCrys}^{\mr{sd},\mu}_\mc{G}(\mf{X})$ the full subcategories spanned by objects of type $\mu$. 
\end{defn}

\begin{rem} By Proposition \ref{prop:type-nu-classification} we have a natural functor 
\begin{equation*}
\cat{FFCrys}^{\mu}_\mc{G}(\mf{X})\to BP_\mu(\mf{X}),\qquad (\mc{Q},\varphi,\Fil^\bullet \mc{Q}_\mf{X})\mapsto \Fil^\bullet\mc{Q}_\mf{X}.
\end{equation*}
\end{rem}

We observe that there is an apparent mismatch between two notions of `type $\mu$': As the type $\mu$-condition on $(\mc{Q},\varphi,\Fil^\bullet \mc{Q}_\mf{X})$ is only about the filtration $\Fil^\bullet \mc{Q}_\mf{X}$ it's not evident that the object $(\mc{Q},\varphi)$ of $\cat{FCrys}_\mc{G}(\mf{X}_k)$ is of type $\mu$ (or rather $\sigma(\mu)$). That said, the following shows that this mismatch is illusory when the objects are strongly divisible.

\begin{prop}\label{prop:sd-implies-type-mu} Suppose that $(\mc{Q},\varphi,\Fil^\bullet \mc{Q}_\mf{X})$ is an object of $\cat{FFCrys}_\mc{G}^{\mr{sd},\mu}(\mf{X})$. Then, the object $(\mc{Q},\varphi)$ of $\cat{FCrys}_\mc{G}(\mf{X}_k)$ is of type $\sigma(\mu)$.
\end{prop}
\begin{proof} By Proposition \ref{prop:type-mu-equiv-char-p} we see that it suffices to assume that $\mf{X}=\Spf(W)$ and that $k$ is algebraically closed. In this case, as $\Fil^\bullet\mc{Q}_\mf{X}$ is a $P_\mu$-torsor, it must by Footnote \ref{footnote:Hensel} be trivial. Choosing a trivialization of $\sigma_{}^\ast\mc{G}\simeq \mc{G}$ we may write $\varphi=g$ for some $g$ in $\mc{G}(K)$. Moreover, as $\mu$ is minuscule it suffices to show that $g$ is bounded by $\sigma(\mu)$ in the Bruhat order. By \cite[Lemma 2.2 (iv)]{RapoportRichartz} it suffices to check that $\rho(g)$ is bounded by $\rho\circ \sigma(\mu)$ for all representations $\rho\colon \mc{G}\to\GL(\Lambda)$. But, by strong divisibility we have that
\begin{equation*}
\sum_{r\in\bb{Z}}p^{-r}\rho(g)(\Fil^r_{\sigma(\mu)}(\Lambda))=\Lambda.
\end{equation*}
This implies that $\rho(g)\rho(\sigma(\mu)(p))^{-1}\Lambda=\Lambda$, and thus that $\rho(g)$ is bounded by $\rho\circ \sigma(\mu)$, as desired.
\end{proof}

Finally, we record for future use that strong divisibility of an object of $\cat{FFCrys}_\mc{G}(\mf{X})$ can be tested on one faithful representation, at least if the object is of type $\mu$.

\begin{lem}\label{lem:one-for-all}  Let $\mf{X}$ be a base formal $W$-scheme, $(\mc{Q},\varphi,\Fil^\bullet\mc{Q}_\mf{X})$ an object of $\cat{FFCrys}^\mu_\mc{G}(\mf{X})$, and $\rho\colon \mc{G}\to\mr{GL}(\Lambda)$ a faithful representation. Then, $(\mc{Q},\varphi,\Fil^\bullet\mc{Q}_{\mf{X}})$ belongs to $\cat{FFCrys}^{\mr{sd},\mu}_\mc{G}(\mf{X})$ if and only if $(\mc{Q},\varphi,\Fil^\bullet\mc{Q}_\mf{X})^{(\Lambda)}$ is strongly divisible. \end{lem}
\begin{proof}We may work \'etale locally on $\mf{X}$ and thus assume that $\Fil^\bullet\mc{Q}_\mf{X}\simeq \Fil^\bullet_\mu|_{\mf{X}}$ and that $\mf{X}=\Spf(R)$ for a base formal $W$-algebra $R$ with formal framing $w$. Choose an identification of $\phi_w^\ast \mc{G}_{\mf{X}}\simeq \mc{G}_\mf{X}$ so that we may think of $\varphi$ as an element $g$ in $\mc{G}(R[\nicefrac{1}{p}])$. As $\rho$ is a closed immersion and $R\to R[\nicefrac{1}{p}]$ is injective we deduce that $g\sigma(\mu)(p)^{-1}$ belongs to $\mc{G}(R)$, as desired. \end{proof}

\subsubsection{The filtered crystalline-aperture equivalence} We are now prepared to state the equivalence between filtered $F$-crystals with $\mc{G}$-structure and $(\mc{G},\mu)$-apertures. 

\begin{thm}\label{thm:BT-over-base-ring} Suppose $p>2$ and that $\mf{X}$ is a base formal $W$-scheme. Then, the functor
\begin{equation*}
\bb{D}_\mr{crys}\colon \mr{BT}^{\mc{G},\mu}_\infty(\mf{X})\to \cat{FFCrys}^{\mr{sd},\mu}_\mc{G}(\mf{X}),\quad \mf{Q}\mapsto \bigg(\mr{R}_{\mf{X}_k}(\mf{Q}_{\mf{X}_k}),\Fil^\bullet_\mr{Hdg} T_\mr{dR}(\mf{Q})\bigg)
\end{equation*}
is an equivalence compatible with the structure maps to $BP_\mu(\mf{X})$.
\end{thm}

\begin{rem} For $\bb{D}_\mr{crys}$ to make sense we need that the underlying $\mc{G}$-bundle of $\Fil^\bullet_\mr{Hdg} T_\mr{dR}(\mf{Q})$ agrees with $\mr{R}_{\mf{X}_k}(\mf{Q}_{\mf{X}_k})_\mf{X}$: this follows from \cite[Theorem 1.19]{IKY3}. That $\bb{D}_\mr{crys}(\mf{Q})$ belongs to $\cat{FFCrys}^{\mr{sd},\mu}_\mc{G}(\mf{X})$ is also quite non-trivial; see \cite[Proposition 1.41]{IKY2} or \cite[Remark 3.7.6 and Construction 3.7.7]{MY}.
\end{rem}

To prove Theorem \ref{thm:BT-over-base-ring}, we first need some setup.

\begin{nota} Let us write $X_\mu(\mf{X})$ for the following groupoid:
\begin{equation*}
X_\mu(\mf{X})\defeq B\mc{G}(\mf{X})\times_{B\mc{G}(\mf{X}_k)} BP_\mu(\mf{X}_k).
\end{equation*}
\end{nota}

Key to our proof are the following constructions. 

\begin{constr} There is a natural map
\begin{equation*}
\pi_1\colon BP_\mu(\mf{X})\to X_\mu(\mf{X}),\qquad P\mapsto (P\times^{P_\mu}\mc{G},P_{\mf{X}_k},\mr{nat.}),
\end{equation*}
where $\mr{nat.}$ is the obvious map. On the other hand, there is a natural map
\begin{equation*}
\pi_2\colon \cat{FCrys}^{\sigma(\mu)}_\mc{G}(\mf{X}_k)\to X_\mu(\mf{X}),\qquad (\mc{Q},\varphi)\mapsto (\mc{Q}_\mf{X},\Fil^\bullet_\mr{MNO}\mc{Q}_{\mf{X}_k},\mr{nat.})
\end{equation*}
where $\Fil^\bullet_\mr{MNO}\mc{Q}$ is as in Definition \ref{defn:MNO-filtration}, and again $\mr{nat.}$ denotes the natural map.
\end{constr}

\begin{prop}\label{prop:j-map} Suppose $p>2$ and that $\mf{X}$ is a base formal $W$-scheme. Then, 
\begin{equation*}
\j\colon \cat{FFCrys}^{\mr{sd},\mu}_\mc{G}(\mf{X})\to \cat{FCrys}_\mc{G}(\mf{X}_k)\times BP_\mu(\mf{X}),\qquad (\mc{Q},\varphi,\Fil^\bullet {\mc{Q}_\mf{X}})\mapsto \big((\mc{Q},\varphi),\Fil^\bullet {\mc{Q}_\mf{X}}\big)
\end{equation*}
factorizes naturally through $\cat{FCrys}^{\sigma(\mu)}_\mc{G}(\mf{X}_k)\times BP_\mu(\mf{X})$. Moreover, there is a natural equivalence $\pi_1\circ \mr{pr}_2\circ\j\simeq \pi_2\circ \mr{pr}_1\circ\j$, and thus a natural map
\begin{equation*}
\j\colon \cat{FFCrys}^{\mr{sd},\mu}_\mc{G}(\mf{X})\to \cat{FCrys}^{\sigma(\mu)}_\mc{G}(\mf{X}_k)\times_{X_\mu(\mf{X})} BP_\mu(\mf{X}).
\end{equation*}
\end{prop}
\begin{proof}  The factorization claim follows from Proposition \ref{prop:sd-implies-type-mu}. To prove the claim about the existence of a natural equivalence $\pi_1\circ \mr{pr}_2\circ\j\simeq \pi_2\circ \mr{pr}_1\circ\j$ we proceed as follows. Unraveling the definition, the content of the claim is that $\Fil^\bullet {\mc{Q}_\mf{X}}|_{\mf{X}_k}$ is naturally isomorphic to $\Fil^\bullet_\mr{MNO}\mc{Q}_{\mf{X}_k}$. As both of these are naturally filtrations on $\mc{Q}_{\mf{X}_k}$, it suffices to check that they are equal, and in fact we can check it representation-by-representation. Let us write $(\mc{E},\varphi,\Fil^\bullet\mc{E}_\mf{X})$ for $(\mc{Q},\varphi,\Fil^\bullet{\mc{Q}_\mf{X}})^{(\Lambda)}$. 

Now, to check this equality we are free to assume that $\mf{X}=\Spf(R)$ where $R$ is a framed base $W$-algebra with framing $w$. As $\phi_w$ is faithfully flat it suffices to check further that this equality holds after pullback along $\phi_w$. As our filtration is locally split, from \cite[Proposition 2.7]{IKY3} we see that $(\mc{E},\varphi,\Fil^\bullet\mc{E})$ being strongly divisible implies the following equality of filtered modules
\begin{equation*}
\phi_w^\ast\Fil^\bullet \mc{E}\otimes_{(R,\Fil^\bullet_\mr{triv})}(R,\Fil^\bullet_p)=\Fil^\bullet_\mr{Nyg}\mc{E},
\end{equation*}
with notation as in loc.\@ cit. Reducing this mod $p$ (which is the same as taking the image in the mod $p$ reduction as the filtration is locally split) we see that 
\begin{equation*}
\phi_w^\ast \Fil^\bullet \mc{E}= \Fil^\bullet_\mr{Nyg}\mc{E}\mod p,
\end{equation*}
and the latter is equal to $\phi_w^\ast \Fil^\bullet_\mr{MNO}\mc{E}_{\mf{X}_k}$ by construction.\end{proof}

\begin{proof}[Proof of Theorem \ref{thm:BT-over-base-ring}] Consider the diagram
\begin{equation*}
\begin{tikzcd}[cramped,column sep=5em]
	{\mr{BT}^{\mc{G},\mu}_\infty(\mf{X})} & {\cat{FFCrys}^{\mr{sd},\mu}_\mc{G}(\mf{X})} \\
	{\mr{BT}_\infty^{\mc{G},\mu}(\mf{X}_k)\times_{X_\mu(\mf{X})} BP_\mu(\mf{X})} & {\cat{FCrys}_\mc{G}^{\sigma(\mu)}(\mf{X}_k)\times_{X_\mu(\mf{X})} BP_\mu(\mf{X}),}
	\arrow["{\bb{D}_\mr{crys}}", from=1-1, to=1-2]
	\arrow["\wr"', from=1-1, to=2-1]
	\arrow["\j",from=1-2, to=2-2]
	\arrow["{\mr{R}_{\mf{X}_k}\times \mr{id}}"', from=2-1, to=2-2]
\end{tikzcd}
\end{equation*}
where the right vertical is the fully faithful map $\j$ from Proposition \ref{prop:j-map} and the left-vertical arrow is the equivalence by Grothendieck--Messing theory (see Theorem \ref{thm:GMM}). This diagram commutes which, unraveling the definitions, reduces quickly to Proposition \ref{prop:Hodge-filtered-bundle-descp-char-p}. As the right vertical arrow is fully faithful, we will then deduce the top arrow is an equivalence as $\mr{R}_{\mf{X}_k}$, and thus the bottom-horizontal arrow, is an equivalence by Proposition \ref{prop:BT-in-char-p}.
\end{proof}

\subsubsection{Griffiths transversality and quintuples} We would now like to give a more concrete description of $\cat{FFCrys}_\mc{G}^{\mr{sd},\mu}(\mf{X})$ in the vein of Proposition \ref{prop:F-crystal-quad-equiv}.

\begin{rem}\label{rem:Kodaira-Spencer} Let $\Fil^\bullet Q$ be a filtered $\mc{G}$-bundle on $\mf{X}$ of type $\mu$. By Proposition \ref{prop:type-nu-classification} we obtain a reduction of structure group $P\subseteq Q$, where $P$ is a $P_\mu$ torsor. It is simple to show that
\begin{equation*}
\frac{\mr{At}_{\mf{X}/W}(Q)}{\mr{At}_{\mf{X}/W}(P)}\simeq \frac{\mr{ad}(Q)}{\mr{ad}(P)}=\frac{\mr{ad}(Q)}{\Fil^0 \mr{ad}(Q)},
\end{equation*}
where the identification of $\mr{ad}(P)=\mr{Fil}^0 \mr{ad}(Q)$ is by unraveling the definitions. 

Suppose then that we have a $\mc{G}$-connection $\nabla$ on $Q$. This then induces a morphism
\begin{equation*}
\kappa\colon T_{\mf{X}/W}\to \mr{At}_{\mf{X}/W}(Q)\twoheadrightarrow \frac{\mr{ad}(Q)}{\Fil^0 \mr{ad}(Q)},
\end{equation*}
called the \emph{Kodaira--Spencer morphism}. 
\end{rem}

\begin{defn}\label{defn:Griffiths} Let $(Q,\Fil^\bullet Q,\nabla)$ be as in Remark \ref{rem:Kodaira-Spencer}. We say that the triple $(Q,\Fil^\bullet Q,\nabla)$ satisfies \emph{Griffiths transversality} if $\kappa$ factorizes through $\Fil^{-1}\mr{ad}(Q)/\Fil^0\mr{ad}(Q)$. 
\end{defn}

\begin{nota} For the rest of this subsubsection, fix a formal framing $w$ of $\mf{X}$, and write $\phi_w$ for its associated Frobenius lift.
\end{nota}

\begin{defn}\label{defn:FFCrys-on-cover}We denote by $\cat{FilTors}^{\varphi,\nabla,\mu}_\mc{G}(\mf{X})$ the category of quadruples $(Q,\varphi,\nabla,\Fil^\bullet Q)$ where $Q$ is a $\mc{G}$-torsor on $\mf{X}$, $\nabla$ is an {integrable and topologically quasi-nilpotent} $\mc{G}$-connection, $\varphi\colon \phi_w^\ast Q[\nicefrac{1}{p}]\to Q[\nicefrac{1}{p}]$ is a $\nabla$-horizontal quasi-isogeny of $\mc{G}$-torsors on $\mf{X}$, and $\Fil^\bullet Q$ is a filtration on $Q$ of type $\mu$ satisfying Griffiths transversality.

We write $\cat{FilTors}^{\varphi,\mr{sd},\nabla,\mu}_\mc{G}(\mf{X})$ for the full subcategory of $\cat{FilTors}^{\varphi,\nabla,\mu}_\mc{G}(\mf{X})$ spanned by those objects such that $\varphi$ satisfies the strong divisibility condition with respect to $\Fil^\bullet Q$. \footnote{We are again suppressing the dependency on $w$.}
\end{defn}

The following is proven in much the same way as Proposition \ref{prop:torsors-with-connection-equiv}.

\begin{prop}\label{prop:torsors-with-connection-equiv-filtered} There is a natural equivalence
\begin{equation*}
 \cat{FFCrys}^{\mu}_\mc{G}(\mf{X})\isomto \cat{FilTors}^{\varphi,\nabla,\mu}_\mc{G}(\mf{X}),\qquad (\mc{Q},\varphi,\Fil^\bullet \mc{Q}_\mf{X})\mapsto (\mc{Q}_\mf{X},\varphi,\nabla_\mc{Q},\Fil^\bullet \mc{Q}_\mf{X}),
\end{equation*}
restricting to an equivalence $\cat{FFCrys}^{\mr{sd},\mu}_\mc{G}(\mf{X})\isomto \cat{FilTors}^{\varphi,\mr{sd},\nabla,\mu}_\mc{G}(\mf{X})$.
\end{prop}
We now define a category of quintuples which gives a more concrete description of $\cat{FFCrys}^{\mr{sd},\mu}_\mc{G}(\mf{X})$.

\begin{defn}\label{defn:quint} Fix a tensor package $\{s_{\alpha,0}\}\subseteq \Lambda^\otimes$ as in Notation \ref{nota:tensor-package}. Consider quintuples $\big(M,\nabla,\varphi,\{s_\alpha\},\Fil^\bullet M\big)$ where the quadruple $(M,\nabla,\varphi,\{s_\alpha\})$ is the same data as in Definition \ref{defn:quad}, and $\Fil^\bullet M$ is a filtration of type $\mu$ on $M$. The quintuple $(M,\nabla,\varphi,\{s_\alpha\},\Fil^\bullet M)$ must satisfy those conditions in Definition \ref{defn:quad}, and:
\begin{itemize}
\item[(4)] $\{s_\alpha\}$ lies in the $0^\text{th}$ part of the filtration that $\Fil^\bullet M$ induces on $M^\otimes$;
\item[(5)] the filtration satisfies Griffiths transversality with respect to $\nabla$;
\item[(6)] $\varphi$ is strongly divisible relative to $\Fil^\bullet M$'
\item[(7)] and the isomorphism $M_{\mf{X}'}\isomto \Lambda\otimes\mc{O}_{\mf{X}'}$ additionally carries $\Fil^\bullet M|_{\mf{X}'}$ to $\Fil^\bullet_\mu(\Lambda)|_{\mf{X}'}$. 
\end{itemize}
With the obvious notion of morphism, we denote the category of such quintuples by $\cat{Quint}^\mu_\mc{G}(\mf{X})$.\footnote{We are again suppressing the dependency on $w$.}
\end{defn}

The following follows easily from Proposition \ref{prop:torsors-with-connection-equiv-filtered} and Remark \ref{rem:tannakian-connection}.

\begin{prop}\label{prop:F-crystal-quint-equiv} There is a natural equivalence
\begin{equation*}
\cat{FFCrys}^{\mr{sd},\mu}_\mc{G}(\mf{X})\isomto \cat{Quint}_\mc{G}^\mu(\mf{X}).
\end{equation*}
\end{prop}

\subsection{Relation to $p$-divisible groups} We next recall the exact relationship between $(\mc{G},\mu)$-apertures, $p$-divisible groups, and classical (crystalline) Dieudonn\'e theory.

\begin{nota} For $1\leqslant n\leqslant \infty$, we denote by $\mr{BT}_{p,n}$ the moduli stack of $n$-truncated $p$-divisible groups (e.g., the stack of $p$-divisible groups when $n=\infty$).
\end{nota}
\begin{nota} For a $p$-adic formal scheme $\mf{X}$ and $1\leqslant n\leqslant \infty$, denote by $\cat{Vect}_{[0,1]}(\mf{X}^\syn\dotimes \bb{Z}/p^n)$ the category of vector bundles $\mc{V}$ on $\mf{X}^\mr{syn}\dotimes\bb{Z}/p^n$ such that for every geometric point $x\colon \Spec(\kappa)\to \mf{X}$, the filtered $\kappa$-vector space $(x_\mr{dR}^\mr{syn})^\ast\mc{V}|_\kappa$ has zero graded pieces outside degrees $0$ and $1$.
\end{nota}

\begin{thm}\label{thm:syntomic-DD} For any $p$-adic formal scheme $\mf{X}$, there are equivalences
\begin{equation*}
\mc{M}\colon \mr{BT}_{p,n}(\mf{X})\isomto \cat{Vect}_{[0,1]}(\mf{X}^\syn\dotimes \bb{Z}/p^n)
\end{equation*}
compatible with truncation maps in $n$ and Cartier duality. Moreover, one has the identifications
\begin{equation*}
\mr{R}_\mf{X}\circ \mc{M}\simeq \bb{D}, \qquad\qquad \bb{D}_\mr{crys}\circ \mc{M}\simeq \bb{D}
\end{equation*}
when $\mf{X}$ is a regular $\bb{F}_p$-scheme and a base formal $W$-scheme, respectively, assuming $p>2$ in the latter case. Here $\bb{D}$ denotes the usual contravariant Dieudonn\'e crystal and filtered Dieudonn\'e crystal, respectively; see \cite[\S4.2]{IKY3}.
\end{thm}
\begin{proof} The first equivalence (the real content of this result) is \cite[Theorem 1.16]{AnschutzLeBrasDD} when $\mf{X}$ is quasi-syntomic and \cite[Theorem A]{GMM} in general. For the comparison when $\mf{X}$ is a regular $\bb{F}_p$-scheme see \cite[Lemma 4.45]{AnschutzLeBrasDD}, and for the filtered comparison over base formal $W$-schemes see \cite[Theorem 4.8]{IKY3}. Note that here we use the contravariant convention of \cite{AnschutzLeBrasDD} and \cite{IKY3}, obtained by taking ordinary linear duals of the
covariant construction in \cite{GMM}. \end{proof}

\begin{rem}\label{rem:BT-tate-module-etale-realization}
Write $\mc{M}_{\mr{cov}}(H)\defeq\mc{M}(H)^\vee$ for the covariant
construction. For $H$ in $\mr{BT}_{p,n}(\mf{X})$, there is a natural
identification on $\mf{X}_\eta$
\begin{equation*}
T_\et(\mc{M}_{\mr{cov}}(H))\simeq
\begin{cases}
H_\eta,&n<\infty,\\
T_p(H),&n=\infty.
\end{cases}
\end{equation*}
The finite-level assertion is \cite[Proposition 11.8.2]{GMM}, and the case $n=\infty$ follows by passage to the inverse limit.
\end{rem}

\subsection{Deformation theory}\label{ss:deformation-theory} We now would like to give a concrete description of the universal deformations for points of $\mr{BT}^{\mc{G},\mu}_\infty$. As we will only need this case, we assume $k$ is a finite field. In particular, we use Notation \ref{nota:finite-field} and the identifications in Proposition \ref{prop:isom-classes-in-BT} without comment.

\begin{nota} Throughout the following we fix $b$ in $\mc{G}(W_j)\sigma(\mu)(p)\mc{G}(W_j)$.
\end{nota}

\subsubsection{An explicit $(\mc{G},\mu)$-aperture}\label{sss:explicit-aperture} We begin by describing very explicitly a $(\mc{G},\mu)$-aperture, or more precisely an object of the category in Definition \ref{defn:FFCrys-on-cover}.

We begin by defining the base formal $W$-scheme on which this $(\mc{G},\mu)$-aperture will live. 

\begin{defn}\label{defn:opposite-unipotent} We define $\mc{U}^-_{\mu,j}$ to be the opposite unipotent over $W_j$ associated to $\mu$:
\begin{equation*}
\displaystyle \mc{U}^-_{\mu,j}\defeq  \left\{g\in \mc{G}_{W_j}:\lim_{t\to 0}\mu(t)^{-1}g\mu(t)=1\right\}.
\end{equation*}
We further define $\widehat{\mc{U}}^-_{\mu,j}$ to be the completion of $\mc{U}^-_{\mu,j}$ at the identity section.
\end{defn}

We will often work ring-theoretically with the objects above, and so we give such notation.

\begin{nota}\label{nota:deformation-ring} Set $R_{\mc{G},\mu,j}\defeq \mc{O}(\wh{\mc{U}}^-_{\mu,j})$ and let $\mf{m}_{\mc{G},\mu,j}$ denote its maximal ideal. Fix an isomorphism
\begin{equation*}
w\colon R_{\mc{G},\mu,j}\simeq W_j\llbracket t_1,\ldots,t_d\rrbracket
\end{equation*}
as in \cite[Lemma 4.2.6]{Ito1}. We write $\sigma$ for the Frobenius lift on $R_{\mc{G},\mu,j}$ restricting to Witt Frobenius on $W_j$ and satisfying $\sigma(t_i)=t_i^p$. Finally, we write $\wh{\Omega}^q_{R_{\mc{G},\mu,j}}$ for the $p$-adic completion of the sheaf of $q$-differentials.
 \end{nota}
 
Finally, important in our construction is the following element of the formal group $\wh{\mc{U}}^-_{\mu,j}$.
 
\begin{nota} We let $U_{\mu,j}^\mr{univ}$ denote the tautological element in $\wh{\mc{U}}^-_{\mu,j}(R_{\mc{G},\mu,j})\subseteq \mc{G}(R_{\mc{G},\mu,j})$.
\end{nota}

\begin{constr}\label{const:deformation} Write $b=k_1\sigma(\mu)(p)k_2$ where $k_1$ and $k_2$ belong to $\mc{G}(W_j)$. Set  $c\defeq \sigma^{-1}(k_2)$ and write $\wt{U}_{\mu,j}^\univ\defeq c^{-1}U_{\mu,j}^\univ c$.\footnote{If $k_2=1$, then $c=1$ and $\wt{U}_{\mu,j}^\univ=U_{\mu,j}^\univ$.} Set:
\begin{equation*}
\mc{Q}_b^\mr{univ}= \mc{G},\qquad \varphi_b^\univ=(\wt{U}^\univ_{\mu,j})^{-1} b,\qquad \Fil^\bullet \mc{Q}_b^\mr{univ}=c^{-1}P_\mu\in (\mc{G}/P_\mu)(R_{\mc{G},\mu,j}),
\end{equation*}
where the last term gives a filtration by Proposition \ref{prop:type-nu-classification} and the fact that $\mc{G}/P_\mu$ parameterizes reductions of structure group of $\mc{G}$ to $P_\mu$.
\end{constr}

\begin{rem}\label{rem:filtration-desc}One can more concretely describe the filtered $\mc{G}$-bundle $\Fil^\bullet \mc{Q}_b^\univ$ in terms of Proposition \ref{prop:Filtered-Rees-equiv}. Namely, it is the filtered Tannakian $\mc{G}$-bundle associating to a representation $\rho\colon\mc{G}\to \GL(\Lambda)$ the module $\Lambda\otimes_{\bb{Z}_p}R_{\mc{G},\mu,j}$ filtered by $\rho(c)^{-1}\cdot (\Fil^\bullet_{\rho\circ\mu}\otimes_W R_{\mc{G},\mu,j})$.
\end{rem}

We now show that there is a unique topologically quasi-nilpotent integrable $\mc{G}$-connection $\nabla$ on $\mc{Q}_b^\mr{univ}$ such that the quadruple $(\mc{Q}_b^\mr{univ},\varphi_b^\mr{univ},\nabla,\Fil^\bullet \mc{Q}_b^\mr{univ})$ is an object of $\cat{FilTors}^{\varphi,\nabla,\mu}_\mc{G}(R_{\mc{G},\mu,j})$.   

To state this cleanly, we first make the following notational definitions.

\begin{nota} We define the operator
\begin{equation*}
\Delta_b\colon \mf{g}[\nicefrac{1}{p}]\otimes_{\Z_p}\wh{\Omega}^1_{R_{\mc{G},\mu,j}/W_j}\to \mf{g}[\nicefrac{1}{p}]\otimes_{\Z_p}\wh{\Omega}^1_{R_{\mc{G},\mu,j}/W_j}
\end{equation*}
by the rule
\begin{equation*}
\Delta_b\left(x\otimes \eta\right)= \varphi_b^\mr{univ}(x)\otimes \,\sigma^\ast(\eta)=\mr{Ad}\left((\wt{U}_{\mu,j}^\univ)^{-1} b\right)(x)\otimes \,\sigma^\ast(\eta).\footnote{Where here $\sigma^\ast$ is literally the pullback map on differentials for the map $\sigma\colon R_{\mc{G},\mu,j}\to R_{\mc{G},\mu,j}$.}
\end{equation*}
\end{nota}

\begin{defn}[{\cite[\S3.1]{Wakabayashi}}] Recall that the \emph{left-invariant Maurer--Cartan form} $\omega_\mc{G}^l$ is the element $\mf{g}\otimes_{\Z_p}\Omega^1_{\mc{G}/\Z_p}$ corresponding to the map 
\begin{equation*}
\omega_\mc{G}^l\colon T_{\mc{G}/\Z_p}\to \mf{g}\otimes_{\Z_p}\mc{O}_\mc{G},\qquad \bigg(v\in T_g(\mc{G})\bigg)\mapsto \bigg((L_{g^{-1}})_\ast(v)\in\mf{g}\bigg),
\end{equation*}
 where $L_{g^{-1}}$ is the left-translation operator. One defines the \emph{right-invariant Maurer--Cartan form} $\omega_\mc{G}^r$ similarly, but by using the right-translation operator $R_{g^{-1}}$.
\end{defn}

\begin{nota}\label{nota:Maurer-Cartan} For a morphism $g\colon \Spf(R)\to \mc{G}$, where $R$ is any $p$-adically complete $W_j$-algebra, we write $g^{-1}dg\defeq g^\ast(\omega_\mc{G}^l)$ and $dgg^{-1}=g^\ast(\omega_\mc{G}^r)$, elements of $\mf{g}\otimes_{\Z_p}\wh{\Omega}^1_{R/W_j}$. 
\end{nota}

\begin{nota} Set $\eta_b^\mr{univ}=(\wt{U}_{\mu,j}^\univ)^{-1}d\wt{U}_{\mu,j}^\univ$, an element of $\mf{g}\otimes_{\Z_p}\wh{\Omega}^1_{R_{\mc{G},\mu,j}/W_j}$.
\end{nota}

We now give an explicit construction of our sought-after connection.

\begin{prop}\label{prop:G-connection-converges} The operator $\Delta_b$ stabilizes $\mf{g}\otimes_{\Z_p}\wh{\Omega}^1_{R_{\mc{G},\mu,j}/W_j}$, and the sum
\begin{equation}\label{eq:univ-G-connection-series}
\nabla_b^\mr{univ}
=\sum_{n\geqslant 0}\Delta_b^n(\eta_b^\univ)
\end{equation}
converges $(t_1,\ldots,t_d)$-adically. Moreover, $\nabla_b^\mr{univ}$ is the unique
solution to the equation
\begin{equation}\label{eq:univ-G-connection-fixed-point}
\eta=\eta_b^\univ+\Delta_b(\eta),\qquad
\eta\in
\mf{g}\otimes_{\Z_p}\wh{\Omega}^1_{R_{\mc{G},\mu,j}/W_j}.
\end{equation}
\end{prop}

\begin{proof}
For notational simplicity, write
\begin{equation*}
R=R_{\mc{G},\mu,j},\qquad I=(t_1,\ldots,t_d),\qquad
U=\wt{U}_{\mu,j}^\mr{univ},\qquad a=U^{-1}b,
\end{equation*}
and put $\mf{g}_R=\mf{g}\otimes_{\Z_p}R$ and $M^q=\mf{g}_R\otimes_R\wh{\Omega}^q_{R/W_j}$. We write $\Delta_{b,q}$ for the operator on $M^q[\nicefrac{1}{p}]$ defined by the same formula as $\Delta_b$.

Write $b=k_1\sigma(\mu)(p)k_2$ as in Construction \ref{const:deformation}. Since $\mu$ is minuscule, the weights of $\mr{Ad}\circ\sigma(\mu)$ on $\mf{g}_{W_j}$ lie in $\{-1,0,1\}$. So, if $x$ is in $\mf{g}_R$ we may write $\mr{Ad}(k_2)(x)=x_{-1}+x_0+x_1$ for the corresponding weight decomposition. We then see that
\begin{equation*}
p\,\mr{Ad}(a)(x)=\mr{Ad}(U^{-1}k_1)(x_{-1}+px_0+p^2x_1),
\end{equation*}
and so $p \mr{Ad}(a)(\mf{g}_R)\subseteq\mf{g}_R$ as $U^{-1}k_1$ is in $\mc{G}(R)$. On the other hand, as $\sigma(t_i)=t_i^p$ one quickly calculates that for every $r\geqslant 0$ and $q\geqslant 1$ one has the containment
\begin{equation*}
\sigma^\ast\bigl(I^r\wh{\Omega}^q_{R/W_j}\bigr)
\subseteq
p^qI^{pr+q(p-1)}\wh{\Omega}^q_{R/W_j}.
\end{equation*}
From this and the containment $p \mr{Ad}(a)(\mf{g}_R)\subseteq\mf{g}_R$ we deduce that
\begin{equation}\label{eq:q-form-contraction}
\Delta_{b,q}(I^rM^q)
\subseteq
p^{q-1}I^{pr+q(p-1)}M^q.
\end{equation}
Taking $r=0$ and $q=1$ shows that $\Delta_b$ stabilizes $M^1=\mf{g}\otimes_{\Z_p}\wh{\Omega}^1_{R/W_j}$, and that  $\Delta_b^n(M^1)\subseteq I^{p^n-1}M^1$ for every $n\geqslant 0$. Since $M^1$ is $I$-adically complete, we deduce that $\nabla_b^\mr{univ}$ converges. Moreover, by the $I$-adic continuity of $\Delta_b$ we see that we have the equality
\begin{equation}\label{eq:Deltab-applied-to-nabla-univ}
\eta_b^\mr{univ}+\Delta_b(\nabla_b^\mr{univ})=\eta_b^\univ+\sum_{n\geqslant 0}\Delta_b^{n+1}(\eta_b^\univ)=\sum_{n\geqslant 0}\Delta_b^n(\eta_b^\univ)=\nabla_b^\univ,
\end{equation}
verifying that \eqref{eq:univ-G-connection-fixed-point} holds.

For the unicity of a solution to \eqref{eq:univ-G-connection-fixed-point}, take two solutions $\eta$ and $\eta'$. Setting $\delta=\eta-\eta'$ we have $\delta=\Delta_b(\delta)$. The discussion in the previous paragraph then shows that $\delta=\Delta_b^n(\delta)$ and thus belongs to $I^{p^n-1}M^1$ for all $n$. As $M^1$ is $I$-adically complete we deduce that $\delta=0$, as desired. 
\end{proof}

Note that using {Example \ref{eg:trivial-torsor-connections}}, we may interpret $\nabla_b^\univ$ as a $\mc{G}$-connection on the trivial $\mc{G}$-torsor over $R_{\mc{G},\mu,j}$. With this interpretation, we have the following.

\begin{prop}\label{prop:integrable-nilpotent-horizontal} The $\mc{G}$-connection $\nabla_b^\univ$ is integrable and topologically quasi-nilpotent, and the map $\varphi_b^\univ$ is $\nabla_b^\univ$-horizontal.
\end{prop}
\begin{proof} We use the same notational shortcuts as at the beginning of the proof of Proposition \ref{prop:G-connection-converges}.

Let us first verify that $\varphi_b^\univ$ is $\nabla_b^\univ$-horizontal. First, note that for $\nu$ in $\mf{g}\otimes_{\Z_p}\wh{\Omega}^1_{R_{\mc{G},\mu,j}/W_j}$, $\nabla_\nu$-horizontality of $\varphi_b^\univ$ is equivalent to the identity
\begin{equation}\label{eq:horizontality-gauge-identity}
\nu=-da\,a^{-1}+\mr{Ad}(a)(\sigma^\ast\nu).
\end{equation}
 Indeed, after applying a representation $\rho$, and writing $A=\rho(a)$ and $\Omega=\rho_\ast(\nu)$, this becomes the usual identity for horizontality:
\begin{equation*}
dA+\Omega A=A\sigma^\ast(\Omega).
\end{equation*}
Now, since $a=U^{-1}b$ and the element $b$ is constant over $W_j$, one can simplify the above to:
\begin{equation*}
-da\,a^{-1}=(\wt{U}_{\mu,j}^\univ)^{-1}d\wt{U}_{\mu,j}^\univ=\eta_b^\univ.
\end{equation*}
Thus, for $\nu=\nabla_b^\mr{univ}$, \eqref{eq:horizontality-gauge-identity} reduces to the previously verified \eqref{eq:univ-G-connection-fixed-point}.

Now, we prove that $\nabla_b^\univ$ is integrable. Let $K$ be the curvature tensor of $\nabla_b^\univ$:
\begin{equation*}
K=d\nabla_b^\mr{univ}
+\frac{1}{2}[\nabla_b^\mr{univ},\nabla_b^\mr{univ}].
\in M^2
\end{equation*}
Using \eqref{eq:Deltab-applied-to-nabla-univ} and the integrability of $\eta_b^\univ$, one checks that $K=\Delta_{b,2}(K)$. Thus, the same argument in the last paragraph of the proof of Proposition \ref{prop:G-connection-converges} shows that $K=0$, as desired.

It remains to check topological quasi-nilpotence. Fix a representation $\rho\colon\mc{G}\to\GL(\Lambda)$, and set
\begin{equation*}
F=\rho(a)\colon
\sigma^\ast \Lambda_{R}[\nicefrac{1}{p}]
\xrightarrow{\ \sim\ }\Lambda_{R}[\nicefrac{1}{p}].
\end{equation*}
Let $D_i$ be the covariant derivative on ${\Lambda_R}$ in the direction $\tfrac{\partial}{\partial t_i}$, and let $D_i^\sigma$ be the corresponding operator on $\sigma^\ast {\Lambda_R}$. Choose $h\geqslant 0$ such that both $p^hF$ and $p^hF^{-1}$ preserve the natural integral lattices. Since $\sigma^\ast(dt_i)$ is divisible by $p$, the pullback connection on $\sigma^\ast {\Lambda_R}$ reduces modulo $p$ to the trivial connection. As $\left(\tfrac{\partial}{\partial t_i}\right)^p=0$ on $R/p$, it follows that for each $s\geqslant 0$:
\begin{equation*}
(D_i^\sigma)^p(\sigma^\ast {\Lambda_R})\subseteq p\sigma^\ast {\Lambda_R},
\qquad
(D_i^\sigma)^{ps}(\sigma^\ast {\Lambda_R})\subseteq p^s\sigma^\ast {\Lambda_R}.
\end{equation*}
Horizontality of $F$ gives $D_i^nF=F(D_i^\sigma)^n$, and so for $x$ in ${\Lambda_R}$ and $y=p^hF^{-1}(x)$ in $\sigma^\ast {\Lambda_R}$, one has
\begin{equation*}
D_i^{ps}(x)
=p^{-h}F\bigl((D_i^\sigma)^{ps}(y)\bigr)
\in p^{s-2h}{\Lambda_R}.
\end{equation*}
Given $N\geqslant 1$, choose $s\geqslant N+2h$, so then $D_i^{ps}({\Lambda_R})\subseteq p^N{\Lambda_R}$ for every $i$. Since the connection is integrable, the operators $D_1,\ldots,D_d$ commute. Thus, if $\alpha=(\alpha_1,\ldots,\alpha_d)$ satisfies $|\alpha|>d(ps-1)$, then some $\alpha_i\geqslant ps$, and hence $D_1^{\alpha_1}\cdots D_d^{\alpha_d}({\Lambda_R})\subseteq p^N{\Lambda_R}$. Since $N$ and $\rho$ were arbitrary, the conclusion follows.\end{proof}

We finally arrive at our desired statement about quadruples in $\cat{FilTors}^{\varphi,\nabla,\mu}_\mc{G}(R_{\mc{G},\mu,j})$ filling out the triple $(\mc{Q}_b^\univ,\varphi_b^\univ,\Fil^\bullet \mc{Q}_b^\univ)$.

\begin{cor}\label{cor:universal-G-connection-characterization}
The unique $\mc{G}$-connection $\nabla$ on $\mc{Q}_b^\mr{univ}$ such that $\big(\mc{Q}_b^\mr{univ},\varphi_b^\mr{univ},\nabla,\Fil^\bullet\mc{Q}_b^\mr{univ}\big)$ is an object of $\cat{FilTors}^{\varphi,\nabla,\mu}_\mc{G}(R_{\mc{G},\mu,j})$ is $\nabla_b^\univ$. Moreover, this quadruple is strongly divisible.
\end{cor}
\begin{proof}
By Propositions \ref{prop:G-connection-converges} and \ref{prop:integrable-nilpotent-horizontal}, and the fact that as $\mu$ is minuscule Griffiths transversality is automatic, we know that $\big(\mc{Q}_b^\mr{univ},\varphi_b^\mr{univ},\nabla_b^\univ,\Fil^\bullet\mc{Q}_b^\mr{univ}\big)$ is an object of $\cat{FilTors}^{\varphi,\nabla,\mu}_\mc{G}(R_{\mc{G},\mu,j})$.

Conversely, suppose that $\big(\mc{Q}_b^\mr{univ},\varphi_b^\mr{univ},\nabla,\Fil^\bullet\mc{Q}_b^\mr{univ}\big)$ is an object of $\cat{FilTors}^{\varphi,\nabla,\mu}_\mc{G}(R_{\mc{G},\mu,j})$. Then, writing $\nabla=\nabla_\eta$ for some $\eta$ in $\mf{g}\otimes_{\Z_p}\wh{\Omega}^1_{R_{\mc{G},\mu,j}/W_j}$ we deduce from the first paragraph of the proof of Proposition \ref{prop:integrable-nilpotent-horizontal} that $\varphi_b^\univ$ being $\nabla$-horizontal implies that $\eta=\eta_b^\univ+\Delta_b(\eta)$. Thus, we are done by the unicity assertion in Proposition \ref{prop:G-connection-converges}.

To prove strong divisibility, take a representation
$\rho\colon\mc{G}\to\GL(\Lambda)$ and set $M=\Lambda\otimes_{\Z_p}R_{\mc{G},\mu,j}$.
Give $M$ the filtration induced by
$\Fil^\bullet\mc{Q}_b^{\mr{univ}}$.
Using the identification $\sigma^*M\simeq M$, we have
\begin{equation*}
 \sum_{r\in\Z}p^{-r}\sigma^*(\Fil^r M) = \rho\bigl(k_2^{-1}\sigma(\mu)(p)^{-1}\bigr)M.
\end{equation*}
Writing $F=\rho((\wt{U}_{\mu,j}^\univ)^{-1}b)$
for the induced Frobenius under this identification, we obtain
\begin{equation}\label{eq:str-div-group-proof}
 \begin{aligned}
 F\left(\sum_{r\in\Z}p^{-r}\sigma^*(\Fil^r M)\right)
 &=
 \rho\bigl((\wt{U}_{\mu,j}^\univ)^{-1}b
             k_2^{-1}\sigma(\mu)(p)^{-1}\bigr)M\\
 &=\rho\bigl((\wt{U}_{\mu,j}^\univ)^{-1}k_1\bigr)M\\
 &=M,
 \end{aligned}
\end{equation}
since $(\wt{U}_{\mu,j}^\univ)^{-1}k_1$ belongs to
$\mc{G}(R_{\mc{G},\mu,j})$, as desired.
\end{proof}

We finally give notation for the aperture associated to the above quadruple.

\begin{nota}\label{nota:universal-aperture} We let $\mf{Q}_b^\univ$ be the $(\mc{G},\mu)$-aperture over $\wh{\mc{U}}^-_{\mu,j}$ corresponding to the quadruple $\big(\mc{Q}_b^\mr{univ},\varphi_b^\mr{univ},\nabla_b^\univ,\Fil^\bullet\mc{Q}_b^\mr{univ}\big)$ under the equivalences from Proposition \ref{prop:torsors-with-connection-equiv-filtered} and Theorem \ref{thm:BT-over-base-ring}.
\end{nota}

\subsubsection{A description of the universal deformation} We now use the explicitly described $(\mc{G},\mu)$-aperture from above to describe the universal deformation of $b$. 

\begin{defn} Let $\mc{Y}$ be a formal prestack over $W$. Let $y$ be a point of $\mc{Y}(k')$ for some perfect extension $k'/k$. Set $W'=W(k')$. Define the \emph{deformation prestack} associated to the pair $(\mc{Y},y)$ to be the functor on $\mc{C}_{W'}$ (the completed base category)\footnote{More precisely, $\mc{C}_{W'}$ is the category of complete Noetherian local rings $R$ equipped with a local map $W'\to R$ satisfying $k'\isomto R/\mf{m}_R$, with morphisms local $W'$-algebra maps.} given by:
\begin{equation*}
    \mr{Def}_{y}\colon \mc{C}_{W'}\to\cat{Grpd},\quad R\mapsto \mr{fib}(\mc{Y}(R)\to \mc{Y}(k');y).
\end{equation*}
If $\mr{Def}_y$ is representable, we call the representing space the \emph{deformation space} of the pair $(\mc{Y},y)$, and the tautological object over that space the \emph{universal deformation}.
\end{defn}

\begin{nota}\label{nota:defm-space} For $b$, thought of as an object of $\mr{BT}^{\mc{G},\mu}_\infty(k_j)$, we write $\mf{D}(\mc{G},b,\mu)$ for the deformation functor.
\end{nota}

\begin{obs} Pulling back $\mf{Q}_b^\univ$ to $k_j$ gives us a $(\mc{G},\mu)$-aperture over $k_j$. Under the equivalence from Setup \ref{setup:simplifying-over-finite-fields}, this $(\mc{G},\mu)$-aperture corresponds to the object of $\cat{Tors}_\mc{G}^{\varphi,\sigma(\mu)}(W_j)$ given by 
\begin{equation*}
(\mc{Q}_b^\univ,\varphi_b^\univ)\otimes_{R_{\mc{G},\mu,j}}W_j=(\mc{G}_{W_j},b).
\end{equation*}
In other words, it is the $(\mc{G},\mu)$-aperture associated to $b$. Thus, we obtain a natural isomorphism
\begin{equation*}
\iota_b^\univ\colon \mf{Q}_b^\univ\otimes_{R_{\mc{G},\mu,j}}k_j\isomto b,
\end{equation*}
and thus the pair $(\mf{Q}_b^\univ,\iota_b^\univ)$ gives a morphism $\rho_{b}^\mr{univ}\colon \wh{\mc{U}}^-_{\mu,j}\to \mf{D}(\mc{G},b,\mu)$.
\end{obs}

The following shows $(\mf{Q}_b^\univ,\iota_b^\univ)$ is indeed the universal deformation of $b$ in $\mr{BT}^{\mc{G},\mu}_\infty$.

\begin{thm}\label{thm:explicit-deformation} The map $\rho_{b}^\mr{univ}\colon \wh{\mc{U}}^-_{\mu,j}\to \mf{D}(\mc{G},b,\mu)$
is an isomorphism.
\end{thm}
\begin{proof} This result is essentially already contained (with some duplication between the references) in \cite[Theorem 4.4.2]{ItoDeformation}, \cite[Proposition 3.32]{IKY2}, and \cite[Proposition 10.2.9]{GMM}. But, let us spell this out a bit more. In \cite[Proposition 3.32]{IKY2} it is explained that the $F$-gauge with $\mc{G}$-structure over $R_{\mc{G},\mu,j}$ constructed in \cite[Theorem 4.4.2]{ItoDeformation}, which we shall denote $(\mf{Q}^\mr{Ito},\iota^\mr{Ito})$, is actually universal (as opposed to just universal on power series rings as verified in loc.\@ cit.\@). Thus, it suffices to verify that $(\mf{Q}^\mr{Ito},\iota^\mr{Ito})$ is isomorphic to $(\mf{Q}_b^\univ,\iota_b^\univ)$.  

By Proposition \ref{prop:torsors-with-connection-equiv-filtered} and Theorem \ref{thm:BT-over-base-ring}, it suffices to show that the object of $\cat{FilTors}^{\varphi,\nabla,\mu}_\mc{G}(R_{\mc{G},\mu,j})$ associated to $(\mf{Q}^\mr{Ito},\iota^\mr{Ito})$ is $\big(\mc{Q}_b^\mr{univ},\varphi_b^\mr{univ},\nabla_b^\univ,\Fil^\bullet\mc{Q}_b^\mr{univ}\big)$. But, given the presentation of $(\mf{Q}^\mr{Ito},\iota^\mr{Ito})$ given in \cite[Theorem 4.4.2]{ItoDeformation}, one easily checks that this is true for all parts of this tuple except the $\mc{G}$-connection $\nabla_b^\univ$. But, we are then done by the unicity from Corollary \ref{cor:universal-G-connection-characterization}.
\end{proof}

\begin{rem}\label{rem:base-change-deformation}Theorem \ref{thm:explicit-deformation} does not require that $b$ is defined over a finite field, and so from this it's simple to see the following: For any perfect extension $\kappa/k_j$, let $b|_\kappa$ be the induced element of $\mr{BT}^{\mc{G},\mu}_\infty(\kappa)$. Then the pullback of $(\mf{Q}_b^\mr{univ},\iota^\mr{univ})$ to $\wh{\mc{U}}^-_{\mu,j}\otimes W(\kappa)$ is the deformation of $b|_\kappa$. 
\end{rem}

\begin{rem}\label{rem:deformation-space-description-on-power-series-rings} To understand the full deformation theory contained in $\mf{D}(\mc{G},b,\mu)$ and the universal deformation living over it one requires the notion of $(\mc{G},\mu)$-apertures. But, one can still characterize it in a more concrete way: for an integer $m\geqslant 0$, one has a bijection
\begin{equation*}
\begin{tikzcd}[cramped]
	{\Hom\bigg(\Spf\big(W_j\ll x_1,\ldots,x_m\rr, (p,x_1,\ldots,x_m)\big),\mf{D}(\mc{G},b,\mu)\bigg)} \\
	\begin{array}{c} \mr{fib}\bigg(\cat{FilTors}^{\varphi,\mr{sd},\nabla,\mu}_\mc{G}(W_j\ll x_1,\ldots,x_m\rr)\to \cat{Tors}^{\varphi,\sigma(\mu)}_\mc{G}(W_j);(\mc{G}_{W_j},b)\bigg), \end{array}
	\arrow["\wr", from=1-1, to=2-1]
\end{tikzcd}
\end{equation*}
given by $f\mapsto f^\ast\big(\mc{Q}_b^\mr{univ},\varphi_b^\mr{univ},\nabla_b^\univ,\Fil^\bullet\mc{Q}_b^\mr{univ}\big)$. This description can be made even more concrete by choosing a tensor package and phrasing things in terms of the quintuples as in Definition \ref{defn:quint}.

That this bijection indeed characterizes the universal deformation follows from the fact that $R_{\mc{G},\mu,j}$ is itself a power series ring. This description is in line with the original work of Faltings in \cite{FaltingsVeryRamified}, and the approach taken by Ito in \cite{ItoDeformation}. 
\end{rem}

We end this section by giving a more direct description of $\mf{D}(\mc{G},b,\mu)(R)$ for an arbitrary Noetherian adic $W_j$-algebra $R$.

\begin{prop}\label{prop:functorial-description-deformation} Let $R$ be a Noetherian adic $W_j$-algebra. There is a functorial identification 
\begin{equation*}
\mf{D}(\mc{G},b,\mu)(R)\simeq \mr{fib}\left(\mr{BT}^{\mc{G},\mu}_\infty(R)\to \mr{BT}^{\mc{G},\mu}_\infty(R_\mr{red});b|_{R_{\mr{red}}}\right),\footnote{That this fiber groupoid is discrete follows as in the proof of \cite[Proposition 3.31]{IKY2}.}
\end{equation*}
given by pulling back the object $(\mf{Q}_b^\mr{univ},\iota_b^\univ)$. 
\end{prop}
\begin{proof} If $R_n=R/(R^{\circ\circ})^n$, the map under consideration is the limit of the maps
\begin{equation*}
\mf{D}(\mc{G},b,\mu)(R_n)\simeq \mr{fib}\left(\mr{BT}^{\mc{G},\mu}_\infty(R_n)\to \mr{BT}^{\mc{G},\mu}_\infty(R_\mr{red});b|_{R_{\mr{red}}}\right),\footnote{Here we are implicitly using the integrability of $\mr{BT}^{\mc{G},\mu}_n$ for each $n$; see the discussion in Remark \ref{rem:integrability}.}
\end{equation*}
and so it suffices to show that each of these is a bijection. We can prove this by induction on $n$. When $n=1$ both sides are singletons, and so the claim is trivial. Suppose we have proven the claim for $n$. Then, it suffices to show that for each $y$ in $\mf{D}(\mc{G},b,\mu)(R_n)$, with the same symbol for its images in the fiber set, the map
\begin{equation*}
\mr{fib}\bigg(\mf{D}(\mc{G},b,\mu)(R_{n+1})\to \mf{D}(\mc{G},b,\mu)(R_n);y\bigg)\to \mr{fib}\bigg(\mr{BT}^{\mc{G},\mu}_\infty(R_{n+1})\to \mr{BT}^{\mc{G},\mu}_\infty(R_n);y\bigg)
\end{equation*}
is bijective. Set $J=R^{\circ\circ}$ and $I_n=J^n/J^{n+1}$, so that $JI_n=0$. By the usual description of infinitesimal liftings and
Grothendieck--Messing theory (see Theorem~\ref{thm:GMM}), respectively,
these fibers are torsors respectively under
\begin{equation*}
 I_n\otimes_{k_j}T_0\mf{D}(\mc{G},b,\mu)
 \qquad\text{and}\qquad
 I_n\otimes_{k_j}T_b\mr{BT}^{\mc{G},\mu}_\infty,
\end{equation*}
where the tangent spaces are relative to $W_j$. The map between the fibers is equivariant for the tangent map induced by the universal aperture, which is an isomorphism by Theorem \ref{thm:explicit-deformation}. The inductive step, and thus the full claim, follows.
\end{proof}

\begin{sectionappendices}

\sectionappendix{Appendix: Tannakian theory of (filtered) torsors}
\label{sec:filtered-Tannakian-appendix}

In this appendix we recall the theory of filtered torsors in the context that is frequently used in this section. See \cite[Appendix A]{WedhornTannakian} and \cite{Ziegler} for generalizations of some of these results to more general group schemes and over more general bases.

\begin{nota}\label{nota:Filtered-torsors} Throughout this appendix we use the following notation:
\begin{itemize}
\item $R$ is a Dedekind domain;
\item $\mc{H}$ is a finite type affine flat group $R$-scheme;
\item $\cat{Rep}_R(\mc{H})$ denotes the exact $R$-linear $\otimes$-category (see Definition \ref{defn:exact-R-linear} below) of (left) representations $\rho\colon \mc{H}\to \mr{GL}(\Lambda)$ where $\Lambda$ is a finite projective $R$-module;
\item $(\ms{X},\ms{O})$ is a ringed topos of $R$-algebras;
\item $\mc{H}_\ms{O}$ is the group object of $\ms{X}$ defined by $\mc{H}_\ms{O}(T)\defeq \mc{H}(\ms{O}(T))$.
\end{itemize}
\end{nota}

\begin{eg} Suppose that $\mc{X}$ is a prestack over $R$ (e.g., in increasing level of generality, a scheme, a formal scheme, or a formal stack). Then, we can naturally build a ringed topos $(\mc{X}_\mr{fppf},\mc{O}_\mc{X})$ whose objects are morphisms $\Spec(A)\to \mc{X}$ for a ring $A$, whose morphisms are maps over $\mc{X}$, and whose covers are given by fppf covers. Then, there is a natural identification
\begin{equation*}
\cat{Vect}(\mc{X}_\mr{fppf},\mc{O}_\mc{X})\simeq \twolim_{\Spec(A)\to\mc{X}}\cat{Vect}(A).
\end{equation*}
\end{eg}

\begin{defn}\label{defn:exact-R-linear} An \emph{exact $R$-linear $\otimes$-category} is a category $\mc{C}$ endowed with the following structure and properties:
\begin{itemize}
\item $\mc{C}$ is equipped with the structure of an exact category; see \cite[Appendix A]{Keller};
\item $\mc{C}$ is equipped with an additive $R$-linear structure (see \stacks{0104} and \stacks{09MI}); 
\item the underlying category of $\mc{C}$ is Karoubian (see \stacks{09SF});
\item there is an $R$-bilinear symmetric monoidal structure $\otimes\colon \mc{C}\times \mc{C}\to \mc{C}$ (see \stacks{0FFJ}).
\end{itemize}

For exact $R$-linear $\otimes$-categories $\mc{C}$ and $\mc{D}$, a functor $F\colon \mc{C}\to\mc{D}$ is an \emph{exact $R$-linear $\otimes$-functor} if it preserves exact sequences and it is $R$-linear (see \stacks{09MK}) and symmetric monoidal (see \stacks{0FFL} and \stacks{0FFY}).
\end{defn}

\begin{defn} Let $\mc{C}$ be an exact $R$-linear $\otimes$-category. Then, the category of \emph{$\mc{H}$-objects in $\mc{C}$} is the category $\mc{H}\text{-}\mc{C}$ consisting of exact $R$-linear $\otimes$-functors $\omega\colon \cat{Rep}_R(\mc{H})\to \mc{C}$ with morphisms given by monoidal natural transformations.
\end{defn}

\begin{defn} We define the following notions.
\begin{itemize} 
\item The category of \emph{$\mc{H}$-bundles on $(\ms{X},\ms{O})$}, denoted $\cat{Tors}_\mc{H}(\ms{X},\ms{O})$ (or $\cat{Tors}_\mc{H}(\ms{X})$ when $\ms{O}$ is clear from context), is the category of $\mc{H}_\mc{O}$-torsors on $\ms{X}$.
\item The category of \emph{Tannakian $\mc{H}$-bundles on $(\ms{X},\ms{O})$}, denoted $\mc{H}\text{-}\cat{Vect}(\ms{X},\ms{O})$ (or $\mc{H}\text{-}\cat{Vect}(\ms{X})$ when $\ms{O}$ is clear from context), is the category of $\mc{H}$-objects in $\cat{Vect}(\ms{X},\ms{O})$.
\item A Tannakian $\mc{H}$-bundle $\omega$ is \emph{locally trivial} if for every object $T$ of $\ms{X}$ there exists a cover $\{T_i\to T\}$ such that $\omega|_{T_i}\simeq \omega_\mr{triv}|_{T_i}$ for all $i$, where $\omega_\mr{triv}(\Lambda)\defeq \Lambda\otimes_R \ms{O}$. We denote the full subcategory of locally trivial objects of $\mc{H}\text{-}\cat{Vect}(\ms{X},\ms{O})$ by $\mc{H}\text{-}\cat{Vect}^\mr{lt}(\ms{X},\ms{O})$. 
\end{itemize}
\end{defn}

\begin{constr}There is a natural functor 
\begin{equation*}
\cat{Tors}_\mc{H}(\ms{X},\ms{O})\to\mc{H}\text{-}\cat{Vect}(\ms{X},\ms{O}),\qquad Q\mapsto \omega_Q
\end{equation*}
where 
\begin{equation*}
\omega_Q(\Lambda)\defeq Q^{(\Lambda)}\defeq Q\wedge^{\mc{H}_\ms{O}}(\Lambda\otimes_R\ms{O}),
\end{equation*}
or, equivalently, $\omega_Q(\Lambda)$ is the natural vector bundle associated to the $\mr{GL}(\Lambda)_\mc{O}$-torsor $\rho_\ast Q$, thought of as a vector bundle in the usual way (e.g., see \cite[Proposition A.4]{IKY1}).
\end{constr}

\begin{prop}[{see \cite[Theorem A.19 and Remark A.20]{IKY1}}]\label{prop:Tannakian-equiv} The functor
\begin{equation*}
\cat{Tors}_\mc{H}(\ms{X},\ms{O})\to\mc{H}\text{-}\cat{Vect}^\mr{lt}(\ms{X},\ms{O}),\qquad Q\mapsto \omega_Q
\end{equation*}
is an equivalence with quasi-inverse given by
\begin{equation*}
\omega\mapsto \bigg(\underline{\mr{Isom}}(\omega_\mr{triv},\omega)\colon T\mapsto \underline{\mr{Isom}}(\omega_\mr{triv}|_T,\omega|_T)\bigg),
\end{equation*}
and every Tannakian $\mc{H}$-bundle is locally trivial if $(\ms{X},\ms{O})=(\mc{X}_\mr{fppf},\mc{O}_\mc{X})$ for a prestack $\mc{X}$ over $R$.
\end{prop}

We would now like to discuss certain filtered upgrades of the above story. 

\begin{defn} Let $\mc{X}$ be a prestack over $R$. A \emph{filtered $\mc{H}$-bundle on $\mc{X}$} is an $\mc{H}$-bundle on the product prestack $\mc{X}\times\bb{A}^1/\bb{G}_m$.\footnote{Here we have $\bb{G}_m$ acting on $\bb{A}^1$ with the usual action: $\lambda\cdot x:=\lambda x$. Equivalently, we are endowing $\bb{Z}[t]$ with the grading where $t$ has degree $1$.} Pulling back along $\mc{X}=\mc{X}\times\bb{G}_m/\bb{G}_m\to \mc{X}\times \bb{A}^1/\bb{G}_m$ gives a $\mc{H}$-bundle $Q$ on $\mc{X}$ which we call the \emph{underlying $\mc{H}$-bundle}. We often then denote the filtered $\mc{H}$-bundle on $\mc{X}$ by $\Fil^\bullet\mc{Q}$. We denote the category of filtered $\mc{H}$-bundles by $\cat{FilTors}_\mc{H}(\mc{X})$.
\end{defn}

\begin{rem} Putting a filtration on an $\mc{H}$-bundle $Q$ on $\mc{X}$, i.e., giving a filtered $\mc{H}$-bundle of the form $\Fil^\bullet Q$, can be succinctly understood as extending the morphism $\mc{X}\to B\mc{H}$ corresponding to $\mc{Q}$ along the natural map $\mc{X}=\mc{X}\times\bb{G}_m/\bb{G}_m\to \mc{X}\times\bb{A}^1/\bb{G}_m$.
\end{rem}

To justify the name `filtered' $\mc{H}$-bundle, let us make the following definition.

\begin{defn} Let $(\ms{X},\ms{O})$ be a ringed topos over $R$.
\begin{itemize} 
\item A \emph{filtered vector bundle on $(\ms{X},\ms{O})$} consists of a vector bundle $\mc{V}$ on $(\ms{X},\ms{O})$ together with an exhaustive separated decreasing $\bb{Z}$-filtered sequence of $\ms{O}$-submodules $\Fil^\bullet\mc{V}$. 
\item A filtered vector bundle $\Fil^\bullet\mc{V}$ is \emph{locally split} if for all $i$ the $\ms{O}$-module 
\begin{equation*}
\mr{Gr}^i\mc{V}\defeq \Fil^i\mc{V}/\Fil^{i+1}\mc{V}
\end{equation*}
 is a vector bundle. 
\end{itemize}
We denote the category of filtered vector bundles on $(\ms{X},\ms{O})$ by $\cat{FilVect}(\ms{X},\ms{O})$ (or just $\cat{FilVect}(\ms{X})$ when $\ms{O}$ is clear from context), and we write $\cat{FilVect}^\mr{ls}(\ms{X},\ms{O})$ for the full subcategory of locally split filtered vector bundles.
\end{defn}

\begin{defn}Endow $\cat{FilVect}(\ms{X},\ms{O})$ with the structure of an exact $R$-linear $\otimes$-category where we define $\Fil^\bullet \mc{V}\otimes \Fil^\bullet\mc{W}$ to have underlying vector bundle $\mc{V}\otimes_{\ms{O}}\mc{W}$ and filtration
\begin{equation*}
\Fil^k(\mc{V}\otimes\mc{W})=\sum_{i+j=k}\mr{im}\bigg(\Fil^i\mc{V}\otimes_{\ms{O}}\Fil^j\mc{W}\to \mc{V}\otimes_\ms{O}\mc{W}\bigg)
\end{equation*}
and where a sequence
\begin{equation*}
0\to \Fil^\bullet \mc{V}_1\to\Fil^\bullet\mc{V}_2\to\Fil^\bullet\mc{V}_3\to 0
\end{equation*}
is exact if the sequence of $\ms{O}$-modules
\begin{equation*}
0\to \Fil^i \mc{V}_1\to\Fil^i\mc{V}_2\to\Fil^i\mc{V}_3\to 0
\end{equation*}
is exact for all $i$. Note that $\cat{FilVect}^\mr{ls}(\ms{X},\ms{O})$ is an exact $R$-linear $\otimes$-subcategory of $\cat{FilVect}(\ms{X},\ms{O})$. 
\end{defn}

\begin{defn} The category of $\mc{H}$-objects in $\cat{FilVect}(\ms{X},\ms{O})$ is called the category of \emph{filtered Tannakian $\mc{H}$-bundles} and denoted $\mc{H}\text{-}\cat{FilVect}(\ms{X},\ms{O})$.
\end{defn}

\begin{rem} This name is apt as every such object gives rise to an underlying $\mc{H}$-bundle $\omega$, so we may then think of the filtered $\mc{H}$-bundle $\Fil^\bullet\omega$ as giving a compatible system of filtrations $\Fil^\bullet \omega(\Lambda)$ on $\omega(\Lambda)$ for all $\Lambda$. We note that, as an exact $R$-linear $\otimes$ functor sends dualizable objects to dualizable objects, in fact every object of $\mc{H}\text{-}\cat{FilVect}(\ms{X},\ms{O})$ automatically factorizes through $\mc{H}\text{-}\cat{FilVect}^\mr{ls}(\ms{X},\ms{O})$, i.e., we have the equality $\mc{H}\text{-}\cat{FilVect}^\mr{ls}(\ms{X},\ms{O})=\mc{H}\text{-}\cat{FilVect}(\ms{X},\ms{O})$.
\end{rem}

\begin{constr}For a prestack $\mc{X}$, there is a natural functor
\begin{equation}\label{eq:Rees}
 \mr{Rees}\colon\cat{Vect}(\mc{X}\times \bb{A}^1/\bb{G}_m)\to \cat{FilVect}(\mc{X})
\end{equation}
To describe this functor it suffices, by passing to the limit, to deal with the case when $\mc{X}=\Spec(A)$. In this case, we see that $\cat{Vect}(\mc{X}\times\bb{A}^1/\bb{G}_m)$ corresponds exactly to the category of finite projective graded $A[t]$-modules $N=\bigoplus_{i\in\bb{Z}}N_i$ where $t$ has weight $1$.\footnote{Here we are using the very well-known correspondence between graded modules over a graded ring $A$ and $\bb{G}_m$-equivariant modules on $\Spec(A)$; see \stacks{03LE}.} Then, we associate to $N$ the $A$-module $M=N/(t-1)$ with the filtration $\Fil^r M=(N_{-r}+(t-1)N)/(t-1)N$.
\end{constr}

We then have the following classical result, sometimes called the \emph{Rees equivalence}.

\begin{prop}[{see \cite[Proposition 1.5]{IKY2} and the surrounding references}]\label{prop:Filtered-Rees-equiv} For a prestack $\mc{X}$ over $R$, the functor $\mr{Rees}$ from \eqref{eq:Rees} induces an exact $R$-linear $\otimes$-equivalence
\begin{equation*}
\mr{Rees}\colon \cat{Vect}(\mc{X}\times\bb{A}^1/\bb{G}_m)\isomto \cat{FilVect}^\mr{ls}(\mc{X}),
\end{equation*}
and thus induces a commuting diagram 
\begin{equation*}
\begin{tikzcd}[cramped,column sep=13em]
	{\cat{FilTors}_\mc{H}(\mc{X})} & {\mc{H}\text{-}\cat{FilVect}(\mc{X})} \\
	{\cat{Tors}_\mc{H}(\mc{X})} & {\mc{H}\text{-}\cat{Vect}(\mc{X})}
	\arrow["{\Fil^\bullet Q\mapsto \mr{Rees}(\Fil^\bullet Q)=\Fil^\bullet \omega_Q}", from=1-1, to=1-2]
	\arrow[from=1-1, to=2-1]
	\arrow[from=1-2, to=2-2]
	\arrow["{Q\mapsto \omega_Q}"', from=2-1, to=2-2]
\end{tikzcd}
\end{equation*}
where the horizontal maps are equivalences.
\end{prop}

We end this section by singling out a specific type of filtered torsor. 

\begin{constr} Let $\nu\colon \bb{G}_{m,R}\to \mc{H}$ be a cocharacter. By pushing out the tautological $\bb{G}_{m,R}$-torsor on $B\bb{G}_{m,R}$ along $\nu^{-1}$ we obtain a $\mc{H}$-torsor $\mr{Gr}^\bullet_\nu\mc{H}$ on $B\bb{G}_{m,R}=\Spec(R)\times B\bb{G}_m$. Pulling $\mr{Gr}^\bullet_\nu\mc{H}$ back along the projection $\Spec(R)\times \bb{A}^1/\bb{G}_m\to \Spec(R)\times B\bb{G}_m$ we obtain a filtered $\mc{H}$-bundle on $\Spec(R)$ which we denote $\Fil^\bullet_\nu \mc{H}$.
\end{constr}

\begin{defn}\label{defn:type-nu} Let $\mc{X}$ be a prestack over $R$ and $\nu\colon \bb{G}_{m,R}\to\mc{H}$ a cocharacter. Then, a filtered $\mc{H}$-bundle $\Fil^\bullet\mc{Q}$ on $\mc{X}$ is \emph{of type $\nu$} if it is locally on $\mc{X}_{\mr{fppf}}$ isomorphic to $\Fil^\bullet_\nu \mc{H}$.\footnote{Note that this is not equivalent to asking that $\Fil^\bullet Q$ and $\Fil^\bullet_\nu \mc{H}$ are locally isomorphic on $(\mc{X}\times\bb{A}^1/\bb{G}_m)_\mr{fppf}$ which is, in fact, automatic as both are locally trivial on that site.} We denote the full subcategory of $\cat{FilTors}_\mc{H}(\mc{X})$ consisting of objects of type $\nu$ by $\cat{FilTors}^\nu_\mc{H}(\mc{X})$.
\end{defn}

\begin{rem}\label{rem:filtered-Tannakian-torsor-type-mu} We can more concretely understand $\Fil^\bullet_\nu \mc{H}$ using the equivalence given in Proposition \ref{prop:Filtered-Rees-equiv}. Namely, $\Fil^\bullet_\nu\mc{H}$ corresponds to the filtered Tannakian $\mc{H}$-bundle $\Fil^\bullet_\nu\omega_\mr{triv}$ on $\Spec(R)$ sending a representation $\rho\colon \mc{H}\to\GL(\Lambda)$ to $\Lambda$ with the filtration
\begin{equation*}
\Fil^r_\nu \Lambda\defeq \bigoplus_{j\geqslant r}\Lambda(j),
\end{equation*}
where $\Lambda(j)$ is the weight $j$-space of the $\bb{G}_{m,R}$ action on $\Lambda$ given by $\rho\circ \nu$.

Thus, we can think of a filtered $\mc{H}$-torsor $\Fil^\bullet Q$ as being of type $\nu$ if locally on $\mc{X}_\mr{fppf}$ one has that $\Fil^\bullet \omega_Q$ is isomorphic as a filtered Tannakian $\mc{H}$-torsor to $\Fil^\bullet_\nu\omega_\mr{triv}$. This gives rise to a full subcategory $\mc{H}\text{-}\cat{FilVect}^\nu(\mc{X})$ of Tannakian $\mc{H}$-bundles of type $\nu$.
\end{rem}

\begin{nota} We define $P_\nu$ to be the group $R$-scheme $\Aut(\Fil^\bullet_\nu\mc{H})$. 
\end{nota}

\begin{rem}The group $R$-scheme $P_\nu$ is evidently a closed subgroup $R$-scheme of $\mc{H}$. Moreover, when $\mc{H}$ is reductive this agrees with the parabolic subgroup $R$-scheme of $\mc{H}$ associated to $\nu$ by the dynamic method as in \cite[Theorem 4.1.7]{ConradReductive}; cf.\@ \cite[Theorem 4.39 (iii)]{Ziegler15}.
\end{rem}

\begin{rem}As filtered $\mc{H}$-torsors satisfy descent, the general theory of twists (see \cite[Chapitre III, Théorème 2.5.1]{Giraud}) implies that for any prestack $\mc{X}$ over $R$ one may construct a functor
\begin{equation*}
\cat{Tors}_{P_\nu}(\mc{X})\to \cat{FilTors}^\nu_\mc{H}(\mc{X}),\qquad P\mapsto \bigg(\omega_P\colon \Lambda\mapsto P^{(\Lambda)}\defeq P\wedge^{(P_\nu)_{\ms{O}}}(\Lambda,\Fil^\bullet_\nu\Lambda)\bigg).
\end{equation*}
We may observe here that the underlying Tannakian $\mc{H}$-bundle of $\omega_P$ is $\omega_Q$ where $Q$ is the $\mc{H}$-torsor on $\mc{X}$ obtained by pushing out along the inclusion $P_\nu\subseteq \mc{H}$. We may then write $\omega_P=\Fil^\bullet_P \omega_Q$.
\end{rem}

The following then again follows from the general theory of twists.

\begin{prop}\label{prop:type-nu-classification} For a prestack $\mc{X}$ over $R$ and cocharacter $\nu\colon \bb{G}_{m,R}\to \mc{H}$, the functor 
\begin{equation*}
\cat{Tors}_{P_\nu}(\mc{X})\to \cat{FilTors}^\nu_\mc{H}(\mc{X}),\qquad P\mapsto  \omega_P=\Fil^\bullet_P\omega_Q.
\end{equation*}
is an equivalence.
\end{prop}

\begin{rem}\label{rem:derived-prestack} Definition \ref{defn:type-nu} and Proposition \ref{prop:type-nu-classification} naturally extend to the derived setting. Indeed, the definitions go through verbatim if one takes a $\mc{H}$-torsor on a derived stack $\mc{Y}$ to be an object of the $\infty$-category $\mr{Map}(\mc{Y},B\mc{H})$. 
\end{rem}

\begin{rem} While we will not need it in this article, we point out that assuming $\mc{H}$ is $R$-smooth, essentially every filtered $\mc{H}$-torsor on a prestack $\mc{X}$ over $R$ is locally of type $\nu$ for some $\nu$; see \cite[Theorem 1.1]{Ziegler} for a precise statement.
\end{rem}

\begin{rem}[Remark on sign conventions]\label{rem:sign-conventions} When comparing the material in this article with that in other articles on related topics, it is helpful to understand the relationship between various sign conventions concerning filtrations. Those of most import are: how $\bb{G}_m$ acts on $\bb{A}^1$ (i.e., whether $\lambda\cdot x$ is $\lambda x$ or $\lambda^{-1}x$) and what the notation $P_\mu$ or $P_\mu^{-}$ that appears there means. 

For the convenience of the reader we summarize these relationships in the following table. For clarity in this table we use the following notation for a cocharacter $\mu$ of an affine group scheme $\mc{G}$:
\begin{equation*}
P(\mu)\defeq \left\{g\in\mc{G}: \lim_{t\to 0}\mu(t)g\mu(t)^{-1}\text{ exists}\right\}
\end{equation*}
as in \cite[Theorem 4.1.7]{ConradReductive}.

\begin{table}[htbp]
  \centering
  \small
  \renewcommand{\arraystretch}{1.2}
  \begin{tabular}{
    @{}
    l
    r@{\;}c@{\;}l
    r@{\;}c@{\;}l
    @{}
  }
    \toprule
    Source
      & \multicolumn{3}{c}{$\mathbb G_m$-action on $\mathbb A^1$}
      & \multicolumn{3}{c}{Parabolic convention}
      \\
    \cmidrule(lr){2-4}
    \cmidrule(lr){5-7}

    This paper
      & $\lambda\cdot x$ & $=$ & $\lambda x$
      & $P_\mu$          & $=$ & $P(\mu)$
      \\

    \cite{IKY2}
      & $\lambda\cdot x$ & $=$ & $\lambda x$
      & $P_\mu^-$        & $=$ & $P(-\mu)$
      \\

    \cite{Ziegler15}
      &                  & \textemdash & 
      & $P(\Fil_\mu^\bullet)$ & $=$ & $P(\mu)$
      \\

    \cite{ItoDescent}/\cite{ItoDeformation}
      & $\lambda\cdot x$ & $=$ & $\lambda^{-1}x$
      & $P_\mu$          & $=$ & $P(\mu)$
      \\

    \cite{GMM}
      & $\lambda\cdot x$ & $=$ & $\lambda^{-1}x$
      & $P_\mu^-$        & $=$ & $P(\mu)$
      \\

    \cite{MY}
      & $\lambda\cdot x$ & $=$ & $\lambda^{-1}x$
      & $P_\mu^-$        & $=$ & $P(\mu)$
      \\

    \cite{BhattFGaugesNotes}
      & $\lambda\cdot x$ & $=$ & $\lambda x$
      &                  & \textemdash &
      \\

    \bottomrule
  \end{tabular}

  \label{tab:sign-conventions}
\end{table}
\end{rem}

\end{sectionappendices}

\section{The test functions $\phi_{\tau,h}$} 
In this section we define the functions $\phi_{\tau,h}^{\mc{G},\mu}$ as alluded to in the introduction, and show they have reasonable properties. We will then compare these with previously defined functions.

\begin{nota} Throughout the following we maintain the notation from Notation \ref{nota:deformation-spaces-setup}, and freely use notions and notations from \S\ref{s:apertures-and-crystalline-theory}, and particularly from \S\ref{ss:deformation-theory}. Furthermore we assume 
\begin{itemize} 
\item $p>2$;
\item $k$ is finite;
\item $j$ belongs to $\bb{N}\cup\{\infty\}$.
\end{itemize}
\end{nota}

\subsection{Aperture tubes with level structure} We begin by defining the spaces whose cohomology will be used to define the functions $\phi_{\tau,h}$. 

\begin{setup}\label{setup:test-functions} Fix $b$ in $\mc{G}(W_j)\sigma(\mu)(p)\mc{G}(W_j)$ thought of as an element of $\mr{BT}^{\mc{G},\mu}_\infty(k_j)$ as in \S\ref{ss:case-of-finite-fields}.
\end{setup}

\begin{defn}\label{defn:aperture-tubes} We define the \emph{aperture tube at $b$} to be the rigid $K_j$-space
\begin{equation*}
\mathcal{D}(\mc{G},b,\mu) \defeq \Spf(R_{\mc{G},\mu,j},\mathfrak{m}_{\mc{G},\mu,j})_\eta.
\end{equation*}
By \'etale realization we obtain a map $\mc{D}(\mc{G},b,\mu)\to B\mc{G}(\bb{Z}_p)$ and we define the \emph{infinite level aperture tube at $b$}, denoted $\mc{D}_\infty(\mc{G},b,\mu)$ to be $\mc{D}(\mc{G},b,\mu)\times_{B\mc{G}(\Z_p)} \ast$, i.e., the frame space for $T_\et(\mf{Q}_b^\univ)$.

For any compact open subgroup $\mathsf{K}\subseteq \mc{G}(\bb{Z}_p)$, we define the \emph{level $\mathsf{K}$ \emph{(}aperture\emph{)} tube at $b$} as 
\begin{equation*}
\mc{D}_\mathsf{K}(\mc{G},b,\mu)\defeq \mc{D}_\infty(\mc{G},b,\mu)/\mathsf{K}.
\end{equation*}
When $\mathsf{K}=\mathsf{K}_n\defeq\ker\left(\mc{G}(\bb{Z}_p)\to \mc{G}(\bb{Z}/p^n)\right)$ we abbreviate $\mc{D}_{\mathsf{K}_n}(\mc{G},b,\mu)$ to $\mc{D}_n(\mc{G},b,\mu)$.
\end{defn}

\begin{nota} Suppose that $\mathsf{K},\mathsf{K}'\subseteq\mc{G}(\Z_p)$ are compact open subgroups, and $g$ is an element of $\mc{G}(\Z_p)$ satisfying $g^{-1}\mathsf{K}g\subseteq \mathsf{K}'$. Right multiplication by $g$ on $\mc{D}_\infty(\mc{G},b,\mu)$ descends to a morphism
\begin{equation*}
t_{\mathsf{K},\mathsf{K}'}(g)\colon \mc{D}_\mathsf{K}(\mc{G},b,\mu)\to \mc{D}_{\mathsf{K}'}(\mc{G},b,\mu).
\end{equation*}
For $\mathsf{K}\subseteq\mathsf{K}'$, put $\pi_{\mathsf{K},\mathsf{K}'}=t_{\mathsf{K},\mathsf{K}'}(1)$, and write $\pi_\mathsf{K}=\pi_{\mathsf{K},\mc{G}(\Z_p)}$.
\end{nota}

We would like to give a functorial description of these spaces $\mc{D}_\mathsf{K}(\mc{G},b,\mu)$.

\begin{defn} Let $R$ be an adic $W_j$-algebra. An object in the essential image of
\begin{equation*}
T_\et\colon \mr{BT}^{\mc{G},\mu}_\infty(R)\to B\mc{G}(\bb{Z}_p)(R[\nicefrac{1}{p}]),
\end{equation*}
is called \emph{apertile of type $\mu$ \emph{(}relative to $R$\emph{)}}. If $R$ is a discrete $W_j$-algebra, we define an object of $B\mc{G}(\bb{Z}_p)(R[\nicefrac{1}{p}])$ to be an \emph{apertile object of type $\mu$} if it lies in the essential image of
\begin{equation*}
T_\et\colon \mr{BT}^{\mc{G},\mu,\mr{alg}}_\infty(R)\to B\mc{G}(\bb{Z}_p)(R[\nicefrac{1}{p}]).
\end{equation*}
These notions readily globalize.
\end{defn}

We have the following concrete, purely generic-fiber description of apertile $\mc{G}(\bb{Z}_p)$-local systems of type $\mu$ in the good reduction situation.

\begin{prop}[{\cite[Theorem 5.1.4]{MY}}] Let $R$ be a smooth adic $W_j$-algebra and assume $p\geqslant 5$. Then, an object $\ms{Q}$ of $B\mc{G}(\bb{Z}_p)(R[\nicefrac{1}{p}])$ is apertile of type $\mu$ if and only if:
\begin{enumerate}
\item $\ms{Q}$ is pointwise crystalline, i.e., for every classical point $x$ of $\Spa(R[\nicefrac{1}{p}])$ and for every representation $\rho$ of $\mc{G}$ the $\bb{Z}_p$-representation of $\Gal(\ov{k(x)}/k(x))$ given by 
$\rho_\ast \ms{Q}_x$ is crystalline.
\item the filtered Tannakian $\mc{G}$-bundle $D_\mr{dR}\circ \omega_\ms{Q}$ is of type $\mu$ in the sense of \emph{Remark \ref{rem:filtered-Tannakian-torsor-type-mu}}.
\end{enumerate}
\end{prop}

\begin{defn} An adic $W_j$-algebra $R$ is \emph{$\eta$-normal} if the map $R\to R[\nicefrac{1}{p}]^\circ$ is an isomorphism.
\end{defn}

\begin{thm}[{\cite[Theorem 4.1.5]{MY}}] Suppose that $R$ is a Noetherian $\eta$-normal adic $W_j$-algebra. Then, $T_\et\colon \mr{BT}^{\mc{G},\mu}_\infty(R)\to B\mc{G}(\bb{Z}_p)(R[\nicefrac{1}{p}])$ is fully faithful. 
\end{thm}

\begin{nota} Let $R$ be a Noetherian $\eta$-normal adic $W_j$-algebra. Then, for an apertile object $\ms{Q}$ of type $\mu$ in $B\mc{G}(\bb{Z}_p)(R[\nicefrac{1}{p}])$ we denote $T_\et^{-1}(\ms{Q})|_{R_\mr{red}}$ by $\mr{sp}(\ms{Q})$ and call it the \emph{specialization} of $\ms{Q}$.
\end{nota}

We then have the following more functorial description of the values of $\mc{D}(\mc{G},b,\mu)(A)$ for a quite general class of inputs, e.g., all affinoid $K_j$-algebras. 

\begin{prop}\label{prop:functorial-description} Let $A$ be a Tate $K_j$-algebra such that $A^\circ$ is Noetherian. Then $\mathcal{D}_n(\mc{G},b,\mu)(A)$ is naturally identified with isomorphism classes of triples $(\ms{Q},\iota,\alpha)$ where:
\begin{itemize}
\item $\ms{Q}$ is an object of $B\mc{G}(\bb{Z}_p)(A)$ which is apertile of type $\mu$ \emph{(}relative to $A^\circ$\emph{)} over $A$;
\item $\iota\colon \mr{sp}(\ms{Q})\isomto b|_{A^\circ_\mr{red}}$ is an isomorphism in $\mr{BT}^{\mc{G},\mu}_\infty(A^\circ_\mr{red})$;
\item $\alpha$ is a trivialization of $\ms{Q}\otimes \bb{Z}/p^n$.\footnote{By definition $\ms{Q}\otimes \bb{Z}/p^n$ denotes the pushout of $\ms{Q}$ along $\mc{G}(\bb{Z}_p)\to \mc{G}(\bb{Z}/p^n)$.}
\end{itemize}
\end{prop}
\begin{proof} The general case easily follows from the $n=0$ case. In the $n=0$ case we see by the full faithfulness of the \'etale realization functor on $A^\circ$ and Proposition \ref{prop:functorial-description-deformation} that the set of isomorphism classes of such triples naturally identifies with $\mf{D}(\mc{G},b,\mu)(A^\circ)$. We then must show that the natural map $\mf{D}(\mc{G},b,\mu)(A^\circ)\to \mc{D}(\mc{G},b,\mu)(A)$ is a bijection. But, choosing an isomorphism $\mf{D}(\mc{G},b,\mu)\simeq \Spf(W_j\llbracket t_1,\ldots,t_d\rrbracket)$ identifies this map with the map $((A^\circ)^{\circ\circ})^d\to (A^{\circ\circ})^d$. But, as $(A^\circ)^{\circ\circ}=A^{\circ\circ}$ the claim follows.
\end{proof}

\subsection{Action of the groups \texorpdfstring{$J_b^\mr{int}$}{J_b}}

We now wish to discuss the action of the group $J_b^\mr{int}$ of automorphisms of $b$ on tubes with level structure, and establish some workhorse results concerning this action that will be used in the finer study of our local test functions. 

\begin{nota}\label{nota:integral-sigma-centralizer} Write $J_b^\mr{int}\defeq \left\{g\in \mc{G}(W_j): gb\sigma(g)^{-1}=b\right\}=\Aut(b)$.
\end{nota}

\begin{rem}There is a natural group action
\begin{equation*}
J_{b}^\mr{int}\times \mf{D}(\mc{G},b,\mu)\to \mf{D}(\mc{G},b,\mu),\qquad (g,(\mf{Q},\iota))\mapsto (\mf{Q},g\circ \iota),
\end{equation*}
where this last composition uses the discussion in \S\ref{ss:case-of-finite-fields}. This naturally induces an action 
\begin{equation*}
J_{b}^\mr{int}\times \mc{D}_\mathsf{K}(\mc{G},b,\mu)\to \mc{D}_\mathsf{K}(\mc{G},b,\mu),
\end{equation*}
which when $\mathsf{K}=\mathsf{K}_n$ can be described in terms of Proposition \ref{prop:functorial-description}: This action associates to the pair $(g,(\ms{Q},\iota,\alpha))$ the triple $(\ms{Q},g\circ \iota,\alpha)$.
\end{rem}

\begin{rem}\label{rem:sigma-centralizer} For $j<\infty$, one has the classical group scheme over $\bb{Q}_p$ given on $R$-points by
\begin{equation*}
{J}_b(R)\defeq\left\{g\in G(R\otimes_{\Q_p}K_j):gb\sigma(g)^{-1}=b\right\}.
\end{equation*}
This admits a closed embedding into $\mr{Res}_{K_j/\mathbb{Q}_p}\,\mc{G}_{K_j}$, and the Zariski closure of ${J}_b$ in $\mr{Res}_{W_j/\bb{Z}_p}\,\mc{G}_{W_j}$ has $\Z_p$-points $J_b^\mr{int}$.
\end{rem}

In the rest of this subsection we aim to establish two continuity properties for the action of $J_b^\mr{int}$ on the spaces $\mf{D}(\mc{G},b,\mu)$ and $\mc{D}_\mathsf{K}(\mc{G},b,\mu)$. 

\begin{obs}The action of $J_b^\mr{int}$ on $\mf{D}(\mc{G},b,\mu)$ gives a map $J_b^\mr{int}\to \mr{Aut}_\mr{cont.}(R_{\mc{G},\mu,j})$ where the subscript $\mr{cont.}$ means the group of ring automorphisms $f$ of $R_{\mc{G},\mu,j}$ which are continuous for the $\mf{m}_{\mc{G},\mu,j}$-adic topology; in particular $f(\mf{m}_{\mc{G},\mu,j})\subseteq \mf{m}_{\mc{G},\mu,j}$.
\end{obs}

\begin{defn} The \emph{$N$-congruence subgroups} of $J_b^\mr{int}$ and $\mr{Aut}_\mr{cont.}(R_{\mc{G},\mu,j})$ are given by
\begin{equation*}
J_b^\mr{int}(N)\defeq J_b^\mr{int}\cap \ker\bigg(\mc{G}(W_j)\to \mc{G}(W_j/p^N)\bigg),\qquad \ker\bigg(\mr{Aut}_\mr{cont.}(R_{\mc{G},\mu,j})\to \mr{Aut}(R_{\mc{G},\mu,j}/\mf{m}_{\mc{G},\mu,j}^N)\bigg).
\end{equation*}
respectively. The \emph{congruence topology} on $J_b^\mr{int}$ and $\mr{Aut}_\mr{cont.}(R_{\mc{G},\mu,j})$ is the unique group topology such that the $N$-congruence subgroups form a neighborhood basis of the identity.
\end{defn}

The groups $J_b^\mr{int}(N)$ can be given the following useful alternative description.

\begin{nota} Let $\mf{Q}$ be a $(\mc{G},\mu)$-aperture over a $p$-nilpotent ring $R$. For each $N\geqslant 1$, denote by $\mf{Q}_N$ its $N$-truncation, i.e., its image in $\mr{BT}^{\mc{G},\mu}_N(R)$. We then write
\begin{equation*}
K_N(\mf{Q})\defeq\ker\left(\Aut(\mf{Q})\to \mr{Aut}(\mf{Q}_N)\right).
\end{equation*}
\end{nota}

\begin{prop}\label{prop:congruence-agrees-with-truncation} For $N\gg 0$ one has 
\begin{equation*}
J_b^\mr{int}(N)=K_N(b)=\ker\bigg(\mr{Aut}(b)\to \mr{Aut}(b_N)\bigg).
\end{equation*}
\end{prop}
\begin{proof} Without loss of generality we may assume that the Hodge filtration of $b$ is $P_\mu$. By \cite[Remark 3.2.20]{MY} (and the references therein) one has
\begin{equation*}
\mr{Aut}(b_N)\subseteq \pi_0(H_{N}^{(\mu)})\defeq \pi_0\left(\mc{G}(W_j/p^N)\times_{\mc{G}(k_j\dquot p^N)}P_\mu(k_j\dquot p^N)\right).
\end{equation*}
Now, a quick calculation shows that there is a short exact sequence
\begin{equation*}
0\to (\mf{g}/\mf{p}_\mu)(k_j)\to \pi_0(H_{N}^{(\mu)})\to \mc{G}(W_j/p^N)\times_{\mc{G}(k_j)}{P}_\mu(k_j)\to 1
\end{equation*}
We claim that the image of $\mr{Aut}(b)$ in $\pi_0(H_N^{(\mu)})$ intersects $(\mf{g}/\mf{p}_\mu)(k_j)$ trivially. This will then prove our claim as we will then have that the image of $\mr{Aut}(b)$ in $\mr{Aut}(b_N)$ embeds into $\mc{G}(W_j/p^N)$. 

To see trivial intersection, take $g$ in $\mr{Aut}(b)\subseteq \mc{G}(W_j)$ whose image $\ov{g}$ in $\pi_0(H_N^{(\mu)})$ lands in $(\mf{g}/\mf{p}_\mu)(k_j)$. In particular, we see that $g$ maps to the identity in $\mc{G}(W_j/p^N)$. From the bijection 
\begin{equation*}
\exp\colon p^N\mf{g}_{W_j}\isomto \left\{g\in \mc{G}(W_j): g=1\mod p^N\right\} ,
\end{equation*}
which exists for $N\gg 0$, we deduce that $g=\exp(p^N x)$ for $x$ in $\mf{g}_{W_j}$. Unraveling the definitions shows that $x=\ov{g}\mod p$ as elements of $(\mf{g}/\mf{p}_\mu)(k_j)$. Thus, we must show that the reduction mod $p$ of $x$ lies in $\mf{p}_\mu(k_j)$. Using Proposition \ref{prop:Hodge-filtered-bundle-descp-char-p} it suffices to show that the reduction of $x$ lies in the $0^\text{th}$ filtered piece of the Mazur--Nygaard--Ogus filtration on $\mf{g}_{k_j}$. But, by taking logarithms of the equality $gb\sigma(g)^{-1}=b$ we deduce that $x=\mr{Ad}(b)(\sigma(x))$ which implies the claim. 
\end{proof}

We now aim to prove the following (cf.\@ the proof of \cite[Proposition 3.15]{ScholzeLK}) which, in particular, says that the map $J_b^\mr{int}\to \mr{Aut}_\mr{cont.}(R_{\mc{G},\mu,j})$ is continuous for the congruence topology. 
 
 \begin{prop}\label{prop:uniform-continuity-base-level} Let $r$ and $t$ be positive integers. Then, there exists an $N\geqslant 1$ such that every $g$ in $J_b^\mr{int}(N)$ lifts to an element $\wt{\alpha}$ of $K_t(\mf{Q}(r))$, where $\mf{Q}(r)=\mf{Q}_b^\univ|_{R_{\mc{G},\mu,j}/\mf{m}^r_{\mc{G},\mu,j}}$.
 \end{prop}
 Before we prove this, we establish the following modified Grothendieck--Messing lifting lemma.
 \begin{lem}\label{lem:GM-variant} Suppose that $s\geqslant 1$ is an integer and that $R'\twoheadrightarrow R=R'/I$ is a square-zero thickening of $W$-algebras such that $p^sR'=0$ and $pI=0$. Fix $\mf{Q}'$ in $\mr{BT}^{\mc{G},\mu}_\infty(R')$ and set $\mf{Q}\defeq \mf{Q}'|_R$. Then, for any $\alpha$ in $K_s(\mf{Q})$, one has that $\alpha^p$ lies in the image of $K_s(\mf{Q}')\to K_s(\mf{Q})$. 
 \end{lem}
\begin{proof} By the Grothendieck--Messing theorem (see Theorem \ref{thm:GMM}) one has a Cartesian diagram
\begin{equation*}
    \begin{tikzcd}
	{\mr{BT}^{\mc{G},\mu}_\infty(R')} & {BP_{\mu}(R')} \\
	{\mr{BT}^{\mc{G},\mu}_\infty(R)} & B\mc{G}(R')\times_{B\mc{G}(R)} BP_{\mu}(R).
	\arrow[from=1-1, to=1-2]
	\arrow[from=1-1, to=2-1]
	\arrow[from=1-2, to=2-2]
	\arrow[from=2-1, to=2-2]
\end{tikzcd}
\end{equation*}
So, let $\alpha$ be in $K_s(\mf{Q})$. To lift $\alpha$ to an automorphism $\wt{\alpha}$ of $\mf{Q}'$ it suffices to give an automorphism $\beta$ of $\Fil^\bullet_\mr{Hdg} T_\mr{dR}(\mf{Q}')$ in $BP_\mu(R')$ and an identification of the images of $\alpha$ and $\beta$ in the groupoid $B\mc{G}(R')\times_{B\mc{G}(R)} BP_{\mu}(R)$. We may take $\beta=\mr{id}$, as long as we can justify that $\alpha$ induces the identity on the pullback along the map $\Spec(R')\to R^\mr{syn}$. But, as $R'$ is $p^s$-torsion, this factorizes through $R^\mr{syn}\dotimes \bb{Z}/p^s$, and thus the claim follows as $\alpha$ pulled back to $R^\mr{syn}\dotimes \bb{Z}/p^s$ is equivalent to $\alpha_s$ which is the identity by assumption.

Now, a priori $\wt{\alpha}$ does not lie in $K_s(\mf{Q}')$. But, if $H=\underline{\mr{Aut}}(\mf{Q}'_s)$, a finite type affine group $R'$-scheme, we have that $\wt{\alpha}_s$ lies in
\begin{equation*}
\ker\left(H(R')\to H(R)\right)\simeq \Hom_{R'}(e^\ast\Omega^1_{H/R'},I),
\end{equation*}
where $e$ is the identity section. But, by assumption this is a $p$-torsion group, and thus we see that $\wt{\alpha}^p_s$ is trivial. Thus, $\wt{\alpha}^p$ is in $K_s(\mf{Q}')$ and lifts $\alpha^p$.
\end{proof}
 
 \begin{proof}[Proof of Proposition \ref{prop:uniform-continuity-base-level}] Choose $m\geqslant1$ so that we have an isomorphism
\begin{equation*}
\exp\colon p^m\mf{g}_{W_j}\isomto \left\{g\in \mc{G}(W_j): g=1\mod p^m\right\} ,
\end{equation*}
with inverse given by $\log$. Set $q=\max\{m,r,t\}$. Then, we claim that $N=q+r$ is sufficient.

So, for $g$ in $J_b^\mr{int}(N)$, set $h=\exp\big(p^{-(r-1)}\log(g)\big)$; note that this makes sense as $\log(g)$ belongs to $p^N\mf{g}_{W_j}$ and $N\geqslant r-1$. We claim that $h$ is in $J_b^\mr{int}(q+1)$ and $h^{p^{r-1}}=g$.
Indeed, the latter is clear by setup, and the equation $gb\sigma(g)^{-1}=b$ gives
$\log(g)=\mr{Ad}(b)\sigma(\log(g))$ which is preserved by the division and exponentiation.
Proposition \ref{prop:congruence-agrees-with-truncation} then gives that $h$ is in $K_q(b)$.

We now build $\wt{\alpha}$ inductively. Start with $\beta_1=h$.
For $1\leqslant i<r$, one has that $\mf{m}_{\mc{G},\mu,j}^i/\mf{m}_{\mc{G},\mu,j}^{i+1}$ is square-zero and killed by $p$, and moreover we have that $p^qR_{\mc{G},\mu,j}/\mf{m}_{\mc{G},\mu,j}^{i+1}=0$. Thus, Lemma \ref{lem:GM-variant} gives inductively $\beta_{i+1}$ in $K_q(\mf Q(i+1))$ whose restriction to $R_{\mc{G},\mu,j}/\mf{m}_{\mc{G},\mu,j}^i$ is $\beta_i^p$. Thus $\beta_r|_{k_j}=h^{p^{r-1}}=g$ and $(\beta_r)_t=1$.
Thus, we may take $\wt{\alpha}=\beta_r$.
 \end{proof}

Finally, we use Proposition \ref{prop:congruence-agrees-with-truncation} to establish the following uniformity condition for the action of the topological group $J_b^\mr{int}$ on the space $\mc{D}_\mathsf{K}(\mc{G},b,\mu)$.

\begin{nota}\label{nota:exhaustion} For any $m\geqslant 1$, write $X_m\subseteq \mc{D}(\mc{G},b,\mu)$ for the affinoid open subspace given by
\begin{equation*}
X_m=\left\{x\in \mc{D}(\mc{G},b,\mu):|f(x)|\leqslant |p|\text{ for all }f\in\mf{m}_{\mc{G},\mu,j}^m\right\}.
\end{equation*}
For any compact open subgroup $\mathsf{K}\subseteq \mc{G}(\bb{Z}_p)$ set $Y_{\mathsf{K},m}\defeq \pi_\mathsf{K}^{-1}(X_m)$, which form an increasing cover of $J_b^\mr{int}$-stable affinoid open subspaces of $\mc{D}_\mathsf{K}(\mc{G},b,\mu)$.
\end{nota}

\begin{prop}\label{prop:uniform-continuity-K-level}
Let $\mathsf{K}\subseteq \mathcal G(\bb{Z}_p)$ be a compact open subgroup. Then, for every open neighborhood $\Delta_{Y_{\mathsf{K},m}}\subseteq U$ there exists $N\geqslant 0$ such that $\Gamma_g^{Y_{\mathsf{K},m}}\subseteq U$ for every $g$ in $J_b^{\mathrm{int}}(N)$.
\end{prop}

In the statement of Proposition \ref{prop:uniform-continuity-K-level}, and its proof below, we write $\Gamma_g^{(-)}$ to indicate we are considering the graph of the action of $g$ on $(-)$.

\begin{proof}
Without loss of generality, $\mathsf{K}=\mathsf{K}_n$ for some $n$. Denote by $\mr{pr}_i$ for $i=1,2$ the projection maps $\mf{D}(\mc{G},b,\mu)\times_{\Spf(W_j)}\mf{D}(\mc{G},b,\mu)\to\mf{D}(\mc{G},b,\mu)$, and further abbreviate $\mf{Q}_{b,n}^\mr{univ}$ to $\mf{Q}_n$. Over the closed point of $\mf{D}(\mc{G},b,\mu)\times_{\Spf(W_j)}\mf{D}(\mc{G},b,\mu)$ there are tautological identifications of $\mr{pr}_1^*\mf{Q}_n$ and $\mr{pr}_2^*\mf{Q}_n$: Set $\mf{I}_n$ to be the completion of $\underline{\mr{Isom}}(\mr{pr}_1^*\mf{Q}_n,\mr{pr}_2^*\mf{Q}_n)$ at this point. This formal scheme is affine and Noetherian by Theorem \ref{thm:GMM} and Remark \ref{rem:integrability}; let us write $\mf{I}_n=\Spf(S,\mf{n})$.

For $g$ in $K_n(b)$, the action on the universal deformation induces
an isomorphism $\mf{Q}_n\isomto g^\ast\mf{Q}_n$ whose restriction to the closed point is the identity of the truncation $b_n$. Together with $(\mr{id},g)$, this thus defines a morphism $\mf{D}(\mc{G},b,\mu)\to\mf{I}_n$ which we denote $s_g$. By Proposition \ref{prop:congruence-agrees-with-truncation}, these
morphisms are defined for every $g$ in $J_b^\mr{int}(N)$ for $N\gg 0$.

Now, fix $r\geqslant 1$ an integer. Applying Proposition \ref{prop:uniform-continuity-base-level} with $t=n$ gives, for $N\gg 0$, a lift
$\wt\alpha$ in $K_n(\mf{Q}(r))$ of an element $g$ in $J_b^\mr{int}(N)$ whose inverse identifies the universal framed deformation modulo
$\mf{m}_{\mc{G},\mu,j}^r$ with its $g$-translate. Representability (and thus discreteness) of the deformation functor implies that $g$ acts trivially on $R_{\mc{G},\mu,j}/\mf{m}_{\mc{G},\mu,j}^r$. Since $\wt\alpha$ has trivial $n$-truncation, uniqueness of the framed comparison gives us an identification $s_g^*\equiv s_1^*\mod \mf{m}_{\mc{G},\mu,j}^r$.

Consider then the affinoid
\begin{equation*}
 Z=\{z\in(\mf I_n)_\eta:
       |f(z)|\leqslant|p|\text{ for all }f\in\mf n^m\}.
\end{equation*}
Note that both projections $\mr{pr}_i$ send $Z$ into $X_m$; let us abbreviate them as the maps $p_1,p_2:Z\to X_m$.
Moreover, observe that $s_g(X_m)\subseteq Z$ for every $g$ in $J_b^\mr{int}(N)$ for $N\gg 0$. Thus, since $\mf m_{\mc G,\mu,j}^r$ is bounded by
$|p|^{\lfloor r/m\rfloor}$ on $X_m$, the preceding congruences applied to affinoid generators of $Z$ imply that for every $a$ in $\mathcal O(Z)$ and $\varepsilon>0$ one has, up to increasing $N$, the bound $\|s_g^*a-s_1^*a\|_{\mathrm{sup},X_m}<\varepsilon$ for every $g$ in $J_b^\mr{int}(N)$.

Applying $T_\et$ to the tautological isomorphism $u\colon \mr{pr}_1^\ast\mf{Q}_n\to \mr{pr}_2^\ast\mf{Q}_n$ gives an isomorphism between the two
pullbacks of the $\mathsf{K}_n$-level cover to $Z$. This isomorphism commutes with the right action on framings, and so descends to a $Z$-isomorphism
\begin{equation*}
 \theta:p_1^*Y_{\mathsf{K}_n,m}\isomto p_2^*Y_{\mathsf{K}_n,m},
 \qquad [\alpha]\longmapsto[T_\et(u)\circ\alpha].
\end{equation*}
Consequently there is a morphism
\begin{equation*}
 e:p_1^* Y_{\mathsf{K}_n,m}\to  Y_{\mathsf{K}_n,m}\times_{\Spa(K_j)} Y_{\mathsf{K}_n,m},
 \qquad
 e=(\operatorname{pr}_{Y_{\mathsf{K}_n,m}},\operatorname{pr}_{Y_{\mathsf{K}_n,m}}\circ\theta),
\end{equation*}
where the second projection is taken from $p_2^*Y_{\mathsf{K}_n,m}$. Along $s_g$, the map $e$ pulls back to $y\mapsto(y,gy)$, and along $s_1$ it pulls back to the diagonal map.

Let $\pi\colon p_1^*Y_{\mathsf{K}_n,m}\to Z$ be the finite projection map.
Since $\pi$ is closed, the subset
\begin{equation*}
V=Z-\pi\big((p_1^*Y_{\mathsf{K}_n,m})-e^{-1}(U)\big)
\end{equation*}
is open, contains $s_1(X_m)$, and satisfies
$\pi^{-1}(V)\subseteq e^{-1}(U)$.
As $X_m$ is quasi-compact, choose finitely many rational open subspaces
\begin{equation*}
W_i=\left\{z\in Z:
|a_{ij}(z)|\leqslant |d_i(z)|\ne 0
\text{ for }1\leqslant j\leqslant r_i\right\}
\subseteq V,
\end{equation*}
with $a_{ij}$ and $d_i$ in $\mathcal O(Z)$, whose union contains $s_1(X_m)$.
Set $T_i=s_1^{-1}(W_i)$.
On each rational affinoid $T_i$, the function $s_1^*d_i$ is invertible
and its inverse has bounded supremum norm.
Thus, the preceding uniform estimate gives, for $N\gg 0$,
\begin{equation*}
\left\|
\frac{s_g^*d_i-s_1^*d_i}{s_1^*d_i}
\right\|_{\mathrm{sup},T_i}<1,
\qquad
\left\|
\frac{s_g^*a_{ij}-s_1^*a_{ij}}{s_1^*d_i}
\right\|_{\mathrm{sup},T_i}<1
\end{equation*}
for every $i,j$ and every $g$ in $J_b^\mr{int}(N)$. Thus, we deduce that $s_g(T_i)\subseteq W_i$ for every $i$, and so $s_g(X_m)\subseteq V$. From this it follows that $\Gamma_g^{Y_{\mathsf{K}_n,m}}\subseteq U$, as desired.
\end{proof}

\subsection{The local test functions \texorpdfstring{$\phi_{\tau,h}$}{phi(tau,h)}}

We finally come to the definition of the local test functions $\phi_{\tau,h}^{\mc{G},\mu}$ which serve as the deeper-level contribution to the cohomology of Shimura varieties.

\subsubsection{Hecke algebras and representations} We first establish some basic notation and results concerning Hecke algebras and their relationship to representations.

\begin{setup} Throughout the following we fix
\begin{itemize}
\item $\mathsf{G}$ to be a unimodular locally profinite topological group;
\item a Haar measure $dg$ giving rational volumes to compact open subgroups;
\item $F$ to be an extension of $\bb{Q}$. 
\end{itemize}
\end{setup}

\begin{defn}\label{defn:Hecke-algebra}The $F$-Hecke algebra for $\mathsf{G}$ is
\begin{equation*}
\mc{H}_F(\mathsf{G})\defeq \left\{f\colon 
\mathsf{G}\to F\text{ locally constant and compactly supported}\right\}.
\end{equation*}
For a compact open subgroup $\mathsf{K}\subseteq\mathsf{G}$ we write $\mc{H}_F(\mathsf{G},\mathsf{K})$ for the subset consisting of functions $f$ which are $\mathsf{K}$-bi-invariant, i.e., $f(k_1 gk_2)=f(g)$ for all $k_1$ and $k_2$ in $\mathsf{K}$. The $F$-Hecke algebra $\mc{H}_F(\mathsf{G})$ is an associative $F$-algebra with convolution product $\ast$:
\begin{equation*}
(f_1 \ast f_2)(g)\defeq \int_{\mathsf{G}}f_1(k)f_2(k^{-1}\cdot g)\, dk.
\end{equation*}
When $F=\Q$ we drop it from the notation, yielding the notation $\mc{H}(\mathsf{G})$ and $\mc{H}(\mathsf{G},\mathsf{K})$.
\end{defn}

\begin{nota}\label{nota:normalized-Hecke-functions} For $g$ in $\mathsf{G}$ and $\mathsf{K}\subseteq\mathsf{G}$ compact open we write $e(g,\mathsf{K})$ for the indicator function $\tfrac{1}{\mr{vol}(\mathsf{K})}\mathds{1}_{\mathsf{K}g\mathsf{K}}$. We shorten $e(1,\mathsf{K})$ to $e(\mathsf{K})$. 
\end{nota}

\begin{defn} An $F[\mathsf{G}]$-module $V$ is called \emph{smooth} if $V=\bigcup_\mathsf{K} V^\mathsf{K}$ as $\mathsf{K}$ ranges over compact open subgroups of $\mathsf{G}$. We say that $V$ is \emph{admissible} if each $V^\mathsf{K}$ is finite-dimensional as an $F$-space.
\end{defn} 

\begin{nota}\label{nota:Hecke-action-repn} If $V$ is a smooth $F[\mathsf{G}]$-module then it naturally inherits the structure of a $\mc{H}_F(\mathsf{G})$-module via the action:
\begin{equation*}
h\ast x=\int_\mathsf{G} h(k)(k\cdot x)\,dk.
\end{equation*}
\end{nota}

\begin{nota}\label{nota:two-group-act} Suppose $\Gamma$ is a group acting $F$-linearly on $V$ in a way that commutes with the action of $\mathsf{G}$. Then $V$ naturally obtains the structure of a $F[\Gamma]\otimes_F \mc{H}_F(\mathsf{G})$-module. We abuse notation and write $\gamma\times h$ for the operator given by the action of $\gamma\otimes h$. 
\end{nota}

\begin{defn} Suppose we are in the situation of Notation~\ref{nota:two-group-act} and that $V$ is admissible. For $\gamma$ in $\Gamma$ and $h$ in $\mc{H}_F(\mathsf{G})$, the operator $\gamma\times h$ has finite-dimensional image.
We define $\tr(\gamma\times h\mid V)$ as its trace on any finite-dimensional $(\gamma\times h)$-stable subspace containing its image.
\end{defn}

\begin{obs}\label{obs:Hecke-action-desc} Suppose we are in the situation of Notation \ref{nota:two-group-act}, and that $\mathsf{K}$ is a compact normal open subgroup of $\mathsf{G}$ and $g$ is an element of $\mathsf{G}$. Then, $\gamma\times e(g,\mathsf{K})$ is precisely the operator
\begin{equation*}
\gamma \times e(g,\mathsf{K})\colon V\to V,\qquad v\mapsto (\gamma\times \ov{g})\cdot v^{\mathsf{K}}
\end{equation*}
where $v^{\mathsf{K}}$ is the projection of $v$ to $V^\mathsf{K}$ and $\ov{g}$ is the image of $g$ in $\mathsf{G}/\mathsf{K}$.

More generally, if $\mathsf{K}\subseteq\mathsf{G}$ is an arbitrary compact open subgroup and $g$ is in $\mathsf{G}$, set $\mathsf{K}^g\defeq \mathsf{K}\cap g\mathsf{K}g^{-1}$. Then, $e(g,\mathsf{K})\ast v$ is precisely $\mr{Tr}_{\mathsf{K}^g}^{\mathsf{K}}(g\cdot v^\mathsf{K})$, where for an element $w$ of $V^{\mathsf{K}^g}$ we write
\begin{equation*}
\mr{Tr}_{\mathsf{K}^g}^{\mathsf{K}}(w)\defeq\sum_{k\in \mathsf{K}/\mathsf{K}^g}k\cdot w.
\end{equation*}
Thus, the action of $e(g,\mathsf{K})$ factors as
\begin{equation*}
V\to V^\mathsf{K}\xrightarrow{g}V^{g\mathsf{K}g^{-1}}\to V^{\mathsf{K}^g}\xrightarrow{\mr{Tr}_{\mathsf{K}^g}^{\mathsf{K}}}V^\mathsf{K}\to V.
\end{equation*}
If $\Gamma$ acts as above, then the same description holds for $\gamma\times e(g,\mathsf{K})$ after composing with $\gamma$.
\end{obs}

\subsubsection{The local test functions and their main properties} We now precisely state the definition of the functions $\phi_{\tau,h}^{\mc{G},\mu}$ and their basic properties. 

\begin{nota} For $g$ in $\mc{G}(\Z_p)$ and $n\geqslant1$, we shorten $e(g,\mathsf{K}_n)$ to $e(g,n)$.
\end{nota}

\begin{nota} Fix an algebraic closure $\ov{K}$ of $K$ and let $C$ denote its $p$-adic completion. Write $\Gamma_K$ for the absolute Galois group of $K$ which agrees with the group of $p$-adically continuous automorphisms of $C$ over $K$. For $1\leqslant j\leqslant \infty$ write $\Gamma_{K_j}$ for the subgroup of elements of $\Gamma_K$ fixing $K_j$. For $\tau$ in $W_K$ with $\mathsf{v}(\tau)=j$, one has $\tau|_{K_\infty}=\sigma^{-rj}$ and $\tau^{-1}|_{K_\infty}=\sigma^{rj}$.
\end{nota}

\begin{defn}\label{defn:coh-for-test-functions} For an integer $1\leqslant j\leqslant \infty$, an element $b$ in $G(K_j)$, a compact open subgroup $\mathsf{K}\subseteq\mc{G}(\bb{Z}_p)$, and a prime $\ell\ne p$, we define 
\begin{equation*}
H^i(\mc{G},b,\mu,\mathsf{K})\defeq \begin{cases}
H^i_\et(\mc{D}_\mathsf{K}(\mc{G},b,\mu)_C,\ov{\Q}_\ell) & \mbox{if}\quad b\in\mc{G}(W_j)\sigma(\mu)(p)\mc{G}(W_j),\\ 0 & \mbox{if}\quad \text{otherwise.}
\end{cases}
\end{equation*}
We then define
\begin{equation*}
H^i(\mc{G},b,\mu)\defeq \varinjlim_{\mathsf{K}\subseteq\mc{G}(\bb{Z}_p)}H^i(\mc{G},b,\mu,\mathsf{K}),
\end{equation*}
which carries an action of $\mc{G}(\bb{Z}_p)$ by modifying the level structure which commutes with the natural action of the Galois group $\Gamma_{K_j}$.
\end{defn}

\begin{rem} Let $j'$ be $\infty$ or an integer divisible by $j$, $b$ an element of $G(K_j)$, and write $b_{j'}$ for $b$ viewed as an element of $G(K_{j'})$. Then, we may observe the following:
\begin{itemize}
\item As $G(K_j)\cap \mc{G}(W_{j'})\sigma(\mu)(p)\mc{G}(W_{j'})=\mc{G}(W_j)\sigma(\mu)(p)\mc{G}(W_j)$, we have that $b$ belongs to $\mc{G}(W_j)\sigma(\mu)(p)\mc{G}(W_j)$ if and only if $b_{j'}$ belongs to $\mc{G}(W_{j'})\sigma(\mu)(p)\mc{G}(W_{j'})$. 
 \item Remark \ref{rem:base-change-deformation} shows that $\mc{D}_\mathsf{K}(\mc{G},b,\mu)_C\simeq \mc{D}_\mathsf{K}(\mc{G},b_{j'},\mu)_C$.
 \end{itemize}
From this we deduce that the underlying $\ov{\bb{Q}}_\ell$-modules $H^i(\mc{G},b,\mu,\mathsf{K})$ and $H^i(\mc{G},b,\mu)$ are independent of $j$, as the notation suggests; in fact so is the $\ov{\Q}_\ell[\mc{G}(\bb{Z}_p)]$-module structure on $H^i(\mc{G},b,\mu)$.
\end{rem}

\begin{rem} By the Hochschild--Serre spectral sequence we have that
\begin{equation*}
    H^i(\mc{G},b,\mu)^{\mathsf{K}}\simeq H^i(\mc{G},b,\mu,\mathsf{K}),
\end{equation*}
and thus $H^i(\mc{G},b,\mu)$ is a smooth $\mc{G}(\bb{Z}_p)$-representation over $\ov{\Q}_\ell$. 
\end{rem}
The following proposition implies that $H^i(\mc{G},b,\mu)$ is, in fact, an admissible representation.

\begin{prop}\label{prop:finite-dimensional} For any $i\geqslant 0$ the $\ov{\Q}_\ell$-spaces $H^i(\mc{G},b,\mu,\mathsf{K})$ are finite-dimensional, and vanish for $i> d=\dim(\mf{D}(\mc{G},b,\mu))$.
\end{prop}
\begin{proof} We may assume that $\mathsf{K}=\mathsf{K}_n$ for some $n$. As $\pi_{\mathsf{K}_n}\colon \mc{D}_n(\mc{G},b,\mu)\to \mc{D}(\mc{G},b,\mu)$ is a $\mc{G}(\Z/p^{n})$-cover we have a canonical isomorphism 
\begin{equation*}
   H^i(\mc{G},b,\mu,\mathsf{K}_n)\simeq H^i_\et(\mc{D}(\mc{G},b,\mu)_C,(\pi_{\mathsf{K}_n})_\ast(\ov{{\Q}}_\ell)), 
\end{equation*}
and $(\pi_{\mathsf{K}_n})_\ast(\ov{{\Q}}_\ell)$ is a $\ov{\Q}_\ell$-local system. As $\mc{D}(\mc{G},b,\mu)\simeq \mf{D}(\mc{G},b,\mu)_\eta$, and $\mf{D}(\mc{G},b,\mu)\to \Spf(W_j)$ is affine and formally of finite type, the claim then follows from Proposition \ref{prop:coh-properties}.
\end{proof}

\begin{nota} Let $\mathsf{K}\subseteq \mc{G}(\Z_p)$ be a compact open subgroup. Set 
\begin{equation*}
    H^\ast(\mc{G},b,\mu)=\sum_{i=0}^d (-1)^i H^i(\mc{G},b,\mu),\qquad  H^\ast(\mc{G},b,\mu,\mathsf{K})=\sum_{i=0}^d (-1)^i H^i(\mc{G},b,\mu,\mathsf{K}),
\end{equation*}
considered as objects of the Grothendieck group of $\mc{G}(\bb{Z}_p)$-admissible modules over the rings $\ov{\Q}_\ell[\Gamma_{K_j}]\otimes_{\ov{\Q}_\ell} \mc{H}_{\ov{\Q}_\ell}(\mc{G}(\bb{Z}_p))$ and $\ov{\Q}_\ell[\Gamma_{K_j}]\otimes_{\ov{\Q}_\ell} \mc{H}_{\ov{\Q}_\ell}(\mc{G}(\bb{Z}_p),\mathsf{K})$, respectively.
\end{nota}

\begin{defn} Let $\tau$ be an element of $W_K$ with $\mathsf{v}(\tau)=j\geqslant 1$ (and so contained in $\Gamma_{K_j}$) and $h$ an element of $\mc{H}(\mc{G}(\Z_p))$. We then define the \emph{test function} 
\begin{equation*}
    {^\ell}\phi_{\tau,h}^{\mc{G},\mu}\colon G(K_j)\to \ov{\Q}_\ell,\quad b\mapsto \tr\left(\tau\times h| H^\ast(\mc{G},b,\mu)\right).
    \end{equation*}
  \end{defn}

Our main result in this section is the following. 

\begin{thm}\label{thm:phi-nice-properties} The function ${^\ell}\phi_{\tau,h}^{\mc{G},\mu}$ is an element of $\mc{H}(G(K_j))$ which is independent of $\ell$.
\end{thm}

\begin{proof} By Observation \ref{obs:Hecke-action-desc} we need only prove the claimed properties for
\begin{equation*}
    {^\ell}\phi^{\mc{G},\mu}_{\tau,e(g,n)}(b)={^\ell}\phi^{\mc{G},\mu}_{\tau,g,n}(b)\defeq \tr\left(\tau\times \ov{g}| H^\ast(\mc{G},b,\mu,\mathsf{K}_n)\right),
\end{equation*}
as $g$ and $n$ vary. We do this in steps, with \textbf{Step 2} inspired by \cite[Proposition 4.3]{ScholzeLK}.

\medskip

\paragraph*{Step 1: Rationality and $\ell$-independence.} Write $\ms{Q}_b^\mr{univ}\defeq T_\et(\mf{Q}_b^\mr{univ})$ and $\ms{Q}_b^\mr{univ}\otimes\bb{Z}/p^n$ for the induced $\mc{G}(\bb{Z}/p^n)$-torsor. Then, $\ms{Q}_b^\mr{univ}\otimes\bb{Z}/p^n$ is represented by $\pi_{\mathsf{K}_n}\colon \mc{D}_n(\mc{G},b,\mu)\to \mc{D}(\mc{G},b,\mu)$. Thus, by Proposition \ref{prop:independence-of-l-bounded}, our claim will follow if we can show that $\ms{Q}_b^\mr{univ}\otimes \bb{Z}/p^n$ extends to the $p$-bounded generic fiber $\Spf(R_{\mc{G},\mu,j},(p))_\eta$; see Definition \ref{defn:p-bounded-generic-fiber}. 

But, this follows from integrability as in Remark \ref{rem:bounded-deformation-space-BT}. From loc.\@ cit.\@ we see that the $(\mc{G},\mu)$-aperture $\mf{Q}_b^\mr{univ}$ on $\Spf(R_{\mc{G},\mu,j},\mf{m}_{\mc{G},\mu,j})$ extends to a $(\mc{G},\mu)$-aperture on $\Spf(R_{\mc{G},\mu,j},(p))$, and applying \'etale realization provides the desired extension of $\ms{Q}_b^\mr{univ}$, and thus $\ms{Q}_b^\mr{univ}\otimes\bb{Z}/p^n$.

\medskip

\noindent\textbf{Step 2: Local constancy.} Fix $b$ in $\mc{G}(W_j)\sigma(\mu)(p)\mc{G}(W_j)$ and consider
\begin{equation*}
F_b\colon G(K_\infty)\to G(K_\infty),\qquad c\mapsto cb\sigma(c)^{-1}.
\end{equation*}
This map is open by \cite[Lemma 4.4]{ScholzeLK}, and thus so is the following neighborhood of $b$
\begin{equation*}
U_b^N\defeq F_b\bigg(\ker\left(\mc{G}(W_\infty)\to \mc{G}\left(W_\infty/p^N\right)\right)\bigg)\cap \mc{G}(W_j)\sigma(\mu)(p)\mc{G}(W_j).
\end{equation*}
Suppose that $b'$ lies in $U_b^N$, and fix $c$ in $\ker(\mc{G}(W_\infty)\to\mc{G}(W_\infty/p^N))$ with $b'=cb\sigma(c)^{-1}$. Letting $b_\infty$ and $b'_\infty$ denote the elements $b$ and $b'$ thought of as elements of $\mc{G}(W_\infty)\sigma(\mu)(p)\mc{G}(W_\infty)$, we see that $c$ gives an isomorphism $\mf{D}(\mc{G},b_\infty,\mu)\isomto \mf{D}(\mc{G},b'_\infty,\mu)$ and so by Remark \ref{rem:base-change-deformation} an isomorphism 
\begin{equation*}
\vartheta_c\colon \mc{D}_n(\mc{G},b,\mu)_C\isomto \mc{D}_n(\mc{G},b',\mu)_C,
\end{equation*}
which is evidently $\mc{G}(\bb{Z}/p^n)$-equivariant but, a priori, not $\tau$-equivariant. In particular, let us denote by $\tau_b$ and $\tau_{b'}$ the action of $\tau$ on the left-hand and right-hand sides, respectively.

Set $a_c=c^{-1}\tau^{-1}(c)=c^{-1}\sigma^{rj}(c)$. Since $b$ and $b'$ are fixed by $\sigma^{rj}$, comparison with $b'=cb\sigma(c)^{-1}$ gives $a_cb\sigma(a_c)^{-1}=b$. Moreover, $a_c=1\mod{p^N}$, so $a_c$ is in $J_{b_\infty}^{\mr{int}}(N)$.
Naturality for the scalar Galois action then gives $\tau_{b'}\circ\vartheta_c=\vartheta_{\tau^{-1}(c)}\circ\tau_b$, and thus
\begin{equation*}
 \vartheta_c^{-1}\circ(\tau_{b'}\times\ov{g})\circ\vartheta_c
 =a_c\circ(\tau_b\times\ov{g}).
\end{equation*}
So, assume that for $N\gg0$ any such $a_c$ acts trivially on cohomology. Then, we deduce that, in fact, the map induced by $\vartheta_c$ on cohomology is $\tau\times\ov{g}$-equivariant and thus
\begin{equation*}
 {^\ell}\phi^{\mc{G},\mu}_{\tau,g,n}(b)= \tr\left(\tau\times \ov{g}| H^\ast(\mc{G},b,\mu,\mathsf{K}_n)\right)=\tr\left(\tau\times \ov{g}| H^\ast(\mc{G},b',\mu,\mathsf{K}_n)\right)= {^\ell}\phi^{\mc{G},\mu}_{\tau,g,n}(b')
\end{equation*}
Constancy of ${^\ell}\phi^{\mc{G},\mu}_{\tau,g,n}$ on $U_b^N$ for an element $b$ in $\mc{G}(W_j)\sigma(\mu)(p)\mc{G}(W_j)$ handles local constancy on $\mc{G}(W_j)\sigma(\mu)(p)\mc{G}(W_j)$. But, as this set is closed and ${^\ell}\phi^{\mc{G},\mu}_{\tau,g,n}$ is identically zero on the complement, it actually fully handles local constancy.

To prove that for $N\gg 0$ any such $a_c$ acts trivially on cohomology it suffices to prove the following stronger claim: There exists some $N\geqslant 0$ such that $J_{b_\infty}^\mr{int}(N)$ acts trivially on $H^i(\mc{G},b_\infty,\mu,\mathsf{K}_n)$ for all $i$. To prove this, observe that from \cite[Lemma 3.9.2]{HuberEC} and the rising open cover $\mc{D}_n(\mc{G},b_\infty,\mu)=\bigcup_m Y_{\mathsf{K}_n,m}$, we deduce that
\begin{equation*}
H^i(\mc{G},b_\infty,\mu,\mathsf{K}_n)=\varprojlim_m H^i_\et((Y_{\mathsf{K}_n,m})_C,\ov{\bb{Q}}_\ell).
\end{equation*}
As the left-hand side is finite-dimensional by Proposition \ref{prop:finite-dimensional} we deduce that for some $m$
\begin{equation*}
H^i(\mc{G},b_\infty,\mu,\mathsf{K}_n)\hookrightarrow H^i_\et((Y_{\mathsf{K}_n,m})_C,\ov{\bb{Q}}_\ell).
\end{equation*}
Thus, for fixed $i$ the claim follows by combining Proposition \ref{prop:uniform-continuity-K-level} with \cite[Proposition 2.6]{ScholzeLK}. Taking the maximum of $N$ over $i=0,\ldots,d$ then gives the desired claim.

\medskip

\paragraph*{Step 3: Compact support.} As we know that ${^\ell}\phi^{\mc{G},\mu}_{\tau,g,n}$ is locally constant, it suffices to show that it is supported in a compact subset of $G(K_j)$. But, by definition, we have that the support is contained in the compact set $\mc{G}(W_j)\sigma(\mu)(p)\mc{G}(W_j)$, from where the claim follows.
\end{proof}

\begin{defn}\label{defn:phitauh} Let $\tau$ be an element of $W_K$ with Frobenius-valuation $\mathsf{v}(\tau)=j$ (and so contained in $\Gamma_{K_j}$) and $h$ an element of $\mc{H}(\mc{G}(\Z_p))$. We then write $\phi_{\tau,h}^{\mc{G},\mu}$ for the unique element of $\mc{H}(G(K_j))$ equal to ${^\ell}\phi_{\tau,h}^{\mc{G},\mu}$ for any $\ell\ne p$. When $\mc{G}$ and $\mu$ are clear from context we shorten this to $\phi_{\tau,h}$. 
\end{defn}

\subsubsection{A formula for $\phi_{\tau,h}$ when $\mu$ is central} We end this section by giving a description of $\phi_{\tau,h}$ in the simplest case: when $\mu$ is central (e.g., $\mc{G}$ is a torus). The reason that this centrality assumption is so simplifying is that it forces $\mf{D}(\mc{G},b,\mu)$ to be a point.

Our ultimate calculation will use $(\mc{G},\mu)$-apertures with $\mu$ central in two specific cases: (a) the setting where $\mu$ is trivial but $b$ is arbitrary; (b) the setting where $\mu$ is arbitrary but $b=\sigma(\mu)(p)$, the so-called \emph{Lubin--Tate case}.

\medskip

\paragraph*{The $\mu$-trivial case.} We first discuss some features of $\mr{BT}^{\mc{G},\mu}_\infty$ when $\mu$ is trivial.

\begin{setup}
Let $u$ be an element of $\mc{G}(W_j)=\mc{G}(W_j)\sigma(\mathbf{1})(p)\mc{G}(W_j)$, where $\mathbf{1}\colon \bb{G}_{m,\Z_p}\to \mc{G}$ is the trivial cocharacter. For notational clarity, we depart from the convention of Setup \ref{setup:test-functions} and write $\mf{Q}_{u}$ for the $(\mc{G},\mathbf{1})$-aperture over $k_j$ corresponding to $(\mc{G}_{W_j},u)$ as in Setup \ref{setup:simplifying-over-finite-fields}.
\end{setup}

\begin{rem} By \cite[Proposition 9.4.1]{GMM} one has an equivalence $F\colon \mr{BT}^{\mc{G},\mathbf{1}}_\infty\simeq B\mc{G}(\bb{Z}_p)$; in particular $\mr{BT}^{\mc{G},\mathbf{1}}_\infty\to \Spf(\Z_p)$ is formally \'etale. Under this equivalence $T_\et(\mf{Q})\simeq F(\mf{Q})_\eta$.
\end{rem}

\begin{nota}\label{nota:aperture-lift} Since $\mr{BT}^{\mc{G},\mathbf{1}}_\infty\to \Spf(\bb{Z}_p)$ is formally \'etale, the $(\mc{G},\mathbf{1})$-aperture $\mf{Q}_{u}$ admits a unique deformation to $W_j$ which we write $\wt{\mf{Q}}_u$.
\end{nota} 

\begin{constr}\label{constr:integral-Shintani-class} Note that $T_\et(\wt{\mf{Q}}_u)$ is a $\mc{G}(\Z_p)$-bundle over $K_j$, and so trivial over $C$. Choose an element $\alpha$ of $T_\et(\wt{\mf{Q}}_u)(C)$, i.e., a trivialization $\mc{G}(\Z_p)_C\isomto T_\et(\wt{\mf{Q}}_u)_C$. 

For any $\tau$ in $W_K$ with $\mathsf{v}(\tau)=j$ there is a unique element $s_{u,\alpha,\tau}$ of $\mc{G}(\Z_p)$ such that $\tau(\alpha)=\alpha s_{u,\alpha,\tau}$. It is simple to check that, up to $\mc{G}(\Z_p)$-conjugacy, $s_{u,\alpha,\tau}$ is independent of $\alpha$ and $\tau$, and only depends on the $\sigma$-conjugacy class of $u$. 
\end{constr}

\begin{defn} Construction \ref{constr:integral-Shintani-class}  gives us a map 
\begin{equation*}
\Sh_j\colon \mc{G}(W_j)/\sigma\text{-conj.}\to\mc{G}(\Z_p)/\text{conj.}
\end{equation*}
which we call the \emph{integral Shintani class} map. 
\end{defn}

We would now like to explain the terminology for the above map.

\begin{constr}
To explain this terminology, set $d$ to be $[K_j:\Q_p]=rj$ and, for $n\geqslant 1$, let $\mathsf{G}_n$ be the connected group $\ov{\bb{F}}_p$-scheme obtained by the Greenberg construction applied to $\mc{G}_{W_n(\ov{\bb{F}}_p)}$; see \cite[Definition 2.1, Proposition 2.2, and Proposition 4.3]{StasinskiGreenberg}. If $F$ denotes arithmetic Frobenius, then there are natural identifications:
\begin{equation}\label{eq:Greenberg-fixed-points}
 \mathsf{G}_n(\ov{\bb{F}}_p)^{F^d}\simeq \mc{G}(W_j/p^n),
 \qquad
 \mathsf{G}_n(\ov{\bb{F}}_p)^F\simeq \mc{G}(\bb{Z}/p^n).
\end{equation}
Thus, taking $F_1=F^d$ and $F_2=F$ in the
usual Shintani norm gives a bijection
\begin{equation}\label{eq:finite-level-Shintani}
 \mc{G}(W_j/p^n)/\!\sim_\sigma \isomto \mc{G}(\bb{Z}/p^n)/\mr{conj.};
\end{equation}
see \cite[\S1.1]{ShojiSorlin}. We have
identified the twisted-conjugacy convention of loc.\@ cit.\@ with ours
by replacing the conjugating element with its inverse.

More explicitly, let $u_n$ be an element of $\mc{G}(W_j/p^n)$, and choose $c$ in $\mathsf{G}_n(\ov{\bb{F}}_p)$ with $u_n=cF(c)^{-1}$ (which exists by Lang's theorem). Setting $\alpha=c^{-1}$ gives
$u_n=\alpha^{-1}F(\alpha)$, and the definition of the Shintani norm
shows that the image of $u_n$ in the right-hand side of
\eqref{eq:finite-level-Shintani} is represented by
\begin{equation}\label{eq:finite-level-Shintani-representative}
 a_n\defeq F^d(c)^{-1}c\in \mathsf{G}_n(\ov{\bb{F}}_p)^F=\mc{G}(\Z/p^n);
\end{equation}
cf.\@ \cite[\S1.4 and (1.5.1)]{ShojiSorlin}.
\end{constr}

We now compare this construction with the integral Shintani class defined above.
 
\begin{prop} For $u$ in $\mc{G}(W_j)$ one has that the image of $\Sh_j(u)$ in $\mc{G}(\Z/p^n)/\mr{conj.}$ agrees with the Shintani construction applied to $u_n=u\mod p^n$ in $\mc{G}(W_j/p^n)$.
\end{prop}
\begin{proof}Write $\widetilde{\mf{Q}}_{u,n}$ for the $n$-truncation of $\widetilde{\mf{Q}}_u$. The equivalence $\mr{BT}^{\mc{G},\mathbf{1}}_n\isomto B\mc{G}(\bb{Z}/p^n)$ from \cite[Proposition 9.4.1]{GMM} carries
$\widetilde{\mf{Q}}_{u,n}$ to the torsor $\Gamma_\syn(\widetilde{\mf{Q}}_{u,n}) =\mr{Isom}(\mc{G},\widetilde{\mf{Q}}_{u,n})$ of syntomic trivializations. Unraveling the quotient description in
the proof of loc.\@ cit.\@, and using our convention for
$\sigma$-conjugacy, identifies the geometric points of this torsor
with the elements $c$ satisfying $u_n=cF(c)^{-1}$. 

Applying \'etale realization to the universal syntomic
trivialization gives a $\mc{G}(\bb{Z}/p^n)$-equivariant map from the
generic fiber of $\Gamma_\syn(\widetilde{\mf{Q}}_{u,n})$ to (the
frame torsor of) $T_\et(\widetilde{\mf{Q}}_{u,n})$ which must be an isomorphism.\footnote{This is also immediate from the description of the
\'etale realization by Frobenius-fixed frames in
\cite[Construction 3.8.1]{MY}.} Finally, an element $\tau$ of $W_K$ with $\mathsf{v}(\tau)=j$ acts on
this unramified torsor through geometric Frobenius on $k_j$, which is
$F^{-d}$. Since $a_n$ belongs to $\mathsf{G}_n(\ov{\bb{F}}_p)^F$, one sees from
\eqref{eq:finite-level-Shintani-representative} that $F^{-d}(c)=ca_n$. Thus, $\tau$ acts on the frame $c$ by right multiplication by the
usual Shintani norm of $u_n$, as desired.
\end{proof}

\begin{rem}We finally remark that one can view the integral Shintani class construction as an integral analogue of Kottwitz's rational norm construction; see Definition \ref{defn:norm-twisted-transfer} below.
\end{rem}

\begin{eg}\label{eg:Shintani-for-tori} Suppose that $\mc{G}=\mc{T}$ is a torus. Then, the integral Shintani class is something more classical; it is the norm map:
\begin{equation*}
\Sh_j(u)=\Nr_{K_j/\Q_p}(u)=u\sigma(u)\cdots \sigma^{rj-1}(u).
\end{equation*}
Indeed, after passing to an unramified extension which splits $\mc{T}$, it suffices to treat $\mc{T}=\bb{G}_m$, where one can make the same calculation as in the discussion following \cite[Definition 3.3]{ScholzeLK}.
\end{eg}

\medskip

\paragraph*{The Lubin--Tate case.} We now turn our attention to the aforementioned \emph{Lubin--Tate case}.

\begin{setup} We set the following notation:
\begin{itemize}
\item $\mc{Z}\defeq Z(\mc{G})^\circ$, and observe that $\mu$ factorizes through $\mc{Z}_W$;
\item $\mc{T}_0\defeq \mr{Res}_{W/\bb{Z}_p}\bb{G}_{m,W}$;
\item normalize the Artin reciprocity map
\begin{equation*}
\Art_K\colon W_K^\mr{ab}\isomto \mc{T}_0(\Q_p)=K^\times
\end{equation*}
so that $\mr{Art}_K$ sends geometric Frobenius to $p$. 
\end{itemize}
\end{setup}

\begin{rem}Using the decomposition $W\otimes_{\bb{Z}_p}W\simeq W^r$, let $\mu_0\colon\bb{G}_{m,W}\to (\mc{T}_0)_W$ be the inclusion in the first factor. Define the \emph{reflex norm} attached to $\mu$ to be the morphism
\begin{equation}\label{eq:central-reflex-norm}
r_\mu\colon \mc{T}_0
\xrightarrow{\mr{Res}_{W/\bb{Z}_p}(\mu)}
\mr{Res}_{W/\bb{Z}_p}\mc{Z}_W
\xrightarrow{\Nr_{W/\bb{Z}_p}}
\mc{Z}.
\end{equation}
The construction in \cite[\S10.4.2]{GMM} gives the identity
\begin{equation}\label{eq:reflex-norm-recovers-mu}
(r_\mu)_W\circ\mu_0=\mu.
\end{equation}
\end{rem}

\begin{nota} For $\tau$ in $W_K$ with $\mathsf{v}(\tau)=j$, set $c_{\mu,\tau}\defeq r_\mu\left(p^{j}\Art_K(\tau)^{-1}\right)$ in $\mc{Z}(\bb{Z}_p)$.
\end{nota}

\begin{setup}
Let $\mf{G}_0$ be the Lubin--Tate formal $W$-module associated with $pX+X^q$, and let $\mf{L}_0^{-1}$ be its $(\mc{T}_0,\mu_0^{-1})$-aperture under the covariant construction of \cite[Proposition 11.7.6]{GMM}, with the type written in our
filtration convention. Writing $\iota\colon\mc{T}_0\to\mc{T}_0$ for inversion, set
$\mf{L}_0\defeq\iota_*\mf{L}_0^{-1}$. Thus $\mf{L}_0$ has type $\mu_0$, and its filtered $F$-crystal has Frobenius factor $\sigma(\mu_0)(p)$.
\end{setup}

\begin{defn}\label{defn:Lmu} Set $\mf{L}_\mu$ to be the pushout of $\mf{L}_0$ along $r_\mu$. The filtered $F$-crystal underlying $\mf{L}_\mu$ then has Frobenius factor $\sigma(\mu)(p)$ by \eqref{eq:reflex-norm-recovers-mu}.
\end{defn}

\begin{rem} By \cite[Proposition 11.8.2]{GMM} and compatibility with duality,
the rank-one $W$-local system associated with $T_\et(\mf{L}_0)$ is the ordinary $W$-linear dual of $T_p(\mf{G}_0)$.
\end{rem}

\begin{nota} Write $\chi_\mr{LT}\colon W_K\to W^\times$ for the character of the rank-one $W$-module $T_p(\mf{G}_0)$.
\end{nota}

\begin{rem} From Lubin--Tate theory one has that if $\mathsf{v}(\tau)=j$, then $\chi_\mr{LT}(\tau)=p^{-j}\Art_K(\tau)$.
\end{rem}

\medskip

\paragraph*{A description of $\phi_{\tau,h}$ for central $\mu$} With the above discussion, we are now ready to state and prove the formula for $\phi_{\tau,h}$ when $\mu$ is central, generalizing \cite[Proposition 4.10]{ScholzeLK}.

\begin{defn}For a conjugacy class $c$ in $\mc{G}(\Z_p)$, and an element $h$ of $\mc{H}(\mc{G}(\Z_p))$, we write
\begin{equation*}
\mr{Av}(h)(c)\defeq \int_{\mc{G}(\Z_p)}h(x^{-1}ax)\,dx,
\end{equation*}
where $a$ is any representative of $c$; this is independent of the choice of $a$.
\end{defn}

\begin{setup}\label{setup:setup-for-central-calc} We suppose that $\mu$ is central and that $b$ is in $\mc{G}(W_j)\sigma(\mu)(p)\mc{G}(W_j)$. Set
\begin{equation*}
u_b\defeq b\sigma(\mu)(p)^{-1}\in \mc{G}(W_j)
\end{equation*}
and define $a_{b,\tau}$, up to conjugacy in $\mc{G}(\Z_p)$, by the formula $[a_{b,\tau}]=\left[c_{\mu,\tau}\operatorname{Sh}_j(u_b)\right]$; the product is unambiguous because $c_{\mu,\tau}$ is central.
\end{setup}

\begin{prop}\label{prop:central-character-calc} With setup as in \emph{Setup \ref{setup:setup-for-central-calc}} one has:
\begin{equation}\label{eq:central-test-function}
\phi_{\tau,h}^{\mc{G},\mu}(b)= \begin{cases}
\mr{Av}(h)([a_{b,\tau}]) & \emph{if}\quad b\in \mc{G}(W_j)\sigma(\mu)(p)\mc{G}(W_j),\\
0 & \emph{otherwise.}
\end{cases}
\end{equation}
In particular, if $\mc{G}=\mc{T}$ is a torus, then
\begin{equation}\label{eq:central-torus-test-function}
\phi_{\tau,h}^{\mc{G},\mu}(b)= \begin{cases}
h\left(r_\mu\left(p^j\Art_K(\tau)^{-1}\right)
 \Nr_{K_j/\Q_p}\left(b\sigma(\mu)(p)^{-1}\right)\right) & \emph{if}\quad b\in \mc{G}(W_j)\sigma(\mu)(p)\mc{G}(W_j),\\
0 & \emph{otherwise.}
\end{cases}
\end{equation} 
\end{prop}

\begin{proof} The formula for tori follows from the general formula via Example \ref{eg:Shintani-for-tori} and a simple calculation. Thus, we focus on the first formula. Moreover, as both sides vanish when $b$ is not in $\mc{G}(W_j)\sigma(\mu)(p)\mc{G}(W_j)$, we may assume that $b$ lies in this double coset.

As $\mu$ is central, one has that $\wh{\mc{U}}^-_{\mu,j}=\Spf(W_j)$ and thus by Theorem \ref{thm:explicit-deformation} we have $\mc{D}(\mc{G},b,\mu)=\Spa(K_j)$. Let $P_b\defeq \mc{D}_\infty(\mc{G},b,\mu)(C)$ which is a trivial right $\mc{G}(\Z_p)$-torsor, and fix a base point $\alpha$ in $P_b$. There is a unique $m_{b,\tau,\alpha}$ in $\mc{G}(\Z_p)$ such that $\tau(\alpha)=\alpha m_{b,\tau,\alpha}$. 
Changing $\alpha$ conjugates $m_{b,\tau,\alpha}$, so the conjugacy class $[m_{b,\tau,\alpha}]$ is independent of all choices.

We first compute $\phi_{\tau,h}^{\mc{G},\mu}(b)$ in terms of this class. Assume first that $h=e(g,n)$ for some $g$ and some $n$. Then, as $\mc{D}_n(\mc{G},b,\mu)_C=P_b/\mathsf{K}_n$ is a finite disjoint union of points, we have the following:
\begin{equation*}
H^\ast(\mc{G},b,\mu,\mathsf{K}_n)=H^0(\mc{G},b,\mu,\mathsf{K}_n)=\ov{\Q}_\ell[\mc{G}(\Z_p)/\mathsf{K}_n].
\end{equation*}
Under this identification, $\tau$ carries $x\mathsf{K}_n$ to $m_{b,\tau,\alpha}x\mathsf{K}_n$, and $g$ carries $x\mathsf{K}_n$ to $xg^{-1}\mathsf{K}_n$. Thus, $\tau\times\ov{g}$ carries $x\mathsf{K}_n$ to $m_{b,\tau,\alpha}xg^{-1}\mathsf{K}_n$. Consequently,
\begin{align*}
\tr\left(\tau\times e(g,n)\mid H^0(\mc{G},b,\mu,\mathsf{K}_n)\right)
&=
\#\left\{x\mathsf{K}_n:
 x^{-1}m_{b,\tau,\alpha}x\in g\mathsf{K}_n
\right\}\\
&=
\int_{\mc{G}(\Z_p)}
 e(g,n)(x^{-1}m_{b,\tau,\alpha}x)
 \,dx\\ &=\mr{Av}(e(g,n))([m_{b,\tau,\alpha}]).
\end{align*}
By summing we deduce $\phi_{\tau,h}^{\mc{G},\mu}(b)=\mr{Av}(h)([m_{b,\tau,\alpha}])$ in general. It remains to show $[m_{b,\tau,\alpha}]=[a_{b,\tau}]$.

To this end, set $z_\mu=\sigma(\mu)(p)$ and $u_b=bz_\mu^{-1}$; note that $z_\mu$ belongs to $\mc{Z}(K)$ and that $u_b$ lies in $\mc{G}(W_j)$. Let $\wt{\mf{Q}}_{u_b}$ be as in Notation \ref{nota:aperture-lift}, and let $\mf{L}_{\mu,j}$ be the base change of $\mf{L}_\mu$ from Definition \ref{defn:Lmu} to $W_j$. Consider the multiplication map
\begin{equation*}
m\colon \mc{G}\times\mc{Z}\to\mc{G},\qquad (g,z)\longmapsto gz,
\end{equation*}
which gives rise to a $(\mc{G},\mu)$-aperture over $W_j$ by $\mf{Q}=m_\ast(\wt{\mf{Q}}_{u_b},\mf{L}_{\mu,j})$. Under the equivalences in Setup \ref{setup:simplifying-over-finite-fields}, $\mf{Q}_{k_j}$ corresponds to $\mc{G}_{W_j}$ with Frobenius $u_bz_\mu=b$; thus $\mf{Q}$ is a deformation of the $(\mc{G},\mu)$-aperture over $k_j$ associated to $b$. As $\wh{\mc{U}}_{\mu,j}^-=\Spf(W_j)$, we then conclude from Theorem \ref{thm:explicit-deformation} that $\mf{Q}\simeq \mf{Q}_b^\mr{univ}$. Note that he two factors definining $\mf{Q}$ have \'etale monodromy $\mr{Sh}_j(u_b)$ and $r_\mu(\chi_\mr{LT}(\tau)^{-1})=c_{\mu,\tau}$, respectively. Compatibility with products and pushouts therefore gives us the desired equality $[m_{b,\tau,\alpha}]=[c_{\mu,\tau}\operatorname{Sh}_j(u_b)]=[a_{b,\tau}]$.
\end{proof}

\begin{rem} Let us observe that for $K=\Q_p$, $\mc{G}=\bb{G}_m$, and $\mu(t)=t^{-1}$, the formula \eqref{eq:central-torus-test-function} is supported on
$p^{-1}W_j^\times$ and there it is equal to
\begin{equation*}
h\left(\Art_{\Q_p}(\tau)\Nr_{\Q_{p^j}/\Q_p}(b)\right).
\end{equation*}
This then is precisely the formula of \cite[Proposition 4.10]{ScholzeLK}, with $b=\delta_{\bb{G}_m}$.
\end{rem}

\subsection{A shtuka-theoretic description of the test functions}
\label{ss:shtuka-description-test-functions} In this section we give an alternative definition of $\phi_{\tau,h}$ in terms of local Shimura varieties. Throughout we assume $j<\infty$.

\begin{nota} Denote by $\mathbf{Perf}_{A}$, for a $p$-adically complete ring $A$, the category of characteristic $p$ perfectoid spaces $S$ equipped with a map $S\to \Spd(A)$. For any other undefined piece of notation or terminology, we direct the reader to the comprehensive discussion in \cite[\S2-3]{PappasRapoportI}. 
\end{nota}

\begin{defn} A \emph{tame local Shimura datum} is a triple $(\mc{H},c,\{\nu\})$ where $\mc{H}$ is a parahoric group $\mathbb{Z}_p$-scheme with generic fiber $H$, $\{\nu\}$ is a conjugacy class of minuscule cocharacters of $H_{\overline{\mathbb{Q}}_p}$,  and $c$ is an element of $H(\breve{\mathbb{Q}}_p)$ inducing an element of $B(H,-\nu)$.
\end{defn}

\begin{defn}\label{defn:integral-local-shimura-variety} For a tame local Shimura datum $(\mc{H},c,\{\nu\})$ with reflex field $F$, define
\begin{equation*}
    \mc{M}_{\mathcal{H},c,\{\nu\}}^\mr{int}\colon \mathbf{Perf}_{\mc{O}_{\breve{F}}}\to \mathbf{Set},
\end{equation*}
by assigning to $S\to \Spd(\mc{O}_{\breve{F}})$ the set of isomorphism classes of tuples $(S^\sharp, \mathscr{P}, \phi_\mathscr{P}, i_r)$, with:
\begin{itemize}
    \item $S^\sharp$ an untilt of $S$ over $\mc{O}_{\breve{F}}$ associated with $S \to \Spd(\mc{O}_{\breve{F}})$;
    \item $(\mathscr{P},\phi_\mathscr{P})$ an $\mc{H}$-shtuka on $S$ with leg along $S^\sharp$ bounded by $\{\nu\}$; see \cite[Definition 2.4.3]{PappasRapoportI};
    \item and $i_r$ a framing (see loc.\@ cit.\@).
\end{itemize}
By \cite[\textsection 25.1]{ScWeBLp}, the presheaf $\mc{M}_{\mc{H},c,\{\nu\}}^\mr{int}$ is a small $v$-sheaf (in the sense of \cite[Definition 12.1]{ScholzeECD}), called the \emph{integral local Shimura variety} associated with $(\mc{H},c,\{\nu\})$.
\end{defn}

\begin{defn}\label{defn:local-shimura-variety} We write $\mc{M}_{\mc{H},c,\{\nu\}}$ for the generic fiber of $\mc{M}^{\mr{int}}_{\mc{H},c,\{\nu\}}$ and call it the \emph{local Shimura variety} associated with $(H,c,\{\nu\})$. To the universal $\mc{H}$-shtuka $(\mathscr{P}^\univ, \phi_{\mathscr{P}^\univ})$ over $\mc{M}_{\mc{H},c,\{\nu\}}$ there is associated a natural $\mc{H}(\Z_p)$-local system $\bb{P}^\univ$; see \cite[\S2.5.1]{PappasRapoportI}. We write $\mc{M}_{\mc{H},c,\{\nu\},\infty}$ for the covering parameterizing trivializations of $\bb{P}^\univ$, and for any compact open subgroup $\mathsf{K}\subseteq\mc{H}(\Z_p)$ we write $\mc{M}_{\mc{H},c,\{\nu\},\mathsf{K}}$ for $\mc{M}_{\mc{H},c,\{\nu\},\infty}/\mathsf{K}$.
\end{defn}

To our fixed triple $(\mc{G},b,\mu)$ we may associate a tame local Shimura datum and a natural point on the special fiber of the associated integral local Shimura variety.

\begin{constr}
Set $\beta=\sigma^{-1}(b_\infty)$. Then, the prismatic
realization of the aperture $b_\infty$ has Frobenius $\beta$; see Remark~\ref{rem:crys-prism} and \cite[Proposition 1.1]{IKY3}. Moreover, itts associated shtuka is bounded by $\{-\mu\}$; see \cite[Example 5.1.9]{ItoDeformation}. We then obtain a local Shimura datum $(\mc{G},b_\infty,\{-\mu\})$, with reflex field $F\subseteq K$.\footnote{In the following we implicitly identify $\mc{O}_{\breve F}$ with $W_\infty$.}

The quasi-isogeny represented by $\beta^{-1}$ identifies the (crystalline realization of the) prismatic $F$-crystal with the isocrystal associated to $b_\infty$: $(\beta^{-1})^{-1}b_\infty\sigma(\beta^{-1})=\beta$. Thus $\beta^{-1}\mc{G}(W_\infty)$ defines a point $x_0$ of $X_{\mc{G}}(b_\infty,\{\mu\})(\ov{k})$ and hence a point of $\mc{M}^{\mr{int}}_{\mc{G},b_\infty,\{-\mu\}}$; see \cite[Remark 5.3.2]{ItoDeformation}.
\end{constr}

We can then define the shtuka-theoretic analogue of the deformation space and Berthelot tube at this point, using the material (and notation) from \cite{GleasonSpecialization}.

\begin{defn} We define $(\mc{M}^{\mr{int}}_{\mc{G},b_\infty,\{-\mu\}})^\wedge_{/x_0}$ to be the formal neighborhood of $\mc{M}^{\mr{int}}_{\mc{G},b_\infty,\{-\mu\}}$ at $x_0$ in the sense of \cite[Definition 4.18]{GleasonSpecialization}.
\end{defn}

\begin{defn}\label{defn:distinguished-shtuka-tube-v2} The \emph{shtuka tube at $b$ with $\mathsf{K}$-level structure} is
\begin{equation*}
(\mc{M}_{\mc{G},b_\infty,\{-\mu\}}^\mr{int})^\circledcirc_{/x_0,\mathsf{K}}\defeq (\mc{M}^{\mr{int}}_{\mc{G},b_\infty,\{-\mu\}})^\wedge_{/x_0}\times_{\mc{M}^\mr{int}_{\mc{G},b_\infty,\{-\mu\}}}\mc{M}_{\mc{G},b_\infty,\{-\mu\},\mathsf{K}}.
\end{equation*}
\end{defn}

\begin{rem}\label{rem:shtuka-tube-descent}
Suppose $j<\infty$, set $d=rj$, and write $B_d=b\sigma(b)\cdots\sigma^{d-1}(b)$. Let $\Phi$ be the $[K:F]^\text{th}$-power of the Weil descent datum on $\mc{M}^\mr{int}_{\mc{G},b_\infty,\{-\mu\}}$ as in \cite[\S3.1.1, equation (3.1.6)]{PappasRapoportI}.\footnote{Recall that we are only assuming that the reflex field $F$ of $\mu$ is contained in $K$.} Since $B_db=b\sigma(B_d)$, we have an endomorphism $[B_d]$ of $\mc{M}^\mr{int}_{\mc{G},b_\infty,\{-\mu\}}$ given by precomposing the framing $i_r$ by $B_d^{-1}$, and we then define $\Psi_{b,j}=\Phi^j\circ[B_d]$. One can check that by setup $\Psi_{b,j}$ gives a Weil descent datum to $K_j$ on $(\mc{M}^{\mr{int}}_{\mc{G},b_\infty,\{-\mu\}})^\wedge_{/x_0}$ and $(\mc{M}_{\mc{G},b_\infty,\{-\mu\}}^\mr{int})^\circledcirc_{/x_0,\mathsf{K}}$ for all $\mathsf{K}$.
\end{rem}

We then have the following comparison between aperture tubes and shtuka tubes.

\begin{thm}[{\cite[Theorem 5.3.5]{ItoDeformation}}]\label{thm:Ito-aperture-shtuka-comparison-v2} There is a Weil-descent-datum equivariant isomorphism 
\begin{equation*}
\Spd(R_{\mc{G},\mu,j}\wh{\otimes}_{W_j} W_\infty)
 \isomto (\mc{M}^{\mr{int}}_{\mc{G},b_\infty,\{-\mu\}})^\wedge_{/x_0},
\end{equation*}
of $v$-sheaves over $\Spd(W_\infty)$. Moreover, the pullback of the universal $\mc{G}$-shtuka is isomorphic to $T_\mr{sht}(\mf{Q}_b^\univ)$, where $T_\mr{sht}$ is the shtuka-realization functor from \cite[Construction 3.20]{IKY1}. 
\end{thm}
\begin{proof} This follows by applying the construction in the proof of
\cite[Theorem 5.3.5]{ItoDeformation} to the prismatic realization
of $\mf{Q}_b^\univ$, with marking $\beta^{-1}$, using
\cite[Proposition 7.1.1]{ItoDeformation}. One then explicitly calculates that this identification is equivariant for the Weil descent data we have defined in Remark \ref{rem:shtuka-tube-descent}.
\end{proof}

The following then is an immediate corollary given \cite[Theorem 3.21]{IKY1} which says that, as $T_\mr{sht}(\mf{Q}_b^\univ)$ is the universal shtuka $(\ms{P}^\univ,\phi_{\ms{P}^\univ})$, $T_\et(\mf{Q}_b^\univ)$ is $\bb{P}^\univ$.

\begin{cor} For every neat compact open subgroup $\mathsf{K}\subseteq\mc{G}(\Z_p)$ one has a Weil-descent-data equivariant isomorphism of $v$-sheaves over $\Spd(K_\infty)$
\begin{equation*}
 \mc{D}_\mathsf{K}(\mc{G},b,\mu)_{K_\infty}^\lozenge\isomto (\mc{M}_{\mc{G},b_\infty,\{-\mu\}}^{\mr{int}})
 ^\circledcirc_{/x_0,\mathsf{K}}.
\end{equation*}
\end{cor}

Using this, we can now immediately rewrite our test functions in purely shtuka-theoretic terms.

\begin{defn}\label{defn:shtuka-tube-cohomology-v2}
For $\ell\ne p$ and $\mathsf{K}\subseteq\mc{G}(\Z_p)$, put
\begin{equation*}
H^i_{\mr{sht}}(\mc{G},b,\mu,\mathsf K) \defeq H^i_\et\bigg((\mc{M}_{\mc{G},b_\infty,\{-\mu\}}^\mr{int})^\circledcirc_{/x_0,\mathsf{K}}\times_{\Spd(K_\infty)}\Spd(C),\ov{\Q}_\ell\bigg).
\end{equation*}
We further set
\begin{itemize}
\item $H^i_{\mr{sht}}(\mc{G},b,\mu)=\varinjlim_{\mathsf K}H^i_{\mr{sht}}(\mc{G},b,\mu,\mathsf K)$;
\item $H^*_{\mr{sht}}(-)$ to be the associated alternating sum of cohomology groups.
\end{itemize}
Each $H^i_\mr{sht}(\mc{G},b,\mu,\mathsf{K})$ carries its natural $\Gamma_{K_j}$-action, using the $K_j$-descent of Remark \ref{rem:shtuka-tube-descent}. Moreover, modifying the level gives $H^i_\mr{sht}(\mc{G},b,\mu)$ a commuting smooth $\mc{G}(\Z_p)$-action.
\end{defn}

\begin{prop}\label{prop:shtuka-formula-phi-v2}
Fix $\tau$ in $W_K$ with $\mathsf v(\tau)=j\geqslant1$, and $h$ in $\mc{H}(\mc{G}(\Z_p))$. Then:
\begin{equation}\label{eq:phitauh-shtuka}
\phi_{\tau,h}^{\mc{G},\mu}(b)
 =
\begin{cases}
\tr\big(\tau\times h\mid H^*_{\mr{sht}}(\mc{G},b,\mu)\big)
 & b\in\mc{G}(W_j)\sigma(\mu)(p)\mc{G}(W_j),\\
0 &\text{otherwise.}
\end{cases}
\end{equation}
\end{prop}

\begin{rem}[A conjectural parahoric extension of test functions]\label{rem:parahoric-extension} While we only defined $H^\ast_\mr{sth}(\mc{G},b,\mu)$ for reductive $\mc{G}$, the same definition makes sense when $\mc{G}$ is parahoric using Definitions \ref{defn:integral-local-shimura-variety} and \ref{defn:local-shimura-variety}. Then, \eqref{eq:phitauh-shtuka} provides a definition of $\phi_{\tau,h}^{\mc{G},\mu}$ for arbitrary parahoric $\mc{G}$. Such a definition is attractive as every reductive group over $\Q_p$ admits a parahoric model.

That said, many of the finer arguments about $\phi_{\tau,h}^{\mc{G},\mu}$ made above require deformation theory, something inaccessible to the shtuka perspective (which is nil-invariant). Thus, while one can currently define the local test functions $\phi_{\tau,h}^{\mc{G},\mu}$, and likely prove a version of the primitive trace formula for Shimura varieties at bad level (see Theorem \ref{thm:prim-trace-formula-II}), their finer properties are not clear. It would be better to define them with $(\mc{G},\mu)$-apertures, when such objects have been defined.
\end{rem}

\subsection{Comparison to previously constructed functions}\label{ss:comparison}

In this final subsection, we compare the local test functions $\phi_{\tau,h}$ constructed above with the functions of \cite{ScholzeLK} and \cite{Youcis}. 

\subsubsection{Comparison to \cite{ScholzeLK}} On \cite[pp.\@ 236-247]{ScholzeLK}, Scholze considers certain $(\mc{G},b,\mu)$ of what he calls (quasi-)EL or PEL type. 

\begin{setup}\label{setup:Scholze-comparison}
Suppose that $\mc{G}$ arises from local (quasi-)EL or PEL data $\mc{D}$ as on \cite[pp.\@ 236-247]{ScholzeLK}, with the center $F$ of $B$ unramified over $\Q_p$ and $p>2$. For clarity, let us write $\mu_\mr{Sch}$ for Scholze's cocharacter and $\Lambda$ for his chosen faithful representation of $\mc{G}$. We further write $w(t)$ for multiplication-by-$t$ map on $\Lambda$, and set $\mu=\mu_\mr{Sch}-w$.
\end{setup}

\begin{rem} For ease of comparison, we comment on the dictionary between the notations from this paper and those used in \cite{ScholzeLK}:
\begin{itemize} 
\item in our notation, Scholze's $E$ and $\kappa_E$ are our $K$ and $k$; respectively,\footnote{Although Scholze further assumes that $K$ is the minimal field of definition of $\mu$, i.e., the reflex field.}
\item Scholze's $r$ equals our $rj=[K_j:\Q_p]$; thus $k_j\simeq\F_{p^{rj}}$ and $K_j\simeq\Q_{p^{rj}}$;
\item Scholze's $\delta$ is our $b$, and his cocharacter is $\mu_\mr{Sch}=\mu+w$.
\end{itemize}
\end{rem}

\begin{setup} Let $\underline{\ov{H}}_\delta$ be a $p$-divisible group with $\mc{D}$-structure over $k_j$ as defined in \cite[Definition 3.3]{ScholzeLK}.
\end{setup}

\begin{nota}
For visual clarity, let us write $\phi_{\tau,h}^{\mr{Sch}}$ for the function in \cite[Definition 4.1]{ScholzeLK}.
\end{nota}

\begin{prop}\label{prop:Scholze-test-function-comparison}
With setup as in \emph{Setup \ref{setup:Scholze-comparison}}, and for a compact open subgroup $\mathsf{K}\subseteq\mc{G}(\Z_p)$, there are natural identifications compatible in $\mathsf{K}$ and equivariant for the Weil descent data:
\begin{equation*}
\mf{D}(\mc{G},b,\mu)\simeq\Spf(R_{\underline{\ov{H}}_\delta}),\qquad \mc{D}_\mathsf{K}(\mc{G},b,\mu)\simeq X_{\underline{\ov{H}}_\delta,\mathsf{K}},
\end{equation*}
where $b=\delta$. Consequently, we have the equality $\phi_{\tau,h}^{\mc{G},\mu}(b)=\phi_{\tau,h}^{\mr{Sch}}(\delta)$.
\end{prop}

\begin{proof} To give an isomorphism $\mc{D}_{\mathsf{K}}(\mc{G},b,\mu)\simeq X_{\underline{\ov{H}}_\delta,\mathsf{K}}$ with the desired compatibilities, it suffices to give an isomorphism $\mf{D}(\mc{G},b,\mu)\simeq\Spf(R_{\underline{\ov{H}}_\delta})$ equivariant for the Weil descent data and matching $T_\et(\mf{Q}_b^\univ)$ with the $\mc{G}(\Z_p)$-bundle obtained via $T_p(\underline{\ov{H}}_\delta^\mr{univ})$ as on \cite[p.\@ 240]{ScholzeLK}.

So, let us write $\mc{M}_{\mr{cov}}$ as in Remark \ref{rem:BT-tate-module-etale-realization}. Then, the filtered crystalline realization (see Theorem \ref{thm:syntomic-DD}) of $\mc{M}_{\mr{cov}}(H)$ is given by $\bb{D}(H)^\vee\simeq\bb{D}(H^D)(1)$, where we write $H^D$ for the Cartier dual. From this, we deduce that Scholze's covariant Frobenius $p\delta\sigma$ becomes $\delta\sigma$, and the filtration has type $\mu_\mr{Sch}-w=\mu$. Thus, \cite[Theorem 3.6]{KimRZ} and \cite[Proposition 3.4]{IKY2} identify the universal filtered crystal with extra structure with the quintuple associated to
$\mf{Q}_b^\univ$ as in Notation \ref{nota:universal-aperture}. 

Note moreover that the covariant construction respects the endomorphisms and twisted polarizations; see \cite[\S11.6]{GMM}. From this, we deduce that its \'etale realization is the universal Tate module by Remark \ref{rem:BT-tate-module-etale-realization}, giving the required identification of the $\mc{G}(\Z_p)$-torsors.
\end{proof}

\begin{rem} In Proposition \ref{prop:Scholze-test-function-comparison} we only verified the agreeance of the functions in the case when both are in the not-identically-zero case of their definition, but equality holds in general. 

Namely, Scholze defines $\phi_{\tau,h}^{\mr{Sch}}(\delta)$ to be zero unless both of the following conditions hold:
\begin{itemize}
\item $\delta$ is associated, by \cite[Proposition 3.10]{ScholzeLK}, to a $p$-divisible group with $\mc{D}$-structure;
\item that $p$-divisible group has controlled cohomology.
\end{itemize}
But, controlled cohomology is automatic in the relevant cases by
\cite[Proposition 3.12]{ScholzeLK}. Moreover, by
\cite[Proposition 4.7]{ScholzeLK} and its proof, the first condition is equivalent to $\delta$ belonging to $\mc{G}(W_j)p^{-1}\sigma(\mu_\mr{Sch})(p)\mc{G}(W_j)$. Since $p^{-1}\sigma(\mu_\mr{Sch})(p)=\sigma(\mu)(p)$ and $b=\delta$, this is precisely our support condition.
\end{rem}

\begin{rem} Let us emphasize the symmetric difference between the cases considered in this paper and those considered in \cite{ScholzeLK}. Scholze's work can only handle EL or PEL situations of type A and C, excluding all orthogonal PEL-type situations, many abelian-type situations, and all exceptional cases. That said, Scholze's setup allows him to deal with some ramified cases, i.e., when $\mc{G}$ is not reductive. \end{rem}

\subsubsection{Comparison with \cite{Youcis}}\label{sss:comparison-with-thesis}
In \cite[\S\S3.2--3.3 and 4.1--4.3]{Youcis}, there is defined local test
functions using towers attached to deformation data of Hodge type and of two
special classes, so-called good data of Type 1 or Type 2, of abelian type. We compare those functions with ours, but in the Type $1$ case, the comparison requires the extra datum introduced below.

\begin{rem}
We use the following translation of conventions.
\begin{itemize}
\item The fields denoted $E$, $\kappa_E$, and $E_j$ in \cite{Youcis} are denoted $K$, $k$, and $K_j$ here.
\item Theorem \ref{thm:syntomic-DD} and \cite{Youcis} both use contravariant Dieudonn\'e theory; the same Hodge filtration is labelled by $\mu$ here and by $\mu^{-1}$ in the convention of \cite[\S3.2]{Youcis}. 
\item The notation $[b]$ from \cite{Youcis} (which we temporarily adopt here) denotes a class of $C_j(\mc{G})$.
\item We will express the constructions of \cite{Youcis} in our Hodge-filtration convention, including the Frobenius pullback in \cite[Lemma 2.5.8]{KimRZ}. Thus, with this, the notion of $[b]$ being \emph{adapted} to $(\mc{G},\mu)$ as in op.\@ cit.\@, means that $[b]$ lies in $C_j(\mc{G},\sigma(\mu))$.
\end{itemize}
\end{rem}

\begin{defn}
A triple
$\mf{d}=(\mc{G},[b],\mu)$ is a \emph{deformation datum of Hodge type} if $[b]$ is adapted to $(\mc{G},\mu)$ and there is a faithful representation
$\rho\colon\mc{G}\to\GL(\Lambda)$ with $(\rho\circ\mu)$-weights in
$\{0,1\}$.

A triple $\mf{d}=(\mc{G},[b],\mu)$ is a \emph{deformation datum of abelian type} if there is a Hodge-type datum
$\mf{d}_1=(\mc{G}_1,[b_1],\mu_1)$ and a surjective central morphism $ f\colon\mc{G}_1\twoheadrightarrow\mc{G}$ identifying $\mf{d}_1^\ad$ with $\mf{d}^\ad$. We call $\mf{d}_1$ an \emph{associated
Hodge-type datum}; it is \emph{good} if $Z(\mc{G}_1)$ has connected fibers
and $\mu_1$ and $\mu$ have the same reflex field. The datum $\mf{d}$ is good
if it admits a good associated datum.

A good abelian-type datum is \emph{of Type~$1$} if an associated datum may be
chosen so that $f^\der\colon\mc{G}_1^\der\to\mc{G}^\der$ is an isomorphism, and is \emph{of Type~$2$} if $\mc{G}$ is adjoint.
\end{defn}

\begin{nota}
We decorate the objects of \cite{Youcis} by a superscript
$\dagger$. Thus, $\mc{D}_\infty^\dagger(\mf{d}_1)$ denotes the Hodge-type
tower attached to $\mf{d}_1$, and $\phi_{\tau,h}^\dagger$ denotes the
Hodge-type or Type~$2$ test function. In Type~$1$, we use the notation
$\phi_{\tau,h}^\dagger$ only after a Type $1$ comparison datum (see Definition \ref{defn:Youcis-Type1-comparison} below) has been fixed;
it then denotes the function obtained from the contracted-product tower in
\cite[(3.48)--(3.49)]{Youcis}.
\end{nota}

\begin{defn}\label{defn:Youcis-Type1-comparison}
Let $\mf{d}=(\mc{G},[b],\mu)$ be a good Type~$1$ datum, and fix a good
associated Hodge-type datum $\mf{d}_1$ realizing Type~$1$. Set $\mc{H}=\mc{G}_1^\der\simeq\mc{G}^\der$. Using Theorem \ref{thm:explicit-deformation} and
\cite[Lemma 3.1.5]{Youcis}, we may identify the generic deformation spaces of $\mf{d}$ and $\mf{d}_1$ through their common adjoint datum: This is just the simple observation that $R_{\mc{G},\mu,j}$ only depends on adjoint data. 

A \emph{Type $1$ comparison datum} consists of a
pro\'etale $\mc{H}(\Z_p)$-torsor $\mc{X}_\infty$ over this generic deformation space, equipped with Weil descent data, and Weil-equivariant isomorphisms of torsors
\begin{equation}\label{eq:Youcis-Type1-comparison}
\mc{X}_\infty\times^{\mc{H}(\Z_p)}\mc{G}_1(\Z_p)
    \isomto \mc{D}_\infty^\dagger(\mf{d}_1),
\qquad
\mc{X}_\infty\times^{\mc{H}(\Z_p)}\mc{G}(\Z_p)
    \isomto \mc{D}_\infty(\mc{G},b,\mu).
\end{equation}
\end{defn}

\begin{prop}\label{prop:abelian-type-test-function-comparison}
Let $\mf{d}=(\mc{G},[b],\mu)$ be of Hodge type, of good Type $2$, or of
good Type $1$ equipped with a Type $1$ comparison datum. Then $\phi_{\tau,h}^{\mc{G},\mu}(b)=\phi_{\tau,h}^\dagger(b)$. 
\end{prop}

\begin{proof}
In the Hodge-type case, choose $\rho\colon\mc{G}\to\GL(\Lambda)$ with $(\rho\circ\mu)$-weights in $\{0,1\}$, and use $\Lambda^\vee$ as the chosen faithful representation in \cite{Youcis}. Then, by Theorem \ref{thm:syntomic-DD}, the $\Lambda$-realization of $\mf{Q}_b^\univ$ is $\mc{M}(H^\univ)$ for the universal $p$-divisible group. Dualizing and applying Remark \ref{rem:BT-tate-module-etale-realization} identifies its $\Lambda^\vee$-realization with $T_p(H^\univ)$, respecting the tensors. Thus the underlying $\mc{G}(\Z_p)$-torsor is the
Tate-module frame torsor used in \cite{Youcis}. The construction and descent of the finite-level covers in \cite[Theorem 3.2.13, Lemma 3.2.15, Proposition 3.2.18, and Corollary 3.2.20]{Youcis} then identify the towers and their Hecke actions. The equality of functions then follows.

Suppose now that $\mf{d}$ is of Type~$1$. The first isomorphism in
\eqref{eq:Youcis-Type1-comparison} supplies the reduction to the derived subgroup required in \cite[(3.46)]{Youcis}. By \cite[(3.48)--(3.49)]{Youcis}, the Type~$1$ tower is then given by $\mc{X}_\infty\times^{\mc{H}(\Z_p)}\mc{G}(\Z_p)$. The second isomorphism in \eqref{eq:Youcis-Type1-comparison} identifies this
tower, equivariantly for all relevant data, with
$\mc{D}_\infty(\mc{G},b,\mu)$. The equality of functions then follows.

Finally, suppose that $\mf{d}$ is of Type~$2$, and put
$\mc{Z}_1=Z(\mc{G}_1)$. As $f([b_1])=[b]$, there is a
$g$ in $\mc{G}(W_j)$ with $b=g f(b_1)\sigma(g)^{-1}$. The connectedness of $\mc{Z}_1$ and \cite[Lemma 3.1.10]{Youcis} furnish us with a lift
$g_1$ in $\mc{G}_1(W_j)$ of $g$. Replacing $b_1$ by
$g_1b_1\sigma(g_1)^{-1}$, we may thus assume that $f(b_1)=b$. The
construction in \cite[(3.52)--(3.54)]{Youcis} gives the Type~$2$ tower as $\mc{D}_\infty^\dagger(\mf{d}_1)/\mc{Z}_1(\Z_p)$. The Hodge-type comparison and the compatibility of $T_\et$ with pushout
along $\mc{G}_1\to\mc{G}_1^\ad\simeq\mc{G}$ identify this quotient with
$\mc{D}_\infty(\mc{G},b,\mu)$, compatibly with Weil descent and Hecke
correspondences, as desired.
\end{proof}

\begin{rem}\label{rem:Youcis-Type1} When a Type 1 comparison datum does not exist, then the construction of the tower in \cite{Youcis} is not correct; one should instead use a Lovering-yoga type construction as in \cite[\S\S3.4.4--3.4.8]{LoveringFCrystals} and \cite{IKY2}, which would then agree with our constructions here.
\end{rem}

\begin{rem}
\label{rem:examples-outside-Youcis-types} Even in their maximal generality, there are many examples of abelian-type Shimura varieties giving rise to $(\mc{G},b,\mu)$  (in the sense of Construction \ref{constr:corr-delta}) which are not good of Type 1 or Type 2. Thus, even in the abelian-type case the local test functions of this paper are substantially more general than those considered in \cite{Youcis}.

We give one simple example of such an abelian-type datum. Let $H$ be a four-dimensional symplectic $\Q$-vector space, and let $V=\bigwedge\nolimits^2_0 H$ be its primitive exterior square. The exterior-square representation
induces a central isogeny $\mr{GSp}_4\to\mr{GSO}(V)$ with kernel $\mu_2$. The image of the Siegel Shimura datum is the orthogonal-similitude Shimura datum of signature $(3,2)$. Geometrically, this is the datum carried by the primitive second cohomology $H^2_{\mr{prim}}(A)= \bigwedge\nolimits^2_0 H^1(A)$ of a principally polarized abelian surface $A$. At every odd hyperspecial prime, this gives a good abelian-type local datum which is not of Type 1 or 2.
\end{rem}

\section{A trace formula for Shimura varieties at bad level}

In this section we establish the primitive version of the trace formula for the cohomology of Shimura varieties at subhyperspecial level using the test functions $\phi_{\tau,h}^{\mc{G},\mu}$. We then explain how one can, assuming an appropriate version of the Langlands--Rapoport conjecture (which is known for abelian-type Shimura varieties), give a more refined version of this trace formula. 

\subsection{Integral canonical models of Shimura varieties}\label{ss:ICMs} 

We first record the setup for the theory of (syntomic) integral canonical models of Shimura varieties as in \cite{IKY2} and \cite{MY}. 

Throughout we use the standard theory of Shimura varieties as in \cite{DeligneModulaire} or \cite{MilneShimura}, and we refer the reader to these references for basic details. 

\begin{setup}\label{setup:sv} Throughout we fix the following unless indicated otherwise:
\begin{itemize}
\item $(\mb{G},\mb{X})$ is a Shimura datum with Hodge cocharacter $\mu_h$ and reflex field $\mb{E}$;
\item $\mb{Z}$ is the center of $\mb{G}$;
\item $\mathsf{K}\subseteq \mb{G}(\bb{A}_f)$ is a neat compact open subgroup which we assume can be factorized as $\mathsf{K}=\mathsf{K}_p\mathsf{K}^p$ where $\mathsf{K}_p\subseteq\mb{G}(\Q_p)$ and $\mathsf{K}^p\subseteq\mb{G}(\bb{A}_f^p)$;
\item $\mr{Sh}_\mathsf{K}(\mb{G},\mb{X})$ is the canonical model of the associated Shimura variety over $\mb{E}$;
\item we fix a prime $p>2$ where $\mb{G}$ is unramified, and a place $v$ of $\mb{E}$ lying over $p$;
\item we set $E\defeq \mb{E}_v$ and write $\mc{O}_E=W=W(k)$ where $k$ is the residue field of $E$;
\item $q=|k|=p^r$;
\item $\Sh_\mathsf{K}$ is used as shorthand for $\Sh_\mathsf{K}(\mb{G},\mb{X})_E$;

\item we write $G=\mb{G}_{\bb{Q}_p}$ and choose a reductive model $\mc{G}$ of $G$ writing then $\mathsf{K}_0\defeq \mc{G}(\bb{Z}_p)$;
\item for $\mathsf{K},\mathsf{K}'\subseteq\mb{G}(\bb{A}_f)$ neat compact open and $g$ in $\mb{G}(\bb{A}_f)$ with $g^{-1}\mathsf{K}g\subseteq \mathsf{K}'$ we denote by 
\begin{equation*}
t_{\mathsf{K},\mathsf{K}'}(g)\colon \Sh_\mathsf{K}(\mb{G},\mb{X})\to \Sh_{\mathsf{K}'}(\mb{G},\mb{X})
\end{equation*}
the unique finite \'etale morphism of these $\mb{E}$-schemes given on $\C$-points by
 \begin{equation*}
 t_{\mathsf{K},\mathsf{K}'}(g)(\mb{G}(\bb{Q})(x,g')\mathsf{K})=\mb{G}(\Q)(x,g'g)\mathsf{K}'
 \end{equation*}
  using the standard double-quotient description of the $\bb{C}$-points of Shimura varieties. We abbreviate $t_{\mathsf{K},\mathsf{K'}}(\mr{id})$ to $\pi_{\mathsf{K},\mathsf{K}'}$;
\item $\ell \ne p$ is a prime.
\end{itemize}
\end{setup}

\subsubsection{\'Etale realization and automorphic \'etale sheaves}\label{sss:etale-realization} We now give setup that will not only be required to formulate the notion of (syntomic) integral canonical model, but which will also be central for the formulation of our desired trace formulae. 

\begin{defn}
[{\cite[Definition 1.5.4 and Lemma 1.5.5]{KSZ}}]
A torus $\mb{T}$ over $\Q$ is \emph{cuspidal} if it is isogenous to a product of a $\Q$-split torus and a $\Q$-torus that is anisotropic over $\R$. For a general $\Q$-torus $\mb{T}$, the \emph{anti-cuspidal part} $\mb{T}_{\mr{ac}}\subset \mb{T}$ is the minimal $\Q$-subtorus with $\mb{T}/\mb{T}_{\mr{ac}}$ cuspidal.
\end{defn}

\begin{rem}
\label{rem:anti-cuspidal}
One may also characterize $\mb{T}_\mr{ac}$ as the smallest subtorus of the maximal anisotropic subtorus $\mb{T}_\mr{a}\subset \mb{T}$ that contains the maximal $\R$-split subtorus of $\mb{T}_{\mr{a},\R}$, i.e., it is the $\Q$-Zariski closure in $\mb{T}$ of the maximal $\R$-split subtorus of $\mb{T}_{\mr{a},\R}$.
\end{rem}

\begin{defn}[{\cite[\S 1.5.6]{KSZ}}]\label{defn:cuspidal-quotient}
We define the \emph{cuspidal quotient} $\mb{G}^c$ of $\mb{G}$ to be $\mb{G}^c = \mb{G}/(\mb{Z}^\circ)_{\mr{ac}}$. 
\end{defn}

\begin{nota}\label{nota:local-cuspidal-quotient} We will consistently write
\begin{itemize}
\item $G^c$ for $\mb{G}^c_{\bb{Q}_p}$;
\item $\mr{pr}^c\colon \mb{G}\to\mb{G}^c$ for the natural map;
\item $\mc{G}^c$ for the corresponding central quotient of $\mc{G}$; see \cite[Proposition 1.1.4]{KisinPappas};
\item and $\mu_h^c$ for the cocharacter of $\mc{G}^c$ induced by $\mu_h$.
\end{itemize}
\end{nota}

\begin{rem} If $(\mb{G},\mb{X})$ is of Hodge type, then $\mb{G}=\mb{G}^c$ (e.g., see \cite[Remark 2.6]{IKY2}).
\end{rem}

\begin{eg}  Let $(\mb{G},\mb{X})=\big(\Res_{F/\Q}\,\mr{GL}_{2,F},(\mf{h}^{\pm})^{[F:\Q]}\big)$ be a Hilbert modular Shimura datum where $F/\Q$ is a totally real extension and $\mf{h}^{\pm}=\C-\R$. Then, $(\mb{Z}^\circ)_\mr{ac}=\Res^1_{F/\Q}\,\bb{G}_{m,F}$: the group of norm-$1$ units in $F^\times$. Thus, $\mb{G}\ne\mb{G}^c$ if $F\ne\Q$.
\end{eg}

We now wish to define the sheaves which form the coefficient systems for our trace formulae.

\begin{defn} Writing $\mathsf{K}=\mathsf{K}_\ell\mathsf{K}^\ell\subseteq\mb{G}(\A_f)$ (where we temporarily allow $\ell=p$), the map
\begin{equation}\label{eq:l-tower}
    \varprojlim_{\scriptscriptstyle\mathsf{K}_\ell'\subseteq \mathsf{K}_\ell}{\Sh}_{\mathsf{K}_\ell'\mathsf{K}^\ell}(\mb{G},\mb{X})\to {\Sh}_{\mathsf{K}}(\mb{G},\mb{X})
\end{equation}
is a $\mathsf{K}_\ell/\mb{Z}(\Q)^{-}_{\mathsf{K}}$-torsor on the pro\'etale site of ${\Sh}_\mathsf{K}(\mb{G},\mb{X})$, where $\mb{Z}(\Q)^{-}_{\mathsf{K}}$ is the closure of $\mb{Z}(\Q)\cap \mathsf{K}$ in $\mathsf{K}$; see \cite[\S1.5.8]{KSZ}. 
The result in loc.\@ cit.\@ shows that $\mathsf{K}_\ell\to \mb{G}^c(\Q_\ell)$ factorizes through $\mathsf{K}_\ell/\mb{Z}(\Q)_{\mathsf{K}}^{-}$. We define $\mb{Et}_{\mathsf{K},\ell}$, the \emph{rational $\ell$-adic \'etale realization \emph{(}at level $\mathsf{K}$\emph{)}}, to be the $\mb{G}^c(\bb{Q}_\ell)$-torsor on the pro\'etale site of $\Sh_\mathsf{K}(\mb{G},\mb{X})$ obtained by the pushout of \eqref{eq:l-tower} along $\mathsf{K}_\ell/\mb{Z}(\bb{Q})_\mathsf{K}^{-}\to \mb{G}^c(\bb{Q}_\ell)$.
\end{defn}

\begin{nota} We write $\mr{Et}_{\mathsf{K},\ell}$ for the pullback of $\mb{Et}_{\mathsf{K},\ell}$ to ${\Sh}_\mathsf{K}$. 
\end{nota} 

\begin{remdef}\label{remdef:automorphic-sheaves}From general Tannakian formalism (see \cite[\S2.1.2]{IKY1}), using $\mb{Et}_{\mathsf{K},\ell}$ we obtain an exact $\bb{Q}_\ell$-linear $\otimes$-functor
\begin{equation*}
\omega_{\mathsf{K},\ell}\colon \cat{Rep}_{\Q_\ell}(\mb{G}^c_{\Q_\ell})\to \cat{Loc}_{\bb{Q}_\ell}(\Sh_\mathsf{K}(\mb{G},\mb{X})).
\end{equation*}
For notational simplicity, for $\xi\colon \mb{G}_{\bb{Q}_\ell}^c\to\GL_{{\Q}_\ell}(V)$ in $\cat{Rep}_{\bb{Q}_\ell}(\mb{G}^c)$ we write
\begin{equation*}
\mc{F}_{\xi,\mathsf{K}}\defeq \omega_{\mathsf{K},\ell}(\xi)=\mb{Et}_{\mathsf{K},\ell}\wedge^{\mb{G}^c(\bb{Q}_\ell)}V,
\end{equation*}
a $\Q_\ell$-local system on $\Sh_\mathsf{K}(\mb{G},\mb{X})$. We call a sheaf of this form an \emph{automorphic local system}.
\end{remdef}

\begin{rem} Tracing through the definitions shows that if $g=g_\ell g^\ell$ in $\mb{G}(\bb{A}_f)$, and $\Int(g_\ell^c)$ is the inner automorphism of $\mb{G}^c(\bb{Q}_\ell)$ associated to the image $g_\ell^c$ of $g_\ell$ in $\mb{G}^c(\Q_\ell)$, then
\begin{equation}\label{eq:Hecke-action-local-system-compat-rat}
  t_{\mathsf{K},\mathsf{K}^{'}}(g)^\ast(\omega_{\mathsf{K}',\ell}(\xi))=\omega_{{\mathsf{K}},\ell}(\xi\circ \mathrm{Int}((g_\ell^c)^{-1})),
\end{equation}
for any neat compact open levels $\mathsf{K}$ and $\mathsf{K}'$ with $g^{-1}\mathsf{K}g\subseteq \mathsf{K}'$.
\end{rem}

While automorphic local systems will be the coefficients for our trace formulae, for the definition of \emph{integral canonical model} of $\Sh_\mathsf{K}$ we need to construct certain $\mc{G}^c(\bb{Z}_p)$-lattices $\mb{Et}^\circ_{\mathsf{K},p}\subseteq \mb{Et}_{\mathsf{K},p}$. 

\begin{defn} Assume $\mathsf{K}=\mathsf{K}_0\mathsf{K}^p$. By \cite[\S1.5.8]{KSZ}, the map $\mathsf{K}_0\to\mc{G}^c(\Z_p)$ factorizes through $\mathsf{K}_0/\mb{Z}(\Q)^{-}_{\mathsf{K}}$. The \emph{integral $p$-adic \'etale realization \emph{(}of level $\mathsf{K}$\emph{)}}, denoted $\mb{Et}_{\mathsf{K},p}^\circ$, is the pro\'etale $\mc{G}^c(\bb{Z}_p)$-torsor on $\Sh_\mathsf{K}(\mb{G},\mb{X})$ given by the pushforward of \eqref{eq:l-tower} along this factorization.
\end{defn}

\begin{nota}We write $\mr{Et}^\circ_{\mathsf{K},p}$ for the pullback of $\mb{Et}^\circ_{\mathsf{K},p}$ to ${\Sh}_\mathsf{K}$. 
\end{nota}

\begin{rem} Again using Tannakian formalism (see \cite[\S2.1.2]{IKY1}), we may alternatively think of $\mr{Et}_{\mathsf{K},p}^\circ$ as an exact $\bb{Z}_p$-linear $\otimes$-functor 
\begin{equation*}
\omega_{\mathsf{K},p}^\circ\colon \cat{Rep}_{\bb{Z}_p}(\mc{G}^c)\to \cat{Loc}_{\Z_p}(\Sh_\mathsf{K}).
\end{equation*}Again tracing through the definitions shows that if $g=g_p g^p$ is in $\mathsf{K}_0\mb{G}(\bb{A}_f^p)$, then
\begin{equation}\label{eq:Hecke-action-local-system-compat-int}
  t_{\mathsf{K},\mathsf{K}^{'}}(g)^\ast(\omega^\circ_{\mathsf{K}',p}(\xi)))=\omega_{{\mathsf{K}},p}^\circ(\xi\circ \mathrm{Int}((g_p^c)^{-1})),
\end{equation}
for any neat compact open levels $\mathsf{K}=\mathsf{K}_0\mathsf{K}^p$ and $\mathsf{K}'=\mathsf{K}_0\mathsf{K}^{'p}$ with $g^{-1}\mathsf{K}g\subseteq \mathsf{K}'$.
\end{rem}

\subsubsection{Integral canonical models} We now discuss the so-called integral canonical models of Shimura varieties with $p$-level given by $\mathsf{K}_0$. While our eventual trace formula will involve deeper level than $\mathsf{K}_0$, we will only require these integral models at level $\mathsf{K}_0$.

\begin{defn} Let $x\colon \Spa(F)\to \Sh_\mathsf{K}^\an$ be a classical point. We say that $\mr{Et}^\circ_{\mathsf{K},p}$ is \emph{\emph{(}potentially\emph{)} crystalline at $x$} if for all representations $\Lambda$ of $\mc{G}^c$ one has that the $\Q_p$-representation of $\Gal(\ov{F}/F)$ given by $\omega_{\mathsf{K},p}^\circ(\Lambda)_x[\nicefrac{1}{p}]$ is (potentially) crystalline.
\end{defn}

\begin{rem}\label{rem:one-rep-enough} By \cite[Proposition 2.21]{IKY1} to check that $\mr{Et}^\circ_{\mathsf{K},p}$ is (potentially) crystalline at $x$ it suffices to check (potential) crystallinity of $\omega_{\mathsf{K},p}^\circ(\Lambda)_x[\nicefrac{1}{p}]$ for a single faithful representation $\Lambda$.
\end{rem}

\begin{defn}[{Integral canonical models, \cite[Definition 3.39]{IKY2} and \cite[Definition 6.5.2]{MY}}]\label{defn:ICM} A separated $W$-model $\ms{S}_\mathsf{K}=\ms{S}_{\mathsf{K}}(\mb{G},\mb{X})$ of $\Sh_\mathsf{K}$ is an \emph{integral canonical model} if:
\begin{enumerate}
\item The $\mc{G}^c(\bb{Z}_p)$-local system $\mr{Et}_{\mathsf{K},p}^\circ$ is apertile of type $-\mu_h^c$ and the associated map of $W$-stacks $\ms{S}_\mathsf{K}\to\mr{BT}^{\mc{G}^c,-\mu_h^c,\mr{alg}}_\infty$ is $p$-adically formally \'etale.
\item One has an equality of classical points
\begin{equation*}
|(\wh{\ms{S}}_{\mathsf{K}})_\eta|^\mr{cl}=\left\{x\in |\Sh_\mathsf{K}^\mr{an}|^\mr{cl}:\mb{Et}_{\mathsf{K},p}^\circ\text{ is potentially crystalline at }x\right\}.
\end{equation*}
\end{enumerate}
An integral canonical model $\ms{S}_\mathsf{K}$ is furthermore called \emph{limpid} if $\mr{Et}_{\mathsf{K},\ell}$ extends to $\ms{S}_{\mathsf{K}}$ for all $\ell\ne p$.
\end{defn}

It will be helpful to name the $(\mc{G}^c,-\mu_h)$-aperture arising on the completion $\wh{\ms{S}}_\mathsf{K}$ of an integral canonical model from the map $\ms{S}_\mathsf{K}\to\mr{BT}^{\mc{G}^c,-\mu_h^c,\mr{alg}}_\infty$

\begin{defn}\label{defn:syntomic-realization} Suppose that $\ms{S}_\mathsf{K}$ is an integral canonical model of $\Sh_\mathsf{K}$. We then call the induced object $\mf{Q}_\mathsf{K}$ of $\mr{BT}^{\mc{G}^c,-\mu_h^c}_\infty(\wh{\ms{S}}_{\mathsf{K}})$ the \emph{syntomic realization \emph{(}at level $\mathsf{K}$\emph{)}}. 
\end{defn}

\begin{rem} By Proposition \ref{prop:Tannakian-equiv}, we may think of $\mf{Q}_\mathsf{K}$ as an exact $\bb{Z}_p$-linear $\otimes$-functor
\begin{equation*}
\omega_{\mathsf{K},p}^\mr{syn}\colon \cat{Rep}_{\bb{Z}_p}(\mc{G}^c)\to \cat{Vect}(\wh{\ms{S}}_\mathsf{K}{}^\mr{syn}).
\end{equation*}
\end{rem}

Now, it is not clear from the definition that integral canonical models, if they exist, are unique in any way. The following mapping property verifies this is the case as it implies that integral canonical models are in fact functorial in $(\mb{G},\mb{X},\mc{G})$.

\begin{prop}[{\cite[Corollary 6.5.7]{MY}}]\label{prop:mapping-prop} Suppose that $\mc{Y}$ is an excellent $\eta$-normal $\mc{O}_E$-scheme with generic fiber $Y$. Then, the following are equivalent for a map $f\colon Y\to \Sh_\mathsf{K}$:
\begin{enumerate}
\item $f$ extends to a map $\mc{Y}\to\ms{S}_\mathsf{K}$;
\item the composition 
\begin{equation*}
Y\to \Sh_\mathsf{K}\xrightarrow{\mr{Et}_{\mathsf{K},p}^\circ} B\mc{G}^c(\bb{Z}_p)\otimes E=\mr{BT}^{\mc{G}^c,-\mu_h^c,\mr{alg}}_\infty \otimes E
\end{equation*}
extends to a map $\mc{Y}\to\mr{BT}^{\mc{G}^c,-\mu_h^c,\mr{alg}}_\infty$.
\end{enumerate} 
\end{prop}

Next, it is not clear that integral canonical models, at least in the way we have defined them here, exist; it's not even entirely obvious that the usual Siegel-type integral models of Mumford are canonical. That said, the following shows that limpid integral canonical models exist for essentially all Shimura varieties.

\begin{thm}\label{thm:ICMs-exist} A limpid integral canonical model of $\Sh_\mathsf{K}$ exists if:
\begin{enumerate}
\item \emph{(}\cite{IKY2}, \cite{MY}\emph{)} $(\mb{G},\mb{X})$ is of abelian type.
\item \emph{(}\cite{MY}\emph{)} $(\mb{G},\mb{X})$ is of pre-abelian type and $p>2$.
\item \emph{(}\cite{MY}\emph{)} $(\mb{G},\mb{X})$ is arbitrary and $p\gg 0$ \emph{(}depending only on $(\mb{G},\mb{X})$\emph{)}.
\end{enumerate}
\end{thm}

\begin{rem}\label{rem:ICMs-exist} Let us slightly expound upon the above results:

\begin{itemize} 
\item When $(\mb{G},\mb{X})$ is of abelian type and $p>2$, it is shown in \cite{IKY2} that the integral models constructed by Kisin in \cite{KisIntShab} are limpid integral canonical models. 

\item When $(\mb{G},\mb{X})$ is of pre-abelian type and $p>2$, it is shown in \cite{MY} how to construct limpid integral canonical models; necessarily these agree with Kisin's models in the abelian-type case by the previous paragraph and Proposition \ref{prop:mapping-prop}. 
\item In \cite{MY} limpid integral canonical models are shown to exist with no restriction on $(\mb{G},\mb{X})$ for $p\gg 0$ (depending only on $(\mb{G}^\ad,\mb{X}^\ad)$). It is further shown that, up to enlarging $p$, these agree with the integral models of Bakker--Shankar--Tsimerman from \cite{BST}.
\end{itemize}
\end{rem}

While limpid integral canonical models are defined at fixed level they propagate to a projective system of integral canonical models at smaller level, a fact leveraged in our trace formulae.

\begin{nota} Let $\Omega_{\mathsf{K}}$ denote the set of levels $\mathsf{K}'\subseteq \mathsf{K}$ of the form $\mathsf{K}_0\mathsf{K}^{'p}$. 
\end{nota}

\begin{prop}[Existence of prime-to-$p$ Hecke towers, see {\cite[Corollary 6.5.11]{MY}}]\label{prop:prime-to-p-Hecke} Suppose that $\Sh_\mathsf{K}$ admits a limpid integral canonical model $\ms{S}_\mathsf{K}$. Then,
\begin{enumerate}
\item If $\mathsf{K}'$ belongs to $\Omega_\mathsf{K}$, then there exists a limpid integral canonical model $\ms{S}_{\mathsf{K}'}$. 
\item If $\mathsf{K}'$ and $\mathsf{K}''$ belong to
$\Omega_\mathsf{K}$ and $g=g_pg^p$ in $\mathsf{K}_0\mb{G}(\bb{A}_f^p)$ satisfies
$g^{-1}\mathsf{K}'g\subseteq\mathsf{K}''$, then $t_{\mathsf{K}',\mathsf{K}''}(g^p)$ extends uniquely to a finite \'etale morphism $\ms{S}_{\mathsf{K}'}\to\ms{S}_{\mathsf{K}''}$ modeling $t_{\mathsf{K}',\mathsf{K}''}(g)$.
\end{enumerate}
\end{prop}

\begin{rem}\label{rem:ICM-propagation} With further assumptions on $(\mb{G},\mb{X})$ and $\mc{G}$, one obtains the same results as in Proposition \ref{prop:prime-to-p-Hecke} but with only the condition that $\mathsf{K}'=\mathsf{K}_0\mathsf{K}^{'p}$; see \cite[Theorem 7.6.25]{MY}. 
\end{rem}

\begin{rem}\label{rem:etale-realization-Hecke} Suppose that $\Sh_\mathsf{K}$ admits a limpid integral canonical model. Using Proposition \ref{prop:prime-to-p-Hecke} we see from \eqref{eq:Hecke-action-local-system-compat-int}, and \cite[Theorem 4.1.5]{MY}, that if $g=g_p g^p$ in $\mathsf{K}_0\mb{G}(\bb{A}_f^p)$, then
\begin{equation}\label{eq:Hecke-action-local-system-compat-syn}
  t_{\mathsf{K}',\mathsf{K}^{''}}(g^p)^\ast(\omega^\mr{syn}_{\mathsf{K}'',p}(\xi)))=\omega_{{\mathsf{K}'},p}^\syn(\xi),
\end{equation}
for any neat compact open levels $\mathsf{K}'$ and $\mathsf{K}''$ in $\Omega_{\mathsf{K}}$ with $(g^p)^{-1}\mathsf{K}'g^p\subseteq \mathsf{K}''$.
\end{rem}

\subsubsection{Alternative proof of (parts of) Theorem \ref{thm:ICMs-exist}}\label{sss:alternative-proofs}

In the final part of this subsection we lay out a proof of cases (1) and (3) of Theorem \ref{thm:ICMs-exist} when $p>2$, different from those in \cite{IKY2} and \cite{MY}. Instead we use Theorem \ref{thm:BT-over-base-ring} to bootstrap from the $F$-crystals with $\mc{G}$-structure constructed in \cite{LoveringFCrystals} and \cite{BST}. The major advantage of this is that one avoids the usage of any major results from \cite{IKY1} and, consequently, of \cite[Theorem A]{GuoReinecke}.

\begin{rem}For the sake of space, we will be quite brief in our explanations below. But, as the reader will see, many arguments are quite standard.
\end{rem}

\begin{assumption} Throughout the following we assume that $p>2$.
\end{assumption}

\paragraph*{Integral canonical models in the abelian-type case} We begin by showing that the results of \cite{LoveringFCrystals} upgrade to prove that Kisin's integral models are canonical.

\begin{setup} Suppose that $(\mb{G},\mb{X})$ is an abelian-type Shimura datum, and let $\ms{S}_{\mathsf{K}}$ be the integral models constructed by Kisin in \cite{KisIntShab}.
\end{setup}

\begin{constr} Using \cite[\S\S3.4.4--3.4.8]{LoveringFCrystals} and \cite[Proposition 4.7.11]{LoveringModels}, one constructs an object of $\cat{Quint}_{\mc{G}^c}^{-\mu_h^c}(\wh{\ms{S}}_\mathsf{K})$ as in Definition \ref{defn:quint}.\footnote{In \cite{LoveringModels} it is assumed that $\mb{Z}^\circ$ splits over a CM extension, but this is not necessary with the definition of $\mb{G}^c$ given here, and his arguments go through essentially verbatim.} Applying Theorem \ref{thm:BT-over-base-ring} we obtain a $(\mc{G}^c,-\mu_h^c)$-aperture on $\wh{\ms{S}}_{\mathsf{K}}$ which we denote $\mf{Q}_\mathsf{K}^\heartsuit$.
\end{constr}

\begin{lem}\label{lem:realization-compat} There is an identification of $T_\et(\mf{Q}_\mathsf{K}^\heartsuit)$ and the restriction of $\mr{Et}^\circ_{\mathsf{K},p}$ to $(\wh{\ms{S}}_\mathsf{K})_\eta$.
\end{lem}
\begin{proof} This essentially follows from the proof of \cite[Theorem 3.5.1]{LoveringFCrystals}. But, for the sake of the reader we comment that in the Hodge type case this follows from the discussion in \cite[\S2.6]{IKY2}, and one reduces from the abelian-type case to the Hodge-type case by carefully tracing through the construction. 
\end{proof}

\begin{constr} Using Lemma \ref{lem:realization-compat} we may promote the map $\wh{\ms{S}}_{\mathsf{K}}\to\mr{BT}^{\mc{G}^c,-\mu_h^c}_\infty$ corresponding to $\mf{Q}_\mathsf{K}^\heartsuit$  to a map $\rho\colon \ms{S}_{\mathsf{K}}\to \mr{BT}^{\mc{G}^c,-\mu_h^c,\mr{alg}}_\infty$. 
\end{constr}

\begin{prop}\label{prop:rho-is-formally-etale} The map $\rho$ is $p$-adically formally \'etale. 
\end{prop}

\begin{proof}In the Hodge-type case this follows from \cite[Corollary 3.3.13]{LoveringFCrystals} using \cite[Proposition 3.4]{IKY2}. One then reduces to the Hodge-type case via standard Hodge-type-to-abelian-type yoga; cf.\@ \cite[Theorem 3.6]{IKY2}.
\end{proof}

\begin{thm} The integral model $\ms{S}_\mathsf{K}$ of $\Sh_\mathsf{K}$ is canonical.
\end{thm}
\begin{proof} Via the construction of the map $\rho$ above, and Proposition \ref{prop:rho-is-formally-etale}, all that remains to be shown is that Condition (2) of Definition \ref{defn:ICM} holds. But, this is \cite[Proposition 2.15]{IKY2}.
\end{proof}

\paragraph*{Integral canonical models for sufficiently large primes} We now show that the integral models in \cite[\S5.1]{BST} are canonical. 

\begin{setup} Suppose that $(\mb{G},\mb{X})$ is an arbitrary Shimura datum, and let $\ms{S}_{\mathsf{K}}$ be the integral models constructed in \cite{BST}.
\end{setup}

\begin{constr} From \cite[\S7.1]{BST} it can be deduced that one has an object of $\cat{Quint}_{\mc{G}^c}^{-\mu_h^c}(\wh{\ms{S}}_\mathsf{K})$ as in Definition \ref{defn:quint} save one thing: loc.\@ cit.\@ does not verify the strong divisibility condition outside of the one fixed faithful representation. But, this is automatic by Lemma \ref{lem:one-for-all}. Applying Theorem \ref{thm:BT-over-base-ring} we obtain a $(\mc{G}^c,-\mu_h^c)$-aperture on $\wh{\ms{S}}_{\mathsf{K}}$ which we denote $\mf{Q}_\mathsf{K}^\mr{BST}$.
\end{constr}

\begin{lem}\label{lem:realization-compat-BST} There is an identification of $T_\et(\mf{Q}_\mathsf{K}^\mr{BST})$ and the restriction of $\mr{Et}^\circ_{\mathsf{K},p}$ to $(\wh{\ms{S}}_\mathsf{K})_\eta$.
\end{lem}
\begin{proof} This follows via the construction of the tensors $\{s_\mr{crys}\}$ in \cite[\S7.1.1]{BST} which show that \'etale realization yields the $\bb{Z}_p$-local system with tensors corresponding to $\mr{Et}_{\mathsf{K},p}^\circ$.  
\end{proof}

\begin{constr} Using Lemma \ref{lem:realization-compat-BST} we may promote the map $\wh{\ms{S}}_{\mathsf{K}}\to\mr{BT}^{\mc{G}^c,-\mu_h^c}_\infty$ corresponding to $\mf{Q}_\mathsf{K}^\mr{BST}$  to a map $\rho\colon \ms{S}_{\mathsf{K}}\to \mr{BT}^{\mc{G}^c,-\mu_h^c,\mr{alg}}_\infty$. 
\end{constr}

\begin{prop}\label{prop:rho-is-formally-etale-BST} The map $\rho$ is $p$-adically formally \'etale. 
\end{prop}

\begin{proof} By \cite[Lemma 6.4.7]{MY} it suffices to show that the Kodaira--Spencer morphism from Remark \ref{rem:Kodaira-Spencer} is an isomorphism. But, this follows from \cite[Theorem 7.9]{BST} and its proof.
\end{proof}

\begin{thm} The integral model $\ms{S}_\mathsf{K}$ of $\Sh_\mathsf{K}$ is canonical.
\end{thm}
\begin{proof} Via the construction of the map $\rho$ above, and Proposition \ref{prop:rho-is-formally-etale-BST}, all that remains is to observe is that Condition (2) of Definition \ref{defn:ICM} holds by the proof in \cite[\S7.8]{MY}.
\end{proof}

\subsection{The primitive trace formula} We now state and prove the primitive trace formula for the Galois--Hecke action on the cohomology of Shimura varieties with subhyperspecial-at-$p$ level. This formula will be in terms of a weighted point count on the fixed points of the mod-$p$ Hecke correspondence, with the weighting given by the local test functions $\phi_{\tau,h}$ from the last section.

\begin{setup}\label{setup:prim-trace-formula} We continue to freely use notation from Setup \ref{setup:sv}. We further let:
\begin{itemize}
\item $\ov{\mb{E}}$ and $\ov{E}$ be algebraic closures of $\mb{E}$ and $E$, respectively;
\item  $\ov{\mb{E}}\hookrightarrow \ov{E}$ be an embedding corresponding to a place $\ov{v}$ of $\ov{\mb{E}}$ extending $v$;
\item $C$ be the $p$-adic completion of $\ov{E}$;
\item $\Gamma_{\mb{E}}$ and $\Gamma_E$ be the absolute Galois groups of $\mb{E}$ and $E$, respectively;
\item $\ms{S}_\mathsf{K}$ be a limpid integral canonical model of $\Sh_\mathsf{K}$.
\end{itemize}
Our choice of embedding gives an isomorphism $\Gamma_E\isomto D_{\ov{v}}\subseteq \Gamma_{\mb{E}}$ where $D_{\ov{v}}$ is the decomposition group of $\ov{v}$. We may identify $\Gamma_E$ with the group of continuous automorphisms of $C$ over $E$.

We will also make frequent use of notation/terminology concerning (cohomological) correspondences and their actions on cohomology; we refer the reader to \S\ref{s:prelims} for such things.
\end{setup}

\subsubsection{Hecke correspondences and cohomology} We now give the basic setup for the Galois--Hecke action on the cohomology of Shimura varieties with automorphic \'etale sheaf coefficients, as well as a rigid-analytic variant that will be useful to clarify the flow of the results.

\begin{defn}\label{defn:pot-crys-loci} Suppose that $\mathsf{L}\subseteq\mb{G}(\bb{A}_f)$ is a neat compact open. A \emph{potentially crystalline locus at level $\mathsf{L}$}, if it exists, is a quasi-compact open subset $\mc{U}_\mathsf{L}$  of $\Sh_\mathsf{L}^\mr{an}$ such that
\begin{equation*}
|\mc{U}_\mathsf{L}|^\mr{cl}=\left\{x\in |\Sh_\mathsf{L}^\mr{an}|^\mr{cl}:\mb{Et}_{\mathsf{L},p}\text{ is potentially crystalline at }x\right\}.
\end{equation*}
If a potentially crystalline locus exists for all neat compact open subgroups $\mathsf{L}$, then we say that $(\mb{G},\mb{X})$ \emph{admits potentially crystalline loci}.
\end{defn}

\begin{rem} By \cite[Corollary 4.3]{HuberCV}, a potentially crystalline locus at level $\mathsf{L}$ is unique.
\end{rem}

The following shows that potentially crystalline loci are preserved by the Hecke action, and that $(\mb{G},\mb{X})$ admitting potentially crystalline loci follows from the existence at one level.

\begin{prop}\label{prop:pot-crys-loci-exist} If $(\mb{G},\mb{X})$ admits a potentially crystalline locus at one level, then it admits potentially crystalline loci. Moreover, if $(\mb{G},\mb{X})$ admits potentially crystalline loci, then the morphisms $t_{\mathsf{L},\mathsf{L}'}(g)\colon \Sh_\mathsf{L}\to \mr{Sh}_{\mathsf{L}'}$ pull back to morphisms $t_{\mathsf{L},\mathsf{L}'}(g)\colon \mc{U}_\mathsf{L}\to \mc{U}_{\mathsf{L}'}$.
\end{prop}
\begin{proof} We begin with an elementary observation. Let $f\colon X\to Y$ be a finite \'etale surjection of rigid $E$-spaces, and $\bb{V}$ a $\Q_p$-local system on $Y$. As $(f^\ast\bb{V})_x=\bb{V}_{f(x)}$, if $V\subseteq X$ is a potentially crystalline locus of $f^\ast\bb{V}$ then $f(V)\subseteq Y$ is a potentially crystalline locus of $\bb{V}$, and if $U\subseteq Y$ is a potentially crystalline locus for $\bb{V}$ then $f^{-1}(U)$ is a potentially crystalline locus for $f^\ast\bb{V}$.

So, suppose $\mathsf{M}$ is a neat level and $\mc{U}_\mathsf{M}\subseteq \Sh_\mathsf{M}^\an$ a potentially crystalline locus. Then, for any other neat level $\mathsf{L}$, take a neat level $\mathsf{N}\subseteq \mathsf{M}\cap \mathsf{L}$. Fix a faithful representation $\xi\colon \mb{G}^c_{\Q_p}\hookrightarrow \GL_{\Q_p}(V)$ and abbreviate $(\mc{F}_{\xi,\mathsf{P}}|_{\Sh_{\mathsf{P}}})^\an=\bb{V}_\mathsf{P}$ for any neat level $\mathsf{P}$. Then, we have that 
\begin{equation*}
\pi_{\mathsf{N},\mathsf{M}}^\ast \bb{V}_{\mathsf{M}}\simeq \bb{V}_{\mathsf{N}}\simeq \pi_{\mathsf{N},\mathsf{L}}^\ast\bb{V}_\mathsf{L},
\end{equation*}
by \eqref{eq:Hecke-action-local-system-compat-rat}. Thus, from the previous paragraph, we deduce that $\mc{U}_\mathsf{L}\defeq \pi_{\mathsf{N},\mathsf{L}}(\pi_{\mathsf{N},\mathsf{M}}^{-1}(\mc{U}_\mathsf{M}))\subseteq \Sh_\mathsf{L}^\mr{an}$ is a potentially crystalline locus for $\bb{V}_\mathsf{L}$, and thus for $\mr{Et}_{\mathsf{L},p}$ by Remark \ref{rem:one-rep-enough}. 
\end{proof}

\begin{rem} If  $\ms{S}_\mathsf{K}$ is an integral canonical model, then $(\wh{\ms{S}}_{\mathsf{K}})_\eta$ is a potentially crystalline locus at level $\mathsf{K}$. Thus, by Proposition \ref{prop:pot-crys-loci-exist}, $(\mb{G},\mb{X})$ admits potentially crystalline loci. 
\end{rem}

For visual clarity, we give shorthand notation for the cohomology spaces that we study below.

\begin{nota}\label{nota:Shimura-cohomology} Suppose that $(\mb{G},\mb{X})$ admits potentially crystalline loci. For a neat compact open subgroup $\mathsf{K}\subseteq\mb{G}(\bb{A}_f)$, a prime $\ell\ne p$, and a $\ov{\Q}_\ell$-representation $\xi$ of $\mb{G}_{\Q_\ell}^c$, we define 
\begin{equation*}
H_c^i(\mb{G},\mb{X},\mathsf{K},\xi)\defeq
H^i_{\et,c}(\Sh_\mathsf{K}(\mb{G},\mb{X})_{\ov{\mb{E}}},\mc{F}_\xi),\qquad \ovc{H}_c^i(\mb{G},\mb{X},\mathsf{K},\xi)\defeq H_{\et,c}^i((\mc{U}_\mathsf{K})_C,\mc{F}_{\xi,\mathsf{K}}^\mr{an}).
\end{equation*}
We then define
\begin{equation*}
H^i_c(\mb{G},\mb{X},\xi)\defeq \varinjlim_{\mathsf{K}}H_c^i(\mb{G},\mb{X},\mathsf{K},\xi),\qquad \ovc{H}^i_c(\mb{G},\mb{X},\xi)\defeq \varinjlim_{\mathsf{K}}\ovc{H}_c^i(\mb{G},\mb{X},\mathsf{K},\xi).
\end{equation*}
Note that $H^i_c(\mb{G},\mb{X},\xi)$ and $\ovc{H}^i_c(\mb{G},\mb{X},\xi)$ inherit an action of $\mb{G}(\bb{A}_f)$ by modifying the level structure, and this action commutes with the natural actions of $\Gamma_{\mb{E}}$ and $\Gamma_E$, respectively.
\end{nota}

\begin{rem} Using \cite[Corollary 7.1.4]{BerkovichEtale}, one has a natural $\Gamma_E$-equivariant identification 
\begin{equation*}
H^i_{\et,c}\left(\Sh_\mathsf{K}(\mb{G},\mb{X})_{\ov{E}},\mc{F}_{\xi,\mathsf{K}}\right)\simeq H^i_{\et,c}\left(\Sh_\mathsf{K}(\mb{G},\mb{X})_C^\mr{an},\mc{F}_{\xi,\mathsf{K}}^\an\right).
\end{equation*}
Via this isomorphism and extension by zero, we obtain natural comparison maps 
\begin{equation}\label{eq:comparison-maps}
\ovc{H}_c^i(\mb{G},\mb{X},\mathsf{K},\xi)\to H_c^i(\mb{G},\mb{X},\mathsf{K},\xi), \qquad \ovc{H}_c^i(\mb{G},\mb{X},\xi) \to H_c^i(\mb{G},\mb{X},\xi),
\end{equation}
which are linear over $\ov{\Q}_\ell[\Gamma_{{E}}]\otimes_{\ov{\Q}_\ell} \mc{H}_{\ov{\Q}_\ell}(\mb{G}(\bb{A}_f),\mathsf{K})$ and $\ov{\Q}_\ell[\Gamma_{{E}}]\otimes_{\ov{\Q}_\ell} \mc{H}_{\ov{\Q}_\ell}(\mb{G}(\bb{A}_f))$, respectively.
\end{rem}

\begin{obs} For any neat level $\mathsf{K}$, we have that $H^i_c(\mb{G},\mb{X},\mathsf{K},\xi)$ and $\ovc{H}^i_c(\mb{G},\mb{X},\mathsf{K},\xi)$ are finite-dimensional and vanish for $i>2\dim_\C(X)$. As $\dim_\C(X)=\dim\Sh_\mathsf{K}(\mb{G},\mb{X})$, the former case is classical, and in the latter case it follows from \cite[Proposition 2.3]{ScholzeLK}.

From the following equalities
\begin{equation*}
H^i_c(\mb{G},\mb{X},\xi)^\mathsf{K}=H^i_c(\mb{G},\mb{X},\mathsf{K},\xi),\qquad \ovc{H}^i_c(\mb{G},\mb{X},\xi)^\mathsf{K}=\ovc{H}^i_c(\mb{G},\mb{X},\mathsf{K},\xi),
\end{equation*}
we deduce that $H^i_c(\mb{G},\mb{X},\xi)$ and $\ovc{H}^i_c(\mb{G},\mb{X},\xi)$ are admissible modules over $\ov{\bb{Q}}_\ell[\mb{G}(\bb{A}_f)]$.
\end{obs}

\begin{nota} For a level $\mathsf{K}$, set 
\begin{equation*}
 H_c^\ast(\mb{G},\mb{X},\mathsf{K},\xi) =\sum_{i=0}^{2\dim_\C(\mb{X})} (-1)^i H^i_c(\mb{G},\mb{X},\mathsf{K},\xi), \qquad H_c^\ast(\mb{G},\mb{X},\xi)=\sum_{i=0}^{2\dim_\C(\mb{X})} (-1)^i H^i_c(\mb{G},\mb{X},\xi),
\end{equation*}
considered as objects of the Grothendieck group of modules over $\ov{\Q}_\ell[\Gamma_\mb{E}]\otimes_{\ov{\Q}_\ell}\mc{H}_{\ov{\Q}_\ell}(\mb{G}(\A_f),\mathsf{K})$ and $\ov{\Q}_\ell[\Gamma_\mb{E}]\otimes_{\ov{\Q}_\ell}\mc{H}_{\ov{\Q}_\ell}(\mb{G}(\A_f))$, respectively. We make the identical definitions for the cohomology spaces $\ovc{H}^i_c(-)$, but the output now only has an action of $\Gamma_E$. 
\end{nota}

We now wish to describe a more geometric understanding of the morphism $\tau\times e(g^{-1},\mathsf{K})$ on these cohomology spaces.

\begin{defn}\label{defn:geometric-Hecke-corr}
Let $\mathsf{K}$ be a neat level, and $g$ an element of $\mb{G}(\bb{A}_f)$. Set $\mathsf{K}^g\defeq \mathsf{K}\cap g\mathsf{K}g^{-1}$. Then, the \emph{Hecke correspondence} is given as follows:
\begin{equation}\label{eq:Hecke-corr}
 c(g,\mathsf{K})=(\pi_{\mathsf{K}^g,\mathsf{K}},t_{\mathsf{K}^g,\mathsf{K}}(g))\colon \Sh_{\mathsf{K}^g}(\mb{G},\mb{X})\rightrightarrows \Sh_\mathsf{K}(\mb{G},\mb{X}).
\end{equation}
The Hecke correspondence restricts to an analytic correspondence $c(g,\mathsf{K})$ on
\begin{equation}\label{eq:Hecke-corr-gr}
 c(g,\mathsf{K})=(\pi_{\mathsf{K}^g,\mathsf{K}},t_{\mathsf{K}^g,\mathsf{K}}(g))\colon \mc{U}_{\mathsf{K}^g}\rightrightarrows\mc{U}_\mathsf{K}.
\end{equation}
\end{defn}

We may extend this to a cohomological correspondence on automorphic \'etale sheaves as follows.

\begin{defn}\label{defn:cohomological-Hecke-correspondence} Let $c(g,\mathsf{K})$ be as in \eqref{eq:Hecke-corr} and \eqref{eq:Hecke-corr-gr}, respectively. Then, for any automorphic \'etale sheaf $\mc{F}_{\xi,\mathsf{K}}$ we define the \emph{cohomological Hecke correspondence} to be the element $u(g,\mathsf{K})$ of $\mr{Coh}_{c(g,\mathsf{K})}(\mc{F}_{\xi,\mathsf{K}},\mc{F}_{\xi,\mathsf{K}})$ and $\mr{Coh}_{c(g,\mathsf{K})}(\mc{F}_{\xi,\mathsf{K}}^\an,\mc{F}_{\xi,\mathsf{K}}^\an)$, respectively, defined as the morphism
\begin{equation*}
    u(g,\mathsf{K})\colon (c(g,\mathsf{K})_2)_!\, c(g,\mathsf{K})_1^\ast\,\mc{F}_{\xi,\mathsf{K}}\to\mc{F}_{\xi,\mathsf{K}},
\qquad  u(g,\mathsf{K})\colon (c(g,\mathsf{K})_2)_!\, c(g,\mathsf{K})_1^\ast\,\mc{F}_{\xi,\mathsf{K}}^\an\to\mc{F}_{\xi,\mathsf{K}}^\an,
\end{equation*}
respectively, each coming via adjunction from \eqref{eq:Hecke-action-local-system-compat-rat}, and the fact that $c(g,\mathsf{K})$ is finite \'etale.
\end{defn}

We then have the following workhorse lemma, which follows from Observation \ref{obs:Hecke-action-desc} and the natural inverse that's introduced when we pass to cohomology.

\begin{lem}\label{lem:hecke-corr-coh} For $\tau$ an element of $\Gamma_\mb{E}$, one has
\begin{equation*}
\tr\left(\tau\times e(g^{-1},\mathsf{K})\mid H^i_c(\mb{G},\mb{X},\xi)\right)=\tr\bigg(\tau\times u(g,\mathsf{K})\mid H^i_{\et,c}(\Sh_\mathsf{K}(\mb{G},\mb{X})_{\ov{\mb{E}}},{\mc{F}}_{\xi,\mathsf{K}})\bigg).
\end{equation*}
For $\tau$ an element of $\Gamma_E$, one has
\begin{equation*}
\tr\left(\tau\times e(g^{-1},\mathsf{K})\mid \ovc{H}^i_c(\mb{G},\mb{X},\xi)\right)=\tr\bigg(\tau\times u(g,\mathsf{K})\mid H^i_{\et,c}((\mc{U}_\mathsf{K})_{C},{\mc{F}}_{\xi,\mathsf{K}}^\an)\bigg).
\end{equation*}
\end{lem}

We end by observing that these Hecke correspondences can be lifted integrally if the $p$-part of $g$ is contained in our fixed hyperspecial.

\begin{defn}\label{defn:integral-Hecke-corr} For $g=g_0g^p$ in $\mathsf{K}_0\mb{G}(\A_f^p)$, observe that $c(g,\mathsf{K})=c(g^p,\mathsf{K})$ as $g_0$ is in $\mathsf{K}_0$. Thus, by Proposition \ref{prop:prime-to-p-Hecke} we obtain \emph{integral geometric Hecke correspondences}
\begin{equation*}
c(g^p,\mathsf{K})=(\pi_{\mathsf{K}^g,\mathsf{K}},t_{\mathsf{K}^g,\mathsf{K}}(g^p))\colon\ms{S}_{\mathsf{K}^g}\rightrightarrows \ms{S}_\mathsf{K},\quad  c(g^p,\mathsf{K})=(\pi_{\mathsf{K}^g,\mathsf{K}},t_{\mathsf{K}^g,\mathsf{K}}(g^p))\colon
\wh{\ms{S}}_{\mathsf{K}^g}\rightrightarrows \wh{\ms{S}}_\mathsf{K},
\end{equation*}
modeling \eqref{eq:Hecke-corr} and \eqref{eq:Hecke-corr-gr}, respectively. Moreover, it's clear that the cohomological Hecke correspondence $u(g^p,\mathsf{K})$ lifts integrally relative to these integral geometric Hecke correspondences.

Both of these correspondences give rise to the same correspondence on the special fiber
\begin{equation*}
c(g^p,\mathsf{K})=(\pi_{\mathsf{K}^g,\mathsf{K}},t_{\mathsf{K}^g,\mathsf{K}}(g^p))\colon\ms{S}_{\mathsf{K}^g,k}\rightrightarrows \ms{S}_{\mathsf{K},k},
\end{equation*}
which we denote $c(g^p,\mathsf{K})$, and write $u(g^p,\mathsf{K})$ for the induced cohomological correspondence.
\end{defn}

\subsubsection{A comparison of tubes}\label{ss:comparison-of-tubes} We now wish to show that the aperture tubes with level structure defined in Definition \ref{defn:aperture-tubes} agree with Berthelot tubes inside of Shimura varieties coming from Hecke-fixed points in the special fiber of integral canonical models. 

To do this, we must first explain how to precisely associate to such a Hecke-fixed point $y$ an element $b=\delta(y)$ of $C_j(\mc{G},\sigma(\mu))$.

\begin{constr}\label{constr:delta-functor} Using the equivalence from Proposition \ref{prop:isom-classes-in-BT}, we define a map
\begin{equation*}
\delta^\mr{stk}=\delta^\mr{stk}_\mathsf{K}\colon \ms{S}_{\mathsf{K}}(\ov{k})\to \mc{C}_\infty(\mc{G}^c,\sigma(-\mu_h^c))
\end{equation*}
associating to a point $x$ in the source the object $(\mf{Q}_\mathsf{K})_x$ of the target.
\end{constr}

\begin{rem}\label{rem:delta-compatabilities} The functor $\delta^\mr{stk}_{\mathsf{K}}$ from Construction \ref{constr:delta-functor} is Frobenius-equivariant: if $\Phi$ denotes the $q$-Frobenius on $\ms{S}_{\mathsf{K}}(\ov{k})$, then there is a natural identification
\begin{equation}\label{eq:delta-Frob-equiv}
\delta^\mr{stk}\left(\Phi(x)\right)=(\sigma^r)^\ast(\delta^\mr{stk}(x)).
\end{equation}
Moreover, if $\mathsf{K}'$ and $\mathsf{K}''$ are in $\Omega_\mathsf{K}$, and $g=g_pg^p$ in $\mathsf{K}_0\mb{G}(\bb{A}_f^p)$ is such that $g^{-1}\mathsf{K}'g\subseteq \mathsf{K}''$ then it follows from \eqref{eq:Hecke-action-local-system-compat-syn} that there is a natural identification
\begin{equation}\label{eq:delta-Hecke-compatability}
 \delta^\mr{stk}_{\mathsf{K}''}\circ t_{\mathsf{K}',\mathsf{K}''}(g^p)\simeq \delta^\mr{stk}_{\mathsf{K}'}.
 \end{equation}
\end{rem}

\begin{constr}\label{constr:corr-delta} Suppose that $g=g_pg^p$ belongs to $\mathsf{K}_0\mb{G}(\bb{A}_f^p)$. Then, for an integer $j\geqslant 1$ we may define a functor 
\begin{equation*}
\delta^\mr{stk}\colon \mr{Fix}\left(c(g^p,\mathsf{K})^{(j)}\right)(\ov{k})\to \mc{C}_\infty\left(\mc{G}^c,\sigma(-\mu_h^c)\right),
\end{equation*}
by taking $y$ to $\delta^\mr{stk}_{\mathsf{K}}(x)$ where $x=c(g^p,\mathsf{K})_2(y)=\Phi^j(c(g^p,\mathsf{K})_1(y))$. By this second equality, as well as \eqref{eq:delta-Frob-equiv} and \eqref{eq:delta-Hecke-compatability}, we deduce a natural isomorphism
\begin{equation*}
(\sigma^{rj})^\ast\delta^\mr{stk}(y)\simeq \delta^\mr{stk}(y).
\end{equation*}
Thus, $\delta^\mr{stk}$ admits a natural lifting
\begin{equation*}
\delta^\mr{stk}\colon \mr{Fix}\left(c(g^p,\mathsf{K})^{(j)}\right)(\ov{k})\to \mc{C}_\infty\left(\mc{G}^c,\sigma(-\mu_h^c)\right)^{h\sigma^{rj}},
\end{equation*}
where the target is the category of homotopy $\sigma^{rj}$-fixed points. Passing to the sets of isomorphism classes we deduce the existence of a natural map
\begin{equation*}
\delta\colon \mr{Fix}\left(c(g^p,\mathsf{K})^{(j)}\right)(\ov{k})\to C_j\left(\mc{G}^c,\sigma(-\mu_h^c)\right).
\end{equation*}
\end{constr}

To state our desired comparison of tubes we first give the following setup.

\begin{setup}\label{setup:tube-comparison} Fix 

\begin{itemize}
\item a neat level $\mathsf{K}'=\mathsf{K}_p\mathsf{K}^p\subseteq \mathsf{K}_0\mathsf{K}^p$;
\item $g=g_pg^p$ in $\mathsf{K}_0\mb{G}(\A_f^p)$; 
\item $j\geqslant 1$ an integer;
\item $y$ a point of $\mr{Fix}(c(g,\mathsf{K})^{(j)})(\ov{k})$;
\item $x=c(g^p,\mathsf{K})_2(y)=\Phi^j(c(g^p,\mathsf{K})_1(y))$;
\item $\tau$ an element of $W_E$ with $\mathsf{v}(\tau)=j$.
\end{itemize}
\end{setup}

\begin{constr}As in Remark \ref{tau-twisted-generic-corr} we get an analytic correspondence
\begin{equation*}
    c(g,\mathsf{K})^\tau(y)\colon \pi_{\mathsf{K}^{'g},\mathsf{K}^g}^{-1}(\wh{\ms{S}}_{\mathsf{K}^g}(y)_C)\rightrightarrows \pi_{\mathsf{K}',\mathsf{K}}^{-1}(\wh{\ms{S}}_\mathsf{K}(x)_C),
\end{equation*}
where $\wh{\ms{S}}_{\mathsf{K}^g}(y)$ and $\wh{\ms{S}}_{\mathsf{K}}(x)$ denote the Berthelot tubes as in Notation \ref{nota:rigid-analytic-general}. 
\end{constr}

To state our desired comparison of tubes, we need one more assumption about the level $\mathsf{K}$.

\begin{assumption}\label{ass:pullback-from-cuspidal} In the notation of Setup \ref{setup:tube-comparison}, assume that $\mathsf{K}_p=(\mr{pr}^c)^{-1}(\mathsf{L}_p)$ for a compact open subgroup $\mathsf{L}_p\subseteq \mc{G}^c(\Z_p)$.
\end{assumption}

With this assumption we can now state the correspondence we wish to compare $c(g,\mathsf{K})^\tau$ to.

\begin{constr} As $\delta(y)$ is in $C_j(\mc{G}^c,\sigma(-\mu_h^c))$ we may consider the deformation space with level structure $\mc{D}_{\mathsf{L}_p}(\mc{G}^c,\delta(y),-\mu_h^c)$. We then define the analytic correspondence
\begin{equation*}
\varsigma(g_p^c,\mathsf{L}_p)^\tau(\delta(y))=(\tau\circ \pi_{\mathsf{L}_p^{g_p^c},\mathsf{L}_p},t_{\mathsf{L}_p^{g_p^c},\mathsf{L}_p}(g_p^c))\colon \mc{D}_{\mathsf{L}_p^{g_p^c}}(\mc{G}^c,\delta(y),-\mu_h^c)_C\rightrightarrows \mc{D}_{\mathsf{L}_p}(\mc{G}^c,\delta(y),-\mu_h^c)_C
\end{equation*}
\end{constr}

Given the above setup, our cohomological comparison of tubes comes in the following form.

\begin{prop}\label{prop:tube-comp} There is an identification of correspondences $c(g,\mathsf{K})^\tau(y)\simeq \varsigma(g_p^c,\mathsf{L}_p)^\tau(\delta(y))$.
\end{prop}
\begin{proof}For the sake of visual clarity, we make the following notational definitions:
\begin{equation*}
 X_y\defeq\wh{\ms{S}}_{\mathsf{K}^g}(y)_C,
 \qquad
 X_x\defeq\wh{\ms{S}}_{\mathsf{K}}(x)_C,
 \qquad
 z\defeq\pi_{\mathsf{K}^g,\mathsf{K}}(y),
\end{equation*}
so $x=t_{\mathsf{K}^g,\mathsf{K}}(g^p)(y)=\Phi^j(z)$. The map $t_{\mathsf{K}^g,\mathsf{K}}(g^p)$ is finite \'etale and hence restricts to an isomorphism $\vartheta_y\colon X_y\isomto X_x$. Set $a_y\defeq (\tau\circ\pi_{\mathsf{K}^g,\mathsf{K}})\circ\vartheta_y^{-1}\colon X_x\to X_x$.

We further make the following notational simplifications:
\begin{equation*}
 \mf{Q}_x\defeq(\mf{Q}_{\mathsf{K}})_x,
 \qquad
 \mf{Q}_y\defeq(\mf{Q}_{\mathsf{K}^g})_y,
 \qquad
 \mf{Q}_z\defeq(\mf{Q}_{\mathsf{K}})_z.
\end{equation*}
Then, \eqref{eq:Hecke-action-local-system-compat-syn}, together with the fact that $\mf{Q}_{\mathsf{K}}$ is defined over $W$, gives isomorphisms
\begin{equation*}
 \epsilon_t\colon\mf{Q}_x\isomto\mf{Q}_y,
 \qquad
 \epsilon_\pi\colon\mf{Q}_z\isomto\mf{Q}_y,
 \qquad
 \epsilon_\Phi\colon
 \mf{Q}_x\isomto(\sigma^{rj})^*\mf{Q}_z.
\end{equation*}
Then, in fact, the descent datum used to define $\delta(y)$ in Construction \ref{constr:corr-delta} is
\begin{equation}\label{eq:tube-comp-descent-datum}
 \alpha_y\defeq
 \epsilon_\Phi^{-1}\circ
 (\sigma^{rj})^*(\epsilon_\pi^{-1})\circ
 (\sigma^{rj})^*(\epsilon_t)
 \colon
 (\sigma^{rj})^*\mf{Q}_x\isomto\mf{Q}_x.
\end{equation}
As $\ms{S}_{\mathsf{K}}\to\mr{BT}^{\mc{G}^c,-\mu_h^c,\mr{alg}}_\infty$ is $p$-adically formally \'etale, it identifies the deformation spaces of $x$ and $\mf{Q}_x$. Passing to generic fibers therefore gives an isomorphism $\lambda_y\colon X_x\isomto \mc{D}(\mc{G}^c,\delta(y),-\mu_h^c)_C$. We claim that $\lambda_y\circ a_y=\tau\circ\lambda_y$. To verify this, let $(\mf{Q},\iota)$ be a deformation of $\mf{Q}_x$. Then $a_y$ sends it to $\big(\tau^*\mf{Q},\alpha_y\circ\tau^*\iota\big)$, but this is precisely the action of $\tau$ on the deformation functor for $\delta(y)$. 

Now as $a_y\circ\vartheta_y= \tau\circ\pi_{\mathsf{K}^g,\mathsf{K}}$ and $\vartheta_y= t_{\mathsf{K}^g,\mathsf{K}}(g^p)|_{X_y}$ (the latter by definition), we see that the equality $\lambda_y\circ a_y=\tau\circ\lambda_y$ gives the commutative diagram
\begin{equation}\label{eq:tube-comp-base-comparison}
\begin{tikzcd}[column sep=8em,row sep=3em]
 X_y
 \ar[r,shift left=.8ex,
     "\tau\circ\pi_{\mathsf{K}^g,\mathsf{K}}"]
 \ar[r,shift right=.8ex,
     "t_{\mathsf{K}^g,\mathsf{K}}(g^p)"']
 \ar[d,"\lambda_y\circ\vartheta_y"']
 &
 X_x
 \ar[d,"\lambda_y"]
 \\
 \mc{D}(\mc{G}^c,\delta(y),-\mu_h^c)_C
 \ar[r,shift left=.8ex,"\tau"]
 \ar[r,shift right=.8ex,"\mr{id}"']
 &
 \mc{D}(\mc{G}^c,\delta(y),-\mu_h^c)_C.
\end{tikzcd}
\end{equation}
This is the base-level version of our result, and so it remains to add level structure.

Applying \'etale realization to the isomorphism of universal apertures inducing
$\lambda_y$, and using
\eqref{eq:Hecke-action-local-system-compat-syn} along $\vartheta_y$,
gives $\mc{G}^c(\Z_p)$-equivariant isomorphisms
\begin{equation*}
 \mr{Et}_{\mathsf{K},p}^\circ|_{X_x}
 \isomto \mc{D}_\infty(\mc{G}^c,\delta(y),-\mu_h^c)_C,
 \qquad
 \mr{Et}_{\mathsf{K}^g,p}^\circ|_{X_y}\isomto
 \mc{D}_\infty(\mc{G}^c,\delta(y),-\mu_h^c)_C,
\end{equation*}
lying over $\lambda_y$ and $\lambda_y\circ\vartheta_y$, respectively. Taking the quotients by $\mathsf{L}_p$ and $\mathsf{L}_p^{g_p^c}$ gives isomorphisms
\begin{align*}
 \lambda_{y,\mathsf{L}_p^{g_p^c}}\colon
 \pi_{\mathsf{K}^{'g},\mathsf{K}^g}^{-1}(X_y)
 &\isomto
 \mc{D}_{\mathsf{L}_p^{g_p^c}}
 (\mc{G}^c,\delta(y),-\mu_h^c)_C,
 \\
 \lambda_{x,\mathsf{L}_p}\colon
 \pi_{\mathsf{K}',\mathsf{K}}^{-1}(X_x)
 &\isomto
 \mc{D}_{\mathsf{L}_p}
 (\mc{G}^c,\delta(y),-\mu_h^c)_C.
\end{align*}
These isomorphisms fit into the commutative diagram
\begin{equation}\label{eq:tube-comp-level-comparison}
\begin{tikzcd}[column sep=9em,row sep=4em]
 \pi_{\mathsf{K}^{'g},\mathsf{K}^g}^{-1}(X_y)
 \ar[r,shift left=.9ex,
     "{c(g,\mathsf{K})^\tau(y)_1}"]
 \ar[r,shift right=.9ex,
     "{c(g,\mathsf{K})^\tau(y)_2}"']
 \ar[d,"\lambda_{y,\mathsf{L}_p^{g_p^c}}"']
 &
 \pi_{\mathsf{K}',\mathsf{K}}^{-1}(X_x)
 \ar[d,"\lambda_{x,\mathsf{L}_p}"]
 \\
 \mc{D}_{\mathsf{L}_p^{g_p^c}}
 (\mc{G}^c,\delta(y),-\mu_h^c)_C
 \ar[r,shift left=.9ex,
     "\tau\circ
      \pi_{\mathsf{L}_p^{g_p^c},\mathsf{L}_p}"]
 \ar[r,shift right=.9ex,
     "t_{\mathsf{L}_p^{g_p^c},\mathsf{L}_p}(g_p^c)"']
 &
 \mc{D}_{\mathsf{L}_p}
 (\mc{G}^c,\delta(y),-\mu_h^c)_C.
\end{tikzcd}
\end{equation}
More explicitly, commutativity means
\begin{align*}
 \lambda_{x,\mathsf{L}_p}
 \circ c(g,\mathsf{K})^\tau(y)_1
 &=
 \big(
 \tau\circ
 \pi_{\mathsf{L}_p^{g_p^c},\mathsf{L}_p}
 \big)
 \circ\lambda_{y,\mathsf{L}_p^{g_p^c}},
 \\
 \lambda_{x,\mathsf{L}_p}
 \circ c(g,\mathsf{K})^\tau(y)_2
 &=
 t_{\mathsf{L}_p^{g_p^c},\mathsf{L}_p}(g_p^c)
 \circ\lambda_{y,\mathsf{L}_p^{g_p^c}}.
\end{align*}
The first equality is obtained from
\eqref{eq:tube-comp-base-comparison} by passing to the indicated
level quotients. For the second, \eqref{eq:Hecke-action-local-system-compat-int} says that
the $p$-component $g_p$ acts on a frame by right multiplication by
$g_p^c$; the prime-to-$p$ component is the map $\vartheta_y$. Right
multiplication by $g_p^c$ descends from the
$\mathsf{L}_p^{g_p^c}$-quotient to the $\mathsf{L}_p$-quotient
because $(g_p^c)^{-1} \mathsf{L}_p^{g_p^c}g_p^c\subseteq \mathsf{L}_p$. 

As the lower row of \eqref{eq:tube-comp-level-comparison} is
$\varsigma(g_p^c,\mathsf{L}_p)^\tau(\delta(y))$ by definition, the commutativity of \eqref{eq:tube-comp-level-comparison} means that we are done.
\end{proof}

\subsubsection{The primitive trace formula} We are now ready to state and prove the primitive form of our trace formula.

\begin{setup}\label{setup:prim-trace-form} We maintain the setup from Setup \ref{setup:tube-comparison} and Assumption \ref{ass:pullback-from-cuspidal}. We then specialize the setup from Setup \ref{setup:bad-level-corr} to our current situation as follows:
\begin{itemize}
\item $\mc{O}_K=W$;
\item $\mf{X}=\wh{\ms{S}}_{\mathsf{K}}$;
\item $c= c(g^p,\mathsf{K})$;
\item $u=u(g^p,\mathsf{K})$;
\item $d=c(g,\mathsf{K}')$;
\item $f\colon d\to c_\eta$ is the morphism $(\pi_{\mathsf{K}',\mathsf{K}};\pi_{\mathsf{K}^{'g},\mathsf{K}^g},\pi_{\mathsf{K}',\mathsf{K}})$;
\item $v=u(g,\mathsf{K}')$ which is naturally isomorphic to $f^\ast(u_\eta)$.
\end{itemize}
\end{setup}

\begin{prop} The following equality holds true:
\begin{equation*}
\tr\left(R\Gamma(\tr_{c(g,\mathsf{K})^\tau(y)})\right)=\tr\bigg(\tau\times e((g_p^c)^{-1},\mathsf{L}_p)|H^\ast(\mc{G}^c,\delta(y),-\mu_h^c)\bigg)=\phi_{\tau,e((g_p^c)^{-1},\mathsf{L}_p)}^{\mc{G}^c,-\mu_h^c}(\delta(y)).
\end{equation*}
\end{prop}
\begin{proof} From Proposition \ref{prop:tube-comp} we deduce that 
\begin{equation*}
\tr\left(R\Gamma(\tr_{c(g,\mathsf{K})^\tau(y)})\right)=\tr\left(R\Gamma\tr_{\varsigma(g_p^c,\mathsf{L}_p)^\tau)(\delta(y))}\right).
\end{equation*}
But, by Observation \ref{obs:Hecke-action-desc} one deduces that 
\begin{equation*}
\tr\left(R\Gamma\tr_{\varsigma(g_p^c,\mathsf{L}_p)^\tau)(\delta(y))}\right)=\tr\bigg(\tau\times e((g_p^c)^{-1},\mathsf{L}_p)|H^\ast(\mc{G}^c,\delta(y),-\mu_h^c)\bigg)=\phi_{\tau,e((g_p^c)^{-1},\mathsf{L}_p)}^{\mc{G}^c,-\mu_h^c}(\delta(y)),
\end{equation*}
from where the conclusion follows.
\end{proof}

From this, we deduce the first form of our primitive trace formula.

\begin{thm}\label{thm:prim-trace-formula-I} Let $g=g_pg^p$ and $\mathsf{K}'$ be as in \emph{Setup \ref{setup:tube-comparison}} and \emph{Assumption \ref{ass:pullback-from-cuspidal}}. There exists a constant $j_0\geqslant 0$, depending only on $g^p$, such that for all $\tau$ in $W_E$ with $j:=\mathsf{v}(\tau)>j_0$ one has
 \begin{equation*}
\tr\left(\tau\times e(g^{-1},\mathsf{K}')\mid \ovc{H}^\ast_c(\mb{G},\mb{X},\xi)\right)=\sum_{y\in\mr{Fix}(c^{(j)})(\ov{k})}\tr\left((u^{(j)}_{c_2({y})})^\dashv\mid ({\mc{F}}_{\xi,\mathsf{K}})_{c_2({y})}\right) \phi_{\tau,e((g_p^c)^{-1},\mathsf{L}_p)}^{\mc{G}^c,-\mu_h^c}(\delta(y)),
\end{equation*}
where for visual simplicity we have written $u=u(g^p,\mathsf{K})$ and $c=c(g^p,\mathsf{K})$. Moreover, if $\ms{S}_{\mathsf{K}}$ is proper, then one can take the constant $j_0$ to be $0$.
\end{thm}
\begin{proof} This follows immediately from Proposition \ref{prop:specific-rigid-analytic-trace-formula} given Proposition \ref{prop:tube-comp}.
\end{proof}

Of course, we are more interested in understanding $H^\ast_{c}(\mb{G},\mb{X},\xi)$ instead of $\ovc{H}^\ast_c(\mb{G},\mb{X},\xi)$. Thus, we would like to isolate conditions on $(\mb{G},\mb{X})$ so that the analogue of Theorem \ref{thm:prim-trace-formula-I} applies. 

\begin{defn} If the map $\ovc{H}^i_c(\mb{G},\mb{X},\xi)\to H^i_c(\mb{G},\mb{X},\xi)$ from \eqref{eq:comparison-maps} is an isomorphism for all $i$ and representations $\xi$ of $\mb{G}^c_{\Q_\ell}$, we say that $(\mb{G},\mb{X})$ has \emph{$\ell$-adically contractible boundary}.
\end{defn}

\begin{rem}\label{rem:l-adically-contractible-boundary} To show that $(\mb{G},\mb{X})$ has $\ell$-adically contractible boundary, it suffices to show that the comparison map $\ovc{H}^i_c(\mb{G},\mb{X},\mathsf{K}',\xi)\to{H}^i_c(\mb{G},\mb{X},\mathsf{K}',\xi)$ from \eqref{eq:comparison-maps} is an isomorphism for all $i$ and $\xi$ and for $\mathsf{K}$ sufficiently small; so, assume that $\mathsf{K}'\subseteq \mathsf{K}$. Moreover, up to changing $\mathsf{K}^p$, we can assume that $\mathsf{K}'=\mathsf{K}_p\mathsf{K}^p$ with $\mathsf{K}_p\subseteq \mathsf{K}_0$. Then, by acyclicity of finite maps, we must only show that the map
\begin{equation*}
H^i_{\et,c}((\mc{U}_\mathsf{K})_C,\pi_\ast\mc{F}_{\xi,\mathsf{K}'})\to H^i_{\et,c}(\Sh_\mathsf{K}(\mb{G},\mb{X})_C,\pi_\ast\mc{F}_{\xi,\mathsf{K}'}),
\end{equation*}
is an isomorphism where $\pi=\pi_{\mathsf{K}',\mathsf{K}}$. 

That said, using Proposition \ref{prop:nearby-omnibus} we have a commutative diagram
\begin{equation*}
\begin{tikzcd}
	{H^i_{\et,c}(\Sh_\mathsf{K}(\mb{G},\mb{X})_C,\pi_\ast\mc{F}_{\xi,\mathsf{K}'}) } & {H^i_{\et,c}((\mc{U}_\mathsf{K})_C,\pi_\ast\mc{F}_{\xi,\mathsf{K}'})} \\
	& {H^i_{\et,c}(\mathscr{S}_{\mathsf{K},\overline{k}},R\Psi\pi_\ast\mc{F}_{\xi,\mathsf{K}'}).}
	\arrow[from=1-2, to=1-1]
	\arrow[from=2-2, to=1-1]
	\arrow["\wr", from=2-2, to=1-2]
\end{tikzcd}
\end{equation*}
Thus, $\ell$-adically contractible boundary is equivalent to the maps
\begin{equation*}
 H^i_{\et,c}\bigl(\mathscr{S}_{\mathsf{K},\ov{k}},
                 R\Psi\pi_\ast\mc{F}_{\xi,\mathsf{K}'}\bigr)
 \to
 H^i_{\et,c}\big(\Sh_\mathsf{K}(\mb{G},\mb{X})_{C},
                 \pi_\ast\mc{F}_{\xi,\mathsf{K}'}\big)
\end{equation*}
being isomorphisms for every $i$ and $\xi$ and all sufficiently
small neat levels $\mathsf{K}'$. This condition is expected to hold in general, but is currently known in the following cases:
\begin{enumerate}
\item if $\mathscr{S}_\mathsf{K}$ is proper (by the proper base change theorem);
\item if $(\mb{G},\mb{X})$ is of abelian type by \cite[Proposition 5.21]{Wu}, extending results from \cite{LanStrohII}.
\end{enumerate}
\end{rem}

Thus, we are now prepared to state the more desirable form of Theorem \ref{thm:prim-trace-formula-I}, which follows immediately from this result and the definition of $\ell$-adically contractible boundary.

\begin{thm}\label{thm:prim-trace-formula-II} Let $g=g_pg^p$ and $\mathsf{K}'$ be as in \emph{Setup \ref{setup:tube-comparison}} and \emph{Assumption \ref{ass:pullback-from-cuspidal}}. Furthermore assume that $(\mb{G},\mb{X})$ has $\ell$-adically contractible boundary. Then, there exists a constant $j_0\geqslant 0$, depending only on $g^p$, such that for all $\tau$ in $W_E$ with $j:=\mathsf{v}(\tau)>j_0$ one has
 \begin{equation*}
\tr\left(\tau\times e(g^{-1},\mathsf{K}')\mid {H}^\ast_c(\mb{G},\mb{X},\xi)\right)=\sum_{y\in\mr{Fix}(c^{(j)})(\ov{k})}\tr\left((u^{(j)}_{c_2({y})})^\dashv\mid ({\mc{F}}_{\xi,\mathsf{K}})_{c_2({y})}\right) \phi_{\tau,e((g_p^c)^{-1},\mathsf{L}_p)}^{\mc{G}^c,-\mu_h^c}(\delta(y)),
\end{equation*}
where for visual simplicity we have written $u=u(g^p,\mathsf{K})$ and $c=c(g^p,\mathsf{K})$. Moreover, if $\ms{S}_{\mathsf{K}}$ is proper, then one can take the constant $j_0$ to be $0$.
\end{thm}

Finally, we enunciate the unconditionality of Theorem \ref{thm:prim-trace-formula-II} in the abelian-type case. 

\begin{thm}\label{thm:prim-trace-formula-III} Suppose that $(\mb{G},\mb{X})$ is of abelian type. Let $g=g_pg^p$ and $\mathsf{K}'$ be as in \emph{Setup \ref{setup:tube-comparison}} and \emph{Assumption \ref{ass:pullback-from-cuspidal}}. Then, there exists a constant $j_0\geqslant 0$, depending only on $g^p$, such that for all $\tau$ in $W_E$ with $j:=\mathsf{v}(\tau)>j_0$ one has
 \begin{equation*}
\tr\left(\tau\times e(g^{-1},\mathsf{K}')\mid {H}^\ast_c(\mb{G},\mb{X},\xi)\right)=\sum_{y\in\mr{Fix}(c^{(j)})(\ov{k})}\tr\left((u^{(j)}_{c_2({y})})^\dashv\mid ({\mc{F}}_{\xi,\mathsf{K}})_{c_2({y})}\right) \phi_{\tau,e((g_p^c)^{-1},\mathsf{L}_p)}^{\mc{G}^c,-\mu_h^c}(\delta(y)),
\end{equation*}
where for visual simplicity we have written $u=u(g^p,\mathsf{K})$ and $c=c(g^p,\mathsf{K})$. Moreover, if $\mb{G}^\mr{ad}$ is $\Q$-anisotropic, then one can take the constant $j_0$ to be $0$.
\end{thm}
\begin{proof} This follows from Theorem \ref{thm:prim-trace-formula-II} given Remark \ref{rem:l-adically-contractible-boundary} and that $\ms{S}_\mathsf{K}$ is proper if and only if $\mb{G}^\mr{ad}$ is $\Q$-anisotropic; see \cite[Corollary 4.51]{Wu} and \cite[Corollary 4.1.7]{MadTorHod}.
\end{proof}

\begin{rem} Despite its apparent purely-technical nature, Assumption \ref{ass:pullback-from-cuspidal} is indicative of a substantive issue. Ultimately the problem is a fundamental disconnect between the levels of the Shimura variety, which are indexed by subgroups of $\mb{G}(\A_f)$, and the \'etale realization functors on the Shimura variety which are $\mb{G}^c(\A_f)$-local systems. This indicates that there should be a better replacement for $\Sh_\mathsf{K}(\mb{G},\mb{X})$ when $\mb{G}\ne\mb{G}^c$ which actually carries $\mb{G}$-realization functors as opposed to just $\mb{G}^c$-realization functors. It is not clear what this would look like, and perhaps such an object only exists in an exotic category like those appearing in the work of Clausen--Scholze.
\end{rem}

\subsection{A Langlands--Rapoport-type refinement of the primitive trace formula}
\label{ss:KSZ-point-counting} 

From the perspective of applications, the primitive trace formulae above are not sufficient. Instead, one wishes to have a more `motivic' description of these traces. This is furnished by the so-called \emph{Langlands--Rapoport conjecture}.

In this section, we describe such a conjectural upgrade, and then explain how the work of \cite{KSZ}, together with the syntomic canonical models of \cite{IKY2}, makes this unconditional in the abelian-type case.

\begin{setup}\label{setup:LR-general} We maintain the setup from Setup \ref{setup:prim-trace-formula}, but also assume that $\mathsf{K}=\mathsf{K}_0\mathsf{K}^p$.
\end{setup}

\subsubsection{The Langlands--Rapoport-$\tau$ conjecture}

To formulate and prove this upgraded version of the primitive trace formula, we must first recall the Langlands--Rapoport-$\tau$ conjecture of \cite[\S2]{KSZ}, which gives the `motivic' description of $\ms{S}_{\mathsf{K}_0}(\ov{k})$ where $\ms{S}_{\mathsf{K}_0}=\varprojlim_{\mathsf{K}=\mathsf{K}_0\mathsf{K}^p}\ms{S}_\mathsf{K}$. 

This requires an immense amount of technical setup. Below we recall the barest version of this setup necessary to state our results; essentially an indexed tour of the notation from loc.\@ cit.\@ We strongly encourage the reader to consult loc.\@ cit.\@ to properly engage with this material.

We begin by discussing the objects which form the indexing set for our decomposition; these can be viewed as the `motivic' analogue of isogeny classes.

\begin{nota}\label{nota:KSZ-admissible-morphisms} We write
\begin{itemize}
\item $\mf{Q}$ for the corrected quasi-motivic gerb; see \cite[\S2.2.8]{KSZ};
\item $\mf{G}_{\mb{G}}$ for the neutral $\ov{\Q}/\Q$-Galois gerb associated with $\mb{G}$; see \cite[Definition 2.1.4]{KSZ}. 
\end{itemize}
We regard both as objects of the category $\mr{pro-}\mc{G}\mr{rb}(\ov{\Q}/\Q)$ of pro-$\ov{\Q}/\Q$ Galois gerbs defined in \cite[Definition 2.1.11]{KSZ}, viewing $\mf{G}_{\mb{G}}$ as a constant pro-object.
\end{nota}

\begin{defn}\label{defn:KSZ-admissible-morphism}
A morphism $\varphi\colon\mf{Q}\to\mf{G}_{\mb{G}}$ in $\mr{pro-}\mc{G}\mr{rb}(\ov{\Q}/\Q)$ is \emph{admissible} if it satisfies the conditions of \cite[Definition 2.4.2]{KSZ}. 
\end{defn}

\begin{defn}\label{nota:KSZ-conjugacy-automorphisms}
Two morphisms $\varphi_1,\varphi_2\colon\mf{Q}\to\mf{G}_{\mb{G}}$ are \emph{conjugate} if $\varphi_2=\Int(g)\circ\varphi_1$ for some $g$ in $\mb{G}(\ov{\Q})$, as in \cite[Definition 2.1.13]{KSZ}. 
\end{defn}

\begin{rem} By \cite[Remark 2.4.4]{KSZ}, admissibility is invariant under conjugation.
\end{rem}

\begin{nota} For an admissible morphism $\varphi$, we denote by $I_\varphi$ the reductive $\Q$-group attached to $\varphi$ as in \cite[\S2.1.14]{KSZ}. 
\end{nota}

We next recall the local comparison sets associated to such a $\varphi$ as above; these may be viewed as parameterizing the actual `motivic' objects contained in a given isogeny class $\varphi$; informally, one might think of these objects as `lattices' in $\varphi$. This amounts to combining a parameterization of the away-from-$p$ lattices and the at-$p$ lattices separately. 
 
\begin{nota}\label{nota:KSZ-local-comparison-sets}
For $w$ either $\infty$ or a prime other than $p$, let $\xi_w$ be the reference local morphism of \cite[\S2.4.1]{KSZ}. Then with $\varphi(w)$ and $\zeta_w$ as in \cite[\S2.2.8]{KSZ}, we define
\begin{equation*}
X_w(\varphi)
\defeq
\left\{
 g\in\mb{G}(\ov{\Q}_w):
 \Int(g)\circ\xi_w=\varphi(w)\circ\zeta_w
\right\}.
\end{equation*}
Thus $X_w(\varphi)$ being non-empty says precisely that the localization of $\varphi$ is conjugate to the prescribed local morphism; see \cite[\S2.4.1]{KSZ}.
\end{nota}

\begin{defn}\label{defn:KSZ-local-LR-set}
With $\theta_\varphi\defeq\varphi(p)\circ\zeta_p$, we define the \emph{$p$-local Langlands--Rapoport set} to be the set $X_p(\varphi)\defeq X_{-\mu_h}(\theta_\varphi)$ as in \cite[\S2.2.7 and \S2.4.1]{KSZ}; it is independent of the choice of Hodge cocharacter in its conjugacy class.
\end{defn}

\begin{rem} As in \cite[\S2.2.7 and \S2.4.1]{KSZ}, the set $X_p(\varphi)$ carries a $q$-Frobenius automorphism which we will not need to explicitly notate.
\end{rem}

\begin{rem} We can give a more explicit description of $X_p(\varphi)$. Namely, choose $g_0$ in the set $\mr{UR}(\theta_\varphi)$ as in \cite[Definition 2.2.3 and Lemma 2.2.4(i)]{KSZ}, and let $b_\varphi$ be the element of $G(\Q_p^{\mr{ur}})\subseteq G(\breve{\Q}_p)$ attached in \cite[Definition 2.2.5]{KSZ} to the unramified morphism $\Int(g_0^{-1})\circ\theta_\varphi$. Then, there is an identification
\begin{equation}\label{eq:KSZ-local-LR-set}
X_p(\varphi)\isomto\left\{g\mc{G}(W(\ov{k})):g^{-1}b_\varphi\sigma(g)\in\mc{G}(W(\ov{k}))\sigma(-\mu_h)(p)\mc{G}(W(\ov{k}))
\right\}.\footnote{ We use the equivalent $\sigma(-\mu_h)$-convention obtained by sending the coset of $u$ to the coset of $b_\varphi\sigma(u)$. }
\end{equation}
\end{rem}

We now combine these away-from-$p$ and at-$p$ sets into a set parameterizing the $\A_f^p$-lattices within a given isogeny class.

\begin{defn}\label{defn:KSZ-prime-to-p-LR-set}
For an admissible $\varphi$, the sets $X_w(\varphi)$ contain integral elements for almost all finite $w\ne p$ by \cite[\S2.4.6]{KSZ}. We then define the \emph{prime-to-$p$ Langlands--Rapoport set} to be
\begin{equation*}
X^p(\varphi)
\defeq
\prod_{w\ne p,\infty}'X_w(\varphi),
\end{equation*}
where the restricted product is taken with respect to the integral subsets; see \cite[\S2.4.7]{KSZ}. This set is independent of the choices of places/integral subsets, and is a nonempty right $\mb{G}(\bb{A}_f^p)$-torsor. We then define the \emph{Langlands--Rapoport set} to be 
\begin{equation*}
X(\varphi)\defeq X_p(\varphi)\times X^p(\varphi).
\end{equation*}
The groups $I_\varphi(\bb{A}_f)$ and $\mb{G}(\bb{A}_f^p)$ act on $X(\varphi)$ on the left and right, respectively. Moreover, the $q$-Frobenius acts on $X(\varphi)$ through its action on $X_p(\varphi)$. These actions are all commuting.
\end{defn}

The actual Langlands--Rapoport conjecture involves a decomposition of $\ms{S}_{\mathsf{K}_0}(\ov{k})$ involving the quotient of $X(\varphi)$ by $I_\varphi(\Q)$. To achieve our desired trace formula though, one only requires a weaker version where one allows a twist of this action.

\begin{defn}[{\cite[\S2.4.7]{KSZ}}]\label{defn:KSZ-adelically-twisted-LR-set}
Let $\eta$ be in $I_\varphi^{\mr{ad}}(\bb{A}_f)$, and let $I_\varphi(\Q)^\eta$ be the image of the diagonal embedding $I_\varphi(\Q)\to I_\varphi(\bb{A}_f)$ under the inner automorphism $\Int(\eta)$ of $I_\varphi(\bb{A}_f)$. Define
\begin{equation*}
S_\eta(\varphi)
\defeq
\varprojlim_{\mathsf{L}^p}
I_\varphi(\Q)^\eta\backslash X(\varphi)/\mathsf{L}^p,
\end{equation*}
where $\mathsf{L}^p$ runs through the compact open subgroups of $\mb{G}(\bb{A}_f^p)$.
\end{defn}

To choose these twists coherently as $\varphi$ varies, Kisin--Shin--Zhu introduce two sheaves of abelian groups on the set of admissible morphisms.

\begin{nota}\label{nota:KSZ-twist-sheaves}
Let $\mc{AM}=\mc{AM}(\mb{G},\mb{X},p,\mc{G})$ be the set of admissible morphisms. Following \cite[\S2.6.16]{KSZ}, write $\varphi_1\approx\varphi_2$ if their algebraic parts $\varphi_1^\Delta$ and $\varphi_2^\Delta$ are $\mb{G}(\ov{\Q})$-conjugate. Let $\mc{H}_{\mr{LR}}$ and $\mc{E}_{\mr{LR}}^p$ denote the sheaves written $\mc{H}$ and $\mc{E}^p$, respectively, in \cite[\S2.6]{KSZ}; their fibers are
\begin{equation*}
\begin{aligned}
\mc{H}_{\mr{LR}}(\varphi)
&\defeq
I_\varphi(\bb{A}_f)\backslash
I_\varphi^{\mr{ad}}(\bb{A}_f)/I_\varphi^{\mr{ad}}(\Q),\\
\mc{E}_{\mr{LR}}^p(\varphi)
&\defeq
I_\varphi(\bb{A}_f^p)\backslash
I_\varphi^{\mr{ad}}(\bb{A}_f^p).
\end{aligned}
\end{equation*}
The natural map $\mc{E}_{\mr{LR}}^p(\varphi)\to\mc{H}_{\mr{LR}}(\varphi)$ is surjective, and both fibers carry canonical abelian group structures; see \cite[\S2.6.13 and Lemma 2.6.14]{KSZ}.
\end{nota}

\begin{rem}  The sheaves $\mc{E}_{\mr{LR}}^p$ and $\mc{H}_{\mr{LR}}$, together with their surjection, are the pullbacks of uniquely determined sheaves and a surjection on $\mc{AM}/\!\approx$. Since the quotient map $\mc{AM}\to\mc{AM}/\!\approx$ factors through the set $\mc{AM}/\mr{conj}$ of $\mb{G}(\ov{\Q})$-conjugacy classes, pulling these descended sheaves back to $\mc{AM}/\mr{conj}$ gives the sheaves used to define $\Gamma(-)_1$ below; see \cite[\S2.6.16]{KSZ}.
\end{rem}

\begin{defn}[{\cite[Definition 2.6.17]{KSZ}}]\label{defn:KSZ-sections-with-descent}
For $\mc{F}$ in $\{\mc{H}_{\mr{LR}},\mc{E}_{\mr{LR}}^p\}$, let $\Gamma(\mc{F})$ be its group of global sections on $\mc{AM}$, $\Gamma(\mc{F})_0$ the subgroup of sections descending to $\mc{AM}/\!\approx$, and $\Gamma(\mc{F})_1$ the subgroup of sections descending to the set of $\mb{G}(\ov{\Q})$-conjugacy classes of admissible morphisms.
\end{defn}

\begin{defn}[{\cite[Definition 2.6.19]{KSZ}}]\label{defn:KSZ-tori-rational}
A section $\tau_{\mr{LR}}$ of $\Gamma(\mc{H}_{\mr{LR}})$ is \emph{tori-rational} if, for every admissible $\varphi$ and every maximal $\Q$-torus $T\subseteq I_\varphi$, the image of $\tau_{\mr{LR}}(\varphi)$ under the composite obstruction map \cite[(2.6.19.1)]{KSZ} is trivial. Tori-rationality for a section of $\Gamma(\mc{E}_{\mr{LR}}^p)$ is defined by the analogous condition in \cite[Definition 2.6.19(ii)]{KSZ}.
\end{defn}

\begin{defn}[{\cite[\S2.7.1]{KSZ}}]\label{defn:KSZ-LR-set}
Let $\tau_{\mr{LR}}$ be in $\Gamma(\mc{H}_{\mr{LR}})_1$. For every admissible $\varphi$, choose $\eta_\varphi$ in $I_\varphi^{\mr{ad}}(\bb{A}_f)$ representing the class $\tau_{\mr{LR}}(\varphi)$. Define the \emph{$\tau_\mr{LR}$-twisted Langlands--Rapoport set} to be
\begin{equation}\label{eq:KSZ-LR-set}
S_{\tau_{\mr{LR}}}(\varphi)
\defeq
S_{\eta_\varphi}(\varphi)
=
\varprojlim_{\mathsf{L}^p}
I_\varphi(\Q)^{\eta_\varphi}\backslash X(\varphi)/\mathsf{L}^p.
\end{equation}
Its isomorphism class as a set with commuting $q$-Frobenius and $\mb{G}(\bb{A}_f^p)$-actions is independent of $\eta_\varphi$ and depends on $\varphi$ only through its $\mb{G}(\ov{\Q})$-conjugacy class.
\end{defn}

We can now finally state Kisin--Shin--Zhu's twisted version of the Langlands--Rapoport conjecture, which they write $\mr{LR}(\mb{G},\mb{X},p,\mc{G},\tau_{\mr{LR}})$; see \cite[\S2.7.1]{KSZ}. 

\begin{nota} Let $\ms{S}_{\mathsf{K}_0}=\varprojlim_{\mathsf{K}'}\ms{S}_{\mathsf{K}'}$, i.e., the limit of integral canonical models over the prime-to-$p$ Hecke tower; see Proposition \ref{prop:prime-to-p-Hecke}.
\end{nota}

\begin{conj}[{The Langlands--Rapoport--$\tau$ conjecture, cf.\ \cite[Conjecture 2.7.3]{KSZ}}]\label{conj:LR-tau}
There exist a tori-rational element $\tau_{\mr{LR}}$ in $\Gamma(\mc{H}_{\mr{LR}})_0$ and a $q$-Frobenius/$\mb{G}(\A_f^p)$-equivariant bijection
\begin{equation}\label{eq:LR-tau-bijection}
\ms{S}_{\mathsf{K}_0}(\ov{k})
\isomto
\coprod_{[\varphi]}S_{\tau_{\mr{LR}}}(\varphi),
\end{equation}
where $[\varphi]$ runs through the $\mb{G}(\ov{\Q})$-conjugacy classes of admissible morphisms.
\end{conj}

\begin{rem}\label{rem:LR-versus-LR-tau}
The usual Langlands--Rapoport conjecture (see \cite[Conjecture 2.5.8]{KSZ}) is the special case in which $\tau_{\mr{LR}}$ is trivial.
\end{rem}

\subsubsection{The upgraded trace formula d'apr\`es Kisin--Shin--Zhu} We now explain, using \cite[Theorem 2.7.4]{KSZ}, how the Langlands--Rapoport-$\tau$ conjecture allows one to give an upgraded version of the primitive trace formula. Again, our explanation will need to import a massive amount of material from op.\@ cit.\@ which we will quickly go over, and we again refer the reader to \cite[\S1-2]{KSZ} for more details.

\begin{setup} We use the same setup as in Setup \ref{setup:LR-general}, but further:
\begin{itemize} 
\item fix a lift $\Phi_v$ in $\Gamma_E$ of geometric $q$-Frobenius;
\item fix functions $f^p$ in $\mc{H}(\mb{G}(\bb{A}_f^p),\mathsf{K}^p)$ and $h$ in $\mc{H}(\mc{G}^c(\Z_p))$;
\item write $[\mu]_{\mb{X}}$ for the image in $\pi_1(G)$ of the Hodge cocharacters of $\mb{X}$;
\item normalize $dg$ and $dg^c$ so that $\mr{vol}(\mc{G}(\Z_p))=\mr{vol}(\mc{G}^c(\Z_p))=1$;
\item let $(\pr^c)^*h\defeq h\circ\pr^c$, an element of $\mc{H}(\mc{G}(\Z_p))$;
\item  for $\xi\colon\mb{G}^c\to\mr{GL}(V)$, we abuse notation and write $\xi(x)$ for $x$ in $\mb{G}(\Q)$ to mean $\xi(\mr{pr}^c(x))$;
\item write $\Gamma_\Q$ for $\Gal(\ov{\Q}/\Q)$;
\item for each rational place $v$ write $\Gamma_v$ for $\Gal(\ov{\Q_v}/\Q_v)$;
\item choose embeddings $\ov{\Q}\hookrightarrow \ov{\Q}_v$ giving rise to an embedding $\Gamma_v\hookrightarrow \Gamma_\Q$.
\end{itemize}
\end{setup}

\begin{eg}\label{eg:pullback-Hecke-function} For $\mathsf{L}_p\subseteq\mc{G}^c(\Z_p)$ compact open, $g_p^c$ in $\mc{G}^c(\Z_p)$, and $g_p$ in $\mc{G}(\Z_p)$ lifting $g_p^c$:
\begin{equation}\label{eq:pullback-Hecke-function}
(\pr^c)^*e(g_p^c,\mathsf{L}_p)=e\left(g_p,(\pr^c)^{-1}(\mathsf{L}_p)\right).
\end{equation}
\end{eg}

Our upgraded version of the primitive trace formula seeks to compute the following terms, using notation adapted from \cite{KSZ}.

\begin{defn}\label{defn:weighted-KSZ-trace}
Let $m\geqslant 1$, fix $\tau$ in $\Phi_v^mI_E$,\footnote{Thus, $m$ is the $j$ in the definition of the local test functions; we use this different letter to match \cite{KSZ}.} and set $n=mr$. We define
\begin{equation}\label{eq:weighted-KSZ-trace}
T(\tau,h,f^p)
\defeq
\tr\left(
\tau\times(\pr^c)^*h\,f^p
\mid H_c^*(\mb{G},\mb{X},\xi)
\right).
\end{equation}
\end{defn}

\begin{nota}\label{nota:local-field-Km}
We continue to use the notation $k_m,W_m,K_m$ from Notation~\ref{nota:finite-field}, with $j=m$, and as above let us write $n=mr=[K_m:\Q_p]$.
\end{nota}

We next give modified versions of our local test functions from Definition \ref{defn:phitauh} to account for the potential discrepancy between $\mb{G}$ and $\mb{G}^c$ (cf.\@ Remark \ref{rem:G-equals-Gc-test-function}).

\begin{nota}
Write $\phi_n^{\mr{KSZ}}\defeq\mathds{1}_{\mc{G}(W_m)\sigma(-\mu_h)(p)\mc{G}(W_m)}$.
\end{nota}

\begin{rem} The function $\phi_n^\mr{KSZ}$ is almost that denoted by $\phi_n$ in \cite[\S1.8.2]{KSZ}. More precisely, with the convention in \eqref{eq:KSZ-local-LR-set}, $\phi_n^{\mr{KSZ}}$ is obtained from the function denoted $\phi_n$
in \cite[\S1.8.2]{KSZ} by precomposition with $\sigma^{-1}$. That said, the two functions have the same twisted orbital integrals and so the difference can be ignored.
\end{rem}

\begin{defn}\label{defn:Shimura-test-function} We define the \emph{Shimura test function} associated to $(\mb{G},\mb{X})$ (with our implicitly chosen $p$ and $\mc{G}$) to be the following function on $G(K_m)$:
\begin{equation}\label{eq:Shimura-test-function}
\phi_{\tau,h}^{\mb{G},\mb{X}}(b)
\defeq
\phi_n^{\mr{KSZ}}(b)\,
\phi_{\tau,h}^{\mc{G}^c,-\mu_h^c}\left(\pr^c(b)\right).
\end{equation}
\end{defn}

\begin{rem} When $\mb{G}\ne \mb{G}^c$, the Shimura test function $\phi_{\tau,h}^{\mb{G},\mb{X}}$ really depends on global data as the quotient $\mc{G}^c$ of $\mc{G}$ does; of course if $\mb{G}=\mb{G}^c$ this is not an issue.
\end{rem}

The following is immediate from Theorem \ref{thm:phi-nice-properties}.

\begin{prop}\label{prop:Shimura-test-function-properties}
The function $\phi_{\tau,h}^{\mb{G},\mb{X}}$ is $\Q$-valued, locally constant, and compactly supported. It is invariant under $\sigma$-conjugation by $\mc{G}(W_m)$, and
\begin{equation*}
\mr{supp}(\phi_{\tau,h}^{\mb{G},\mb{X}})
\subseteq
\mc{G}(W_m)\sigma(-\mu_h)(p)\mc{G}(W_m).
\end{equation*}
\end{prop}

\begin{eg}Let $e_0^c=e(\mc{G}^c(\Z_p))$. Then, as $\mc{D}(\mc{G}^c,\mr{pr}^c(b),-\mu_h^c)_C$ is an open polydisk over $C$, and thus has vanishing higher $\ell$-adic cohomology (as $\ell\ne p$), we see that there is an equality $\phi_{\tau,e_0^c}^{\mc{G}^c,-\mu_h^c}=\mathds{1}_{\mc{G}^c(W_m)\sigma(-\mu_h^c)(p)\mc{G}^c(W_m)}$. Consequently, $\phi_{\tau,e_0^c}^{\mb{G},\mb{X}}=\phi_n^{\mr{KSZ}}$.
\end{eg}

\begin{rem}\label{rem:G-equals-Gc-test-function}
Suppose that $\mb{G}=\mb{G}^c$. The support condition in Definition \ref{defn:coh-for-test-functions} makes the first factor in \eqref{eq:Shimura-test-function} redundant. Thus $\phi_{\tau,h}^{\mb{G},\mb{X}}=\phi_{\tau,h}^{\mc{G},-\mu_h}$. 
\end{rem}

We next recall the notion of Kottwitz parameters which will serve as the index for our ultimate point-counting formula; cf.\@ the triples appearing in \cite{KottwitzAnnArbor}. Roughly these form the group-theoretic description of the isogeny classes parameterized by the parameters $\varphi$ from the Langlands--Rapoport-$\tau$ conjecture.

\begin{nota}\label{nota:Kottwitz-centralizer}
Fix $\gamma_0$ in $\mb{G}(\Q)$ semisimple, and write $I_0\defeq I_{\gamma_0}$.
\end{nota}

\begin{defn}\label{defn:Kottwitz-D}
We define
\begin{equation}\label{eq:Kottwitz-D}
\mf{D}(I_0,\mb{G};\bb{A}_f^p)
\defeq
\prod_{w\ne p,\infty}'
\ker\left(
H^1(\Q_w,I_0)\to H^1(\Q_w,\mb{G})
\right),
\end{equation}
where the restricted product is taken with respect to the trivial classes.
\end{defn}

\begin{defn}[{Kottwitz parameters, \cite[Definition 1.6.4]{KSZ}}]\label{defn:KSZ-Kottwitz-parameter}
A \emph{Kottwitz parameter} is a triple $c=(\gamma_0;a,[b])$ with the following properties:
\begin{enumerate}
\item $\gamma_0$ in $\mb{G}(\Q)$ is semisimple and $\R$-elliptic;
\item $a$ belongs to $\mf{D}(I_0,\mb{G};\bb{A}_f^p)$;
\item $[b]$ belongs to $B(I_{0,\Q_p})$ and $\kappa_G([b]_G)=-[\mu]_{\mb{X}}$ where $[b]_G$ is the image of $[b]$ in $B(G)$ and $\kappa_G$ is the Kottwitz map as in \cite[\S7]{KotIsoII}.
\end{enumerate}
We write $\mf{KP}$ for the set of Kottwitz parameters and $\mf{KP}(\gamma_0)$ for those with first component $\gamma_0$.
\end{defn}

\begin{defn}[{\cite[Definition 1.6.5]{KSZ}}]\label{defn:pn-admissible-KP}
A Kottwitz parameter $c=(\gamma_0;a,[b])$ is \emph{$p^n$-admissible} if there is a representative $b$ of $[b]$ and an element $u$ in $G(\breve{\Q}_p)$ such that
\begin{equation}\label{eq:KSZ-pn-admissible}
u^{-1}\gamma_0u= \delta\sigma(\delta)\cdots\sigma^{n-1}(\delta).
\end{equation}
where $\delta=u^{-1}b\sigma(u)$. We denote the set of $p^n$-admissible Kottwitz parameters by $\mf{KP}_a(p^n)$.
\end{defn}

\begin{obs}[{\cite[Lemma 1.6.7]{KSZ}}]\label{obs:delta-from-Kottwitz-parameter}
If $c$ is $p^n$-admissible, then $\delta$ belongs to $G(K_m)$, and its $\sigma$-conjugacy class in $G(K_m)$ depends only on $c$.
\end{obs}

\begin{constr}\label{constr:classical-KP}
Let $c=(\gamma_0;a,[b])$ belong to $\mf{KP}_a(p^n)$, and choose $\delta$ as above. Write $a=(a_w)_{w\ne p,\infty}$. For each $w\ne p,\infty$, choose $g_w$ in $\mb{G}(\ov{\Q}_w)$ such that the cocycle $\rho\longmapsto g_w\rho(g_w)^{-1}$ represents $a_w$; we choose $g_w=1$ at all but finitely many places $w$ which is possible because $a_w$ is trivial for almost all $w$. Define the \emph{classical Kottwitz parameter} associated with $c$ to be
\begin{equation*}
\gamma_w=g_w^{-1}\gamma_0g_w\in\mb{G}(\Q_w),
\qquad
\gamma=(\gamma_w)_{w\ne p,\infty}\in\mb{G}(\bb{A}_f^p).
\end{equation*}
The equivalence class of $(\gamma_0;\gamma,\delta)$ depends only on $c$.
\end{constr}

\begin{nota}\label{nota:Kottwitz-invariant-group} Set $K_c\defeq\ker\left(\pi_1(I_0)\to\pi_1(\mb{G})\right)$.
\end{nota}

\begin{rem}\label{rem:pi1-surj} Since $I_0$ contains a maximal torus of $\mb{G}$, the map $\pi_1(I_0)\to\pi_1(\mb{G})$ is surjective; see \cite[Proposition 1.7.3]{KSZ}.
\end{rem}

\begin{defn}\label{defn:Kottwitz-invariant-group}
We define
\begin{equation}\label{eq:Kottwitz-invariant-group}
\fE(I_0,\mb{G};\bb{A}/\Q)
\defeq
\frac{(K_c)_{\Gamma,\mr{tors}}}
{\displaystyle
\sum_w\mr{im}\left(
\ker\left((K_c)_{\Gamma_w,\mr{tors}}
\to\pi_1(I_0)_{\Gamma_w}\right)
\to(K_c)_{\Gamma,\mr{tors}}
\right)}.
\end{equation}
By \cite[Proposition 1.7.3]{KSZ}, this group is finite.
\end{defn}

\begin{constr}\label{constr:Kottwitz-beta-away}
Let $c=(\gamma_0;a,[b])$ be a Kottwitz parameter. For $w\ne p,\infty$, the component $a_w$ determines $\beta_w(c)$ in $\pi_1(I_0)_{\Gamma_w,\mr{tors}}$ through the abelianized localization map of \cite[\S1.7.5]{KSZ}. Choose $\wt{\beta}_w(c)$ in $\pi_1(I_0)$ whose image in $\pi_1(I_0)_{\Gamma_w}$ is $\beta_w(c)$ and whose image in $\pi_1(\mb{G})$ is $0$, taking $\wt{\beta}_w(c)=0$ for almost all $w$. Such lifts exist by Remark \ref{rem:pi1-surj}; see \cite[\S1.7.5]{KSZ}.
\end{constr}

\begin{constr}\label{constr:Kottwitz-beta-p}
Set $\beta_p(c)\defeq\kappa_{I_0}([b])$ in $\pi_1(I_0)_{\Gamma_{\Q_p}}$. Choose $\wt{\beta}_p(c)$ in $\pi_1(I_0)$ whose image in $\pi_1(I_0)_{\Gamma_{\Q_p}}$ is $\beta_p(c)$ and whose image in $\pi_1(\mb{G})$ is $-[\mu]_{\mb{X}}$. Such a lift again exists by Remark \ref{rem:pi1-surj}.
\end{constr}

\begin{constr}\label{constr:Kottwitz-beta-infinity}
Choose an elliptic maximal torus $T\subseteq\mb{G}_{\R}$ containing $\gamma_0$ and an element $h_\infty$ in $\mb{X}$ which factors through $T$. The Hodge cocharacter of $h_\infty$ determines $\beta_\infty(c)$ in $\pi_1(I_0)_{\Gamma_{\R}}$. Choose $\wt{\beta}_\infty(c)$ in $\pi_1(I_0)$ whose image in $\pi_1(I_0)_{\Gamma_{\R}}$ is $\beta_\infty(c)$ and whose image in $\pi_1(\mb{G})$ is $[\mu]_{\mb{X}}$. Such a lift again exists by Remark \ref{rem:pi1-surj}; see \cite[\S1.7.5]{KSZ}.
\end{constr}

\begin{defn}\label{defn:KSZ-Kottwitz-invariant}
The sum $\wt{\beta}(c)\defeq\sum_w\wt{\beta}_w(c)$ is finite and belongs to $K_c$: The images of $\wt{\beta}_p(c)$ and $\wt{\beta}_\infty(c)$ in $\pi_1(\mb{G})$ cancel, while every other summand maps to $0$. Its image
\begin{equation}\label{eq:Kottwitz-invariant}
\alpha(c)
\defeq
[\wt{\beta}(c)]
\in\fE(I_0,\mb{G};\bb{A}/\Q)
\end{equation}
is independent of all choices and is called the \emph{Kottwitz invariant} of $c$; see \cite[\S1.7.5]{KSZ}.
\end{defn}

\begin{constr}\label{constr:inner-forms-away}
For $w\ne p,\infty$, let $I_w$ be the inner form of $I_{0,\Q_w}$ determined by the image of $a_w$ in $H^1(\Q_w,I_0^{\mr{ad}})$.
\end{constr}

\begin{constr}\label{constr:inner-form-p}
For a $p^n$-admissible parameter, let $I_p$ be the $\sigma$-centralizer of $b$ in $I_0$: 
\begin{equation*}
I_p(R)
=
\left\{
 g\in I_0(R\otimes_{\Q_p}\breve{\Q}_p):
 g^{-1}b\sigma(g)=b
\right\}
\end{equation*}
\end{constr}

\begin{constr}\label{constr:inner-form-infinity}
Let $I_\infty$ be the inner form of $I_{0,\R}$ determined by the Cartan involution of $(I_0/\mb{Z})_{\R}$ induced by $\Int(h_\infty(i))$, with $h_\infty$ chosen as in Construction \ref{constr:Kottwitz-beta-infinity}.
\end{constr}

\begin{defn}\label{defn:KSZ-inner-form}
Assume that $c$ is $p^n$-admissible and $\alpha(c)=0$. There is an inner form $I(c)$ of $I_0$ over $\Q$, unique up to isomorphism of inner forms, whose localization at every place $w$ is $I_w$; see \cite[Proposition 1.7.12]{KSZ}. 
\end{defn}

\begin{nota}\label{nota:KSZ-measures}
Let $c$ in $\mf{KP}_a(p^n)$ satisfy $\alpha(c)=0$. Fix a Haar measure $di_p\,di^p$ on $I(c)(\bb{A}_f)$. We transport $di^p$ and $di_p$ to the identity components of the centralizers of $\gamma$ and $\delta$, respectively, as in \cite[\S1.8.2]{KSZ}. Let $dx^p$ be the Haar measure on $\mb{G}(\bb{A}_f^p)$ used in the Hecke algebra containing $f^p$, and let $dx_p$ be the Haar measure on $G(K_m)$ which gives $\mc{G}(W_m)$ volume $1$.
\end{nota}

We now describe the weighting factors that will be used in our ultimate trace formula, and can be viewed as counting the number of lattices within a given isogeny class.

\begin{defn}\label{defn:prime-to-p-orbital-integral}
For $f^p$ in $\mc{H}(\mb{G}(\bb{A}_f^p))$, define
\begin{equation}\label{eq:prime-to-p-orbital-integral}
\mr{O}_\gamma(f^p)\defeq
\int_{\mb{G}_\gamma^{\circ}(\bb{A}_f^p)\backslash\mb{G}(\bb{A}_f^p)}f^p(x^{-1}\gamma x)\,\frac{dx^p}{di^p}.
\end{equation}
\end{defn}

\begin{defn}\label{defn:twisted-orbital-integral-4-4}
Let $G_\delta^\sigma$ be the $\sigma$-centralizer of $\delta$ in $G$:
\begin{equation*}
G_\delta^\sigma(R)=\left\{ x\in G(R\otimes_{\Q_p}K_m):x^{-1}\delta\sigma(x)=\delta\right\},
\end{equation*}
and put $I_\delta^\sigma=(G_\delta^\sigma)^\circ$. For $\phi$ in $\mc{H}(G(K_m))$, define
\begin{equation}\label{eq:twisted-orbital-integral-4-4}
\mr{TO}_\delta(\phi)
\defeq
\int_{I_\delta^\sigma(\Q_p)\backslash G(K_m)}
\phi(x^{-1}\delta\sigma(x))\,\frac{dx_p}{di_p}.
\end{equation}
\end{defn}

For notational simplicity, we combine the above away-from-$p$ weighting factors and at-$p$ weighting factors into one term.

\begin{defn}\label{defn:weighted-orbital-term}
For $c$ in $\mf{KP}_a(p^n)$ with $\alpha(c)=0$, define
\begin{equation}\label{eq:weighted-orbital-term}
\mr{O}(c,\tau,h,f^p)\defeq \mr{O}_\gamma(f^p)\,
\mr{TO}_\delta\left(\phi_{\tau,h}^{\mb{G},\mb{X}}\right).
\end{equation}
This is independent of the representative $(\gamma_0;\gamma,\delta)$ of the classical Kottwitz parameter, with the measures transported as in Notation \ref{nota:KSZ-measures}; see \cite[\S1.8.2]{KSZ}.
\end{defn}

As in \cite{KSZ}, we will need to modify our weighting factors by a volume term and a term having to do with connected components of the Shimura variety.

\begin{defn}\label{defn:KSZ-c1}
Let $c$ be as in Notation \ref{nota:KSZ-measures}. Set 
\begin{equation*}
\mb{Z}_{\mathsf{K}_0\mathsf{K}^p}\defeq\mb{Z}(\bb{A}_f)\cap\mathsf{K}_0\mathsf{K}^p,\qquad \mb{Z}(\Q)_{\mathsf{K}_0\mathsf{K}^p}\defeq\mb{Z}(\Q)\cap\mathsf{K}_0\mathsf{K}^p.
\end{equation*}
By \cite[Lemma 1.8.5 and \S1.8.6]{KSZ}, the subgroup $I(c)(\Q)\mb{Z}_{\mathsf{K}_0\mathsf{K}^p}$ is closed in $I(c)(\bb{A}_f)$ and the quotient below is compact. We equip $I(c)(\Q)\mb{Z}_{\mathsf{K}_0\mathsf{K}^p}$ with the Haar measure giving its compact open subgroup $\mb{Z}_{\mathsf{K}_0\mathsf{K}^p}$ volume $1$; see the proof of \cite[Lemma 3.7.4]{KSZ}. Relative to this measure and the measure $di_p\,di^p$ fixed in Notation \ref{nota:KSZ-measures}, we define
\begin{equation}\label{eq:KSZ-c1}
c_1(c,\mathsf{K}^p)
\defeq
\mr{vol}\left(
I(c)(\Q)\mb{Z}_{\mathsf{K}_0\mathsf{K}^p}
\backslash I(c)(\bb{A}_f)
\right).
\end{equation}
\end{defn}

\begin{defn}\label{defn:abelianized-Sha}
Let $H$ be a connected reductive group over $\Q$, and define
\begin{equation*}
\Sha_{\mr{ab}}(\Q,H) \defeq \ker\left(H^1_{\mr{ab}}(\Q,H)\to \prod_{w\text{ a place of }\Q}H^1_{\mr{ab}}(\Q_w,H)\right).
\end{equation*}
\end{defn}

\begin{defn}[{ \cite[\S1.2.5]{KSZ}}]\label{defn:relative-Sha}
For the inclusion $I_0\subseteq\mb{G}$, define
\begin{equation}\label{eq:relative-Sha}
\Sha_{\mb{G}}(\Q,I_0)
\defeq
\ker\left(
\Sha_{\mr{ab}}(\Q,I_0)
\to
\Sha_{\mr{ab}}(\Q,\mb{G})
\right).
\end{equation}
\end{defn}

\begin{nota}[{\cite[\S1.8.6]{KSZ}}]\label{nota:KSZ-c2}
Write $c_2(\gamma_0)\defeq\left|\Sha_{\mb{G}}(\Q,I_0)\right|$. 
\end{nota}

\begin{defn}\label{defn:KSZ-iota}
Let $H$ be a connected reductive group over $\Q$ and let $\epsilon$ in $H(\Q)$ be semisimple. We define $\iota_H(\epsilon)\defeq[H_\epsilon(\Q):H_\epsilon^\circ(\Q)]$.
\end{defn}

\begin{defn}\label{defn:KSZ-iota-bar}
With the notation of Definition \ref{defn:KSZ-iota}, define $\ov{\iota}_H(\epsilon)\defeq\left|(H_\epsilon/H_\epsilon^\circ)(\Q)\right|$. 
\end{defn}

\begin{lem}[{cf.\@ \cite[\S1.8.6]{KSZ}}]\label{lem:measure-independence}
The product $c_1(c,\mathsf{K}^p)\mr{O}(c,\tau,h,f^p)$ is independent of the choice of Haar measures.
\end{lem}

\begin{nota}\label{nota:Sigma-Kp}
Let $\Sigma_{\R\text{-}\mr{ell}}(\mb{G})$ be the set of stable conjugacy classes of semisimple $\R$-elliptic elements of $\mb{G}(\Q)$. We choose a collection $\Sigma_{\mathsf{K}^p}\subseteq\mb{G}(\Q)$ whose image in $\Sigma_{\R\text{-}\mr{ell}}(\mb{G})$ contains exactly one representative of each $\mb{Z}(\Q)_{\mathsf{K}_0\mathsf{K}^p}$-translation orbit.
\end{nota}

As indicated before, Kottwitz parameters are related to the parameters $\varphi$ appearing in the Langlands--Rapoport-$\tau$ conjecture. We now make this precise.

\begin{defn}[{\cite[Definitions 3.1.1 and 3.2.1]{KSZ}}]\label{defn:KSZ-gunstig-gelegen}
An \emph{LR pair} is a pair $(\varphi,\epsilon)$ consisting of a morphism $\varphi\colon\mf{Q}\to\mf{G}_{\mb{G}}$ and an element $\epsilon$ in $I_\varphi(\Q)$; it is \emph{semi-admissible} if $\varphi$ is admissible. For each $\rho$ in $\Gamma_\Q$, fix a lift $q_\rho$ in $\mf{Q}^{\mr{top}}$ as in \cite[\S2.6.1]{KSZ}, and write $I_{\varphi,\epsilon}^{\circ}\defeq(I_\varphi)_\epsilon^{\circ}$. We say that $(\varphi,\epsilon)$ is \emph{g\"unstig gelegen} (abbreviated \emph{gg}) if:
\begin{enumerate}
\item the embedding $I_{\varphi,\ov{\Q}}\hookrightarrow\mb{G}_{\ov{\Q}}$ carries $\epsilon$ to an element of $\mb{G}(\Q)$ which is $\R$-elliptic;
\item  $g_\rho$ is in $(\mb{G}_{\ov{\Q},\varphi^\Delta})_\epsilon^{\circ}(\ov{\Q})=I_{\varphi,\epsilon}^{\circ}(\ov{\Q})$ for every $\rho$ in $\Gamma_\Q$, where $\varphi(q_\rho)=g_\rho\rtimes\rho$ with $g_\rho$ in $\mb{G}(\ov{\Q})$.
\end{enumerate}
Here $\mb{G}_{\ov{\Q},\varphi^\Delta}$ denotes the centralizer in $\mb{G}_{\ov{\Q}}$ of the image of $\varphi^\Delta$. The second condition is independent of the choices of $q_\rho$ by \cite[Remark 3.2.2]{KSZ}. We write $\mf{LRP}^{gg}$ for the set of gg LR pairs and $\mf{LRP}^{gg}_a(q^m)$ for its subset of $q^m$-admissible pairs; see \cite[Definition 3.2.13]{KSZ}.
\end{defn}

\begin{nota}\label{nota:controlled-LR-lift}
Fix a tori-rational element $\tau_{\mr{LR}}$ in $\Gamma(\mc{H}_{\mr{LR}})_0$. Let $\mc{E}_{\mr{LR}}^p$ be the twisting sheaf of Notation \ref{nota:KSZ-twist-sheaves}. By \cite[Lemma 2.6.20 and \S3.6]{KSZ}, we may choose a tori-rational lift $\wt{\tau}_{\mr{LR}}$ in $\Gamma(\mc{E}_{\mr{LR}}^p)_0$ of $\tau_{\mr{LR}}$. For a semi-admissible Langlands--Rapoport pair $(\varphi,\epsilon)$ satisfying Definition \ref{defn:KSZ-gunstig-gelegen}, we abbreviate the Kottwitz parameter $t(\varphi,\epsilon,\wt{\tau}_{\mr{LR}}(\varphi))$ of \cite[\S3.5--\S3.6]{KSZ} to $t_{\wt{\tau}_{\mr{LR}}}(\varphi,\epsilon)$.
\end{nota}

Now, up until this point, we have not used that $\ms{S}_{\mathsf{K}}$ is an integral canonical model (in the sense of Definition \ref{defn:ICM}) explicitly. That said, we now impose the following compatibility at $p$ between the map $\delta$ from \ref{constr:corr-delta} and the Langlands--Rapoport-$\tau$ setup.

\begin{nota}\label{nota:LR-local-coordinate}
Fix $m\geqslant 1$, put $n=mr$, and fix a point $y$ of $c(g^p,\mathsf{K})^{(m)}$. Fix $(\varphi,\epsilon)$ in $\mf{LRP}^{gg}_a(q^m)$ representing the corresponding Langlands--Rapoport class, as permitted by \cite[Corollary 3.2.17]{KSZ}, and put $c\defeq t_{\wt{\tau}_{\mr{LR}}}(\varphi,\epsilon)=(\gamma_0;a,[b])$.

 By \cite[Proposition 3.5.8]{KSZ}, the parameter $c$ is $p^n$-admissible. Fix a representative $(\gamma_0;\gamma,\delta)$ of the classical Kottwitz parameter associated with $c$, and fix $g$ in $G(K_m)$ such that $g\mc{G}(W_m)$ is the $p$-coordinate of $y$ in the calculation of \cite[Lemma 3.7.4]{KSZ}. With the convention of Definition \ref{defn:KSZ-local-LR-set}, we then define
\begin{equation}\label{eq:LR-p-coordinate-Hodge-condition}
\delta^{\mr{LR}}(y)
\defeq
g^{-1}\delta\sigma(g)
\in
\mc{G}(W_m)\sigma(-\mu_h)(p)\mc{G}(W_m).
\end{equation}
\end{nota}

\begin{defn}\label{defn:LR-realization-compatible}
A Langlands--Rapoport--$\tau$ bijection is \emph{compatible with syntomic realization} if, for every $m\geqslant 1$ and every fixed point $y$ as in Notation \ref{nota:LR-local-coordinate},  one has $\delta(y)=\left\langle\pr^c\left(\delta^{\mr{LR}}(y)\right)\right\rangle$ in $C_m(\mc{G}^c,\sigma(-\mu_h^c))$, where $\delta(y)$ is as in Construction \ref{constr:corr-delta}.
\end{defn}

\begin{assumption}\label{ass:integral-LR-tau}
We assume that Conjecture \ref{conj:LR-tau} holds for the $\tau_{\mr{LR}}$ fixed in Notation \ref{nota:controlled-LR-lift}, and a Langlands--Rapoport-$\tau$ bijection which is compatible with syntomic realization.
\end{assumption}

\begin{obs}\label{obs:LR-local-weight}
Assume that the Langlands--Rapoport--$\tau$ bijection is compatible with syntomic realization as in Definition \ref{defn:LR-realization-compatible}. For a fixed point with data as in Notation \ref{nota:LR-local-coordinate}, one has $\phi_{\tau,h}^{\mc{G}^c,-\mu_h^c}(\delta(y))=\phi_{\tau,h}^{\mb{G},\mb{X}}\left(\delta^{\mr{LR}}(y)\right)$. 
\end{obs}

We now repeat the point-counting argument of \cite[\S3.7]{KSZ} with the $p$-local weighting factor $\phi_n^\mr{KSZ}$ from loc.\@ cit.\@ replaced by the Shimura test function $\phi_{\tau,h}^{\mb{G},\mb{X}}$.

\begin{nota}\label{nota:KSZ-fixed-point-piece}
Let $(\varphi,\epsilon)$ be a semi-admissible Langlands--Rapoport pair satisfying Definition \ref{defn:KSZ-gunstig-gelegen}, and put $c=t_{\wt{\tau}_{\mr{LR}}}(\varphi,\epsilon)$. For $f^p=e((g^p)^{-1},\mathsf{K}^p)$, we write $\mc{O}(\varphi,\epsilon,m,g^p)$ for the corresponding subset of the fixed-point set defined in \cite[Lemma 3.7.3]{KSZ}.
\end{nota}

\begin{defn}\label{defn:weighted-fixed-point-sum}
For the set in Notation \ref{nota:KSZ-fixed-point-piece}, define
\begin{equation}\label{eq:weighted-fixed-point-sum}
\mr{wt}_{\tau,h}
\left(\mc{O}(\varphi,\epsilon,m,g^p)\right)
\defeq
\sum_{y\in\mc{O}(\varphi,\epsilon,m,g^p)}
\phi_{\tau,h}^{\mc{G}^c,-\mu_h^c}(\delta(y)).
\end{equation}
\end{defn}

\begin{prop}\label{prop:weighted-KSZ-count}
Assume that the Langlands--Rapoport bijection is compatible with syntomic realization. Let $(\varphi,\epsilon)$ and $c$ be as in \emph{Notation \ref{nota:KSZ-fixed-point-piece}}, and let $f^p=e((g^p)^{-1},\mathsf{K}^p)$; by \cite[Proposition 3.6.3]{KSZ}, one has $\alpha(c)=0$. If $(\varphi,\epsilon)$ is $q^m$-admissible, then $c$ is $p^n$-admissible by \cite[Proposition 3.5.8]{KSZ}, and one has the equality
\begin{equation}\label{eq:weighted-KSZ-count}
\mr{wt}_{\tau,h}
\left(\mc{O}(\varphi,\epsilon,m,g^p)\right)= \iota_{I_\varphi}(\epsilon)^{-1}
 c_1(c,\mathsf{K}^p)\mr{O}(c,\tau,h,f^p).
\end{equation}
Otherwise, the fixed-point set on the left is empty.
\end{prop}
\begin{proof}
The last assertion is the first assertion of \cite[Lemma 3.7.4]{KSZ}. We therefore assume that $(\varphi,\epsilon)$ is $q^m$-admissible. Then $c$ is $p^n$-admissible by \cite[Proposition 3.5.8]{KSZ}, and we follow the rest of the proof of \cite[Lemma 3.7.4]{KSZ}, with the following observations.

The weighted sum over the prime-to-$p$ coordinate in loc.\@ cit.\@ gives the orbital integral $\mr{O}_\gamma(f^p)$. After the change of variables in \cite[(3.7.4.2)]{KSZ} and with our conventions, the set of $p$-coordinates can be written
\begin{equation}\label{eq:KSZ-Wp}
W_p
=
\left\{
 g\mc{G}(W_m)\in G(K_m)/\mc{G}(W_m):
 g^{-1}\delta\sigma(g)
 \in\mc{G}(W_m)\sigma(-\mu_h)(p)\mc{G}(W_m)
\right\}.
\end{equation}
By Observation \ref{obs:LR-local-weight}, the summand attached to $g\mc{G}(W_m)$ is the value of $\phi_{\tau,h}^{\mb{G},\mb{X}}$ at $g^{-1}\delta\sigma(g)$, which vanishes outside the set in \eqref{eq:KSZ-Wp}. So, the weighted sum over the $p$-coordinate is $\mr{TO}_\delta\left(\phi_{\tau,h}^{\mb{G},\mb{X}}\right)$.

The remaining quotient by the global centralizer and the comparison of measures are unchanged; they give the factor $\iota_{I_\varphi}(\epsilon)^{-1}c_1(c,\mathsf{K}^p)$ exactly as in \cite[Lemma 3.7.4]{KSZ}. Combining the prime-to-$p$ and $p$-adic factors proves \eqref{eq:weighted-KSZ-count}.
\end{proof}

\begin{nota}\label{nota:admissible-LR-pairs}
Write $[\mf{LRP}^{gg}_a(q^m)]$ for the set of $\mb{G}(\ov{\Q})$-conjugacy classes of $q^m$-admissible Langlands--Rapoport pairs, with notation as in \cite[\S3.1]{KSZ}; $\mb{Z}(\Q)_{\mathsf{K}_0\mathsf{K}^p}$ acts by $z\cdot[\varphi,\epsilon]=[\varphi,z\epsilon]$.
\end{nota}

\begin{cor}\label{cor:weighted-LR-sum}
Assume that $(\mb{G},\mb{X})$ has $\ell$-adically contractible boundary, and that \emph{Assumption \ref{ass:integral-LR-tau}} holds. Then, there exists an integer $m_0\geqslant 1$, depending only on $f^p$, such that, for every $m\geqslant m_0$ and every $\tau$ in $\Phi_v^mI_E$, one has
\begin{equation}\label{eq:weighted-LR-sum}
T(\tau,h,f^p)
=
\sum_{[\varphi,\epsilon]}
\iota_{I_\varphi}(\epsilon)^{-1}
 c_1\left(t_{\wt{\tau}_{\mr{LR}}}(\varphi,\epsilon),\mathsf{K}^p\right)\cdot
 O\left(t_{\wt{\tau}_{\mr{LR}}}(\varphi,\epsilon),\tau,h,f^p\right)
 \tr(\xi(\epsilon)),
\end{equation}
where the sum runs over a set of representatives for $[\mf{LRP}^{gg}_a(|k|^m)]
\big/
\mb{Z}(\Q)_{\mathsf{K}_0\mathsf{K}^p}$. If $\ms{S}_\mathsf{K}$ is proper, then one may take $m_0=1$.
\end{cor}
\begin{proof}
By linearity, it suffices to take $h=e((g_p^c)^{-1},\mathsf{L}_p)$ and $f^p=e((g^p)^{-1},\mathsf{K}^p)$. Choose a lift $g_p$ in $\mc{G}(\Z_p)$ of $g_p^c$ and put $\mathsf{K}_p=(\pr^c)^{-1}(\mathsf{L}_p)$. Example \ref{eg:pullback-Hecke-function} identifies the operator in \eqref{eq:weighted-KSZ-trace} with the Hecke operator $e((g_pg^p)^{-1},\mathsf{K}_p\mathsf{K}^p)$ at level $\mathsf{K}_p\mathsf{K}^p$.

We now apply Theorem \ref{thm:prim-trace-formula-II}. More precisely, we use \eqref{eq:LR-tau-bijection} to decompose the fixed-point set as in \cite[Lemma 3.7.3]{KSZ}. Then, the naive local term for the automorphic \'etale sheaf on the piece indexed by $(\varphi,\epsilon)$ is $\tr(\xi(\epsilon))$, and Proposition \ref{prop:weighted-KSZ-count} takes care of the remaining part.
\end{proof}

To pass from \eqref{eq:weighted-LR-sum} to the final Kottwitz sum, we group Langlands--Rapoport pairs by their Kottwitz parameter. The next definition and lemma give the resulting multiplicity.

\begin{defn}[{\cite[(3.7.5.3)]{KSZ}}]\label{defn:KSZ-Ac}
Fix $m\geqslant 1$ and the controlled lift $\wt{\tau}_{\mr{LR}}$ of Notation \ref{nota:controlled-LR-lift}. For a Kottwitz parameter $c$, define
\begin{equation}\label{eq:KSZ-Ac}
\mc{A}(c)
\defeq
\sum_{\substack{[\varphi,\epsilon]\in[\mf{LRP}^{gg}_a(|k|^m)]\\
 t_{\wt{\tau}_{\mr{LR}}}(\varphi,\epsilon)\simeq c}}
\iota_{I_\varphi}(\epsilon)^{-1}.\footnote{The dependence on $m$ and $\wt{\tau}_{\mr{LR}}$ is suppressed from the notation.}
\end{equation}
\end{defn}

\begin{obs}\label{obs:weighted-nonvanishing}
If $\mr{O}(c,\tau,h,f^p)$ is non-zero, then so is $\mr{O}_\gamma(f^p)\,\mr{TO}_\delta(\phi_n^{\mr{KSZ}})$.
\end{obs}

\begin{lem}\label{lem:weighted-KSZ-grouping}
Let $\gamma_0$ be in $\Sigma_{\mathsf{K}^p}$ and let $c$ in $\mf{KP}(\gamma_0)\cap\mf{KP}_a(p^n)$ satisfy $\alpha(c)=0$. Assume that $\mr{O}(c,\tau,h,f^p)\ne 0$, and let $k(c)$ be the number of elements of $\mf{KP}(\gamma_0)\cap\mf{KP}_a(p^n)$ with vanishing Kottwitz invariant which are isomorphic to $c$. Then $k(c)$ is finite, and
\begin{equation}\label{eq:weighted-KSZ-grouping}
\mc{A}(c)=k(c)\,\ov{\iota}_{\mb{G}}(\gamma_0)^{-1}c_2(\gamma_0).
\end{equation}
\end{lem}
\begin{proof}
Observation \ref{obs:weighted-nonvanishing} verifies the nonvanishing hypothesis of \cite[Lemma 3.7.6]{KSZ}; the finiteness of $k(c)$ and Equation \eqref{eq:weighted-KSZ-grouping} are precisely the conclusions of that lemma. Despite the result from loc.\@ cit.\@ not being directly related here, its proof still applies as it only concerns the fibers of the map from Langlands--Rapoport pairs to Kottwitz parameters and is independent of the $p$-local weighting function.
\end{proof}

\begin{thm}[The Langlands--Kottwitz--Scholze formula]\label{thm:Langlands-Kottwitz-Scholze-formula}
Suppose that $(\mb{G},\mb{X})$ has $\ell$-adically contractible boundary, and that \emph{Assumption \ref{ass:integral-LR-tau}} holds. Then, there exists an integer $m_0\geqslant 1$, depending only on $f^p$, such that, for every $m\geqslant m_0$ and every $\tau$ in $\Phi_v^mI_E$, one has, with $n=mr$,
\begin{equation}\label{eq:Langlands-Kottwitz-Scholze-formula}
T(\tau,h,f^p)=
\sum_{\gamma_0\in\Sigma_{\mathsf{K}^p}}\sum_{\substack{c\in\mf{KP}(\gamma_0)\cap\mf{KP}_a(p^n)\\ \alpha(c)=0}}
\ov{\iota}_{\mb{G}}(\gamma_0)^{-1}
 c_2(\gamma_0)\tr(\xi(\gamma_0))
 c_1(c,\mathsf{K}^p)
 \mr{O}_\gamma(f^p)
 \mr{TO}_\delta\left(\phi_{\tau,h}^{\mb{G},\mb{X}}\right),
\end{equation}
where for each $c$ in the inner sum, $(\gamma_0;\gamma,\delta)$ denotes any representative of its associated classical Kottwitz parameter. If $\ms{S}_{\mathsf{K}_0\mathsf{K}^p}$ is proper, then one may take $m_0=1$.
\end{thm}
\begin{proof}
Starting from Corollary \ref{cor:weighted-LR-sum}, we may group the Langlands--Rapoport pairs according to the isomorphism class of their Kottwitz parameter: This gives the analogue here of \cite[(3.7.5.2)]{KSZ}, but with our $\mr{O}(c,\tau,h,f^p)$ in place of the unweighted orbital term.

For every parameter with nonzero weighted orbital term, Observation \ref{obs:weighted-nonvanishing} verifies the support hypothesis of \cite[Proposition 3.6.5]{KSZ}, and thus the parameter occurs in the image of $t_{\wt{\tau}_{\mr{LR}}}$. So, extending the grouped sum to all Kottwitz parameters satisfying the conditions in the index of the sum introduces only zero terms. Lemma \ref{lem:weighted-KSZ-grouping} then replaces the multiplicity $\mc{A}(c)$ of every nonzero term by $k(c)\,\ov{\iota}_{\mb{G}}(\gamma_0)^{-1}c_2(\gamma_0)$. Finally, the construction of $t_{\wt{\tau}_{\mr{LR}}}(\varphi,\epsilon)$ has first component $\gamma_0=\epsilon$, so $\tr(\xi(\epsilon))=\tr(\xi(\gamma_0))$, as desired.
\end{proof}

\subsubsection{The abelian-type case} Finally, we wish to show that Theorem \ref{thm:Langlands-Kottwitz-Scholze-formula} is unconditional in the case when $(\mb{G},\mb{X})$ is of abelian type. For this, we must only show that the Langlands--Rapoport-$\tau$ bijections constructed by Kisin--Shin--Zhu are compatible with syntomic realization as in Definition \ref{defn:LR-realization-compatible}.

We will do this via the standard Lovering-type yoga; see \cite[\S5.1]{DanielsYoucis}. In particular, we shall prove the claim for Hodge-type Shimura varieties, and then for toral-type Shimura varieties, and finally for abelian-type Shimura varieties (the latter two will be done in the same proof).

\begin{prop}\label{prop:Hodge-LR-realization-compatibility}
Suppose that $(\mb{G},\mb{X})$ is of Hodge type. Then, the
Langlands--Rapoport--$\tau$ bijection constructed in
\cite[\S5.10 and \S5.13]{KSZ} is compatible with syntomic realization.
\end{prop}
\begin{proof}
Let $y$, $(\varphi,\epsilon)$, and $n=mr$ be as in
Notation \ref{nota:LR-local-coordinate}, and let $x$ be the point
of Construction \ref{constr:corr-delta}.
Write $\Upsilon_x$ for the integral $F$-isocrystal with
$\mc{G}$-structure of \cite[Lemma 5.2.4]{KSZ}, base changed to
$W_\infty$. By \cite[Theorem 2.18 and Proposition 2.19]{IKY2} and
Proposition \ref{prop:BT-in-char-p}, applying $\pr^c$ to
$\Upsilon_x$ gives the crystalline realization of
$(\mf{Q}_{\mathsf{K}})_x$.
Choose a lift of $x$ to $\ms{S}_{\mathsf{K}_0}$ and the marking
and trivializations of \cite[\S5.10.3]{KSZ} for the matching of
\cite[Theorem 5.13.9]{KSZ}.

Under these trivializations, let us denote by $d\sigma$ the Frobenius on
$\Upsilon_x[\nicefrac{1}{p}]$, and let $u\mc{G}(W_\infty)$ pick out the integral lattice $\Upsilon_x$. Let $g_0$ and $b_\varphi$ be the corresponding
choices in the description \eqref{eq:KSZ-local-LR-set}.
Note that the elements $d$ and $b_\varphi$ are denoted by $\delta_y$ and
$b_{y'}$ in \cite[\S5.10.3]{KSZ}, respectively.

Let us now consider the element $e$ of $\mc{G}(W_\infty)$ chosen in loc.\@ cit.\@, which satisfies $b_\varphi=ed\sigma(e)^{-1}$. Under the lattice description of \cite[Remark 5.6.6]{KSZ} and \eqref{eq:KSZ-local-LR-set}, the map defining the bijection in \cite[Proposition 5.10.4]{KSZ} sends the $p$-coordinate $u\mc{G}(W_\infty)$ to $eu\mc{G}(W_\infty)$ in $X_p(\varphi)$. Indeed, the two maps
of \cite[\S5.10.3]{KSZ}, followed by our convention change, give
\begin{equation*}
 b_\varphi\sigma\bigl(e\sigma^{-1}(d^{-1}u)e^{-1}\bigr)
 =eu\sigma(e)^{-1},
\end{equation*}
as can be quickly calculated. 

Now, let us use the controlled representative $\eta_\varphi=(1,\eta_\varphi^p)$
of Notation \ref{nota:controlled-LR-lift}. Then, by \cite[\S2.4.7]{KSZ},
the action of $\epsilon$ on $X_p(\varphi)$, in these coordinates, is given by
\begin{equation*}
 v\mc{G}(W_\infty)\longmapsto\epsilon_pv\mc{G}(W_\infty),
 \qquad \epsilon_p\defeq g_0^{-1}\epsilon g_0.
\end{equation*}
In particular, we have the equality $\epsilon_p b_\varphi=b_\varphi\sigma(\epsilon_p)$; see \cite[Proposition 2.2.6(ii)]{KSZ}.
In the same coordinates, let us write $v\mapsto c_p^{-1}v$ for the
composite map of \cite[Lemma 3.7.4]{KSZ}, with the choice satisfying \cite[(1.6.5.1)]{KSZ}. Then, we have
\begin{equation*}
 \delta=c_p^{-1}b_\varphi\sigma(c_p),\qquad
 c_p^{-1}\epsilon_p c_p
 =\delta\sigma(\delta)\cdots\sigma^{n-1}(\delta).
\end{equation*}
Note that our convention changes commute with $v\mapsto c_p^{-1}v$, since
\begin{equation*}
 c_p^{-1}\bigl(b_\varphi\sigma(v)\bigr)
 =\delta\sigma(c_p^{-1}v).
\end{equation*}
The finite-field calculation in \cite[Lemma 3.7.4]{KSZ} implies that after
changing $u$ within its right integral coset we may take $g=c_p^{-1}eu$ in $G(K_m)$. Consequently, we have the equality $\delta^{\mr{LR}}(y)=g^{-1}\delta\sigma(g)=u^{-1}d\sigma(u)$.

The comparison maps respect the quasi-isogeny and Hecke identifications by the proof of \cite[Proposition 5.10.4]{KSZ}. Thus the descent in \eqref{eq:tube-comp-descent-datum} is represented by $\epsilon_p^{-1}(b_\varphi\sigma)^n$. On the other hand, the norm identity gives
\begin{equation*}
 c_p^{-1}\epsilon_p^{-1}(b_\varphi\sigma)^n c_p=\sigma^n.
\end{equation*}
Since $\sigma^n(g)=g$, the descended integral crystal over $W_m$
has Frobenius factor $g^{-1}\delta\sigma(g)$. Thus, applying $\pr^c$ and Proposition \ref{prop:isom-classes-in-BT} yields $\delta(y)=\left\langle\pr^c\big(\delta^{\mr{LR}}(y)\big)\right\rangle$, as desired.
\end{proof}

To carry out our Lovering-type yoga to from the Hodge-type case to the abelian-type case, we make use of the following notion of nice Hodge-type liftings.

\begin{nota}[{\cite[Definition 6.2.1]{KSZ} }]\label{nota:nice-Hodge-lifting}
Let $(\mb{G}_2,\mb{X}_2,p,\mc{G}_2)$ be an unramified Shimura datum of abelian type. A \emph{nice lifting} is an unramified Hodge-type datum $(\mb{G}_1,\mb{X}_1,p,\mc{G}_1)$ together with an isomorphism of Shimura data $(\mb{G}_1^{\mr{ad}},\mb{X}_1^{\mr{ad}})\isomto(\mb{G}_2^{\mr{ad}},\mb{X}_2^{\mr{ad}})$ such that $E(\mb{G}_1,\mb{X}_1)_p=E(\mb{G}_1^{\mr{ad}},\mb{X}_1^{\mr{ad}})_p$, the induced isomorphism $\mb{G}_1^{\mr{ad}}\isomto\mb{G}_2^{\mr{ad}}$ lifts to a unique central isogeny $\mb{G}_1^{\mr{der}}\to\mb{G}_2^{\mr{der}}$, and its localization at $p$ extends to an isomorphism $\mc{G}_1^{\mr{ad}}\isomto\mc{G}_2^{\mr{ad}}$. 
\end{nota}

\begin{rem} By \cite[Lemma 6.2.2]{KSZ}, nice liftings of unramified abelian-type Shimura data always exist.
\end{rem}

\begin{prop}\label{prop:KSZ-integral-realization-compatibility}
Suppose that $(\mb{G},\mb{X})$ is of abelian type. The Langlands--Rapoport--$\tau$ bijection furnished by \cite[Theorem 6.3.5]{KSZ} can be chosen to be compatible with the syntomic realization.
\end{prop}

\begin{rem}[Setup for the proof]\label{rem:proof-setup} Before we begin the proof, let us first describe some setup. First, we choose the nice lifting $(\mb{G}_1,\mb{X}_1,p,\mc{G}_1)$, the matching, and the element $\tau_{\mr{LR}}$ used in \cite[Theorem 6.3.5]{KSZ} and its proof, with the representatives $\eta_\varphi=(1,\eta_\varphi^p)$ of Notation \ref{nota:controlled-LR-lift}. Second, in the proof we freely use \cite[Theorem 2.18]{IKY2}, which shows that the crystalline realizations of the syntomic canonical models are the crystalline canonical models from \cite{LoveringFCrystals}; for readability we shall refer to these common objects as `canonical'.

We will then roughly follow the same strategy used in the proof of \cite[Theorem 3.6]{IKY2}, making use of the data from \cite[Lemma 2.9]{IKY2}. For this construction, set $\mb{E}_1=E(\mb{G}_1,\mb{X}_1)$, $\mb{F}=\mb{E}_1\mb{E}$, and
 $\mb{T}=\mr{Res}_{\mb{F}/\Q}\,\bb{G}_{m,\mb{F}}$. The auxiliary datum is formed over the common field $\mb{F}$, as in \cite[Proposition 3.4.2 and \S3.4.3]{LoveringFCrystals}. With $h_T$ and the maps of that construction, we obtain

\begin{equation*}
 \begin{gathered}
 (\mb{G}_1\times\mb{T},\mb{X}_1\times\{h_T\})
 \xleftarrow{\alpha}(\mb{G}_2,\mb{X}_2)
 \xrightarrow{\beta}(\mb{G},\mb{X}),\\
 \mc{G}_1\times\mc{T}^c
 \xleftarrow{\alpha^c}
 \mc{G}_2^c=\mc{G}_1\times_{\mc{G}_1^{\mr{ab}}}\mc{T}^c
 \xrightarrow{\beta^c}\mc{G}^c,
 \end{gathered}
\end{equation*}
as described in loc.\@ cit. Let us write $\pi\colon\mc{G}_1\to\mc{G}_1^{\mr{ab}}$ and $\rho\colon\mc{T}^c\to\mc{G}_1^{\mr{ab}}$ for the maps defining $\mc{G}_2^c$. 

Write $\ms{S}_1,\ms{S}_T,\ms{S}_2,\ms{S}$ for the integral canonical
models over $W_\infty$ at compatible neat prime-to-$p$ levels, which
we suppress. Functoriality, as in \cite[Lemma 2.9]{IKY2}, gives
\begin{equation*}
 \ms{S}_1\times_{W_\infty}\ms{S}_T
 \xleftarrow{\alpha_{\ms{S}}}\ms{S}_2
 \xrightarrow{\beta_{\ms{S}}}\ms{S}.
\end{equation*}
Fix $x_2$ in $\ms{S}_2(\ov{k})$, and set
\begin{equation*}
 x_1\defeq\pr_1\alpha_{\ms{S}}(x_2),\qquad
 x_T\defeq\pr_T\alpha_{\ms{S}}(x_2),\qquad
 x\defeq\beta_{\ms{S}}(x_2).
\end{equation*}
Here $\pr_1$ and $\pr_T$ are the projections from
$\ms{S}_1\times_{W_\infty}\ms{S}_T$. 
We will then verify that the relevant Langlands--Rapoport-$\tau$ conjecture is compatible with syntomic realization at the points $x_1$ (already completed in Proposition \ref{prop:Hodge-LR-realization-compatibility}), $x_T$, $x_2$, and $x$, in that order. 

Let first spell out what that means. In integral trivializations for an $F$-crystal with $\mc{H}^c$-structure (for a group $\mc{H}$ occurring above), write $b$ and $b'$ in $\mc{H}^c(K_\infty)$ for the Langlands--Rapoport and canonical Frobenius factors (really the first after applying $\mr{pr}^c$). The required comparison is an element $a$ in $\mc{H}^c(W_\infty)$ such that $b'=a^{-1}b\sigma(a)$, chosen compatibly with the maps of Shimura data as well as the quotient and Hecke identifications. For a Hecke-fixed point, let us set $n=mr$ as in Notation \ref{nota:LR-local-coordinate}, and write $D\sigma^n$ and $D'\sigma^n$ for the respective descents of \eqref{eq:tube-comp-descent-datum}. The same element must also satisfy the equation $D'=a^{-1}D\sigma^n(a)$. Indeed, these two identities identify the descended integral crystals of $b$ and $b'$, and hence their classes by Proposition \ref{prop:isom-classes-in-BT}. 

Once this comparison at the points are done, we will have shown that the Langlands--Rapoport-$\tau$ conjecture on $\ms{S}$ is compatible at all points in the image of $\ms{S}_2$. We then finally use compatibility with the Hecke action to reduce from a general point to such a point in the image.
\end{rem}

\begin{proof} We continue the setup from Remark \ref{rem:proof-setup}.

\paragraph*{Step 1: the case of $x_T$} We use the Langlands--Rapoport bijection for $(\mb{T},\{h_T\})$, from \cite[Proposition 3.6.7]{KisinModp}. Denote by $b_T$ and $ b_T'$ the Langlands--Rapoport and canonical Frobenius factors, respectively. Using \cite[\S4.3.13 and Proposition 4.3.14]{KSZ} applied to $\mb{T}^c$, and \cite[\S\S3.2.1--3.2.2]{LoveringFCrystals}, we obtain an element $t$ in $\mc{T}^c(K_\infty)$ with $b_T'=t^{-1}b_T\sigma(t)$. Now, both $b_T$ and $b_T'$ have valuation $\sigma(-\mu_{h_T}^c)$, by \eqref{eq:KSZ-local-LR-set} and Proposition \ref{prop:BT-in-char-p}. Thus, choosing a $\sigma$-invariant element $\nu$ of $X_\ast(\mc{T}^c_{K_\infty})$ with  $\langle\chi,\nu\rangle=v_p(\chi(t))$ for $\chi$ in $X^*(\mc{T}^c_{K_\infty})$, we may set $t^\circ\defeq t\nu(p)^{-1}$, an element of $\mc{T}^c(W_\infty)$. Then, by design $(t^\circ)^{-1}b_T\sigma(t^\circ)=b_T'$. But, also this correction $t\mapsto t^\circ$ is multiplicative, commutes with torus morphisms and $\sigma$, and fixes integral elements. It therefore preserves Frobenius compatibility and the integral
identifying maps in \eqref{eq:tube-comp-descent-datum}. Thus this
Langlands--Rapoport bijection is compatible with syntomic realization at $x_T$.

\medskip

\paragraph*{Step 2: the case of $x_2$} For $(\mb{G}_2,\mb{X}_2)$, the Langlands--Rapoport-$\tau$ bijection of
\cite[\S6.2.3 and Theorem 6.2.4]{KSZ} is obtained from the chosen
Hodge-type bijection on a connected component, using
$\mb{G}_2^{\mr{der}}\isomto\mb{G}_1^{\mr{der}}$ and applying the
quotient construction to both sides; see
\cite[Corollary 3.8.6 and (4.6.11)--(4.6.12)]{KisinModp}.

Note that the projections $\ms{S}_2\to\ms{S}_1$ and $\ms{S}_2\to\ms{S}_T$ commute with the relevant Langlands--Rapoport-$\tau$ bijections, in the obvious sense. The first by definition, and the second follows from the component identifications in \cite[Propositions 3.6.7 and 3.6.10]{KisinModp}.

Now, fix integral trivializations of the two underlying
$\mc{G}_2^c$-torsors at $x_2$, and let $b_2$ and $b_2'$ be the
Langlands--Rapoport and canonical Frobenius factors. In the induced trivializations at $x_1$, the integral isomorphism from the canonical crystal to the Langlands--Rapoport crystal supplied by the proof of Proposition \ref{prop:Hodge-LR-realization-compatibility} is represented by an element $a_1$ of $\mc{G}_1(W_\infty)$. 

On the other hand, the Hodge and toral comparisons are chosen from the same
special-point comparison of \cite[\S5.7.4]{KSZ}, so their images
on the common abelianized factor agree, by \cite[\S4.5.1 and Proposition 4.5.2]{KSZ} and \cite[proof of Theorem 3.5.1]{LoveringFCrystals}.
This agreement is preserved by the chosen twists, as checked in
\cite[Theorem 5.12.5]{KSZ} and part (III) of the proof of
\cite[Proposition 5.11.13]{KSZ}. 

From this discussion we deduce $\pi(a_1)=\rho(t)$ for the toral comparison
element $t$ from \textbf{Step 1}. Since $\pi$ is defined over $\Z_p$, this equality implies the following containments/equalities:
\begin{equation*}
 \rho(t)\in\mc{G}_1^{\mr{ab}}(W_\infty),\qquad
 \rho_*(\nu)=0,\qquad \rho(t^\circ)=\pi(a_1),
\end{equation*}
where $\nu$ and $t^\circ$ are those from \textbf{Step 1}. Consequently,
\begin{equation*}
 a_2\defeq(a_1,t^\circ)\in\mc{G}_2^c(W_\infty),\qquad
 b_2'=a_2^{-1}b_2\sigma(a_2).
\end{equation*}
The Frobenius identity and compatibility with the identifying
maps in \eqref{eq:tube-comp-descent-datum} can be checked after
applying the closed immersion $\alpha^c$, where they are exactly
the Hodge and toral comparisons. 

\medskip

\paragraph*{Step 3: the case of $x$} For $(\mb{G},\mb{X})$, the bijection fixed at the beginning of the
proof is obtained from the same connected Hodge-type bijection
by \cite[Corollary 3.8.12 and (4.6.11)--(4.6.12)]{KisinModp},
as used in \cite[proof of Theorem 6.2.4]{KSZ}.
This construction is compatible with the map $\beta$. Let us fix integral trivializations (induced by those at $x_2$) at $x$ with Langlands--Rapoport and canonical Frobenii $b$ and $b'$, respectively.

Then, by construction, one has $b=\beta^c(b_2)$, and by \cite[(2.4.1)]{IKY2} one has $b'=\beta^c(b_2')$. Thus, we may take the comparing element $a\defeq\beta^c(a_2)\in\mc{G}^c(W_\infty)$. 
The Hodge comparisons intertwine the quotient actions, including
the adjoint action, by \cite[proofs of Propositions 4.4.14 and 4.4.17]{KisinModp} and \cite[proof of Theorem 6.2.4]{KSZ}; and the toral comparisons also intertwine by direct inspection. The corresponding canonical-model compatibility is proved in
\cite[proof of Theorem 3.5.1]{LoveringFCrystals}.
Hence the comparisons descend to the image of $\beta_{\ms{S}}$, compatibly with Frobenius and the identifying maps in \eqref{eq:tube-comp-descent-datum}.

\medskip 

\paragraph*{Step 4: Hecke translation} Finally, the Hecke translates of
the image of $\beta_\ms{S}$ cover the target as the prime-to-$p$ level varies;
see \cite[Lemmas 2.9 and 2.4]{IKY2}. The syntomic realization is Hecke-equivariant by \eqref{eq:Hecke-action-local-system-compat-syn} and \eqref{eq:delta-Hecke-compatability}. On the Langlands--Rapoport side, $h^p$ in $\mb{G}(\mathbb{A}_f^p)$
acts by $(z_p,z^p)\mapsto(z_p,z^ph^p)$, leaving $z_p$ unchanged, and
the bijection is equivariant; see \cite[\S\S2.4.7, 2.7.1, and Theorem 6.3.5]{KSZ}.
Thus the comparisons extend equivariantly to all components.
\end{proof}

We then obtain the following unconditional version of the Langlands--Kottwitz--Scholze formula in the abelian-type setting.

\begin{cor}\label{cor:Langlands-Kottwitz-Scholze-abelian-type}
Suppose that $(\mb{G},\mb{X})$ is of abelian type. Then there exists an integer $m_0\geqslant 1$, depending only on $f^p$, such that for every $m\geqslant m_0$ and $\tau$ in $\Phi_v^mI_E$:
\begin{equation*}
T(\tau,h,f^p)=
\sum_{\gamma_0\in\Sigma_{\mathsf{K}^p}}\sum_{\substack{c\in\mf{KP}(\gamma_0)\cap\mf{KP}_a(p^n)\\ \alpha(c)=0}}
\ov{\iota}_{\mb{G}}(\gamma_0)^{-1}
 c_2(\gamma_0)\tr(\xi(\gamma_0))
 c_1(c,\mathsf{K}^p)
 \mr{O}_\gamma(f^p)
 \mr{TO}_\delta\left(\phi_{\tau,h}^{\mb{G},\mb{X}}\right),
\end{equation*}
If $\mb{G}^{\mr{ad}}$ is $\Q$-anisotropic, then one may take $m_0=1$.
\end{cor}
\begin{proof} Given the above discussion, the only things to remark are: (a) by Remark \ref{rem:l-adically-contractible-boundary} we automatically have $\ell$-adic contractibility of the boundary; and (b) if $\mb{G}^{\mr{ad}}$ is $\Q$-anisotropic, then the integral canonical model is proper; see \cite[Corollary 4.51]{Wu} and \cite[Corollary 4.1.7]{MadTorHod}. 
\end{proof}

\subsection{Stabilization}
\label{ss:geometric-stabilization}

We now give a stabilized form of the Langlands--Kottwitz--Scholze formula. This essentially follows from closely imitating \cite[\S8]{KSZ}, and so we discuss the setup and proof fairly briefly. We strongly encourage the reader to refer to loc.\@ cit.\@ for details, and to see \S\ref{s:transfer} for some of the relevant harmonic analysis (although only locally at $p$).

\begin{setup}\label{setup:geometric-stabilization}
Retain Setup \ref{setup:LR-general}, and fix:
\begin{itemize}
\item $f^p$ and $h$ as in \eqref{eq:weighted-KSZ-trace}, with $f^p$
bi-$\mathsf{K}^p$-invariant and $\mathsf{K}_0\mathsf{K}^p$ neat;
\item an irreducible algebraic representation $\xi$ of
$\mb{G}^c_{\ov{\Q}_\ell}$ and an isomorphism
$\iota_\ell\colon\ov{\Q}_\ell\isomto\C$;
\item $m\geqslant1$ and $\tau$ in $\Phi_v^mI_E$.
\end{itemize}
Put $n=mr$ and $\phi=\phi_{\tau,h}^{\mb{G},\mb{X}}$, and use $W_m$ and $K_m$
from Notation \ref{nota:local-field-Km}. 
\end{setup}

\begin{nota}\label{nota:stabilization-measures}
Set:
\begin{itemize}
\item $U_0=\mb{Z}(\A_f)\cap\mathsf{K}_0\mathsf{K}^p=U_{0,p}U_0^p$;
\item $\mf{X}_0=U_0\mb{Z}(\R)$.
\end{itemize}
We use the adelic quotient measures of
\cite[\S8.1.2]{KSZ}, with $U_0$ of volume $1$, and write
$\tau_{\mf{X}_0}$ for the central-quotient Tamagawa volume as defined in loc.\@ cit.
\end{nota}

\begin{nota}\label{nota:stabilization-coefficient}
For a $p^n$-admissible Kottwitz parameter $c=(\gamma_0;a,[b])$ with
$\alpha(c)=0$, use $I_0=I_{\gamma_0}$ of Notation
\ref{nota:Kottwitz-centralizer} and the real inner form $I_\infty$ of
Construction \ref{constr:inner-form-infinity}. Put
\begin{equation*}
c_1^{\mr{KSZ}}(c,\mathsf{K}^p)
\defeq\frac{\tau_{\mf{X}_0}(I_0)}
{\mr{vol}(\mb{Z}(\R)\backslash I_\infty(\R))},
\end{equation*}
with centralizer measures as in \cite[\S8.1.2 and Lemma 8.1.3]{KSZ}.
\end{nota}

\begin{rem}\label{rem:stabilization-normalization}
One has an identification of $c_1$ of Definition \ref{defn:KSZ-c1} with $c_1^{\mr{KSZ}}$ in the normalization of the Haar measures discussed above.
\end{rem}

We next specify the endoscopic transfers of functions needed for this stabilization.

\begin{nota}\label{nota:stabilization-transfers}
Let $\mc{E}_{\mr{ell}}(\mb{G})$ be representatives for the elliptic
endoscopic data of \cite[Definition 7.2.2]{KSZ}.
Choose the following data:
\begin{itemize}
\item a global $z$-extension $\mb{G}_1\to\mb{G}$ unramified at $p$,
and its associated extensions $\mb{H}_1\to\mb{H}$ with common induced
torus kernel $\mb{Z}_1$, as in \cite[Lemmas 7.2.5--7.2.6]{KSZ};
\item the transfer factors of \cite[\S8.2.1]{KSZ};
\item ordinary transfers $f^{\mf{e},p}$ of $f^p$ and twisted transfers
$\phi^{\mf{e}}$ of $\phi$; see \S\ref{sss:twisted-transfer} for the latter.
\end{itemize}
\end{nota}

\begin{nota}\label{nota:stabilization-central-character}
Let $\mb{Z}_{\mf{e}}$ be the inverse image of
$\mb{Z}\subseteq Z(\mb{H})$ in $\mb{H}_1$, with projection
$q_{\mf{e}}\colon\mb{Z}_{\mf{e}}\to\mb{Z}$, and let
$\lambda_{\mf{e}}$ be the transfer character of
\cite[(8.2.1.1)]{KSZ}. The finite transfer functions of Notation
\ref{nota:stabilization-transfers} transform under $\mb{Z}_1(\A_f)$
by $\lambda_{\mf{e}}|_{\mb{Z}_1(\A_f)}$.
For a compact open subgroup $U_p\subseteq U_{0,p}$, set:
\begin{itemize}
\item $\mf{X}=U_pU_0^p\mb{Z}(\R)$;
\item $\chi|_{U_pU_0^p}=1$ and
$\chi|_{\mb{Z}(\R)}=\omega_\xi^{-1}$, where $\omega_\xi$ is the
central character of $\xi$;
\item $\mf{X}_{\mf{e}}=q_{\mf{e}}^{-1}(\mf{X})$ and
$\Omega_{\mf{e}}=(\chi\circ q_{\mf{e}})\lambda_{\mf{e}}^{-1}$.
\end{itemize}
We will furthermore use the quotient measures of \cite[\S8.2.1]{KSZ}.
\end{nota}

\begin{rem}\label{rem:stabilization-central-character-convention}
Below we will use $\Omega_{\mf{e}}$ as a subscript on a trace distribution to specify the character of the representations. Thus the test functions we will consider satisfy
\begin{equation*}
F(zx)=\Omega_{\mf{e}}(z)^{-1}F(x),
\qquad z\in\mf{X}_{\mf{e}},\quad x\in\mb{H}_1(\A);
\end{equation*}
see \cite[\S7.1.4]{KSZ}. At $p$, with the
compatible $z$-pair used below, this gives
\begin{equation*}
\chi_1=\lambda_{\mf{e},p}|_{Z_1(\Q_p)},
\qquad
\Omega_{\mf{e},p}|_{Z_1(\Q_p)}=\chi_1^{-1};
\end{equation*}
see \cite[\S7.4.7]{KSZ}. Note this datum gives a character on
$\mf{X}_{\mf{e}}$, not necessarily on all of
$Z(\mb{H}_1)(\A)$.
\end{rem}

Using the following elementary lemma, we may choose a $U_p$ so that the additional central averaging preserves the transferred orbital integrals. 

\begin{lem}\label{lem:stabilization-central-data}
Set $\mc{A}=\ker(\mc{G}\to\mc{G}^c)$. Then, there is a compact open
$U_p\subseteq U_{0,p}$ containing $\mc{A}(\Z_p)$ for which the
finite transfers, together with the archimedean factor of
\cite[\S8.2.5]{KSZ}, may be chosen to have transformation character $\Omega_{\mf{e}}^{-1}$.
\end{lem}

\begin{nota}\label{nota:stabilization-global-function}
With this central data established, set
\begin{equation}\label{eq:stabilization-test-function}
f_{\tau,h,\xi}^{\mf{e}}\defeq
f^{\mf{e},p}\phi^{\mf{e}}f_\xi^{\mf{e}},
\end{equation}
where $f_\xi^{\mf{e}}$ is the archimedean factor of \cite[\S8.2.5]{KSZ}. We also set
\begin{equation*}
\iota(\mb{G},\mf{e})\defeq
\tau(\mb{G})\tau(\mb{H})^{-1}
|\operatorname{Out}(\mf{e})|^{-1}.
\end{equation*}
Here $\tau(\mb{L})$ denotes Tamagawa number of $\mb{L}$, and $\mr{Out}(\mf{e})$ is the automorphism image as in \cite[\S7.2.3]{KSZ}; compare \cite[(8.3.7.1)]{KSZ}.
\end{nota}

\begin{defn}[{\cite[(8.3.7.2)]{KSZ}}]
\label{defn:stabilization-stable-elliptic-distribution}
The stable elliptic distribution is
\begin{equation*}
\mr{ST}_{\mr{ell},\Omega_{\mf{e}}}^{\mb{H}_1}(F)
\defeq
\tau_{\mf{X}_{\mf{e}}}(\mb{H}_1)
\sum_{\gamma\in\Sigma_{\mr{ell},\mf{X}_{\mf{e}}}(\mb{H}_1)}
|\operatorname{Stab}_{\mf{X}_{\mf{e}}}(\gamma)|^{-1}
\mr{SO}_\gamma(F),
\end{equation*}
where $\Sigma_{\mr{ell},\mf{X}_{\mf{e}}}(\mb{H}_1)$ is the set of
stable conjugacy classes of elliptic semisimple elements of
$\mb{H}_1(\Q)$, modulo multiplication by
$\mf{X}_{\mf{e},\Q}\defeq\mf{X}_{\mf{e}}\cap\mb{Z}_{\mf{e}}(\Q)$,
and $\mr{Stab}_{\mf{X}_{\mf{e}}}(\gamma)$ consists of those
$z$ in $\mf{X}_{\mf{e},\Q}$ for which $z\gamma$ is stably conjugate to
$\gamma$, before taking this quotient. Note that the factor
$\tau_{\mf{X}_{\mf{e}}}(\mb{H}_1)$ makes use of the central quotient measures
of \cite[\S8.1.2 and Lemma 8.2.2]{KSZ}.
\end{defn}

\begin{obs}\label{obs:stabilization-central-refinement}
With $m$ and $\mathsf{K}^p$ fixed, shrinking $U_p$ to an open compact subgroup containing $\mc{A}(\Z_p)$ does not change the stable elliptic values. 
\end{obs}

We next observe the support condition needed in imitating the stabilization argument from \cite[\S8]{KSZ} below.

\begin{obs}\label{obs:stabilization-support}
One has $|\phi|\leqslant B\phi_n^{\mr{KSZ}}$ for some $B\geqslant0$.
Consequently nonvanishing of $\mr{TO}_\delta(\phi)$ implies
nonvanishing of $\mr{TO}_\delta(\phi_n^{\mr{KSZ}})$.
\end{obs}

Finally, we can state and prove our stabilized trace formula.

\begin{thm}\label{thm:geometric-stabilization}
Assume \eqref{eq:Langlands-Kottwitz-Scholze-formula} with the normalization
of \emph{Remark \ref{rem:stabilization-normalization}}. Then
\begin{equation}\label{eq:geometric-stabilization}
T(\tau,h,f^p)=
\sum_{\mf{e}\in\mc{E}_{\mr{ell}}(\mb{G})}
\iota(\mb{G},\mf{e})
\mr{ST}_{\mr{ell},\Omega_{\mf{e}}}^{\mb{H}_1}
(f_{\tau,h,\xi}^{\mf{e}}).
\end{equation}
\end{thm}
\begin{proof} This essentially follows by copying \cite[Theorem 8.3.10]{KSZ}; we now elaborate on this.

First, using \cite[Lemma 8.1.3]{KSZ} and character orthogonality on $\fE(I_0,\mb{G};\A/\Q)$ we may rewrite the Langlands--Kottwitz--Scholze formula as in \cite[(8.1.4.2)]{KSZ}, with $\phi$ replacing $\phi_n$;
note that Observation \ref{obs:stabilization-support} verifies the Kottwitz support
condition as in loc.\@ cit. We may then apply the twisted transfer of \cite[Proposition 7.4.13]{KSZ} at $p$, ordinary transfer away from $p$, and \cite[(8.2.5.2)]{KSZ} at $\infty$. Note that, as in the proofs of
\cite[Lemmas 8.2.8 and 8.3.8]{KSZ}, only relatively regular classes contribute here, and moreover their centralizers in $\mb{H}_1$ are connected. The global image and fiber calculations of \cite[Lemma 8.3.6 and proof of Theorem 8.3.10]{KSZ} now give \eqref{eq:geometric-stabilization}, where for the central measure comparison needed in this argument we make use of Observation \ref{obs:stabilization-central-refinement}. 
\end{proof}

\section{The Scholze--Shin conjecture and applications to Shimura varieties}\label{s:Scholze-Shin}

In this final section we describe the Scholze--Shin conjecture  relating the local Langlands correspondence to the local test functions $\phi_{\tau,h}$; see \cite[Conjectures 7.1 and 7.2]{ScholzeShin}. We then use this to give conditional decompositions of the cohomology of Shimura varieties and formulas for the semisimplified partial Hasse--Weil $\zeta$-functions.

The real major updates to loc.\@ cit.\@ are: (a) the existence of the functions $\phi_{\tau,h}$ allows one to state this conjecture in full generality (at least for unramified groups); (b) to accompany this maximum generality we must use a slightly more robust theory of (twisted) endoscopy as explained in \cite[\S7]{KSZ}; and (c) we update the precise statement to not require dependence on some hypothetical local Langlands correspondence, using instead the work of Fargues--Scholze.

\begin{rem} One will see that we have devoted considerable space to making the statement of the Scholze--Shin conjecture in the generality considered here as self-contained, both mathematically and expositionally, as possible. This is done in contrast to something like \S\ref{ss:KSZ-point-counting} or \S\ref{ss:geometric-stabilization}. In fact, we have given more self-contained treatment here even to some concepts already used in the above sections.

Besides the issues of feasibility (the amount of setup in these other sections is too great), this reflects the following: Unlike the above sections, where we feel as though it is the theorems themselves which are valuable, the real content of the discussion below is in its setup.
\end{rem}

\subsection{Some preliminaries on (twisted) endoscopic transfer}\label{s:transfer}

The Scholze--Shin conjecture relies, as part of its statement, on the (twisted) endoscopic transfer of functions. In this subsection we discuss these constructions in the generality needed below. 

\begin{rem} This subsection is an extension of \cite[\S3.2]{ScholzeShin}. The main difference is the removal of the simplifying assumption that $G^\mr{der}$ is simply connected made in loc.\@ cit.\@ using \cite[\S7]{KSZ}.
\end{rem}

\begin{nota}\label{nota:endoscopic-transfer} Throughout we use Notations \ref{nota:deformation-spaces-setup} and \ref{nota:finite-field}, as well as the following:
\begin{itemize}
\item $1\leqslant j<\infty$;
\item $d\defeq[K_j:\Q_p]=rj$;
\item $\mf{e}=(H,\mc{H},s,\eta)$ is an endoscopic datum for $G$ in the sense
of \cite[Definition 7.2.2]{KSZ};
\item using \cite[Lemma 7.2.5]{KSZ}, we choose an unramified $z$-extension
\begin{equation}\label{eq:z-extension-G-endoscopic-transfer}
    \mf{z}\colon\quad 1\to Z_1\to G_1\to G\to 1
\end{equation}
with $Z_1$ an induced torus, $G_1^\der=G_1^\sc$ and, as in the proof of
\cite[Proposition 7.4.13]{KSZ} and \cite[Lemma 7.2.13]{KSZ}, the required norm maps on its central torus are surjective.
\end{itemize}
\end{nota}

\begin{rem}\label{rem:normalization-transfer}
Below we use the normalization of transfer factors from \cite{KSZ}. The
normalization in \cite{ScholzeShin} differs by the convention explained
in \cite[Remark 3.2]{ScholzeShin}. In particular, the analogue of the
pairing in \eqref{eq:characterization-twisted-endoscopic-transfer} is
written with an inverse in loc.\@ cit.
\end{rem}

\subsubsection{$L$-embeddings associated to $\mf{z}$}

We first recall that our fixed $z$-extension $\mf{z}$ allows us to replace $\mf{e}$ by an
endoscopic datum represented by an honest $L$-embedding.

\begin{obsnota}
\label{obsnota:z-pair-endoscopic-transfer}
Using \cite[Lemma 7.2.6 and Lemma 7.2.9]{KSZ}, we may associate with $\mf{e}$ and $\mf{z}$ a central extension
\begin{equation}\label{eq:z-extension-H-endoscopic-transfer}
    1\to Z_1\to H_1\to H\to 1,
\end{equation}
an $L$-morphism $\zeta_{H_1}\colon\mc{H}\to{}^L\!H_1$, and an $L$-embedding $\eta_1\colon{}^L\!H_1\to{}^L\!G_1$, compatible with $\eta$. If $s_1$ denotes the image of $s$ in $\wh{G}_1$,
then $\mf{e}_1=(H_1,{}^L\!H_1,s_1,\eta_1)$ is an endoscopic datum for $G_1$.
\end{obsnota}

\begin{constr}The content of Observation/Notation \ref{obsnota:z-pair-endoscopic-transfer} also produces a continuous character $\chi_1\colon Z_1(\Q_p)\to\C^\times$; choose a splitting $W_{\Q_p}\to\mc{H}$ and apply the classically normalized Langlands correspondence for tori\footnote{By the classical normalization we mean the one that sends uniformizers to arithmetic Frobenius, i.e., $\chi_1^{-1}$ corresponds to the displayed parameter under the \emph{geometric} normalization.} to the composition:\begin{equation*}
    W_{\Q_p}\to\mc{H}
    \xrightarrow{\zeta_{H_1}}{}^L\!H_1
    \to{}^L\!Z_1.
\end{equation*}
This character is independent of
the chosen splitting by \cite[\S7.2.10]{KSZ}. 
\end{constr}

\begin{nota}\label{nota:fixed-central-character-Hecke-algebra}
We write $\mc{H}_{\chi_1}(H_1(\Q_p))$ for the space of locally constant
functions $f\colon H_1(\Q_p)\to\C$ which are compactly supported modulo
$Z_1(\Q_p)$ and satisfy
\begin{equation*}
    f(zx)=\chi_1(z)f(x),
    \qquad z\in Z_1(\Q_p),\quad x\in H_1(\Q_p).
\end{equation*}
\end{nota}

\subsubsection{Endoscopic transfer}

With notation as in Observation/Notation \ref{obsnota:z-pair-endoscopic-transfer}, we first discuss endoscopic transfer of functions on $G(\Q_p)$ to functions on $H_1(\Q_p)$. 

\begin{constr}
\label{constr:ordinary-endoscopic-transfer}
Let $h$ be an element of $\mc{H}(G(\Q_p))$, and consider the natural averaging map
\begin{equation}\label{eq:ordinary-averaging-G1-to-G}
\begin{aligned}
    \mr{Av}_{Z_1(\Q_p)}\colon
    \mc{H}(G_1(\Q_p))
    &\to \mc{H}(G(\Q_p)),\\
    h_1&\mapsto
    \left(
        g\mapsto
        \int_{Z_1(\Q_p)}
        h_1(zg_1)\,dz
    \right),
\end{aligned}
\end{equation}
where $g_1$ in $G_1(\Q_p)$ is any lift of $g$ (which exists since $Z_1$ is an induced torus). This map is surjective, and
we choose $h_1$ in $\mc{H}(G_1(\Q_p))$ whose image under \eqref{eq:ordinary-averaging-G1-to-G} is $h$.

Choose an ordinary $\mf{e}_1$-transfer $h_1^{\mf{e}_1}$ in $\mc{H}(H_1(\Q_p))$ of $h_1$, whose existence follows from endoscopic transfer for
the datum $\mf{e}_1$; see \cite[Proposition 7.4.9]{KSZ}. We then obtain a function on $H_1(\Q_p)$ with central character $\chi_1$ by the following definition:
\begin{equation}\label{eq:ordinary-central-character-averaging}
    h^{\mf{e}}(x)
    \defeq
    \int_{Z_1(\Q_p)}
    \chi_1(z)^{-1}h_1^{\mf{e}_1}(zx)\,
    dz.
\end{equation}
By construction, $h^{\mf{e}}$ belongs to $\mc{H}_{\chi_1}(H_1(\Q_p))$.
\end{constr}

\begin{defn}\label{defn:ordinary-endoscopic-transfer}
With notation as above, we call any function $h^{\mf{e}}$ in $\mc{H}_{\chi_1}(H_1(\Q_p))$ obtained from $h$ by Construction \ref{constr:ordinary-endoscopic-transfer} an \emph{ordinary $\mf{e}$-transfer of $h$}.
\end{defn}

The individual function $h^{\mf{e}}$ is not canonical. That said, at least after fixing the $z$-pair, Haar measures, and transfer factors, its stable orbital integrals are (largely) canonical.

To state this claim precisely, we first recall some standard definitions and notation.

\begin{defn}[{see {\cite[\S7.4.1 and \S7.4.3]{KSZ}.}}] \label{defn:stable-conjugacy} Let $L$ be a reductive $\Q_p$-group. Two semisimple elements $\gamma$ and $\gamma'$ of $L(\bb{Q}_p)$ are \emph{stably conjugate} if there exists $x$ in $L(\ov{\Q}_p)$ such that $\gamma'=x^{-1}\gamma x$ and the class of the cocycle \begin{equation*} \Gamma_{\Q_p}\to L_\gamma(\ov{\Q}_p), \qquad \rho\mapsto x\rho(x)^{-1}, \end{equation*} belongs to the image of $H^1(\Q_p,I_\gamma)\to H^1(\Q_p,L_\gamma)$.
\end{defn}

\begin{nota}
 We write 
 \begin{itemize}
 \item $\Gamma(L)$ for the set of semisimple $L(\Q_p)$-conjugacy classes in $L(\Q_p)$;
 \item $\Sigma(L)$ for the set of stable semisimple conjugacy classes in $L(\Q_p)$.
 \end{itemize}
 For $\gamma$ in $L(\bb{Q}_p)$ semisimple, we write $[\gamma]$ and $[\gamma]_s$ for its image in $\Gamma(L)$ and $\Sigma(L)$, respectively.
\end{nota}

The following definitions recall the notion of transfers of stable conjugacy classes from $H_1$ to $G_1$ (and thus to $G$), which is necessary to state the formula for the stable orbital integrals of $h^\mf{e}$.

\begin{defn}[{\cite[\S7.4.2]{KSZ}}]
\label{defn:relative-regularity}
Let $\gamma_{H_1}$ in $H_1(\Q_p)$ be a semisimple element, and choose a maximal torus $T_{H_1}\subseteq H_1$ containing $\gamma_{H_1}$. We obtain a canonical
Galois-stable $G_1(\ov{\Q}_p)$-conjugacy class of embeddings $\j\colon (T_{H_1})_{\ov{\Q}_p}\to (G_1)_{\ov{\Q}_p}$ obtained by dualizing the inclusion of dual tori
\begin{equation*}
    \wh{T}_{H_1}\subseteq\wh{H}_1
    \xrightarrow{\eta_1}\wh{G}_1;
\end{equation*}
see \cite[\S7.4.2]{KSZ}. Via $\j$, we regard
$\Phi(T_{H_1},H_1)$ as a subset of
$\Phi(\jmath(T_{H_1}),G_1)$. We say that $\gamma_{H_1}$ is
\emph{$(G_1,H_1)$-regular} if $\alpha(\gamma_{H_1})\ne 1$ for every $\alpha$ in $\Phi(\jmath(T_{H_1}),G_1)$ but not in $\Phi(T_{H_1},H_1)$. 
\end{defn}

\begin{nota} Let $\Sigma(H_1)_{(G_1,H_1)\mr{-reg}}$ be the subset of $(G_1,H_1)$-regular elements of $\Sigma(H_1)$.
\end{nota}

We can now precisely formulate the notion of transfer alluded to above.

\begin{defn}[{\cite[\S7.4.2]{KSZ}}]
Let $\gamma_{H_1}$ in $H_1(\Q_p)$ be $(G_1,H_1)$-regular. Then, with notation as in Definition \ref{defn:relative-regularity}, the element $\jmath(\gamma_{H_1})$ determines a semisimple Galois-stable $G_1(\ov{\Q}_p)$-conjugacy class which is independent of choices. Since $G_1$ is quasi-split and $G_1^\der=G_1^\sc$, this geometric conjugacy class contains an element of $G_1(\Q_p)$; see \cite[Theorem 4.4]{KottwitzConjugacy}. We thus obtain a map
\begin{equation*}
\mr{transf}_{(G_1,H_1)}\colon \Sigma(H_1)_{(G_1,H_1)\mr{-reg}}\to \Sigma(G_1),
\end{equation*}
and by composing it with the natural map $\Sigma(G_1)\to\Sigma(G)$ a natural map
\begin{equation*}
\mr{transf}_{(G,H_1)}\colon \Sigma(H_1)_{(G_1,H_1)\mr{-reg}}\to \Sigma(G).
\end{equation*}
\end{defn}

Finally, we give notation that will allow us to talk about (a slight refinement of) the rational conjugacy classes within a given stable conjugacy class.

\begin{nota}\label{nota:gamma-x} Let $L$ be a reductive $\Q_p$-group and let $\gamma$ in $L(\Q_p)$ be semisimple. Set 
\begin{equation*} 
\mf{D}(I_\gamma,L;\Q_p) \defeq \ker\left( H^1(\Q_p,I_\gamma)\to H^1(\Q_p,L) \right). \end{equation*} 
For $[x]$ in $\mf{D}(I_\gamma,L;\Q_p)$, choose $x$ in $L(\ov{\Q}_p)$ such that $\rho\mapsto x\rho(x)^{-1}$ is an $I_\gamma$-valued cocycle representing $[x]$. As $[x]$ maps trivially to $H^1(\Q_p,L)$ one may choose $x$ so that $\gamma_x\defeq x^{-1}\gamma x$ belongs to $L(\Q_p)$. We write $\gamma[x]$ for the element $[\gamma_x]$ of $\Gamma(L)$. The resulting map
\begin{equation*} 
\mf{D}(I_\gamma,L;\Q_p)\to\Gamma(L), \qquad [x]\mapsto\gamma[x], 
\end{equation*} 
has image equal to the set of conjugacy classes contained in $[\gamma]_s$. 
\end{nota}

We give our last piece of notational setup by setting our notation for (stable) orbital integrals.

\begin{defn}\label{defn:orbital-integrals}
Let $L$ be a reductive $\bb{Q}_p$-group and let
$\gamma$ in $L(\Q_p)$ be semisimple. For $f$ in $\mc{H}(L(\Q_p))$, we define the \emph{orbital integral}
\begin{equation*}
    \mr{O}_\gamma(f)
    \defeq
    \int_{I_\gamma(\Q_p)\backslash L(\Q_p)}
    f(g^{-1}\gamma g)\,
    \frac{dg}{di},
\end{equation*}
where the measure on the quotient is the one induced by the fixed Haar measures. Evidently $\mr{O}_\gamma(f)$ only depends on $[\gamma]$ and thus the notation $\mr{O}_{c}(f)$ makes sense for $c$ in $\Gamma(L)$. 

Similarly, for such a $\gamma$ we define the \emph{stable orbital integral} by
\begin{equation*}
    \mr{SO}_\gamma(f)
    \defeq
    \sum_{[x]\in\mf{D}(I_\gamma,L;\Q_p)}
    e(I_{\gamma[x]})\,\mr{O}_{\gamma[x]}(f),
\end{equation*}
where $e(I_{\gamma[x]})$ is the Kottwitz sign; see \cite{KottwitzSign}. It's clear that $\mr{SO}_\gamma(f)$ only depends on $[\gamma]_s$, and thus the notation $\mr{SO}_c(f)$ is well-defined for $c$ in $\Sigma(L)$.
\end{defn}

We finally come to the following, which says that while the function $h^\mf{e}$ is not unique, its stable orbital integrals at $(G_1,H_1)$-regular stable conjugacy classes are uniquely determined. It is only these stable orbital integrals that will be relevant for the Scholze--Shin conjecture.

\begin{prop}[{\cite[Proposition 7.4.11]{KSZ}}]
\label{prop:characterization-ordinary-endoscopic-transfer}
Let $c$ be an element of $\Sigma(H_1)_{(G_1,H_1)\mr{-reg}}$. Then, 
\begin{equation}\label{eq:ordinary-transfer-identity}
    \displaystyle \mr{SO}_{c}(h^{\mf{e}}) =\sum_{[x]\in\fD(I_\gamma,G;\Q_p)}
    e(I_{\gamma[x]})
    \Delta^\mf{e}(\gamma_{H_1},\gamma[x])
    \mr{O}_{\gamma[x]}(h) \qquad \emph{where }[\gamma]_s=\mr{transf}_{(G,H_1)}(c)
\end{equation}
where $\Delta^\mf{e}(c,\gamma[x])$ is the normalized
endoscopic transfer factor \emph{(}see \cite[\S7.4.4 and \S7.4.7]{KSZ}\emph{)}.
\end{prop}

\begin{proof}
Choose $h_1$ and $h_1^{\mf{e}_1}$ as in Construction
\ref{constr:ordinary-endoscopic-transfer}. By \cite[Proposition 7.4.11]{KSZ} one has that the analogous equality as in \eqref{eq:ordinary-transfer-identity} holds for $h_1$ and $h_1^{\mf{e}_1}$. Given the definition of $h^\mf{e}$, the claim then follows by averaging along $Z_1(\Q_p)$.
\end{proof}

\subsubsection{Twisted endoscopic transfer}\label{sss:twisted-transfer}

We now discuss the case of twisted endoscopic transfer of functions on $G(K_j)$ to functions on $H_1(\Q_p)$. 

Following \cite{KSZ}, we initially assume, for the sake of readability, that $s$ itself is $\Gamma_{\Q_p}$-invariant. We will then reduce from the general case to this case.

\begin{constr}
\label{constr:twisted-base-change-datum}
Assume that $s$ is in $Z(\wh{H})^{\Gamma_{\Q_p}}$. Set $R_{d,1}\defeq\Res_{K_j/\Q_p}G_{1,K_j},$ and let $\theta_1$ be the
automorphism of $R_{d,1}$ induced by the arithmetic Frobenius of
$K_j/\Q_p$. There is an $L$-map
$\iota_d\colon{}^L\!G_1\to{}^L\!R_{d,1}$ which, via $\wh{R}_{d,1}=\wh{G}_1^{\,d}$, is the diagonal embedding on dual groups.

Set $\wt{s}_1=(s_1,1,\ldots,1)$, and let
\begin{equation*}
 a\colon W_{\Q_p}\to
 \mr{Cent}\bigl((\iota_d\circ\eta_1)(\wh{H}_1),\wh{R}_{d,1}\bigr)
\end{equation*}
be the unramified $1$-cocycle for the action induced by
$\iota_d\circ\eta_1$, characterized by
$a(\mr{Frob}_p)=\wt{s}_1$, where $\mr{Frob}_p$ is an arithmetic
Frobenius lift. Writing $\wt{\eta}_1=a\cdot(\iota_d\circ\eta_1)$, the quadruple
$\wt{\mf{e}}_1=(H_1,{}^L\!H_1,\wt{s}_1,\wt{\eta}_1)$
is a twisted endoscopic datum for $(R_{d,1},\theta_1)$;
see \cite[\S7.2.14]{KSZ}.
\end{constr}

\begin{rem}
Construction \ref{constr:twisted-base-change-datum} is the $z$-pair
version of \cite[\S3.2, equation (6)]{ScholzeShin}.
\end{rem}

We next recall the notion of twisted endoscopic transfer of functions
associated with $\wt{\mf{e}}_1$.

\begin{constr}
\label{constr:twisted-endoscopic-transfer}
Let $\phi$ be an element of $\mc{H}(G(K_j))$, and consider the natural
surjective averaging map
\begin{equation}\label{eq:averaging-G1-to-G-endoscopic-transfer}
\begin{aligned}
    \mr{Av}_{Z_1(K_j)}\colon
    \mc{H}(G_1(K_j))
    &\to \mc{H}(G(K_j)),\\
    \phi_1&\mapsto
    \left(
        g\mapsto
        \int_{Z_1(K_j)}
        \phi_1(zg_1)\,
        dz,
    \right),
\end{aligned}
\end{equation}
where $g_1$ in $G_1(K_j)$ is any lift of $g$. Choose $\phi_1$ in
$\mc{H}(G_1(K_j))$ whose image under
\eqref{eq:averaging-G1-to-G-endoscopic-transfer} is $\phi$.

Choose a twisted $\wt{\mf{e}}_1$-transfer
$\phi_1^{\wt{\mf{e}}_1}$ in $\mc{H}(H_1(\Q_p))$ of $\phi_1$; see \cite[\S7.4.12]{KSZ}. We then
impose the central character $\chi_1$ by setting
\begin{equation}\label{eq:twisted-central-character-averaging}
    \phi^\mf{e}(x)
    \defeq
    \int_{Z_1(\Q_p)}
    \chi_1(z)^{-1}
    \phi_1^{\wt{\mf{e}}_1}(zx)\,
    dz.
\end{equation}
By construction, $\phi^{\mf{e}}$ belongs to
$\mc{H}_{\chi_1}(H_1(\Q_p))$.
\end{constr}

\begin{defn}\label{defn:twisted-endoscopic-transfer}
With notation as above, we call any function $\phi^{\mf{e}}$ in
$\mc{H}_{\chi_1}(H_1(\Q_p))$ obtained from $\phi$ by Construction
\ref{constr:twisted-endoscopic-transfer} a
\emph{twisted $\mf{e}$-transfer of $\phi$}.
\end{defn}

Again the individual function $\phi^{\mf{e}}$ is not canonical, but its stable orbital integrals at $(G_1,H_1)$-regular elements are. To state this we again must recall some further definitions and notation in this twisted setting.

\begin{defn}Let $L$ be a reductive $\Q_p$-group and $\sigma$ the arithmetic Frobenius of $K_j/\Q_p$. Two elements $\delta$ and
$\delta'$ of $L(K_j)$ are \emph{$\sigma$-conjugate} if
$\delta'=g^{-1}\delta\sigma(g)$ for some $g$ in $L(K_j)$. For an element $\delta$ of $L(K_j)$, write $L^\sigma_\delta$ for its \emph{$\sigma$-centralizer} which, for a $\bb{Q}_p$-algebra $A$, has points
\begin{equation*}
L_\delta^\sigma(A)=\left\{
g\in L(A\otimes_{\Q_p}K_j):
g^{-1}\delta\sigma(g)=\delta
\right\}.
\end{equation*}
We then set $I_\delta^\sigma\defeq (L_\delta^\sigma)^\circ$.
\end{defn}

\begin{rem} We have previously encountered $\sigma$-conjugacy in \S\ref{ss:case-of-finite-fields}, and $\sigma$-centralizers in Remark \ref{rem:sigma-centralizer}. As the context for these considerations is totally different from the one here, we have chosen to redefine things here and use distinct notation.
\end{rem}

\begin{defn}
\label{defn:twisted-orbital-integrals} For $\delta$ in $L(K_j)$ $\sigma$-semisimple,\footnote{Recall that this means that $\delta\sigma(\delta)\cdots\sigma^{d-1}(\delta)$ is semisimple.} and $\phi$ in $\mc{H}(L(K_j))$, we define its
\emph{$\sigma$-twisted orbital integral} at $\delta$ by
\begin{equation*}
    \mr{TO}_\delta(\phi)
    \defeq
    \int_{I_\delta^\sigma(\Q_p)\backslash L(K_j)}
    \phi(g^{-1}\delta\sigma(g))\,
    \frac{dg}
         {di}.
\end{equation*}
This depends only on the $\sigma$-conjugacy class of $\delta$.
\end{defn}

The twisted analogue of transfer of stable conjugacy classes is the notion of norms.

\begin{defn}[{\cite[\S5]{KottwitzConjugacy}}]
\label{defn:norm-twisted-transfer}
Let $\delta_1$ be an element of
$R_{d,1}(\Q_p)=G_1(K_j)$. A semisimple element $\gamma_{0,1}$ of
$G_1(\Q_p)$ is called a \emph{degree $d$ norm of $\delta_1$} if it is
conjugate in $G_1(\ov{\Q}_p)$ to
\begin{equation*}
N_d(\delta_1)\defeq \delta_1\sigma(\delta_1)\cdots\sigma^{d-1}(\delta_1).
\end{equation*}
This condition depends only on the $\sigma$-conjugacy class of $\delta_1$ and the stable class $[\gamma_{0,1}]_s$.

For $c$ in $\Sigma(H_1)_{(G_1,H_1)\mr{-reg}}$, we say that $c$ is a \emph{norm of $\delta_1$} if there exists a degree $d$ norm $\gamma_{0,1}$ of $\delta_1$ such that $[\gamma_{0,1}]_s=\mr{transf}_{(G_1,H_1)}(c)$. We say that $c$ is a \emph{norm from $R_{d,1}(\Q_p)$} if it is a norm of some element of $R_{d,1}(\Q_p)$.
\end{defn}

\begin{nota} \label{nota:D-r} Let $\gamma_0$ be a semisimple element of $G(\Q_p)$, and as usual let $B(I_{\gamma_0})$ denote the set of $\sigma$-conjugacy classes in $I_{\gamma_0}(K_\infty)$. We define 
\begin{equation*}
\mf{D}_d(\gamma_0,G;\Q_p)=\left\{[b]\in B(I_{\gamma_0}): \text{there exists }a\in G(K_\infty)\text{ with } a^{-1}\gamma_0a = a^{-1}N_d(b)\sigma^d(a)\right\}.
 \end{equation*}
\end{nota}

\begin{constr} For $[b]$ in $\mf{D}_d(\gamma_0,G;\Q_p)$, choose $b$ and $a$ as in the condition for membership in $\mf{D}_d(\gamma_0,G;\Q_p)$ and set $\delta\defeq a^{-1}b\sigma(a)$. By \cite[Lemma 1.6.7]{KSZ}, $\delta$ belongs to $G(K_j)$, its $\sigma$-conjugacy class is independent of the choices of $a$ and $b$, and $\gamma_0$ is a degree $d$ norm of $\delta$. We denote this $\sigma$-conjugacy class by $\delta[b]$, and abusively also use $\delta[b]$ for any representative of it.
\end{constr} 

\begin{nota} For $\gamma_0$ a semisimple element of $G(\bb{Q}_p)$ and $[b]$ an element of $\mf{D}_d(\gamma_0,G;\Q_p)$, write
\begin{equation*} 
\beta_p(\gamma_0,[b]) \defeq \kappa_{I_{\gamma_0}}([b]) \in X^\ast\big(Z(\wh{I}_{\gamma_0})^{\Gamma_{\Q_p}}\big),
 \end{equation*} 
where $\kappa_{I_{\gamma_0}}\colon B(I_{\gamma_0})\to X^\ast(Z(\widehat{I}_{\gamma_0})^{\Gamma_{\Q_p}})$ is the Kottwitz map as in \cite[\S7]{KotIsoII}.
\end{nota}
We finally come to the following, which says that while the function
$\phi^{\mf{e}}$ is not unique, its stable orbital integrals at $(G_1,H_1)$-regular stable conjugacy classes are.

\begin{prop}[{\cite[Proposition 7.4.13]{KSZ}}]
\label{prop:characterization-twisted-endoscopic-transfer} Suppose that $c$ in $\Sigma(H_1)_{(G_1,H_1)\mr{-reg}}$ is exhibited as a norm from $R_{d,1}(\Q_p)$ via the element $\gamma_{0,1}$ in $G_1(\Q_p)$ and let $\gamma_0$ be its image in $G(\Q_p)$. Then,
\begin{equation}\label{eq:characterization-twisted-endoscopic-transfer}
    \mr{SO}_{c}(\phi^\mf{e})
    =
    \Delta_0^{\mf{e}}(c,\gamma_0)
    \sum_{[b]\in\mf{D}_d(\gamma_0,G;\Q_p)}
    e(I^\sigma_{\delta[b]})
    \left\langle
        \beta_p(\gamma_0,[b]),s
    \right\rangle
    \mr{TO}_{\delta[b]}(\phi).
\end{equation}
On the other hand, if $c$ is not a norm then $\mr{SO}_c(\phi^\mf{e})=0$. Here $\Delta_0^{\mf{e}}(c,\gamma_0)$ is the normalized endoscopic transfer factor as in \emph{\cite[Proposition 7.4.13]{KSZ}}.
\end{prop}

\begin{proof}
Choose $\phi_1$ and $\phi_1^{\wt{\mf{e}}_1}$ as in Construction
\ref{constr:twisted-endoscopic-transfer}. Proposition
\cite[Proposition 7.4.13]{KSZ} gives the asserted vanishing and the corresponding orbital-integral identity upstairs on $G_1$. Moreover, \cite[Lemma 7.3.6]{KSZ} identifies $\mf{D}_d(\gamma_{0,1},G_1;\Q_p)$ with $\mf{D}_d(\gamma_0,G;\Q_p)$, and under this identification the associated $\sigma$-conjugacy classes, Kottwitz invariants, and Kottwitz signs map to those appearing above. Finally, the comparison of transfer factors in \cite[(7.4.12.1)]{KSZ} contributes the factor $\langle\beta_p(\gamma_0,[b]),s\rangle$.This establishes the claim upstairs on $G_1$, and the claim on $G$ itself follows by averaging over $Z_1(K_j)$.
\end{proof}

Finally, we recall, following \cite[\S7.4.17]{KSZ}, how to remove the assumption that $s$ belongs to $Z(\wh{H})^{\Gamma_{\Q_p}}$ in specific cases (that are sufficient for the setting of the Scholze--Shin conjecture).

\begin{rem}[{\cite[\S7.4.17 and Corollary 7.4.18]{KSZ}}]
\label{rem:general-s-endoscopic-transfer} To remove the assumption that $s$ lies in $Z(\wh{H})^{\Gamma_{\Q_p}}$ we must assume that our function $\phi$ has the property that the Kottwitz invariant $\kappa_G$ is constant on the support of $\phi$. Let us write $\nu_\phi$ for this common character of $Z(\wh{G})$. 

Suppose now only that $s$ is in $Z(\wh{H})^{\Gamma_{\Q_p}}Z(\wh{G})$ (the general case for an endoscopic datum), and choose a factorization $s=s's''$, where $s'$ belongs to $Z(\wh{H})^{\Gamma_{\Q_p}}$ and $s''$ belongs to $Z(\wh{G})$. If
$\mf{e}'\defeq(H,\mc{H},s',\eta)$, we define
$\phi^{\mf{e}}\defeq\nu_\phi(s'')\phi^{\mf{e}'}$. Equivalently, for every $[b]$ such that $\mr{TO}_{\delta[b]}(\phi)\ne0$, the character $\beta_p(\gamma_0,[b])$ admits the extension 
\begin{equation*}
    \wt{\beta}_p(\gamma_0,[b])
    \in
    X^\ast\big(
        Z(\wh{I}_{\gamma_0})^{\Gamma_{\Q_p}}Z(\wh{G}),
    \big)
\end{equation*}
prescribed in \cite[\S7.4.17]{KSZ} whose restriction to $Z(\wh{G})$ is $\nu_\phi$. 

So then, if $c$ is a norm from $R_{d,1}(\Q_p)$, formula
\eqref{eq:characterization-twisted-endoscopic-transfer} then takes the uniform form
\begin{equation}\label{eq:characterization-twisted-endoscopic-transfer-general-s}
\begin{aligned}
    \mr{SO}_{c}(\phi^{\mf{e}})
    &=
    \Delta_0^{\mf{e}}(c,\gamma_0)
    \sum_{[b]\in\mf{D}_d(\gamma_0,G;\Q_p)}
    e(I^\sigma_{\delta[b]})
    \left\langle
        \wt{\beta}_p(\gamma_0,[b]),s
    \right\rangle
    \mr{TO}_{\delta[b]}(\phi),
\end{aligned}
\end{equation}
where a summand is understood to be zero if
$\mr{TO}_{\delta[b]}(\phi)=0$. If $c$ is not a norm from
$R_{d,1}(\Q_p)$, then also as in the previous case $\mr{SO}_{c}(\phi^{\mf{e}})=0$. \end{rem}

\begin{rem}[Dependence on $z$-extension]
\label{rem:changing-z-extension} In both our discussion of endoscopic transfer and twisted endoscopic transfer we made an auxiliary choice: the $z$-extension $\mf{z}$. While the transferred functions are not literally independent of $\mf{z}$, the dependence is mild.

Namely, the stable orbital integrals of these transfers are compatible with enlarging $\mf{z}$. More precisely, suppose that
$\mf{z}'$ is an enlargement of $\mf{z}$ as in
\cite[Lemma 7.2.13]{KSZ}, and write $H_1'$ and $\chi_1'$ for the
associated objects. There is then a compatible inclusion
$H_1\to H_1'$, and $\chi_1'$ restricts to $\chi_1$. Restriction along $H_1(\Q_p)\hookrightarrow H_1'(\Q_p)$ identifies $\mc{H}_{\chi_1'}(H_1'(\Q_p))\isomto \mc{H}_{\chi_1}(H_1(\Q_p))$. With compatible measures and transfer factors, the transferred stable orbital integrals agree under this identification; see the proof of \cite[Proposition 7.4.13]{KSZ}.
\end{rem}

\subsection{The Bernstein center, the spectral Bernstein center, and the Fargues--Scholze map}
\label{s:bernstein-centers}

The more conceptual version of the Scholze--Shin conjecture is stated in terms of maps between the Bernstein and spectral Bernstein centers. We now briefly review this material, in particular recording the large improvements in our knowledge since the writing of \cite{ScholzeShin}. 

\begin{nota}\label{nota:Bernstein-center}
In contrast to previous sections, throughout this subsection:
\begin{itemize}
\item $F$ is a $p$-adic field;
\item $G$ is a connected reductive group over $F$;
\item $\cat{Rep}^{\mr{sm}}_\C(G(F))$ is the category of smooth complex representations of $G(F)$.
\end{itemize}
\end{nota}

\subsubsection{The Bernstein center}\label{ss:Bernstein-center}

We first briefly recall the theory of the Bernstein center of $G(F)$, in particular recalling the explicit geometric description used later.

\begin{defn}
\label{defn:Bernstein-center}
The \emph{Bernstein center} of $G(F)$ is the ring
\begin{equation*}
\mc{Z}(G)
\defeq
\End\left(
\id_{\cat{Rep}^{\mr{sm}}_\C(G(F))}
\right),
\end{equation*}
i.e., the ring of natural-in-$\pi$ endomorphisms $z(\pi)\colon \pi\to \pi$.
\end{defn}

\begin{rem} If $\pi$ is an irreducible object of $\cat{Rep}_\C^\mr{sm}(G(F))$, then $z(\pi)$ must be the action of a scalar by Schur's lemma. We abuse notation and let $z(\pi)$ also denote this scalar.
\end{rem}

We next recall the Bernstein variety $\Omega(G)$, whose ring of regular functions identifies with $\mc Z(G)$. This requires recalling some basic definitions from the theory of the Bernstein center; see (for example) \cite{RocheBernstein} for a nice overview of this theory.

\begin{setup}
\label{setup:supercuspidal-support}
A \emph{supercuspidal pair} for $G$ consists of a pair $(M,\sigma)$, where 
\begin{itemize}
\item $M\subseteq G$ is a Levi subgroup;
\item  $\sigma$ is an irreducible supercuspidal representation of $M(F)$.
\end{itemize}
An object $\pi$ of $\cat{Rep}_\C^\mr{sm}(G(F))$ has \emph{supercuspidal support} $\mr{sc}(\pi)=(M,\sigma)$ if, for some parabolic subgroup $P\subseteq G$ with Levi quotient $M$, the representation $\pi$ is an irreducible subquotient of the normalized parabolic induction $i_P^G(\sigma)$.

Two supercuspidal pairs $(M,\sigma)$ and $(M',\sigma')$ are \emph{inertially equivalent} if there exist $g$ in $G(F)$ and an unramified character (see \eqref{eq:unrm-char} below) $\chi$ of $M(F)$ such that
\begin{equation*}
M'=gMg^{-1},
\qquad
\sigma'\simeq{}^g\sigma\otimes{}^g\chi,
\end{equation*}
where ${}^g (-)(x)\defeq (-)(g^{-1}xg)$. We write $\mf{s}=[M,\sigma]_G$ for the inertial equivalence class of $(M,\sigma)$ and $\mf{B}(G)$ for the set of such classes.
\end{setup}

We next recall the algebraic tori from which the Bernstein components are formed.

\begin{nota} \label{nota:unramified-character-torus}
For a Levi subgroup $M\subseteq G$, set
\begin{equation*}
M(F)^1
\defeq
\bigcap_{\chi\in X_F^\ast(M)}
\ker\left(|\chi|_F\colon M(F)\to\R_{>0}\right)
\end{equation*}
and define the space of unramified characters of $M$:
\begin{equation}\label{eq:unrm-char}
X_{\mr{nr}}(M) \defeq \Hom\left(M(F)/M(F)^1,\C^\times\right).
\end{equation}
As $M(F)/M(F)^1$ is a free $\Z$-module, $X_{\mr{nr}}(M)$ is naturally an algebraic torus.

For a supercuspidal representation $\sigma$ of $M(F)$, there is a natural finite subgroup of $X_{\mr{nr}}(M)$:
\begin{equation*}
X_{\mr{nr}}(M,\sigma)\defeq\left\{\chi\in X_{\mr{nr}}(M):\sigma\otimes\chi\simeq\sigma\right\},
\end{equation*}
Thus, for $\mf{s}=[M,\sigma]_G$, the quotient $D_{\mf{s}}\defeq X_{\mr{nr}}(M)/X_{\mr{nr}}(M,\sigma)$ is again an algebraic torus.
\end{nota}

\begin{nota}
For $\mf{s}=[M,\sigma]_G$, we define
\begin{equation*}
W_{\mf{s}} \defeq\frac{\big\{g\in N_G(M)(F):{}^g\sigma\simeq\sigma\otimes\chi \text{ for some }\chi\in X_{\mr{nr}}(M) \big\}}{M(F)}.
\end{equation*}
This is a finite group acting algebraically on $D_{\mf{s}}$.
\end{nota}

We can now state Bernstein's description of the center in terms of this variety:

\begin{thm}[{\cite{BernsteinDeligne}}]
\label{thm:Bernstein-decomposition}
Let $\mf{s}=[M,\sigma]_G$ be an element of $\mf{B}(G)$.
\begin{enumerate}
\item The set $\Omega_\mf{s}(G)$ of $G(F)$-conjugacy classes of supercuspidal pairs inertially equivalent to $(M,\sigma)$ has the structure of an affine variety, namely, $\Omega_{\mf{s}}(G)\simeq D_{\mf{s}}/W_{\mf{s}}$.
\item If $\cat{Rep}^{\mf{s}}_\C(G(F))$ denotes the full subcategory of $\cat{Rep}^\mr{sm}_\C(G(F))$ of smooth representations all of whose irreducible subquotients have supercuspidal support in $\mf{s}$, then:
\begin{equation}\label{eq:Bernstein-category-decomposition}
\cat{Rep}^{\mr{sm}}_\C(G(F))= \prod_{\mf{s}\in\mf{B}(G)}\cat{Rep}^{\mf{s}}_\C(G(F)).
\end{equation}
\item Set $\Omega(G)\defeq\coprod_{\mf s\in\mf B(G)}\Omega_{\mf s}(G)$.
Via \eqref{eq:Bernstein-category-decomposition}, there is a canonical
isomorphism
\begin{equation}\label{eq:Bernstein-center-functions}
 \mc Z(G)\isomto
 \Gamma\bigl(\Omega(G),\mc O_{\Omega(G)}\bigr).
\end{equation}
\end{enumerate}
\end{thm}

We next recall the distribution-theoretic description of the Bernstein center. 

\begin{defn}
\label{defn:essentially-compact-invariant-distributions}
A \emph{distribution} on $G(F)$ is a $\C$-linear map $D\colon \mc{H}_\C(G(F))\to\C$. For $f$ in $\mc{H}_\C(G(F))$, set $\check{f}(g)\defeq f(g^{-1})$. If $D$ is a distribution and $f$ is in $\mc{H}_\C(G(F))$, we define \begin{equation}\label{eq:distribution-convolution-function}
D*f\colon G(F)\to\C, \qquad (D*f)(g)\defeq D\bigl(x\mapsto f(x^{-1}g)\bigr).
\end{equation}
We say that $D$ is \emph{essentially compact} if $D*f$ is in $\mc{H}_\C(G(F))$ for all $f$, and is \emph{$G(F)$-invariant} if
\begin{equation*}
D({}^g f)=D(f),
\qquad
{}^g f(x)\defeq f(g^{-1}xg).
\end{equation*}
We write $\mr{Dist}_{\mr{ec}}(G(F))^{G(F)}$ for the space of essentially compact $G(F)$-invariant distributions.
\end{defn}

\begin{rem} If $D_1$ and $D_2$ are essentially compact distributions, then the operation 
\begin{equation}\label{eq:convolution-distributions}
(D_1*D_2)(f)
\defeq
\check{D}_1(D_2*\check{f}),
\qquad
\check{D}_1(h)\defeq D_1(\check{h})
\end{equation}
endows $\mr{Dist}_{\mr{ec}}(G(F))^{G(F)}$ with the structure of a commutative $\C$-algebra.
\end{rem}

\begin{thm}[{\cite[\S1.4--1.7]{BernsteinDeligne}}]
\label{thm:Bernstein-center-as-distributions}
There is a canonical isomorphism of $\C$-algebras
\begin{equation}\label{eq:Bernstein-center-as-distributions}
\mc{Z}(G)
\isomto
\mr{Dist}_{\mr{ec}}(G(F))^{G(F)},
\qquad
z\longmapsto D_z,
\end{equation}
characterized as follows: for every object $\pi$ of $\cat{Rep}_\C^\mr{sm}(G(F))$ and every $f$ in $\mc{H}_\C(G(F))$ one has 
\begin{equation}\label{eq:Bernstein-distribution-action}
\pi(D_z*f) = z(\pi)\circ\pi(f),
\end{equation}
where for $h$ in $\mc{H}_\C(G(F))$, $\pi(h)$ is the natural operator on $\pi$ induced by $h$; see \emph{Notation \ref{nota:Hecke-action-repn}}.
\end{thm}

\begin{obsnota}
\label{obsnota:Bernstein-center-two-realizations}
From now on we confuse the notation for $z$ and $D_z$, using $z$ for both. Thus, for $f$ in $\mc{H}_\C(G(F))$, the notation $z(f)$ denotes evaluation of this distribution at $f$, while $z*f$ denotes the convolution in \eqref{eq:distribution-convolution-function}. 

Now, Theorem \ref{thm:Bernstein-decomposition} gives the Bernstein isomorphism
\begin{equation}\label{eq:Bernstein-isomorphism-functions}
\mc{Z}(G)
\isomto \Gamma\bigl(\Omega(G),\mc{O}_{\Omega(G)}\bigr),\qquad z\longmapsto\widehat{z}.
\end{equation}
If $\pi$ is an irreducible object of $\cat{Rep}_\C^{\mr{sm}}(G(F))$ then $z(\pi)=\wh{z}(\mr{sc}(\pi))$ as elements of $\C$. Thus \eqref{eq:Bernstein-distribution-action} implies that for every $f$ in $\mc{H}_\C(G(F))$ one has
\begin{equation}\label{eq:Bernstein-center-convolution-characterization}
\pi(z*f) = z(\pi)\cdot \pi(f)=\widehat{z}(\mr{sc}(\pi))\cdot \pi(f),
\end{equation}
and thus
\begin{equation}\label{eq:Bernstein-center-trace}
\tr(z*f\mid\pi)= z(\pi)\tr(f\mid \pi)=\widehat{z}\bigl(\mr{sc}(\pi)\bigr) \tr(f\mid\pi).
\end{equation}
\end{obsnota}

This distribution-theoretic description of the Bernstein center allows us to define various subspaces of it of a harmonic-analytic nature which feature prominently below.

\begin{defn}
\label{defn:unstable-functions}
A function $f$ in $\mc{H}_\C(G(F))$ is \emph{unstable} if $\mr{SO}_\gamma(f)=0$ for every strongly regular semisimple\footnote{A semisimple element $\gamma$ of $G(F)$ is \emph{strongly regular} if its centralizer in $G$ is a torus.} element $\gamma$ of $G(F)$. Write $\mc{H}_\C(G(F))^{\mr{unst}}$ for the space of unstable functions. 
\end{defn}

\begin{defn} If $f_1$ and $f_2$ are two elements of $\mc{H}_\C(G(F))$, we write $f_1=_\mr{st} f_2$ if $f_1-f_2$ is unstable, i.e., $f_1$ and $f_2$ have the same strongly regular stable orbital integrals.
\end{defn}

\begin{defn}
\label{defn:stable-very-stable-center}
We make the following definitions:
\begin{enumerate}
\item We say $z$ in $\mc{Z}(G)$ is \emph{stable} if $z\left(\mc{H}_\C(G(F))^\mr{unst}\right)=0$. The set of such stable elements is denoted $\mc{Z}^{\mr{st}}(G)$.

\item We say $z$ in $\mc{Z}(G)$ is \emph{very stable} if $z\ast\mc{H}_\C(G(F))^{\mr{unst}}\subseteq \mc{H}_\C(G(F))^{\mr{unst}}$. The set of such very stable elements is denoted $\mc{Z}^{\mr{vst}}(G)$.
\end{enumerate}
\end{defn}

\begin{prop}[{see \cite[\S1.1]{HansenStableBernstein}}]
\label{prop:stable-very-stable-basic-properties}
One has inclusions
\begin{equation}\label{eq:stable-very-stable-inclusions}
\mc{Z}^{\mr{vst}}(G) \subseteq \mc{Z}^{\mr{st}}(G) \subseteq \mc{Z}(G).
\end{equation}
Moreover:
\begin{enumerate}
\item $\mc{Z}^{\mr{vst}}(G)$ is a $\C$-subalgebra of $\mc{Z}(G)$;
\item $\mc{Z}^{\mr{st}}(G)$ is naturally a module over $\mc{Z}^{\mr{vst}}(G)$.
\end{enumerate}
\end{prop}
\begin{rem} It is not known in general that $\mc{Z}^{\mr{st}}(G)$ is a subalgebra of $\mc{Z}(G)$, but it is conjectured to be so; see \cite[p.~27]{BKV}.
\end{rem}

\subsubsection{The spectral Bernstein center}\label{sss:spectral-Bernstein} We now recall the spectral Bernstein center of $G$. Throughout, we continue to use Notation \ref{nota:Bernstein-center}, and fix a $\Gamma_F$-stable pinning of $\wh{G}$.

\begin{nota} Choose a lift of (arithmetic) Frobenius $\mr{Fr}$ and a topological generator $s$ of the tame inertia group $I_F/P_F$. Define $W_F^0\subseteq W_F$ to be the subgroup obtained as the preimage of $s^{\Z[\nicefrac{1}{q}]}\rtimes \mr{Fr}^\Z$ (where $q$ is the size of the residue field of $F$) under the projection map $W_F\to W_F/P_F$. 
\end{nota}

\begin{defn}[{\cite{DHKM}}] Define the presheaf on $\C$-algebras $R$
\begin{equation*}
Z^1(W_F^0,\wh{G})(R)\defeq \left\{\varphi\colon W_F^0\to \wh{G}(R):\begin{aligned}(1) & \qquad \varphi(w_1 w_2)=\varphi(w_1)(w_1\cdot \varphi(w_2))\\ (2) & \qquad \varphi|_{P_F}\text{ is locally constant.}\end{aligned}\right\}
\end{equation*}
By the discussion on \cite[p.\@ 1831]{DHKM}, $Z^1(W_F^0,\wh{G})$ is a locally of finite type $\C$-scheme called the \emph{framed moduli of $L$-parameters}. 

There is a natural conjugation action of $\wh{G}$ on $Z^1(W_F^0,\wh{G})$ and a $\C$-point $\varphi$ of $Z^1(W^0_F,\wh{G})$ is called \emph{semisimple} if its $\wh{G}$-orbit is closed. 
\end{defn}

\begin{nota}
\label{nota:spectral-parameter-space}
The \emph{coarse moduli space of semisimple $L$-parameters} is the affine GIT quotient\footnote{Note that while $Z^1(W_F^0,\wh{G})$ is not affine, it is a disjoint union of affine $\wh{G}$-stable open subschemes, and thus the affine GIT quotient here makes sense.}
\begin{equation*}
X^{\mr{spec}}_G\defeq Z^1(W_F^0,\wh{G})/\!\!/\wh{G},
\end{equation*}
which, by \cite[Theorem 4.18]{DHKM}, is canonically independent of all choices.
\end{nota}

\begin{rem} Recall that classically a Weil--Deligne $L$-parameter for $G$ is a pair $(r,N)$ where
\begin{itemize}
\item $r\colon W_F\to{}^LG=\wh{G}(\C)\rtimes W_F$ is a
homomorphism of the form $r(w)=(r_0(w),w)$, with $r_0|_{I_F}$ locally constant and of finite image;
\item $N$ is a nilpotent element of $\mr{Lie}(\wh{G})$ satisfying
\begin{equation*}
\mr{Ad}(r(\mr{Fr}))(N)=qN,\qquad
 \mr{Ad}(r(\gamma))(N)=N\end{equation*}
for every $\gamma$ in $I_F$.
\end{itemize}
The set of Weil--Deligne $L$-parameters carries a natural action of $\wh{G}(\C)$ via conjugation. Then, despite the somewhat baroque presentation, $Z^1(W^0_F,\wh{G})(\C)$ is really naturally $\wh{G}(\C)$-equivariantly identified with the set of Weil--Deligne $L$-parameters; see \cite[Proposition 2.7]{DHKM}. 

Under this bijection semisimplicity in $Z^1(W_F^0,\wh{G})(\C)$ corresponds to semisimplicity of $(r,N)$: that $N=0$ and $r(\mr{Fr})$ is semisimple; see \cite[Corollary 4.16 and Remark 6.9(ii)]{DHKM}. We shall write such objects just as homomorphisms $\lambda\colon W_F\to {}^L G$, and call them \emph{infinitesimal parameters}.
\end{rem}

We now come to the Galois-theoretic analogue of the Bernstein center.

\begin{defn}
The \emph{spectral Bernstein center} of $G$ is
\begin{equation}\label{eq:spectral-Bernstein-center}
\mc{Z}^{\mr{spec}}(G)\defeq \mc{O}(X_G^\mr{spec})=\Gamma\big(Z^1(W_F^0,\wh{G}),\mc{O}\big)^{\wh{G}}.
\end{equation}
\end{defn}

Below it will be useful to understand $X_G^\mr{spec}$ in terms of  Haines's stable Bernstein variety $X_G^\mr{Hai}$, the Galois-theoretic analogue of $\Omega(G)$ from Theorem \ref{thm:Bernstein-decomposition}; see \cite{HaistBC}.

\begin{defn}[{\cite[\S3]{BorelCorvallis}}]
\label{defn:L-Levi-minimal-factorization}
A \emph{standard $L$-Levi subgroup} of ${}^L G$ is a subgroup of the form $\mc{M}=\mc{M}^{\circ}\rtimes W_F$, where $\mc{M}^{\circ}$ is a $W_F$-stable standard Levi subgroup of $\wh{G}$. An \emph{$L$-Levi subgroup} is a $\wh{G}$-conjugate of a standard $L$-Levi subgroup. For any such $\mc{M}$ we write $\mc{M}^{\circ}\defeq \mc{M}\cap\wh{G}$.
\end{defn}

\begin{defn}
Let $\mc{M}\subseteq{}^L G$ be an $L$-Levi subgroup, and set
\begin{equation}\label{eq:Haines-unramified-torus}
X(\mc{M}^{\circ})\defeq\left(\left(Z(\mc{M}^{\circ})^{I_F}\right)_{\mr{Fr}}\right)^{\circ};
\end{equation}
this has the natural structure of a complex algebraic torus. 
\end{defn}

\begin{defn}\label{defn:unramified-central-twists} Let $\lambda\colon W_F\to\mc{M}$ be an infinitesimal parameter and $z$ an element of $X(\mc{M}^{\circ})$. Choose a lift $\dot z$ in $Z(\mc{M}^{\circ})^{I_F}$ of $z$, and let $\dot z\lambda\colon W_F\to\mc{M}$ be uniquely characterized by the conditions that $(\dot z\lambda)|_{I_F}=\lambda|_{I_F}$ and $(\dot z\lambda)(\mr{Fr})=\dot z\lambda(\mr{Fr})$. Then the $\mc{M}^{\circ}$-conjugacy class $(z\lambda)_{\mc{M}^{\circ}}$ of $\dot z\lambda$ depends only on $z$. We call it the \emph{unramified central twist} of $\lambda$ by $z$.
\end{defn}

\begin{defn}
\label{defn:inertial-equivalence-infinitesimal-parameters}
Let $\lambda_1,\lambda_2\colon W_F\to\mc{M}$ be infinitesimal parameters which factor minimally through the same $L$-Levi subgroup $\mc{M}$.
\begin{enumerate}
\item We say that $\lambda_1$ and $\lambda_2$ are \emph{inertially equivalent under $\mc{M}^{\circ}$} if, for some $z$ in $X(\mc{M}^{\circ})$, one has $(z\lambda_1)_{\mc{M}^{\circ}}=(\lambda_2)_{\mc{M}^{\circ}}$. We write $[\lambda]_{\mc{M}^{\circ}}$ for the resulting inertial class.

\item Two $\wh{G}(\C)$-conjugacy classes $(\lambda_1)_{\wh{G}}$ and $(\lambda_2)_{\wh{G}}$ are \emph{inertially equivalent} if they admit representatives $\lambda_1$ and $\lambda_2$ factorizing minimally through a common $L$-Levi $\mc{M}$ in which they are inertially equivalent under $\mc{M}^\circ$.  
\end{enumerate}
The second relation is independent of all choices and is an equivalence relation; see \cite[\S5.3]{HaistBC}. 
\end{defn}

\begin{setup}[{\cite[\S5.3]{HaistBC}}]
\label{setup:Haines-unramified-parameter-torus}
Let $\mf{t}_{\mc{M}^\circ}$ be an inertial equivalence class under $\mc{M}^\circ$ and $\mf{t}$ the associated inertial equivalence class for $\wh{G}$; every inertial class for $\wh{G}$ arises in this way, with $(\mc{M},\mf{t}_{\mc{M}^{\circ}})$ determined up to $\wh{G}$-conjugacy. Write $X_{\mf{t}_{\mc{M}^{\circ}}}$ and $X_\mf{t}$ for the sets of $\mc{M}^{\circ}$-conjugacy classes and $\wh{G}$-conjugacy classes of infinitesimal parameters in $\mf{t}_{\mc{M}^\circ}$ and $\mf{t}$ respectively.

The action of the torus $X(\mc{M}^{\circ})$ on $X_{\mf{t}_{\mc{M}^{\circ}}}$ by unramified central twist is transitive. For a base point $[\lambda]_{\mc{M}^{\circ}}$, the stabilizer $\mr{stab}_{\lambda}$ is finite and acts algebraically on $X(\mc{M}^\circ)$. Thus, $X_{\mf{t}_{\mc{M}^\circ}}=X(\mc{M}^\circ)/\mr{stab}_\lambda$ has a natural algebraic structure. If we furthermore set
\begin{equation*}
N_{\wh{G}}\bigl(\mc{M},[\lambda]_{\mc{M}^{\circ}}\bigr)
\defeq
\left\{
 n\in N_{\wh{G}}(\mc{M}):
 [\,{}^n\!\lambda\,]_{\mc{M}^{\circ}}
 =
 [\lambda]_{\mc{M}^{\circ}}
\right\},
\end{equation*}
where ${}^n\!\lambda(w)=n\lambda(w)n^{-1}$, and $W_{\mf{t}}\defeq N_{\wh{G}}\big(\mc{M},[\lambda]_{\mc{M}^{\circ}}\big)/\mc{M}^{\circ}$, then $W_\mf{t}$ is a finite group, which acts algebraically on $X_{\mf{t}_{\mc{M}^{\circ}}}$ and for which the natural map $X_{\mf{t}_{\mc{M}^\circ}}\to X_\mf{t}$ induces a bijection $X_{\mf{t}_{\mc{M}^\circ}}/W_\mf{t}\isomto X_\mf{t}$, endowing $X_\mf{t}$ with natural algebraic structure.
\end{setup}

We now arrive at the Haines variety which, as indicated before, has followed essentially the same construction (translated to the Galois-theoretic setting) as $\Omega(G)$ from Theorem \ref{thm:Bernstein-decomposition}. 

\begin{defn} The \emph{Haines variety} is $X_G^\mr{Hai}=\bigsqcup_{\mf{t}}X_\mf{t}$. 
\end{defn}

Within the analogy of $X_G^\mr{Hai}$ and $\Omega(G)$, the following is an analogue of Theorem \ref{thm:Bernstein-decomposition}.

\begin{thm}[{\cite[Theorem 6.10]{DHKM}}]
\label{thm:Haines-spectral-parameter-space}
There is a canonical isomorphism $X_G^{\mr{Hai}}\isomto X_G^{\mr{spec}}$.
\end{thm}

\begin{rem} We will often treat the identification of Theorem \ref{thm:Haines-spectral-parameter-space} as implicit, only highlighting its usage when we feel it's clarifying.
\end{rem}

\subsubsection{The Fargues--Scholze map} A relationship between $\mc{Z}(G)$ and $\mc{Z}^\mr{spec}(G)$ should be seen as a manifestation of some kind of local Langlands correspondence for $G$; see \cite[Proposition 4.23]{HaistBC}. But, while such a correspondence does not yet exist for arbitrary $G$, this relationship between Bernstein centers does now exist unconditionally.

\begin{thm}[Fargues--Scholze and Scholze, {\cite{FarguesScholze,ScholzeMotivic}}]
\label{thm:Fargues-Scholze-Bernstein-map}
There is a natural $\C$-algebra map 
\begin{equation*}
\Psi^\mr{FS}_G\colon \mc{Z}^{\mr{spec}}(G)\to\mc{Z}(G).
\end{equation*}
\end{thm}

\begin{nota}
We write $\mc{Z}^{\mr{FS}}(G)\defeq\mr{im}\left(\Psi^\mr{FS}_G\colon\mc{Z}^{\mr{spec}}(G)\to\mc{Z}(G)\right)$.
\end{nota}

\begin{rem}\label{rem:parameter-description-of-Psi} In fact, Fargues--Scholze's construction gives something even more precise. For every irreducible object $\pi$ of $\cat{Rep}_\C^\mr{sm}(G(F))$ they define an infinitesimal parameter
\begin{equation*}
\varphi_\pi^{\mr{FS}}\colon W_F\to{}^L G,
\end{equation*}
well-defined up to $\wh{G}(\C)$-conjugacy, such that for every $z$ in $\mc{Z}^{\mr{spec}}(G)$ one has the equality
\begin{equation}\label{eq:Fargues-Scholze-eigenvalue}
\Psi^\mr{FS}_G(z)(\pi) = z\left([\varphi_\pi^{\mr{FS}}]\right).
\end{equation}
Here $[\varphi_\pi^{\mr{FS}}]$ denotes the corresponding closed point of $X_G^{\mr{spec}}=X_G^\mr{Hai}$ defined by $\varphi_\pi^\mr{FS}$. 
\end{rem}

It is useful to push the content of Remark \ref{rem:parameter-description-of-Psi} even further.

\begin{prop}[{\cite[Corollary IX.7.3]{FarguesScholze}}]
\label{prop:Fargues-Scholze-parabolic-induction}
Let $P\subseteq G$ be a parabolic with Levi $M$, $\sigma$ an irreducible object of $\cat{Rep}_\C^\mr{sm}(M(F))$, and $\pi$ an irreducible subquotient of $i_P^G(\sigma)$. Then
\begin{equation}\label{eq:Fargues-Scholze-parabolic-induction}
[\varphi_\pi^{\mr{FS}}]= \left[{}^L \eta_{M,G}\circ
\varphi_\sigma^{\mr{FS}}\right],
\end{equation}
as elements of $X_G^\mr{spec}=X_G^\mr{Hai}$, where ${}^L\eta_{M,G}\colon{}^L M\to{}^L G$ is the usual $L$-embedding.
\end{prop}

\begin{rem} From Proposition \ref{prop:Fargues-Scholze-parabolic-induction}, we see that $[\varphi_\pi^{\mr{FS}}]$ depends only on the supercuspidal support of $\pi$. We therefore obtain a morphism of complex varieties
\begin{equation*}
\lambda_G^{\mr{FS}}\colon
\Omega(G)\to X_G^\mr{spec}=X_G^{\mr{Hai}},\qquad \mr{sc}(\pi)\mapsto[\varphi_\pi^{\mr{FS}}].
\end{equation*}
Under the Bernstein isomorphism \eqref{eq:Bernstein-isomorphism-functions}, we then have the equality $\Psi^\mr{FS}_G=(\lambda_G^\mr{FS})^\ast$. 
\end{rem}

\begin{rem}
\label{rem:FS-finiteness-status} If $\lambda_G^\mr{FS}$ arose from a local Langlands correspondence for $G$ (e.g., as in \cite[Proposition 4.23]{HaistBC}), one would expect that it is a finite morphism. For example, this is part of \cite[Conjecture 6.3]{ScholzeShin}. While the finiteness of $\lambda_G^\mr{FS}$ is not currently known in general it is known that $\lambda_G^\mr{FS}$ is locally finite, i.e., finite on each component of $\Omega(G)$; see \cite[Corollary 3.5]{DHKM2}. That said, the recent work \cite{CotnerFengII} seems quite promising to upgrade this claim to actual finiteness in the near future (at least in many cases).
\end{rem}

If again $\Psi_G^\mr{FS}$ arose from a local Langlands correspondence, then one would expect $\Psi_G^\mr{FS}$ to factorize through $\mc{Z}^\mr{st}$; for example this is part of \cite[Conjecture 6.3]{ScholzeShin}. In fact, an even stronger claim is now known unconditionally.

\begin{thm}[{\cite{HansenStableBernstein}}]
\label{thm:Hansen-very-stable-center}
One has the containment $\mc{Z}^{\mr{FS}}(G)\subseteq\mc{Z}^{\mr{vst}}(G)$.
\end{thm}

\begin{rem}
\label{rem:status-SS-section-6} It is expected that the inclusions
\begin{equation*}
\mc{Z}^{\mr{FS}}(G)\subseteq \mc{Z}^{\mr{vst}}(G)\subseteq \mc{Z}^\mr{st}(G),
\end{equation*}
are all equalities; the second equality is \cite[Conjecture 1.1.4]{Varma}. In fact, if $G$ is quasi-split one expects that $\Psi^\mr{FS}_G\colon
\mc{Z}^{\mr{spec}}(G)\to \mc{Z}^{\mr{vst}}(G)$ is an isomorphism of $\C$-algebras.
\end{rem}

We next explain how Hansen's theorem allows us to unambiguously associate to any atomically stable virtual character of $G(F)$ a point of $X_G^\mr{spec}$. 

\begin{defn}[{cf.\@ \cite[Definition 5.1]{BMYCharacterization}}]
\label{defn:atomically-stable-character}
A \emph{virtual character} of $G(F)$ is a finite $\C$-linear combination $\Theta= \sum_{i=1}^m a_i\Theta_{\pi_i}$ where $\Theta_{\pi_i}$ is the Harish-Chandra character of an irreducible element of $\cat{Rep}_\C^\mr{sm}(G(F))$, which we view as an invariant distribution on $G(F)$. The set $\{\pi_1,
\ldots,\pi_m\}$ is called the \emph{support} of $\Theta$.

We call $\Theta$ \emph{stable} if $\Theta\left(\mc{H}_\C(G(F))^{\mr{unst}}\right)=0$. Assume furthermore that the $\pi_i$ are pairwise nonisomorphic and that every $a_i$ is nonzero. We say that $\Theta$ is \emph{atomically stable} if it is stable and no nonzero virtual character supported on a proper subset of $\{\pi_1,\ldots,\pi_m\}$ is stable.
\end{defn}

\begin{prop}[{cf.\@ \cite[Corollary 2.7]{HansenStableBernstein}}]
\label{prop:FS-atomically-stable-character}
Let $\Theta=\sum_{i=1}^m a_i\Theta_{\pi_i}$ be an atomically stable virtual character. Then for any $i$ and $j$ we have the equality $[\varphi_{\pi_i}^{\mr{FS}}]=[\varphi_{\pi_j}^{\mr{FS}}]$ in $X_G^\mr{spec}$.
\end{prop}
\begin{proof} Take any element $z$ of $\mc{Z}^\mr{spec}(G)$, and observe that 
\begin{equation*}
\Psi_G^\mr{FS}(z)\ast\Theta=\sum_{i=1}^m \Psi_G^\mr{FS}(z)(\pi_i)a_i\Theta_{\pi_i}=\sum_{i=1}^mz\left([\varphi_{\pi_i}^\mr{FS}]\right)a_i\Theta_{\pi_i},
\end{equation*}
with the last equality by \eqref{eq:Fargues-Scholze-eigenvalue}. We claim that $\Psi_G^\mr{FS}(z)\ast\Theta$ is stable. Indeed, by Theorem \ref{thm:Hansen-very-stable-center} we have that $\Psi_G^\mr{FS}(z)$ is in $\mc{Z}^\mr{vst}(G)$ and so for any $f$ in $\mc{H}_\C(G(F))^\mr{unst}$ we have that $\Psi_G^\mr{FS}(z)\ast f$ is in $\mc{H}_\C(G(F))^\mr{unst}$ from where we deduce the following from the stability of $\Theta$
\begin{equation*}
(\Psi^\mr{FS}_G(z)\ast \Theta)(f)=\Theta(\Psi^\mr{FS}_G(z)\ast f)=0,
\end{equation*}
from where the claim follows. But, then
\begin{equation*}
\Psi^\mr{FS}_G(z)\ast\Theta-z\left([\varphi_{\pi_1}^\mr{FS}]\right)\Theta=\sum_{i=2}^m\left(z\left([\varphi_{\pi_i}^\mr{FS}]\right)-z\left([\varphi_{\pi_1}^\mr{FS}]\right)\right)a_i\Theta_{\pi_i}
\end{equation*}
is stable with support contained in $\{\pi_1,\ldots,\pi_m\}$. Thus, by the atomic stability of $\Theta$ we deduce that this distribution is $0$ which, as each $a_i$ is non-zero, implies that $z\left([\varphi_{\pi_i}^\mr{FS}]\right)$ is independent of $i$. As this holds for any $z$ we deduce the claim.
\end{proof}

\begin{nota} Let $\Theta=\sum_{i=1}^m a_i\Theta_{\pi_i}$ be an atomically stable virtual character. By Proposition \ref{prop:FS-atomically-stable-character} one has that $[\varphi_{\pi_i}^\mr{FS}]$ is independent of $i$. We denote their common value by $\varphi_\Theta^{\mr{FS}}$. For every $z$ in $\mc{Z}^{\mr{spec}}(G)$ and every $f$ in $\mc{H}_\C(G(F))$, one has
\begin{equation}\label{eq:stable-character-FS-center}
\Theta\left(\Psi^\mr{FS}_G(z)*f\right)=z\left(\varphi_\Theta^{\mr{FS}}\right)\Theta(f).
\end{equation}
\end{nota} 

\subsubsection{The fixed-central-character setting}

Finally, for the Scholze--Shin conjecture it will be useful to quickly adapt the above discussion to the situation of Section \ref{s:transfer}. We now fix notation as in Notation \ref{nota:endoscopic-transfer} and Observation/Notation \ref{obsnota:z-pair-endoscopic-transfer}. 

\medskip

\paragraph*{The Bernstein center} We begin by giving the appropriate analogues of the various parts of the Bernstein center story from \S\ref{ss:Bernstein-center}.

\begin{defn}
\label{defn:fixed-central-character-Bernstein-center}
Set $\omega_1\defeq\chi_1^{-1}$, and define
\begin{itemize}
\item $\cat{Rep}^{\mr{sm}}_\C(H_1(\Q_p))_{\omega_1}$ to be the full subcategory of objects of $\cat{Rep}_\C^\mr{sm}(H_1(\Q_p))$ on which the group $Z_1(\Q_p)$ acts via $\omega_1$;
\item  $\mc{Z}(H_1,\chi_1)\defeq\End\left(\id_{\cat{Rep}^{\mr{sm}}_\C(H_1(\Q_p))_{\omega_1}}\right)$.
\end{itemize}
\end{defn}

We next isolate the analogues of Definitions \ref{defn:unstable-functions} and \ref{defn:stable-very-stable-center}.

\begin{defn} We make the following definitions:
\begin{enumerate}
\item A function $f$ in $\mc{H}_{\chi_1}(H_1(\Q_p))$ is \emph{unstable} if $\mr{SO}_\gamma(f)=0$ for every strongly regular semisimple element $\gamma$ of $H_1(\Q_p)$; write $\mc{H}_{\chi_1}(H_1(\Q_p))^{\mr{unst}}$ for the set of unstable functions.
\item If $f_1$ and $f_2$ are two elements of $\mc{H}_{\chi_1}(H_1(\Q_p))$, we write $f_1=_\mr{st} f_2$ if $f_1-f_2$ is unstable.
\item We say $z$ in $\mc{Z}(H_1,\chi_1)$ is \emph{stable} if $z\left(\mc{H}_{\chi_1}(H_1(\Q_p))^{\mr{unst}}\right)=0$. We denote the set of such elements by $\mc{Z}^{\mr{st}}(H_1,\chi_1)$.
\item We say $z$ in $\mc{Z}(H_1,\chi_1)$ is \emph{very stable} if $z*\mc{H}_{\chi_1}(H_1(\Q_p))^{\mr{unst}} \subseteq \mc{H}_{\chi_1}(H_1(\Q_p))^{\mr{unst}}$. We denote the set of such elements by $\mc{Z}^{\mr{vst}}(H_1,\chi_1)$.
\end{enumerate}
\end{defn}

\begin{prop}[{cf.\@ Proposition \ref{prop:stable-very-stable-basic-properties}}]
\label{prop:fixed-central-character-stable-center}
One has inclusions
\begin{equation*}
\mc{Z}^{\mr{vst}}(H_1,\chi_1)
\subseteq
\mc{Z}^{\mr{st}}(H_1,\chi_1)
\subseteq
\mc{Z}(H_1,\chi_1).
\end{equation*}
Moreover, $\mc{Z}^{\mr{vst}}(H_1,\chi_1)$ is a $\C$-subalgebra of $\mc{Z}(H_1,\chi_1)$ and $\mc{Z}^{\mr{st}}(H_1,\chi_1)$ is a module over it.
\end{prop}

\paragraph*{The spectral Bernstein center} We now turn to the spectral Bernstein center.

\begin{constr} The inclusion $Z_1\hookrightarrow H_1$ induces an $L$-homomorphism ${}^Lq_{Z_1}\colon{}^L H_1\to{}^L Z_1$. As $Z_1$ is a torus, the geometrically normalized Langlands correspondence for tori identifies $\omega_1$ with a parameter $\lambda_{\omega_1}\colon W_{\Q_p}\to{}^L Z_1$.
\end{constr}

\begin{nota}
\label{nota:fixed-central-character-spectral-space}
The $L$-homomorphism ${}^Lq_{Z_1}$ induces a map $X_{H_1}^{\mr{spec}}\to X_{Z_1}^{\mr{spec}}$. Set:
\begin{itemize}
\item $X_{H_1,\chi_1}^{\mr{spec}}\defeq X_{H_1}^{\mr{spec}}{\times}_{X_{Z_1}^{\mr{spec}}}\{[\lambda_{\omega_1}]\}$.
\item $\mc{Z}^{\mr{spec}}(H_1,\chi_1)\defeq \mc{O}\bigl(X_{H_1,\chi_1}^{\mr{spec}}\bigr)$.
\end{itemize}
\end{nota}

\paragraph*{The Fargues--Scholze morphism} Finally, we give the necessary adaptations of the Fargues--Scholze map to this setting.

\begin{prop}[Fargues--Scholze, Scholze, and Hansen]
\label{prop:FS-fixed-central-character}
There is a natural $\C$-algebra map
\begin{equation*}
\Psi^\mr{FS}_{H_1,\chi_1}\colon
\mc{Z}^{\mr{spec}}(H_1,\chi_1)\to \mc{Z}^{\mr{vst}}(H_1,\chi_1).
\end{equation*}
\end{prop}
\begin{proof} This follows from Theorem \ref{thm:Fargues-Scholze-Bernstein-map} and Theorem \ref{thm:Hansen-very-stable-center} given the compatibility of the Fargues--Scholze construction with central characters; see \cite[Theorem I.9.6(iii)]{FarguesScholze}. \end{proof}

We finally record the fixed-central-character analogue of Proposition \ref{prop:FS-atomically-stable-character}.

\begin{propnota}
\label{prop:Hansen-atomically-stable-character}
Let $\Theta=\sum_{i=1}^m a_i\Theta_{\pi_i}$ be an atomically stable virtual character of $H_1(\Q_p)$, with each $\pi_i$ in $\cat{Rep}_\C^\mr{sm}(H_1(\Q_p))_{\omega_1}$. Then the points $[\varphi_{\pi_i}^{\mr{FS}}]$ of $X_{H_1,\chi_1}^{\mr{spec}}$ do not depend on $i$. 

We denote the common value of $[\varphi_{\pi_i}^\mr{FS}]$ by $\varphi_\Theta^{\mr{FS}}$. Then, for every $z$ in $\mc{Z}^{\mr{spec}}(H_1,\chi_1)$ and every $f$ in $\mc{H}_{\chi_1}(H_1(\Q_p))$, one has the equality
\begin{equation}\label{eq:fixed-central-character-stable-character}
\Theta\left(\Psi^\mr{FS}_{H_1,\chi_1}(z)*f\right)=z\left(\varphi_\Theta^{\mr{FS}}\right)\Theta(f).
\end{equation}
\end{propnota}

\subsection{The Scholze--Shin conjecture}

We are now prepared to properly formulate the analogue of the Scholze--Shin conjecture in the general unramified situation. 

\begin{nota}\label{nota:Scholze-Shin-setup}
Throughout the following we maintain the notation from Notation \ref{nota:deformation-spaces-setup} and that set out in Sections \S\ref{s:transfer} and \S\ref{s:bernstein-centers}. We further fix:
\begin{itemize}
\item an integer $j\geqslant 1$ and an element $\tau$ in $W_K$ with $\mathsf{v}(\tau)=j$;
\item an element $h$ in $\mc{H}(\mc{G}(\Z_p))$, viewed by extension by zero as an element of $\mc{H}(G(\Q_p))$;
\item an endoscopic datum $\mf{e}=(H,\mc{H},s,\eta)$ for $G$, together with the $z$-pair and the notation $H_1$, $G_1$, $s_1$, $\eta_1$, and $\chi_1$ fixed in Observation/Notation \ref{obsnota:z-pair-endoscopic-transfer};
\item a maximal torus $T$ of $G_{\ov{\Q}_p}$ through which $\mu$ factors;
\item a set of positive roots $\Phi^+$ for $T$ with which $\mu$ is dominant;
\item $\rho_G$ to be $\tfrac{1}{2}\sum_{\alpha\in\Phi^+}\alpha$;
\item a lift $\mu_1\colon \bb{G}_{m,K}\to G_1\otimes K$ lifting $\mu$.
\end{itemize}
We abbreviate the functions $\phi_{\tau,h}^{\mc{G},\mu}$ from Definition \ref{defn:phitauh} to $\phi_{\tau,h}$.
\end{nota}

In essence, the Scholze--Shin conjecture describes the values of $\Psi^\mr{FS}_G(z)$ for certain special elements $z$ of $\mc{Z}^\mr{spec}(G)$ in terms of the function $\phi_{\tau,h}$. Thus, to state it properly we need to describe these special elements $z$, and the modifications of $\phi_{\tau,h}$ necessary to describe $\Psi^\mr{FS}_G(z)$.

\medskip

\paragraph*{The functions $z_{\mu,\tau}^{\mf{e},\mr{spec}}$ and $z_{\mu,\tau}^\mf{e}$} We begin by describing the special elements $z$ of $\mc{Z}^\mr{spec}(G)$, and more generally $\mc{Z}^\mr{spec}(H_1,\chi_1)$, that the Scholze--Shin conjecture concerns.

\begin{constr}[{cf.\@ \cite[Lemma 2.1.2]{KottwitzTwisted}}]\label{constr:Kottwitz} We let
\begin{equation*}
r_{\mu_1}\colon \wh{G}_1\rtimes {W_K}\to\GL(V_{\mu_1})
\end{equation*}
be the unique representation such that
\begin{enumerate}
\item $r_{\mu_1}|_{\wh{G}_1}$ is the irreducible representation $V_{\mu_1}$ of extremal weight $\mu_1$;
\item the action of ${W_K}$ via $r_{\mu_1}$ on the extremal weight vector in $V_{\mu_1}$ is trivial.
\end{enumerate}
\end{constr}

\begin{constr}
\label{constr:Scholze-Shin-spectral-element}
Define the map
\begin{equation}\label{eq:Scholze-Shin-spectral-factor}
z_{\mu,\tau}^{\mf{e},\mr{spec}}\colon X_{H_1,\chi_1}^\mr{spec}\to \C,\qquad [\lambda]\mapsto \tr\bigg(s_1\cdot \eta_1(\lambda(\tau))\mid r_{\mu_1}\bigg)q^{j\langle \rho_G,\mu\rangle}.
\end{equation}
\end{constr}

\begin{rem} Here there is the occurrence of $s_1$, rather than $s_1^{-1}$ as in \cite{ScholzeShin}. This difference reflects the difference in transfer-factor normalizations fixed in Remark \ref{rem:normalization-transfer}.
\end{rem}

\begin{prop}[{cf.\@ \cite[Proposition 4.28]{HaistBC}}] The map $z_{\mu,\tau}^{\mf{e},\mr{spec}}$ is algebraic, and thus determines an element of $\mc{Z}^{\mr{spec}}(H_1,\chi_1)$.
\end{prop}

For visual simplicity we give notation to the image of $z_{\mu,\tau}^{\mf{e},\mr{spec}}$ under $\Psi_{H_1,\chi_1}^\mr{FS}$.

\begin{defn}\label{defn:Scholze-Shin-central-element} We define $z_{\mu,\tau}^\mf{e}$ to be the element of $\mc{Z}(H_1,\chi_1)$ given by $\Psi^\mr{FS}_{H_1,\chi_1}(z_{\mu,\tau}^{\mf{e},\mr{spec}})$.
\end{defn}

\begin{rem} For the fixed $z$-pair, the elements $z_{\mu,\tau}^{\mf{e},\mr{spec}}$ and $z_{\mu,\tau}^{\mf{e}}$ are independent of the chosen lift $\mu_1$ of $\mu$. Indeed, changing the lift twists $r_{\mu_1}$ by a character factoring through $\wh{Z}_1$, whose contribution on $s_1\eta_1(\lambda(\tau))$ is trivial by construction.
\end{rem}

Finally, we would like to specialize the above definitions to the case of the principal endoscopic datum $\mf{e}_{\triv}=(G,{}^L\!G,1,\id)$.

\begin{rem}[The principal central elements]
\label{rem:principal-Scholze-Shin-central-elements}
Consider the principal endoscopic datum $\mf{e}=\mf{e}_{\triv}=(G,{}^L\!G,1,\id)$. In this case the character $\chi_1$ is the trivial character $\mathbf{1}$. 

Observe that there are canonical identifications
\begin{equation*}
X_G^{\mr{spec}}\simeq X_{G_1,\mathbf{1}}^{\mr{spec}},
\qquad
\mc{Z}^{\mr{spec}}(G)\simeq\mc{Z}^{\mr{spec}}(G_1,\mathbf{1}),
\qquad
\mc{Z}(G)\simeq\mc{Z}(G_1,\mathbf{1}).
\end{equation*}
We write $z_{\mu,\tau}^\mr{spec}$ in $\mc{Z}^\mr{spec}(G)$ and $z_{\mu,\tau}$ in $\mc{Z}(G)$ for the elements corresponding to $z_{\mu,\tau}^{\mf{e}_{\triv},\mr{spec}}$ and $z_{\mu,\tau}^{\mf{e}_{\triv}}$, respectively. Equivalently, if
\begin{equation*}
r_{\mu}\colon\wh{G}\rtimes W_K\to\GL(V_{\mu})
\end{equation*}
is the representation attached to $\mu$ as in Construction \ref{constr:Kottwitz}, then
\begin{equation*}
z_{\mu,\tau}^{\mr{spec}}([\lambda])=\tr \left(\lambda(\tau)\mid r_{\mu}\right) q^{j\langle\rho_G,\mu\rangle},
\end{equation*}
and $z_{\mu,\tau}=\Psi_G^{\mr{FS}}(z_{\mu,\tau}^{\mr{spec}})$. 
\end{rem}

\medskip

\paragraph*{The functions $f_{\tau,h}^\mf{e}$} We next describe the modifications of $\phi_{\tau,h}$ used to describe $z_{\mu,\tau}^\mf{e}$.

\begin{defn} As the support of $\phi_{\tau,h}$ is contained in $\mc{G}(W_j)\sigma(\mu)(p)\mc{G}(W_j)$ on which the Kottwitz invariant is constant, the discussion from Remark \ref{rem:general-s-endoscopic-transfer} applies. Thus, it makes sense to consider an $\mf{e}$-transfer $f_{\tau,h}^\mf{e}$ of $\phi_{\tau,h}$.
\end{defn}

\begin{rem}[Principal base-change transfers]
\label{rem:principal-base-change-transfer}
Maintain the notation of Remark \ref{rem:principal-Scholze-Shin-central-elements}. We have a natural identification $\mc{H}_\C(G(\Q_p))\simeq \mc{H}_{\mathbf{1}}(G_1(\Q_p))$. For any twisted $\mf{e}_{\triv}$-transfer $f_{\tau,h}^{\mf{e}_{\triv}}$ in $\mc{H}_{\mathbf{1}}(G_1(\Q_p))$ of $\phi_{\tau,h}$, we write $f_{\tau,h}$ for the corresponding element of $\mc{H}_\C(G(\Q_p))$.
\end{rem}

\paragraph*{The Scholze--Shin conjecture} We are now finally ready to state the Scholze--Shin conjecture in full (unramified) generality.

\begin{conj}[Scholze--Shin]
\label{conj:Scholze-Shin}
For any endoscopic datum $\mf{e}$, one has $z_{\mu,\tau}^\mf{e}\ast h^\mf{e}=_\mr{st} f_{\tau,h}^\mf{e}$. 
\end{conj}

\begin{rem} At first glance it appears that the Scholze--Shin conjecture is not well-defined as the transfer of functions is not. That said, this is not the case by the following observations:
\begin{itemize}
\item if $f^{\mf{e},'}_{\tau,h}$ is another $\mf{e}$-transfer of $\phi_{\tau,h}$ then $f^\mf{e}_{\tau,h}=_\mr{st} f^{\mf{e},'}_{\tau,h}$;
\item by Proposition \ref{prop:FS-fixed-central-character} we have that $z_{\mu,\tau}^\mf{e}$ is very stable, and thus if $h^{\mf{e},'}$ is another $\mf{e}$-transfer of $h$ then $z_{\mu,\tau}^\mf{e}\ast h^\mf{e}=_\mr{st} z_{\mu,\tau}^\mf{e}\ast h^{\mf{e},'}$. 
\end{itemize}
\end{rem}

In practice, it is often more useful to state the following version of the Scholze--Shin conjecture in terms of virtual characters.

\begin{rem}
\label{rem:Scholze-Shin-character-form}
Let $\Theta$ be as in the statement of Proposition \ref{prop:Hansen-atomically-stable-character}. Then, if Conjecture \ref{conj:Scholze-Shin} holds, we have the following identity
\begin{equation}\label{eq:Scholze-Shin-character-form}
\Theta\left(f_{\tau,h}^\mf{e}\right)= z_{\mu,\tau}^{\mf{e},\mr{spec}}\left(\varphi_\Theta^{\mr{FS}}\right)
\Theta\left(h^\mf{e}\right).
\end{equation}
\end{rem}

\begin{rem}[The principal case of the Scholze--Shin conjecture]
\label{rem:principal-endoscopic-transfer}
For the principal datum $\mf{e}_{\triv}$ we obtain the following form of the Scholze--Shin conjecture:
\begin{equation}\label{eq:principal-case-Scholze--Shin}
z_{\mu,\tau}\ast h =_\mr{st} f_{\tau,h}.
\end{equation}
Written in terms of stable virtual characters, we obtain the following conjecture
\begin{equation}\label{eq:Scholze-Shin-character-form-principal-case}
\Theta\left(f_{\tau,h}\right)= z_{\mu,\tau}^{\mr{spec}}\left(\varphi_\Theta^{\mr{FS}}\right)
\Theta\left(h\right).
\end{equation}
\end{rem}

\begin{rem}[A conjectural extension of the Scholze--Shin conjecture] In essence the Scholze--Shin conjecture (which, for simplicity, we only consider in the principal case) attempts to describe some values of $\Psi_G^\mr{FS}(z)$ in terms of the geometrically defined functions $\phi_{\tau,h}$. 

While desirable, as it gives a potentially extrinsic geometric way of understanding $\Psi_G^\mr{FS}$, it has a major flaw: The $z_{\mu,\tau}^\mr{spec}$ only give a relatively small number of elements of $\mc{Z}^\mr{spec}(G)$, and so it's not clear that a description of $z_{\mu,\tau}=\Psi_G^\mr{FS}(z_{\mu,\tau}^\mr{spec})$ is sufficient to characterize $\Psi_G^\mr{FS}$. 

In some cases, so-called \emph{good groups} in the sense of \cite[Definition 3.1]{BMYCharacterization}, the Scholze--Shin conjecture should be enough. Namely, op.\@ cit.\@ shows that, at least if one has the Scholze--Shin conjecture for all hyperendoscopic groups, then it (together with the usual desiderata) is enough to characterize a (discrete) local Langlands correspondence.

That said, there is a natural way to conjecturally generalize the setup which should be sufficient in all cases: one should use \emph{excursion operators}. Namely, for each excursion datum $\mc{D}$ in the sense of \cite[Definition VIII.4.2]{FarguesScholze} one can form an element $z^\mr{spec}_\mc{D}$ of $\mc{Z}^\mr{spec}(G)$, and, in fact, almost every element of $\mc{Z}^\mr{spec}(G)$ is of this form.\footnote{More precisely, as we're working over $\C$, every regular function on a locus with bounded wild ramification is the restriction of such an element. The full spectral Bernstein center is recovered by taking the inverse limit over these loci; see \cite[\S VIII.3.2 and \S IX.5]{FarguesScholze} for precise statements.} Such excursion operators are not so far afield from what we have considered. Namely, set $R=\mr{Ind}_{\wh{G}\rtimes W_K}^{\wh{G}\rtimes W_{\Q_p}}r_{\mu}$ and let $P$ be the $\wh{G}$-equivariant projector onto the identity-coset summand. Then $z_{\mu,\tau}^\mr{spec}=z_{\mc{D}}^\mr{spec}$ for
\begin{equation*}
 \mc{D}=\bigl(\{1,2\},R\boxtimes R^\vee,
 (P\otimes1)\mr{coev}_R,
  q^{j\langle\rho_G,\mu\rangle}\mr{ev}_R,(\tau,1)\bigr).
\end{equation*}
where $\mr{ev}_R$ and $\mr{coev}_R$ stand for `evaluation' and `coevaluation', respectively.

It is then natural to predict that there exist for any excursion datum $\mc{D}$ and $h$ in $\mc{H}_\C(\mc{G}(\Z_p))$ a function $f_{\mc{D},h}$ in $\mc{H}_\C(G(\Q_p))$, and a version of the Scholze--Shin conjecture should hold:
\begin{equation*}
z_\mc{D}\ast h=_\mr{st} f_{\mc{D},h},\qquad \Theta(f_{\mc{D},h})=z_\mc{D}^\mr{spec}(\varphi_\Theta^\mr{FS})\Theta(h),
\end{equation*}
where $z_\mc{D}=\Psi_G^\mr{FS}(z_\mc{D}^\mr{spec})$ and $\Theta$ is an atomically stable virtual character.

The definition of $f_{\mc{D},h}$ is currently unclear, but there is a natural guess. Namely, one might imagine that there is a prismatic version of the moduli space of shtukas with an arbitrary number of legs, for which the case of a single minuscule leg $\mu$ should recover something close to $\mr{BT}^{\mc{G},\mu}_\infty$. Then $f_{\mc{D},h}$ should be defined in a way analogous to the single-let setting, but using the deformation spaces of these multi-legged prismatic deformation spaces. In fact, one can probably currently formulate (although not give a fine study of) such functions using shtukas in place of prismatic theory; compare with \S\ref{ss:shtuka-description-test-functions}.
\end{rem}

\subsection{Applications to Shimura varieties}
\label{ss:applications-to-Shimura-varieties}

In this final section we record the application of Theorem \ref{thm:geometric-stabilization} to understand, conditional on several things, the decomposition of the cohomology of Shimura varieties, and the computation of their partial semisimple Hasse--Weil $\zeta$-functions. The ideas here are extensions of those from \cite{HaistBC} and \cite{ScholzeShin}, and we refer the reader to these sources for more discussion on terms below.

\begin{setup}\label{setup:Haines-application}
Retain Setup \ref{setup:sv} and Notations \ref{nota:stabilization-measures} and
\ref{nota:stabilization-coefficient}, and further fix/write/assume:
\begin{itemize}
\item $\mb{G}^{\mr{ad}}$ is $\Q$-anisotropic;
\item $\mb{G}=\mb{G}^c$;
\item an isomorphism $\iota_\ell\colon\ov{\Q}_\ell\isomto\C$;
\item an irreducible algebraic representation $\xi$ of
$\mb{G}_{\ov{\Q}_\ell}$;
\item a geometric Frobenius lift $\Phi_v$ in $W_E$;
\item $d=\dim_\C\mb{X}$;
\item $\mu=-\mu_h$.
\end{itemize}
\end{setup}

\subsubsection{The case of no endoscopy} We begin by describing the results in the simplifying case of no endoscopy as in \cite[Theorem 6.3.2]{HaistBC}.

\begin{assumption}\label{ass:coh-decomp-no-end} We assume that:
\begin{itemize}
\item $\mb{G}$ has \emph{no endoscopy} in the sense of
\cite[Theorem 6.3.2]{HaistBC};\footnote{Equivalently, the group $\mf{E}(I_{\gamma_0},\mb{G};\A/\Q)$ from Definition \ref{defn:Kottwitz-invariant-group} vanishes for every semisimple
$\gamma_0$ in $\mb{G}(\Q)$.}
\item Theorem \ref{thm:geometric-stabilization} holds for every $h$ and $f^p$, and for $\tau$ in $\Phi_v^mI_E$ for sufficiently large $m$;
\item for $\tau$ and $h$ as above \eqref{eq:principal-case-Scholze--Shin}, i.e., the principal case of the Scholze--Shin conjecture, holds.
\end{itemize}
\end{assumption}

We first specify the archimedean data necessary to state our results.

\begin{setup}\label{setup:Haines-archimedean-factor} We fix/write:
\begin{itemize} 
\item $\pi_\infty^0$  an essentially discrete-series representation of $\mb{G}(\R)$ with central and infinitesimal characters equal to those of $\xi^\vee$;
\item $f_{\pi_\infty^0}$ a pseudocoefficient for $\pi_\infty^0$; see \cite[\S6.3]{HaistBC} and \cite[\S5]{ScholzeShin};
\item $f_{\xi,\infty}=(-1)^d f_{\pi_\infty^0}$.
\end{itemize}
\end{setup}

We next specify the multiplicities which will occur in our formulae.

\begin{nota}\label{nota:Haines-automorphic-multiplicity} Let us write:
\begin{itemize}
\item $\Pi_\infty(\xi)$ for the irreducible admissible representations of
$\mb{G}(\R)$ with the same central and infinitesimal characters as $\xi^\vee$;
\item $m_{\mr{disc}}(\pi_f\otimes\pi_\infty)$ for the automorphic multiplicity as in
\cite[\S6.3]{HaistBC}, with coefficient conventions as in
\cite[\S5]{ScholzeShin}.
\end{itemize}
\end{nota}

We combine the above discussion to define the coefficient of the terms in our decomposition.

\begin{nota}[{\cite[Theorem 6.3.2]{HaistBC} and \cite[\S5]{ScholzeShin}}]\label{nota:Haines-coefficients} For an irreducible admissible representation $\pi_f$ of $\mb{G}(\A_f)$, set
\begin{equation}\label{eq:Haines-coefficient}
a_\xi(\pi_f)\defeq\sum_{\pi_\infty\in\Pi_\infty(\xi)}
m_{\mr{disc}}(\pi_f\otimes\pi_\infty)
\tr(f_{\xi,\infty}\mid\pi_\infty).
\end{equation}
\end{nota}

We now discuss the notation needed to describe the Galois-theoretic components of our decompositions.

\begin{nota}\label{nota:ss} If $V=\varinjlim V_i$ is a filtered colimit of $\ov{\Q}_\ell$-representations of $W_E$, write $V^\mr{ss}=\varinjlim V_i^\mr{ss}$ where $V_i^\mr{ss}$ is the usual Jordan--H\"{o}lder semisimplification. These operations naturally pass to the level of Grothendieck groups.
\end{nota}

\begin{nota}\label{nota:Weil-representation} Let $\pi_p$ be an irreducible smooth $G(\Q_p)$-representation. Write
\begin{equation*}
R_{\pi_p}\defeq \left(r_{-\mu_h}\circ
\varphi_{\pi_p}^{\mr{FS}}|_{W_E}\right)^{\mr{ss}}
\otimes|\cdot|_E^{-\tfrac{d}{2}},
\end{equation*}
where $|\tau|_E\defeq |\mr{Art}_E(\tau)|_E$ is normalized so that geometric Frobenius has absolute value $q^{-1}$.
\end{nota}

\begin{prop}[{cf.\@ \cite[Theorem 6.3.2]{HaistBC}}]\label{prop:cohomology-decomp-no-end}
With \emph{Setup \ref{setup:Haines-application}} and \emph{Assumption \ref{ass:coh-decomp-no-end}}, there is a decomposition of virtual
$\mb{G}(\A_f^p)\times\mathsf{K}_0\times W_E$-representations
\begin{equation}\label{eq:Haines-compact-cohomology}
H_c^*(\mb{G},\mb{X},\xi)^\mr{ss}
=\sum_{\pi_f=\pi^p\otimes\pi_p}a_\xi(\pi_f)
\bigg[\pi^p\boxtimes\pi_p|_{\mathsf{K}_0}\boxtimes R_{\pi_p}\bigg],
\end{equation}
where $\pi_f$ ranges over irreducible admissible $\mb{G}(\A_f)$-representations.
\end{prop}
\begin{proof} By linear independence of characters and the Brauer--Nesbitt theorem, it suffices to verify equality of the traces of
$\tau\times h\times f^p$ for all $h$ and $f^p$ as in
Assumption \ref{ass:coh-decomp-no-end} and all $\tau$ in
$\Phi_v^mI_E$ with $m$ sufficiently large; cf.\@ the argument in \cite[proof of Theorem 9.1]{ScholzeShin}. On the one hand, the pseudostabilization argument in the proof of \cite[Theorem 6.3.2]{HaistBC} shows that, given Theorem \ref{thm:geometric-stabilization}, one has that this trace on the left-hand side is
\begin{equation*}
\sum_{\pi_f=\pi^p\otimes \pi_p}a_\xi(\pi_f)
\tr(f^p\mid\pi^p)\tr(f_{\tau,h}\mid \pi_p).
\end{equation*}
On the other hand, the trace of the right-hand side is 
\begin{equation*}
\sum_{\pi_f=\pi^p\otimes\pi_p}a_\xi(\pi_f)
\tr(f^p\mid\pi^p)\tr(h\mid\pi_p)\tr(\tau\mid R_{\pi_p}).
\end{equation*}
Then, the principal case of the Scholze--Shin conjecture says that these quantities are equal. 
\end{proof}

\begin{rem} Of course, it would be desirable to have the conclusion of Proposition \ref{prop:cohomology-decomp-no-end} hold as virtual $\mb{G}(\A_f)\times W_E$-representations; this is closely related to \cite[Question 7.5]{ScholzeShin}. 
\end{rem}

This has an immediate finite-level corollary by passing to invariants.\footnote{We are implicitly using the fact that compact-group invariants are exact in characteristic $0$.}

\begin{cor}\label{cor:cohomology-decomp-no-end-finite-level}
With setup as in \emph{Proposition \ref{prop:cohomology-decomp-no-end}}, if $\mathsf{K}=\mathsf{K}_p\mathsf{K}^p$ is neat and satisfies $\mathsf{K}_p\subseteq\mathsf{K}_0$, then one has a decomposition in the Grothendieck group of $\ov{\Q}_\ell[W_E]$-modules:
\begin{equation}\label{eq:Haines-finite-level-cohomology}
H^\ast_c(\mb{G},\mb{X},\mathsf{K},\xi)^\mr{ss}
=\sum_{\pi_f=\pi^p\otimes\pi_p}a_\xi(\pi_f)\dim(\pi_f^{\mathsf{K}}) R_{\pi_p}.
\end{equation}
\end{cor}

As in \cite{HaistBC}, we obtain implications for the partial semisimple Hasse--Weil $\zeta$-functions of Shimura varieties, the definition of which we now recall.

\begin{nota}\label{nota:ss-local-factor}
Let $V$ be a finite-dimensional virtual $\ov{\Q}_\ell$-representation
of $\Gamma_{\mb{E}}$, and let $S$ be a set of finite places of
$\mb{E}$ not dividing $\ell$.
For $v$ in $S$, put $V_v=V|_{W_{\mb{E}_v}}$ and $q_v=|\kappa(v)|$,
and let $I_v$ and $\Phi_v$ denote inertia and a geometric Frobenius
lift in $W_{\mb{E}_v}$, respectively. Define
\begin{equation*}
 \zeta_S^{\mr{ss}}(V,s)
 =\prod_{v\in S}
 \det\left(1-q_v^{-s}\Phi_v\mid(V_v^{\mr{ss}})^{I_v}\right)^{-1}.
\end{equation*}
\end{nota}

\begin{nota}\label{nota:semisimple-Shimura-zeta} Let us fix a neat level $\mathsf{K}$ and define $S(\mathsf{K})$ to be the set of places $v$ of $\mb{E}$ not dividing $\ell$ lying over a prime $p$ where $\mathsf{K}_p$ is contained in a hyperspecial subgroup of $\mb{G}(\Q_p)$. We then write
\begin{equation*}
\zeta^\mr{ss}_{S(\mathsf{K})}(\mb{G},\mb{X},\mathsf{K},\xi,s)\defeq \zeta_{S(\mathsf{K})}^\mr{ss}\left(H^\ast_c(\mb{G},\mb{X},\mathsf{K},\xi),s\right).
\end{equation*}
\end{nota}

\begin{cor}
Assume the setup of \emph{Corollary \ref{cor:cohomology-decomp-no-end-finite-level}} holds for every place $v$ in $S(\mathsf{K})$, and write $R_{\pi_p,v}$ for $R_{\pi_p}$ formed with $E=\mb{E}_v$. Then, we have an equality
\begin{equation*}
\zeta^\mr{ss}_{S(\mathsf{K})}(\mb{G},\mb{X},\mathsf{K},\xi,s)=\prod_{v\in S(\mathsf{K})}\prod_{\pi_f=\pi^p\otimes\pi_p}\det\left(1-q_v^{-s}\Phi_v\mid R_{\pi_p,v}^{I_v}\right)^{
-a_{\xi}(\pi_f)\dim\pi_f^{\mathsf{K}}}.
\end{equation*}
\end{cor}

\subsubsection{Endoscopic cohomology and local factors}\label{sss:end-decomp}
We now extend the preceding discussion to cases where the no-endoscopy hypothesis does not hold.

\begin{setup}\label{setup:Haines-endoscopic-application}
Retain Setup \ref{setup:Haines-application}, and fix:
\begin{itemize}
\item a compact open level $\mathsf{K}=\mathsf{K}_p\mathsf{K}^p$ with
$\mathsf{K}_p\subseteq\mathsf{K}_0$ and
$\mathsf{K}_0\mathsf{K}^p$ neat;
\item $h=e(\mathsf{K}_p)$ and $f^p=e(\mathsf{K}^p)$;
\item an integer $m_0\geqslant1$ such that
\eqref{eq:geometric-stabilization} holds for every $m\geqslant m_0$
and $\tau$ in $\Phi_v^mI_E$.
\end{itemize}
\end{setup}

\begin{setup}
For each $\mf{e}$ in $\mc{E}_{\mr{ell}}(\mb{G})$, we:
\begin{itemize}
\item write $H_1=\mb{H}_{1,\Q_p}$;
\item fix a transfer $h^{\mf{e}}$ in $\mc{H}_{\chi_1}(H_1(\Q_p))$ of $h$;
\item use the $z$-pairs, measures, and factors $f^{\mf{e},p}$ and $f_\xi^{\mf{e}}$ from Notations \ref{nota:stabilization-transfers} and \ref{nota:stabilization-global-function};
\item and use the central elements $z_{\mu,\tau}^{\mf{e}}$ associated with Construction \ref{constr:Scholze-Shin-spectral-element}.
\end{itemize}
\end{setup}

\begin{nota}\label{nota:Haines-endoscopic-localization}
For each $\mf{e}$, we choose the central data of Notation
\ref{nota:stabilization-central-character} with $U_p$ small enough that
$h^{\mf{e}}$ transforms under $q_{\mf{e}}^{-1}(U_p)$ by
$\Omega_{\mf{e},p}^{-1}$; the existence of a choice follows as in \cite[\S7.1.4]{KSZ}. 
\end{nota}

\begin{assumption}\label{ass:Haines-endoscopic-spectral-expansion}
For each $\mf{e}$, assume there are constants $b_{\mf{e},j}$ in $\C$ and stable
virtual characters $\Theta_{\mf{e},j}$ of $H_1(\Q_p)$, indexed by
$\mc{J}_{\mf{e}}$, such that:
\begin{itemize}
\item The nonzero constituents of $\Theta_{\mf{e},j}$ have
$q_{\mf{e}}^{-1}(U_p)$-character $\Omega_{\mf{e},p}$, and common Fargues--Scholze parameter $[\lambda_{\mf{e},j}]$ in $X_{H_1,\chi_1}^{\mr{spec}}$.
\item For every $\tau$ as above,
\begin{equation}\label{eq:Haines-endoscopic-stable-expansion}
\mr{ST}_{\mr{ell},\Omega_{\mf{e}}}^{\mb{H}_1}
(f_{\tau,h,\xi}^{\mf{e}})
=\mr{ST}_{\mr{disc},\Omega_{\mf{e}}}^{\mb{H}_1}
(f_{\tau,h,\xi}^{\mf{e}})
=\sum_{j\in\mc{J}_{\mf{e}}}b_{\mf{e},j}
\Theta_{\mf{e},j}(f_{\tau,h}^{\mf{e}}),
\end{equation}
where the first term is as in Definition \ref{defn:stabilization-stable-elliptic-distribution} and the second is as in \cite[\S9.1.2]{KSZ}.
\end{itemize}
We further assume that:
\begin{itemize}
\item Only finitely many pairs $(\mf{e},j)$ satisfy
$b_{\mf{e},j}\Theta_{\mf{e},j}(h^{\mf{e}})\ne0$.
\item For every $(\mf{e},j)$ with $b_{\mf{e},j}\ne0$, there are a
finite-dimensional semisimple $\C$-representation $V_{\mf{e},j}$ of $W_E$
with finite inertia image and
$S_{\mf{e},j}$ in $\End_{W_E}(V_{\mf{e},j})$ such that
\begin{equation}\label{eq:Haines-endoscopic-spectral-trace}
z_{\mu,\tau}^{\mf{e},\mr{spec}}([\lambda_{\mf{e},j}])
=\tr(S_{\mf{e},j}\tau\mid V_{\mf{e},j})
\end{equation}
for every $\tau$ as above.
\item Conjecture \ref{conj:Scholze-Shin} holds for every $\mf{e}$ and $\tau$
as above.
\end{itemize}
\end{assumption}

\begin{lem}\label{lem:Haines-endoscopic-local-calculation}
Under the preceding assumptions, for $b_{\mf{e},j}\ne0$ one has
\begin{equation}\label{eq:Haines-endoscopic-local-calculation}
\Theta_{\mf{e},j}(f_{\tau,h}^{\mf{e}})
=\tr(S_{\mf{e},j}\tau\mid V_{\mf{e},j})
\Theta_{\mf{e},j}(h^{\mf{e}}).
\end{equation}
\end{lem}
\begin{proof}
On each constituent, $z_{\mu,\tau}^{\mf{e}}$ acts by the common scalar
$z_{\mu,\tau}^{\mf{e},\mr{spec}}([\lambda_{\mf{e},j}])$, by
\eqref{eq:Fargues-Scholze-eigenvalue} and Proposition
\ref{prop:FS-fixed-central-character}. So the claim follows by applying
\eqref{eq:Haines-endoscopic-spectral-trace} and \cite[Corollary 3.2.1]{HaistBC}.
\end{proof}

\begin{nota}\label{nota:Haines-endoscopic-coefficients}
Let $\Sigma_{\mathsf{K},\xi}$ represent the isomorphism classes of
irreducible $W_E$-constituents of the $V_{\mf{e},j}$ with
$b_{\mf{e},j}\Theta_{\mf{e},j}(h^{\mf{e}})$ nonzero, and set for $\sigma$ in $\Sigma_{\mathsf{K},\xi}$:
\begin{equation}\label{eq:Haines-endoscopic-coefficient}
n_{\mathsf{K},\xi}(\sigma)\defeq
\sum_{\substack{\mf{e},j\\
 b_{\mf{e},j}\Theta_{\mf{e},j}(h^{\mf{e}})\ne0}}\iota(\mb{G},\mf{e})b_{\mf{e},j}
\Theta_{\mf{e},j}(h^{\mf{e}})
\tr\left(S_{\mf{e},j}\mid
\Hom_{W_E}(\sigma,V_{\mf{e},j})\right),
\end{equation}
with $\iota(\mb{G},\mf{e})$ as in Notation \ref{nota:stabilization-global-function}.
\end{nota}

The following is proven in much the same way as Proposition \ref{prop:cohomology-decomp-no-end}.

\begin{thm}\label{thm:Haines-endoscopic-cohomology}
Under \emph{Setup \ref{setup:Haines-endoscopic-application}} and
\emph{Assumption \ref{ass:Haines-endoscopic-spectral-expansion}}, one has
\begin{equation}\label{eq:Haines-endoscopic-cohomology}
H_c^\ast(\mb{G},\mb{X},\mathsf{K},\xi)^{\mr{ss}}
=\sum_{\sigma\in\Sigma_{\mathsf{K},\xi}}
n_{\mathsf{K},\xi}(\sigma)[\sigma]
\end{equation}
in the Grothendieck group of $\ov{\Q}_\ell[W_E]$-modules. 
\end{thm}

\begin{cor}\label{cor:Haines-endoscopic-zeta}
Assume that the setup in \emph{Theorem
\ref{thm:Haines-endoscopic-cohomology}} holds for every place $v$ in $S(\mathsf{K})$. Then, one has an equality
\begin{equation}\label{eq:Haines-endoscopic-zeta}
\zeta^{\mr{ss}}_{S(\mathsf{K})}(\mb{G},\mb{X},\mathsf{K},\xi,s)
=\prod_{v\in S(\mathsf{K})}
\prod_{\sigma\in\Sigma_{\mathsf{K},\xi,v}}
\det\left(1-q_v^{-s}\Phi_v\mid\sigma^{I_v}\right)^{-n_{\mathsf{K},\xi,v}(\sigma)}.
\end{equation}
\end{cor}

\appendix

\section{Some $\ell$-independence results for the cohomology of rigid spaces}\label{s:l-independence-appendix}

\addtocontents{toc}{\protect\setcounter{tocdepth}{1}}

In this section we extend some of the $\ell$-independence results from \cite{Mieda} to certain non-quasi-compact smooth rigid spaces. 

\begin{nota}  Fix $K/\Q_p$ a finite extension. We use
Notation \ref{nota:rigid-analytic-trace-formula}, but furthermore: 
\begin{itemize}
\item fix a uniformizer $\varpi$; 
\item and write $W_K^+=\{\tau\in W_K:\mathsf{v}(\tau)>0\}$.
\end{itemize}
\end{nota}

\subsection{The $p$-bounded generic fiber} To state our theorem precisely we must first develop some basic terminology and results concerning the `bounded generic fiber' of an affine formal scheme formally of finite type over $\mc{O}_K$.

\begin{defn}\label{defn:p-bounded-generic-fiber} Let $\mf{X}=\Spf(R)$ with $R$ an adic ring formally of finite type over $\mc{O}_K$. Then, the \emph{$p$-bounded generic fiber} of $\mf{X}$ is 
\begin{equation*}
\mf{X}^b_\eta\defeq \Spf(R,(p))_\eta\defeq \Spa(R)\times_{\Spa(\mc{O}_K)}\Spa(K)=\Spa(R[\nicefrac{1}{p}],\wt{R}), 
\end{equation*}
where in the third quantity $R$ is endowed with the $p$-adic topology, and where in the fourth $\wt{R}$ is the integral closure of $\mr{im}(R\to R[\nicefrac{1}{p}])$ in $R[\nicefrac{1}{p}]$.
\end{defn}

\begin{rem} While not called by this name, there is precedent for $p$-bounded generic fibers in the literature. For example, in \cite[\S3]{StrDefone}, such objects were studied under the name \emph{quasi-compactifications}. In \cite{Kappen}, Kappen defines appropriate globalizations of $p$-bounded generic fibers under the name \emph{uniformly rigid spaces}. That said, these uniformly rigid spaces are not adic spaces in the sense of Huber, but likely can be interpreted in Clausen--Scholze's language of analytic stacks. The results of this section likely extend to these settings, correctly interpreted. 
\end{rem}

To help conceptualize this notion, we study the topological properties of the $p$-bounded fiber in the simplest case which, in fact, is the only case used in the main body of the paper.

\begin{eg}[The $p$-bounded open unit disk]\label{eg:bounded-unit-open-disk} For $d\geqslant 1$, let us set 
\begin{equation*}
R_d\defeq \mc{O}_K\llbracket T_1,\ldots,T_d\rrbracket,\qquad \mf{m}_d=(\varpi,T_1,\ldots,T_d).
\end{equation*}
We notate the bounded generic fiber of $\Spf(R_d,\mf{m}_d)$ as the \emph{bounded open disk}, and write it $\bb{D}^{d,b}_K$. 

Observe that there is a natural $p$-adically continuous map 
\begin{equation*}
j_d\colon (\mc{O}_K\langle T_1,\ldots,T_d\rangle,(p))\to (R_d,(p)),
\end{equation*}
and so we obtain a natural map of adic $K$-spaces $j_d\colon \bb{D}^{d,b}_K\to \bb{B}^d_K$. We claim that 

\begin{enumerate}
\item the open embedding $\bb{D}_K^d\to \bb{B}_K^d$ factorizes through an open embedding $\bb{D}_K^d\to \bb{D}_K^{d,b}$;
\item the image of $j_d$ is given by 
\begin{equation*}
    \bb{D}_K^d=\mr{sp}^{-1}(0)^\circ\subseteq \mr{im}(j_d)=\mr{sp}^{-1}(\Spec(\mc{O}_{\bb{A}^d_k,0}))\subseteq \bb{B}^d_K,
\end{equation*}
where $\mr{sp}\colon |\bb{B}^d_K|\to |\bb{A}^d_k|$ is the specialization map and $0$ is the origin in $\bb{A}^d_k(0)$;
\item $j_d$ is a homeomorphism onto its image when $d=1$;
\item but $j_d$ is non-injective for $d\geqslant 2$. 
\end{enumerate}

\noindent\textbf{Claim (1).} Let 
\begin{equation*}
\mr{sp}^b\colon |\bb{D}^{d,b}_K|\to |\Spec(k\llbracket T_1,
\ldots,T_d\rrbracket)|,
\end{equation*}
be the specialization map for $\Spf(R_d,(p))$. Then, as per usual, the map 
\begin{equation*}
(R_d,(p))\to (R_d,(p))^\wedge_{\mf{m}_d}=(R_d,\mf{m}_d)
\end{equation*}
induces an isomorphism 
\begin{equation*}
\bb{D}^d_K=\Spf(R_d,\mf{m}_d)_\eta\isomto (\mr{sp}^b)^{-1}(0)^\circ\subseteq \bb{D}^{d,b}_K,
\end{equation*}
That this gives a factorization of the inclusion $\bb{D}^d_K\hookrightarrow \bb{B}^d_K$ through $j_d$ follows as the composition
\begin{equation*}
    (\mc{O}_K\langle T_1,\ldots,T_d\rangle,(p))\xrightarrow{j_d}  (R_d,(p))\to (R_d,\mf{m}_d)
\end{equation*}
is the natural map. This verifies claim (1).

\medskip

\noindent\textbf{Claim (2).} Let us observe that $\mr{im}(j_d)\subseteq \sp^{-1}(\Spec(\mc{O}_{\bb{A}^d_k,0}))$ merely by the functoriality of specialization. To see the reverse inclusion, let $x$ belong to $\sp^{-1}(\Spec(\mc{O}_{\bb{A}^d_k,0}))$ and let $k(x)^+$ be its valuation ring. As the map of discrete rings $\mc{O}_K\langle T_1,\ldots,T_d\rangle_{\mf{m}_d}\to R_d$ is faithfully flat, we deduce that the map
\begin{equation*}
\Spec(R_d\otimes_{\mc{O}_K\langle T_1,\ldots,T_m\rangle_{\mf{m}_d}} k(x)^+)\to\Spec(k(x)^+),
\end{equation*}
is faithfully flat. Take a prime of the source lying over the maximal ideal of the target, and choose a valuation ring dominating it. Restricting this to $R_d[\nicefrac{1}{p}]$ gives an element of $\bb{D}^{d,b}_K$ which, by construction, maps to $x$ in $\bb{B}^d_K$. This verifies claim (2).

\medskip

\noindent\textbf{Claim (3).} In this case one can do a simple analysis, similar to the classification of points in $\bb{B}_K^1$, to deduce that the only points of $\bb{D}^{1,b}_K$ not lying in $\bb{D}^1_K$ are the Gauss point $\eta_G^{1,b}$ where, more generally
\begin{equation*}
   \eta_G^{d,b} \colon R_d[\nicefrac{1}{p}]\to \R,\qquad \sum_{(i_1,\ldots,i_d)\in\bb{N}^d}a_{(i_1,\ldots,i_d)} T_1^{i_1}\cdots T_d^{i_d} \mapsto \sup_{(i_1,\ldots,i_d)}|a_{(i_1,\ldots,i_d)}|,
\end{equation*}
(where $|a_{(i_1,\ldots,i_d)}|$ is using the native absolute value on $\mc{O}_K$), and the point
\begin{equation*}
   \eta_{G,-}^{b,1} \colon R_1[\nicefrac{1}{p}]\to \R\times \{1^-\}^\Z,\qquad \sum_{i\in\bb{N}}a_i T_1^i \mapsto \sup_{i}|a_i|(1^-)^i,
\end{equation*}
where $1<1^-<x$ for every $x$ in $(1,\infty)$. Thus, to show that $j_1$ is injective it suffices to show that $j_1(\eta^{b,1}_G)\ne j_1(\eta^{b,1}_{G,-})$; but this is clear.

\medskip

\noindent\textbf{Claim (4).} We do this for $d=2$, as the generalization for $d>2$ is obvious. Choose an element $f(T_1)$ in $T_1\mc{O}_K\llbracket T_1\rrbracket$ whose image in $T_1 k\llbracket T_1\rrbracket$ is transcendental over $k(T_1)$. There is a natural continuous quotient map $R_2[\nicefrac{1}{p}]\to R_1[\nicefrac{1}{p}]$ sending $T_2$ to $f(T_1)$. We may then pull back the $1$-dimensional Gauss valuation $\eta_G^{1,b}$ from above to obtain a valuation $\eta$. By inspection one has that $\eta\ne \eta_G^{2,b}$. That said, we claim that $j_2(\eta)=j_2(\eta_G^{2,b})$, and that both are equal to the usual Gauss norm $\eta_G^2$ on $\bb{B}^2_K$. It's clear that $j_2(\eta_G^{2,b})=\eta_G^2$, and so it remains to show that $j_2(\eta)=\eta_G^{2}$. But, this is equivalent to the claim that if $g(T_1,T_2)$ in $ K\langle T_1,T_2\rangle$ is non-zero, and we can scale it by powers of the uniformizer $\pi$ so that it lies in $\mc{O}_K\langle T_1,T_2\rangle$ and has non-zero reduction in $k[T_1,T_2]$, then $\eta(g)=1$. But, as $f(T_1)$ is transcendental mod $\pi$ we have that 
\begin{equation*}
g(T_1,f(T_1))\ne 0\mod \pi,
\end{equation*}
 and thus, by definition, we have that $\eta(g)=1$.
\end{eg}

While many of the features of Example \ref{eg:bounded-unit-open-disk} extend more generally, the only one we will need in the sequel is the following.

\begin{prop} Suppose that $\mf{X}=\Spf(R)$ is formally of finite type over $\mc{O}_K$. Then, there exists a canonical open embedding $\mf{X}_\eta\hookrightarrow \mf{X}_\eta^b$.
\end{prop}
\begin{proof} Let $I$ be an ideal of definition of $R$ which, without loss of generality, contains $p$. Then, we have a natural closed embedding $Z=\Spec(R/I)\to \Spec(R/p)$, and we observe that $\mf{X}_\eta = \mf{X}^b_\eta(Z)$. As $\mf{X}^b_\eta(Z)$ is open in $\mf{X}^b_\eta$, the claim follows.
\end{proof}

We will be interested not only in the $p$-bounded generic fibers of formally of finite type formal $\mc{O}_K$-schemes, but also in certain covers of them. But, the following says that, in fact, no new objects are gained by considering such covers.

\begin{prop}\label{prop:finite-over-bounded-is-bounded} Suppose that $\mf{X}=\Spf(R)$ is formally of finite type over $\mc{O}_K$ and write $X=\mf{X}_\eta$ and $X^b=\mf{X}_\eta^b$. Then, for any finite map $f\colon Y^b\to X^b$ there exists a formally of finite type formal $\mc{O}_K$-scheme $\mf{Y}=\Spf(S)$ such that $Y^b=\mf{Y}_\eta^b$ and a finite morphism $\mf{Y}\to \mf{X}$ modeling $f$. Moreover, if $f$ is open, then with $Y=\mf{Y}_\eta$ one has $Y=f^{-1}(X)$.
\end{prop}

To prove this proposition, it's useful to establish the following basic lemma whose proof is simple and thus left to the reader.

\begin{lem}\label{lem:basic-properties-formally-finite-type} Let $R$ be a formally of finite type $\mc{O}_K$-algebra, and $R\to S$ finite. Then, there is a unique topology on $S$ with $R\to S$ a continuous map of formally of finite type $\mc{O}_K$-algebras; namely the $IS$-adic topology where $I$ is an ideal of definition of $R$.
\end{lem}

\begin{proof}[Proof of Proposition \ref{prop:finite-over-bounded-is-bounded}]
Write $Y^b=\Spa(A,A^+)$, where $A$ is finite over $R[\nicefrac{1}{p}]$
and $A^+$ is the integral closure of $\wt{R}$ in $A$;
see \cite[\S1.4.4]{HuberEC}.
Choose an ideal of definition $I$ of $R$
containing $p$, as well as generators $a_1,\ldots,a_t$ of $A$ as an
$R[\nicefrac{1}{p}]$-module. As $A$ is finite over $R[\nicefrac{1}{p}]$,
each $a_i$ is integral over $R[1/p]$. Choose $N$ such that every $b_i=p^Na_i$ is integral over $R$, and set
\begin{equation*}
 S=\operatorname{im}(R\to A)[b_1,\ldots,b_t]\subseteq A;
\end{equation*}
so $S$ is finite over $R$ and $S[\nicefrac{1}{p}]=A$.
Equip $S$ with the $IS$-adic topology and set
$\mf{Y}=\Spf(S)$.
By Lemma \ref{lem:basic-properties-formally-finite-type},
$\mf{Y}$ is formally of finite type over $\mc{O}_K$,
and $\mf{Y}\to\mf{X}$ is finite.

Now, the $p$-adic topology on $S$ induces the given topology on $A$, and by transitivity of integrality, both $A^+$ and the integral closure $\wt{S}$ of $S$ in $A$ are the integral closure of the image of $R$ in $A$.Thus, we conclude that
\begin{equation*}
 \mf{Y}_\eta^b=\Spa(S[1/p],\wt{S})=\Spa(A,A^+)=Y^b,
\end{equation*}
compatibly with the map to $X^b$. 

To see that $f^{-1}(X)=Y$ when $f$ is open, we proceed as follows. Write
\begin{equation*}
    \mr{sp}_1\colon |X^b|\to |\Spec(R/p)|,\qquad \mr{sp}_2\colon |Y^b|\to |\Spec(S/p)|.
\end{equation*}
Thus, we have the following equalities
\begin{equation*}
 \begin{aligned} f^{-1}(X) &=f^{-1}(\mr{sp}_1^{-1}(\Spec(R/I))^\circ)\\ &= f^{-1}(\mr{sp}^{-1}(\Spec(R/I)))^\circ\\ &=\mr{sp}_2^{-1}(\mf{f}^{-1}(\Spec(R/IR))^\circ\\ &=\mr{sp}_2^{-1}(\Spec(S/IS))^\circ\\ &=Y,
 \end{aligned}
\end{equation*}
where in the second equality we have used the fact that $f$ is open.
\end{proof}

\subsection{Finiteness and vanishing of cohomology}

We first observe that results of Berkovich imply that the $\ell$-adic cohomology of $p$-bounded finite covers of $p$-bounded generic fibers of formally of finite type formal $\mc{O}_K$-schemes is finite.

\begin{defn}[{\cite[\S3.1]{BerkovichFiniteness}}]\label{defn:constructible} Let $\mf{X}$ be a locally formally of finite type formal $\mc{O}_K$-scheme. A sheaf of $\Lambda$-modules $\mc{F}$ on $\mf{X}_\eta$ is \emph{$\mf{X}$-constructible} if there exists an affine open cover $\{\Spf(A_i)\}$ of $\mf{X}$ such that $\mc{F}|_{\Spf(A_i)_\eta}$ is the pullback of a constructible $\Lambda$-module on $\Spec(A_i[\nicefrac{1}{p}])$.
\end{defn}

\begin{eg}[{\cite[\S3.1]{BerkovichFiniteness}}] With notation as in Definition \ref{defn:constructible}, any locally constant sheaf of $\Lambda$-modules on $\mf{X}_\eta$ is $\mf{X}$-constructible.
\end{eg}

\begin{prop}\label{prop:coh-properties} Let $\mf{X}$ be a formally of finite type formal $\mc{O}_K$-scheme of dimension $n$. Fix a locally constant $\Lambda$-module $\mc{F}$ on $\mf{X}_\eta$. Then, $H^i_\et(\mf{X}_C,\mc{F}_C)$ is finite as a $\Lambda$-module for each $i\geqslant 0$ and vanishes for $i>n$ \emph{(}resp.\@ $i>2n$\emph{)} if $\mf{X}$ is affine \emph{(}resp.\@ in general\emph{)}.
\end{prop}
\begin{proof} The first claim is \cite[Corollary 3.1.2]{BerkovichFiniteness}. The second claim follows as $\mf{X}_\eta$ is quasi-Stein as in \cite[\S3.3]{Hansen}; see \cite[p.\@ 234]{ScholzeLK}. Indeed, we may then apply \cite[Corollary 3.4]{Hansen} (resp.\@ \cite[Corollary 2.8.3]{HuberEC}) and the method of proof of \cite[Corollary 3.4]{Hansen}, using \cite[Proposition 5.9.4]{Fargues} to pass to general $\Lambda$-coefficients.
\end{proof}

The following is an immediate corollary of Proposition \ref{prop:coh-properties} by pushing forward $\mc{F}$ along $f$.

\begin{cor} Let $\mf{X}$ be a formally of finite type formal $\mc{O}_K$-scheme of dimension $n$ and let $X=\mf{X}_\eta$. Then, for any finite \'etale cover $f\colon Y\to X$ and any locally constant $\Lambda$-module $\mc{F}$ on $Y$ the $\Lambda$-module $H^i_\et(Y_C,\mc{F}_C)$ is finite for each $i\geqslant 0$ and vanishes for $i>n$ \emph{(}resp.\@ $i>2n$\emph{)} if $\mf{X}$ is affine \emph{(}resp.\@ in general\emph{)}.
\end{cor}

\subsection{($p$-bounded) generic fibers and $\ell$-independence} In this subsection we extend some of the $\ell$-independence results from \cite{Mieda} to certain non-quasi-compact rigid spaces.

\begin{thm}\label{thm:l-independence-formally-of-finite-type} Let $\mf{X}$ be a formally of finite type formal $\mc{O}_K$-scheme with $X=\mf{X}_\eta$ smooth over $K$. Then, for any $\tau$ in $W_K^+$ the quantity $\tr\left(\tau| H^\ast_\et(X_C,\ov{\Q}_\ell)\right)$ is an integer independent of $\ell$.
\end{thm}
\begin{proof} It suffices to treat $\mf{X}=\Spf(A)$. Choose an ideal of definition
$I=(f_1,\ldots,f_s)$ containing $p$, and set
\begin{equation*}
 X_m=\{x\in X: |f_a(x)|^m\leqslant |p|\text{ for all }a\},
\end{equation*}
which are smooth affinoid $K$-spaces exhausting $X$.

We claim that, for each $\ell\ne p$, restriction induces
$\Gamma_\eta$-equivariant isomorphisms
\begin{equation*}
 H^i_\et(X_C,\Q_\ell)\isomto H^i_\et((X_m)_C,\Q_\ell),
 \qquad\text{ for all }i\geqslant0,
\end{equation*}
for all sufficiently large $m$. Indeed, start by choosing a strictly semistable compact hypercovering $\mf{Y}_\bullet\to\mf{X}$ of $\mf{X}$ as in \cite[end of \S3.3]{BerkovichFiniteness}, and write $Y_a=(\mf{Y}_a)_\eta$ and $Y_{a,m}=Y_a\times_X X_m$. Note that each $\mf{Y}_a$ is of finite type over $\mf{X}$ and locally algebraizable. Now, by \cite[\S5.7 and Corollary 5.8]{HuberFinitenessII}, Poincar\'e duality, and descent for finite affine open covers, the restriction map
\begin{equation}\label{eq:restr-isom}
 H^j_\et(Y_{a,C},\bb{F}_\ell)
 \to H^j_\et(Y_{a,m,C},\bb{F}_\ell),\qquad \text{For all }j\geqslant 0
\end{equation}
once $m$ is sufficiently large. By \cite[\S1.2]{BerkovichFiniteness}, compact hypercoverings satisfy universal cohomological descent, which allows the comparison of the first-quadrant descent spectral sequences for $Y_\bullet\to X$ and $Y_{\bullet,m}\to X_m$. Writing $d=\dim X$, we observe that only simplicial degrees at most $2d+1$ are needed to compare their abutments through degree $2d$. The claim then follows by choosing $m$ large enough for these finitely many comparisons. Using standard arguments, one then deduces that \eqref{eq:restr-isom} holds with $\Z/\ell^n$-coefficients for every $n$, and thus for $\Q_\ell$-coefficients. 

As the desired independence-of-$\ell$ claim holds for $(X_m)_C$ by \cite[Theorem 7.1.10]{Mieda}, we deduce from the above discussion that it also holds for $X_C$, as desired.
\end{proof}

Using Theorem \ref{thm:l-independence-formally-of-finite-type}, we are able to prove the following proposition which is important in the main body of the text (see Theorem \ref{thm:phi-nice-properties}).

\begin{prop}\label{prop:independence-of-l-bounded} Suppose that $\mf{X}=\Spf(R)$ is formally of finite type over $\mc{O}_K$ and write $X=\mf{X}_\eta$, which we assume smooth, and $X^b=\mf{X}_\eta^b$. Then, for any finite \'etale map $f\colon Y^b\to X^b$ and any finite order element $g$ of $\Aut(Y^b/X^b)$ and any $\tau$ in $W_K^+$, the quantity $\tr\left(\tau\times g| H^\ast_\et(f^{-1}(X))_C,\ov{\Q}_\ell)\right)$ is in $\Z$ and independent of $\ell$.
\end{prop}

We first establish a Galois-descent-like result for formally of finite type $\mc{O}_K$-schemes.

\begin{lem}\label{lem:Galois-descent-formally-of-finite-type}
Suppose $\mf X=\Spf(R)$ is formally of finite type over
$\mc O_K$, and write $X=\mf X_\eta$.
Let $g$ be a finite-order $\mc O_K$-automorphism of $\mf X$,
and let $\tau$ be in $W_K^+$.
Put $r=\mathsf{v}_K(\tau)$, and let $L/K$ be the unramified
extension of degree $r$ inside $\ov{K}$. Then there is a formally of finite type $\mc{O}_L$-algebra $R'$ such that, writing $X'=\Spf(R')_\eta$, there is $f\colon X'_C\isomto X_C$ with $f\circ\tau_{X'_C}=(g\circ\tau_{X_C})\circ f$.
\end{lem}
\begin{proof} With $L$ as in the statement, replacing $(K,R)$ by
$(L,R\otimes_{\mc{O}_K}\mc{O}_L)$, we may assume that $\mathsf v_K(\tau)=1$. Let $m=|g|$. Set $K_m$ to be the unramified extension of $K$ of degree $m$. Then, $\mathcal{O}_K\to \mathcal{O}_{K_m}$ is a Galois cover with Galois group $\Gamma=\{\sigma^\Z\}\simeq \Z/m\Z$, where $\sigma=\tau|_{K_m}$. Let $\sigma_R$ denote the natural action of $\sigma$ on $R_{\mc{O}_{K_m}}$ and let $\Gamma$ act on $R_{\mc{O}_{K_m}}$ by $\sigma\cdot x=(\sigma_R\times g)(x)$. Let $T$ denote $R_{\mc{O}_{K_m}}$ with this action of $\Gamma$. Finally, write $\gamma=g\circ\sigma_R$.

By Galois descent for affine schemes we know that the natural map $T^\Gamma\otimes_{\mc{O}_K}\mc{O}_{K_m}\to T$ is an isomorphism. As $T$ is formally of finite type over $\mc{O}_{K_m}$, and thus formally of finite type over $\mc{O}_K$, we may choose an adic surjection $h\colon \mc{O}_K\langle x_1,\ldots,x_r\rangle\llbracket y_1,\ldots,y_s\rrbracket\to T$. Let us write $P= \mc{O}_K\langle X_{i,j}\rangle_{\substack{\mathsmaller {i=1,\ldots,r}\\ \mathsmaller{j=1\ldots,m}}}\llbracket Y_{\ell,k}\rrbracket_{\substack{\mathsmaller{\ell=1,\ldots,s}\\ \mathsmaller{k=1,\ldots,m}}}$, and define the map 
\begin{equation*}
h'\colon P\to T,\qquad X_{i,j}\mapsto e_j(\gamma(h(x_i)),\ldots,\gamma^m(h(x_i))),\quad Y_{\ell,k}\mapsto e_k(\gamma(h(y_\ell)),\ldots,\gamma^m(h(y_\ell)))
\end{equation*}
where $e_j$ and $e_k$ are the $j^\text{th}$ and $k^\text{th}$ elementary symmetric polynomials, respectively.

 This map is finite and continuous, and thus adic by Lemma \ref{lem:basic-properties-formally-finite-type}. Then, by design the map $h'$ factorizes through $T^\Gamma$ and, as $P$ is Noetherian, we deduce that the map $h'\colon P\to T^\Gamma$ must also be finite. In particular, we see from Lemma \ref{lem:basic-properties-formally-finite-type} that the $(p,Y_{\ell,k})$-adic topology on $T^\Gamma$ makes it into a formally of finite type $\mc{O}_K$-algebra. Moreover, by setup we see that the natural map $T^\Gamma\otimes_{\mc{O}_K}\mc{O}_{K_m}\to T$ is now an isomorphism of formally of finite type $\mc{O}_{K_m}$-algebras. 

It is clear by construction that $R'\defeq T^\Gamma$ satisfies the desired conditions.
\end{proof}
\begin{proof}[Proof of Proposition \ref{prop:independence-of-l-bounded}] To begin, we apply Proposition \ref{prop:finite-over-bounded-is-bounded} to obtain $\mf{Y}=\Spf(S)$ as in the statement of that proposition. Moreover, replacing the generators $b_i$ of $S$ as in the proof of loc.\@ cit.\@  their
$g$-translates we may assume that $g$ stabilizes $S$. We may then apply Lemma \ref{lem:Galois-descent-formally-of-finite-type} to obtain an $S'$ and a $Y'$ as in this statement. Then, by Theorem \ref{thm:l-independence-formally-of-finite-type} we know that $\tr\left(\tau| H^\ast_\et(Y'_C,\ov{\Q}_\ell)\right)$ is an integer independent of $\ell$. But, by design, we have an isomorphism $f\colon Y'_C\to Y_C$ such that $f\circ \tau_{Y'_C}=(g\times\tau_{Y_C})\circ f$, and thus 
\begin{equation*}
    \tr\left(\tau| H^\ast_\et(Y'_C,\ov{\Q}_\ell)\right)=\tr\left(g\times \tau| H^\ast_\et(Y_C,\ov{\Q}_\ell)\right),
\end{equation*}
from where the conclusion follows.
\end{proof}

\bibliographystyle{test2}
\bibliography{reference}

\clearpage
\phantomsection
\addtocontents{toc}{\protect\setcounter{tocdepth}{1}}
\addcontentsline{toc}{section}{\protect\tocsection{}{}{Index of notation}}
\markboth{Index of notation}{Index of notation}

\begin{center}
  \normalfont\scshape Index of notation
\end{center}
\medskip

\begingroup
\small
\setlength{\LTpre}{0pt}
\setlength{\LTpost}{0pt}
\setlength{\LTleft}{\fill}
\setlength{\LTright}{\fill}
\renewcommand{\arraystretch}{1.16}

\begin{longtable}{@{}
  >{\raggedright\arraybackslash}p{.62\linewidth}
  @{\hspace{.04\linewidth}}
  >{\raggedright\arraybackslash}p{.34\linewidth}
  @{}}

\toprule
\textbf{Symbol} & \textbf{Definition location} \\
\midrule
\endfirsthead

\toprule
\textbf{Symbol} & \textbf{Definition location} \\
\midrule
\endhead

\bottomrule
\endfoot

\bottomrule
\endlastfoot

\(\mf X_\eta\), \(\mf X(x)\)
  & Notation~\ref{nota:rigid-analytic-general} \\

\(u^\dashv\)
  & Notation~\ref{nota:adjoint-correspondence} \\

\(\mr{Fix}(c)\)
  & Definition~\ref{defn:fixed-scheme} \\

\(\Phi_Y\), \(c^{(n)}\)
  & Notation~\ref{nota:Frobenius-twist-corr} \\

\(\sigma\)
  & Remark \ref{rem:frobenius-galois-pullback}/
    Notation \ref{nota:deformation-spaces-setup}
     \\

\(\Psi_{\mf X}\), \(R\Psi_{\mf X}\)
  & Definition~\ref{defn:nearby-cycles} \\

\(\Fil^\bullet_{\mr{Hdg}}T_{\mr{dR}}(\mf Q)\)
  & Definition~\ref{defn:BTGmu} \\

\(\mr{BT}^{\mc G,\mu}_n(\mf X)\)
  & Definition~\ref{defn:BTGmu} \\

\(T_\et\)
  & Notation~\ref{nota:etale-realization} \\

\(\mr{BT}^{\mc G,\mu,\mr{alg}}_n(Y)\)
  & Definition~\ref{defn:algebraic-apertures} \\

\(\cat{FCrys}^{\nu}_{\mc G}(S)\)
  & Definition~\ref{defn:crystalline-type} \\

\(\Fil^\bullet_{\mr{MNO}}\mc Q_S\)
  & Definition~\ref{defn:MNO-filtration} \\

\(\cat{FilTors}^{\nu}_{\mc H}(\mc X)\)
  & Definition~\ref{defn:type-nu} \\

\(\cat{Quad}_{\mc G}^{\mu}(\mf X)\)
  & Definition~\ref{defn:quad} \\

\(k_j\), \(W_j\), \(K_j\)
  & Notation~\ref{nota:finite-field} \\

\(\mc C_j(\mc G)\), \(C_j(\mc G)\),
\(C_j(\mc G,\sigma(\mu))\)
  & Notation~\ref{nota:sigma-conjugacy-classes} \\

\(\cat{FFCrys}^{\mu}_{\mc G}(\mf X)\),
\(\cat{FFCrys}^{\mr{sd},\mu}_{\mc G}(\mf X)\)
  & Definition~\ref{defn:filtered-crystalline-type} \\

\(\cat{Quint}_{\mc G}^{\mu}(\mf X)\)
  & Definition~\ref{defn:quint} \\

\(\mc U^-_{\mu,j}\), \(\widehat{\mc U}^-_{\mu,j}\)
  & Definition~\ref{defn:opposite-unipotent} \\

\(R_{\mc G,\mu,j}\), \(\mf m_{\mc G,\mu,j}\)
  & Notation~\ref{nota:deformation-ring} \\

\(\mf Q_b^{\mr{univ}}\)
  & Notation~\ref{nota:universal-aperture} \\

\(\mf D(\mc G,b,\mu)\)
  & Notation \ref{nota:defm-space} \\

\(R_{\mr{red}}\)
  & \hyperref[conv:notation-and-conventions]
      {Notation and Conventions} \\

\(\phi_{\tau,h}^{\mc G,\mu}\), \(\phi_{\tau,h}\)
  & Definition~\ref{defn:phitauh} \\

\(\mc D(\mc G,b,\mu)\), \(\mc D_\infty\),
\(\mc D_{\mathsf K}\), \(\mc D_n\)
  & Definition~\ref{defn:aperture-tubes} \\

\(J_b^{\mr{int}}\)
  & Notation~\ref{nota:integral-sigma-centralizer} \\

\(\mc H_F(\mathsf G)\), \(\mc H(\mathsf G)\)
  & Definition~\ref{defn:Hecke-algebra} \\

\(e(g,\mathsf K)\), \(e(\mathsf K)\)
  & Notation~\ref{nota:normalized-Hecke-functions} \\

\(H^i(\mc G,b,\mu,\mathsf K)\), \(H^i(\mc G,b,\mu)\)
  & Definition~\ref{defn:coh-for-test-functions} \\

\(\Sh_{\mathsf K}(\mb G,\mb X)\), \(\Sh_{\mathsf K}\)
  & Setup~\ref{setup:sv} \\

\(\mb G^c\)
  & Definition~\ref{defn:cuspidal-quotient} \\

\(G^c\), \(\mc G^c\), \(\mu_h^c\), \(\mr{pr}^c\)
  & Notation~\ref{nota:local-cuspidal-quotient} \\

\(\mc F_{\xi,\mathsf K}\)
  & Remark/Definition~\ref{remdef:automorphic-sheaves} \\

\(\ms S_{\mathsf K}\)
  & Definition~\ref{defn:ICM} \\

\(\mf Q_{\mathsf K}\)
  & Definition~\ref{defn:syntomic-realization} \\

\(\mc U_{\mathsf L}\)
  & Definition~\ref{defn:pot-crys-loci} \\

\(H^i_c(\mb G,\mb X,\mathsf K,\xi)\),
\(\ovc H^i_c(\mb G,\mb X,\mathsf K,\xi)\)
  & Notation~\ref{nota:Shimura-cohomology} \\

\(c(g,\mathsf K)\)
  & Definition~\ref{defn:geometric-Hecke-corr} \\*

\(u(g,\mathsf K)\)
  & Definition~\ref{defn:cohomological-Hecke-correspondence} \\

\(\delta(y)\)
  & Construction~\ref{constr:corr-delta} \\

\(X^p(\varphi)\), \(X(\varphi)\)
  & Definition~\ref{defn:KSZ-prime-to-p-LR-set} \\

\(S_{\tau_{\mr{LR}}}(\varphi)\)
  & Definition~\ref{defn:KSZ-LR-set} \\

\(T(\tau,h,f^p)\)
  & Definition~\ref{defn:weighted-KSZ-trace} \\

\(\phi_{\tau,h}^{\mb G,\mb X}\)
  & Definition~\ref{defn:Shimura-test-function} \\

\(L_\gamma\), \(I_\gamma\)
  & \hyperref[conv:notation-and-conventions]
      {Notation and Conventions} \\

\(c=(\gamma_0,a,[b])\), \(\mf{KP}\)
  & Definition~\ref{defn:KSZ-Kottwitz-parameter} \\

\(\alpha(c)\)
  & Definition~\ref{defn:KSZ-Kottwitz-invariant} \\

\(\mr O_\gamma(f)\)
  & Definition~\ref{defn:orbital-integrals} \\

\(\mr{TO}_\delta(\phi)\)
  & Definition~\ref{defn:twisted-orbital-integrals} \\

\(\mr O(c,\tau,h,f^p)\)
  & Definition~\ref{defn:weighted-orbital-term} \\

\(\phi^{\mf e}\)
  & Definition~\ref{defn:twisted-endoscopic-transfer} \\

\(f^{\mf e}_{\tau,h,\xi}\), \(\iota(\mb G,\mf e)\)
  & Notation~\ref{nota:stabilization-global-function} \\

\(\mr{ST}^{\mb H_1}_{\mr{ell},\Omega_{\mf e}}\)
  & Definition~\ref{defn:stabilization-stable-elliptic-distribution} \\

\(\mr{SO}_\gamma(f)\)
  & Definition~\ref{defn:orbital-integrals} \\

\(\mc H_{\chi_1}(H_1(\Q_p))\)
  & Notation~\ref{nota:fixed-central-character-Hecke-algebra} \\

\(h^{\mf e}\)
  & Definition~\ref{defn:ordinary-endoscopic-transfer} \\

\(\mc Z(G)\)
  & Definition~\ref{defn:Bernstein-center} \\

\(\Omega(G)\), \(\Omega_{\mf s}(G)\)
  & Theorem~\ref{thm:Bernstein-decomposition} \\

\(\mc Z^{\mr{st}}(G)\), \(\mc Z^{\mr{vst}}(G)\)
  & Definition~\ref{defn:stable-very-stable-center} \\

\(X_G^{\mr{spec}}\)
  & Notation~\ref{nota:spectral-parameter-space} \\

\(\mc Z^{\mr{spec}}(G)\)
  & Equation~\eqref{eq:spectral-Bernstein-center} \\

\(\Psi_G^{\mr{FS}}\)
  & Theorem~\ref{thm:Fargues-Scholze-Bernstein-map} \\

\(\varphi_\pi^{\mr{FS}}\)
  & Remark~\ref{rem:parameter-description-of-Psi} \\

\(\omega_1\), \(\mc Z(H_1,\chi_1)\)
  & Definition~\ref{defn:fixed-central-character-Bernstein-center} \\

\(z_{\mu,\tau}^{\mf e,\mr{spec}}\)
  & Construction~\ref{constr:Scholze-Shin-spectral-element} \\

\(z_{\mu,\tau}^{\mf e}\)
  & Definition~\ref{defn:Scholze-Shin-central-element} \\

\(z_{\mu,\tau}^{\mr{spec}}\), \(z_{\mu,\tau}\)
  & Remark~\ref{rem:principal-Scholze-Shin-central-elements} \\

\(a_\xi(\pi_f)\)
  & Notation~\ref{nota:Haines-coefficients} \\

\(R_{\pi_p}\)
  & Notation~\ref{nota:Weil-representation} \\

\(\zeta^{\mr{ss}}_{S(\mathsf K)}(\mb G,\mb X,\mathsf K,\xi,s)\)
  & Notation~\ref{nota:semisimple-Shimura-zeta} \\

\(\mf X_\eta^b\)
  & Definition~\ref{defn:p-bounded-generic-fiber} \\

\end{longtable}

\endgroup

\end{document}